\documentclass{article}
\usepackage{arxiv}
\usepackage{comment}
\usepackage[utf8]{inputenc} 
\usepackage[T1]{fontenc}    
\usepackage{hyperref}       
\usepackage{url}            
\usepackage{booktabs}       
\usepackage{amsmath}
\usepackage{amsfonts}       
\usepackage{nicefrac}       
\usepackage{microtype}      
\usepackage{xcolor}
\usepackage{graphicx}
\usepackage{stmaryrd}
\usepackage{bm}
\usepackage[normalem]{ulem}
\usepackage[tickmarkheight=0.1cm]{todonotes}
\usepackage[ruled,vlined]{algorithm2e}

\hypersetup{
    colorlinks=true, 
    linkcolor=blue, 
    urlcolor=red, 
    citecolor=magenta,
    linktoc=all 
}

\newcommand{\re}{\mathbb{R}} 
\newcommand{\im}{\mathfrak{i}} 

\newcommand{\poly}{\mathbb{P}} 

\newcommand{\emh}{{e-\frac{1}{2}}}
\newcommand{\eph}{{e+\frac{1}{2}}}
\newcommand{\half}{\frac{1}{2}}
\newcommand{\ud}{\textrm{d}}
\newcommand{\pd}[2]{\frac{\partial #1}{\partial #2}}
\newcommand{\od}[2]{\frac{\ud #1}{\ud #2}}

\usepackage{caption}
\usepackage{subcaption}

\newcommand{\avg}[1]{\overline{#1}}
\newcommand{\au}{\avg{u}}

\newcommand{\myvector}[1]{\mathsf{#1}}
\newcommand{\vu}{\myvector{u}}
\newcommand{\vA}{\myvector{A}}
\newcommand{\vD}{\myvector{D}}
\newcommand{\vH}{\myvector{H}}
\newcommand{\vT}{\myvector{T}}
\newcommand{\vI}{\myvector{I}}
\newcommand{\vU}{\myvector{U}}
\newcommand{\vf}{\myvector{f}}
\newcommand{\vg}{\myvector{g}}
\newcommand{\vF}{\myvector{F}}
\newcommand{\vG}{\myvector{G}}
\newcommand{\vb}{\myvector{b}}

\newcommand{\vV}{\myvector{V}}
\newcommand{\vp}{\myvector{p}}
\newcommand{\cfl}{\mathrm{CFL}}
\newcommand{\bs}{\boldsymbol}

\newcommand{\Kref}{\hat{K}} 

\newcommand{\rnum}[1]{\left\llbracket #1 \right\rrbracket}

\DeclareMathOperator{\sign}{sign}
\DeclareMathOperator{\minmod}{minmod}

\newtheorem{remark}{Remark}

\title{Lax-Wendroff flux reconstruction method for hyperbolic conservation laws}

\author{
Arpit~Babbar \\
Centre for Applicable Mathematics\\
Tata Institute of Fundamental Research\\
Bangalore -- 560065\\
\texttt{arpit@tifrbng.res.in} \\
\And
Sudarshan~Kumar~Kenettinkara\\
School of Mathematics\\
Indian Institute of Science Education and Research\\
Trivandrum --695551\\
\texttt{sudarshan@iisertvm.ac.in} \\
\And
Praveen~Chandrashekar\thanks{Corresponding author} \\
Centre for Applicable Mathematics\\
Tata Institute of Fundamental Research\\
Bangalore -- 560065\\
\texttt{praveen@math.tifrbng.res.in}
}
\begin{document}
\maketitle
\begin{abstract}
The Lax-Wendroff method is a single step method for evolving time dependent solutions governed by partial differential equations, in contrast to Runge-Kutta methods that need multiple stages per time step. We develop a flux reconstruction version of the method in combination with a Jacobian-free Lax-Wendroff procedure that is applicable to general hyperbolic conservation laws. The method is of collocation type, is quadrature free and can be cast in terms of matrix and vector operations. Special attention is paid to the construction of numerical flux, including for non-linear problems, resulting in higher CFL numbers than existing methods, which is shown through Fourier analysis and yielding uniform performance at all orders. Numerical results up to fifth order of accuracy for linear and non-linear problems are given to demonstrate the performance and accuracy of the method.
\end{abstract}
\section{Introduction}
The numerical solution of hyperbolic PDE is usually based on the Godunov approach where a Riemann problem is solved to estimate the numerical flux which may be used in a finite volume or discontinuous Galerkin (DG) method. The resulting set of ODE are then marched forward in time using a multi-stage Runge-Kutta scheme under some CFL condition. For hyperbolic problems, the spatial and time order are usually taken to be equal due to the CFL condition $\Delta t = O(\Delta x)$, which requires atleast as many stages as the order of the scheme. The need to perform many stages of RK scheme means that costly limiters and data exchange in case of parallel computations must be performed several times for each time step~\cite{Dumbser2018}. There are also some order barriers in the sense that at high orders, we need more stages than the order of the method. In contrast to RK schemes, the Lax-Wendroff method achieves high order accuracy in time with a single step but requires more information to estimate the higher order time derivative terms in a Taylor expansion. It is expected that the single stage nature of LW-type schemes can make them more efficient than RK schemes~\cite{Qiu2005b}.

In the context of hyperbolic conservation laws, the Lax-Wendroff (LW) time discretization in conjunction with a wide range of   spatial schemes were extensively studied in the literature. These temporal schemes are essentially based on the classical second-order Lax-Wendroff method~\cite{Lax1960}. The Lax-Wendroff temporal discretization, originally  referred to as Taylor-Galerkin method, was used in the continuous finite element spatial schemes by Safian et al.~\cite{Safjan1995} and Tabarrok et al.~\cite{Tabarrok1994}, followed by further improvements in~\cite{Youn1995}. The case of  discontinuous finite element spatial schemes was studied in~\cite{Choe1991,Choe1992}. All these methods are confined to a certain order of accuracy in both space and time. In the finite difference framework, the LW time discretization was originally proposed by Qui and Shu~\cite{Qiu2003} with the  WENO approximation of spatial derivatives~\cite{Shu1989}. As an extension to this, a combination with alternative WENO method was developed in~\cite{Jiang2013}.  The discontinuous Galerkin spatial discretization combined with the LW temporal  scheme was originally proposed in~\cite{Qiu2003,Qiu2005b} (abbreviated as LWDG)  with an advantage of having arbitrary order of accuracy in both space and time, in other words, with no theoretical order barrier. It was further studied in~\cite{Qiu2007}, where the performance of various numerical fluxes were analyzed for the Euler equations of compressible flaws.   It is observed that the LWDG schemes are more compact and cost effective for certain problems like the two dimensional Euler system of compressible gas dynamics, especially when nonlinear limiters are applied. In~\cite{Guo2015} it is found that the LWDG method of~\cite{Qiu2005b} need not exhibit the super-convergence property. In order to overcome this issue, a modified version of LWDG was proposed~\cite{Guo2015} using the local DG framework of Cockburn et al.~\cite{Cockburn1998a}. The resulting scheme was found to satisfy the super-convergence property.  For linear conservation laws the stability and accuracy properties of LWDG scheme are explored in~\cite{Sun2017} with the modified LWDG scheme of Guo et al.~\cite{Guo2015}.

 Another significant contribution towards the single step temporal discretization was made by Toro et al., initially for linear equations  in~\cite{Toro2001} and for nonlinear systems in~\cite{Titarev2002}, following the idea of generalized Riemann problem (GRP)~\cite{Artzi2006,Han2010}. These  are widely known as arbitrary high order derivative (ADER) methods. Though its inception was in the finite volume spatial setup, later it was extended to finite difference and discontinuous Galerkin frameworks~\cite{Dumbser2006}. In the sequel, several authors have contributed to this approach with the aim of shaping up  a compact single time step scheme, see~\cite{Titarev2005,Kaser2005,Dumbser2007,Castro2008} and references therein. In the flavour of ADER methods, Dumbser et al. proposed an efficient DG spatial solver in~\cite{Dumbser2008} and a finite difference WENO spatial solver in~\cite{Dumbser2008a}. These are compact schemes that replace the so called Cauchy-Kowalevski procedure in the original ADER scheme with an element local space-time Galerkin  predictor step and a discontinuous Galerkin  corrector step, which are also found to  be suitable for stiff source terms and further studied by Gassner et al.in~\cite{Gassner2011a}.  These methods have been extended to the divergence free MHD problems  with a finite volume WENO spatial scheme in~\cite{Balsara2009a}. Through a modification of the method in~\cite{Dumbser2008,Gassner2011a}, Guthrey et al. in~\cite{Guthrey2019} proposed a regionally implicit ADER discontinuous Galerkin solver which is stable for higher CFL numbers. A simplified Cauchy-Kowalevski procedure is developed in~\cite{Montecinos2020} which is efficient, easier to implement for any system and can be used in ADER type schemes.

The generic versions of LWDG and ADER methods require the computation of high-order flux derivatives for each hyperbolic system and may require the use of symbolic manipulation software to perform the algebra. At higher orders of accuracy, we need higher order derivatives which needs the computation of flux Jacobian and other higher order tensors. This increases the computational task and the process has to be performed for each PDE system. In order to overcome this difficulty,  an approximation procedure was originally developed in~\cite{Zorio2017} in the finite difference scenario and further studied by several other authors~\cite{Lee2021,Burger2017,Carrillo2021, Carrillo2021a}. These approximation procedures for LW type solvers are found to be computationally more efficient and easier in implementation. As a single time step method, the resulting schemes are efficient for solving  hyperbolic conservation laws. Moreover, it is independent of the specific form of flux function in the governing equation as it is free from Jacobian and other higher version of derivatives.

The flux reconstruction (FR) method~\cite{Huynh2007} is a class of  discontinuous Spectral Element Method for the discretization of conservation laws. FR methods utilize a nodal basis which is usually based on some solution points like Gauss points, to approximate the solution with piecewise polynomials. The main idea is to construct a continuous approximation of the flux utilizing a numerical flux at the cell interfaces and a correction function. The solution at the nodes is then updated by a collocation scheme in combination with a Runge-Kutta method. The choice of the correction function affects the accuracy and stability of the method; by properly choosing the correction function and solution points, FR method can be shown to be equivalent to some discontinuous Galerkin and spectral difference schemes, as shown in~\cite{Huynh2007,Trojak2021}. In~\cite{Vincent2011a}, linear stability analysis of FR is performed through a broken Sobolev norm, leading to a 1-parameter family of correction functions which encompasses the stable correction functions found in~\cite{Huynh2007}. The family of stable correction functions has been extended in~\cite{Vincent2015,Trojak2021}, see~\cite{Trojak2021} for a review. For the 1-parameter family of correction functions in~\cite{Vincent2011a}, non-linear stability for E-fluxes was studied in~\cite{Jameson2012} where the significance of solution points was pointed out, with Gauss-Legendre points being the most resistant to aliasing driven instabilities. In another study on accuracy with different choices of solution points~\cite{Witherden2021}, the optimality of Gauss-Legendre points was again observed. The long term error behaviour of FR schemes has been studied in~\cite{Offner2019,Abgrall2020}, while dispersion and dissipation errors have been analyzed in~\cite{Vincent2011,Asthana2015,Vermeire2017}. The computationally efficient performance of FR has been noted in~\cite{Vincent2016,Lopez2014,Vandenhoeck2019}, which is attributed to the structured computation of finite element methods  suitable for modern hardware~\cite{Vincent2016}. The quadrature-free nature of FR methods together with the ability to cast the operations as matrix-vector operations that can be performed efficiently using optimized BLAS kernels makes these methods ideal for use on modern vector processors~\cite{Vincent2016}.

In the present work, we combine the Lax-Wendroff method for time discretization with the FR method for spatial discretization, since each of these two methods have their advantages as discussed above.  The flux reconstruction schemes have been combined with the LW time discretization in~\cite{Lou2020}, where the flux derivatives are computed in the standard way by explicitly using chain rule of differentiation and making use of the PDE. This is a tedious process that has to be performed for each PDE system and leads to large amount of computations. The results in~\cite{Burger2017} reveal that this approach would be computationally more expensive and depends on the nature of flux in the governing equation. In this work we propose to combine the approximate LW procedure~\cite{Zorio2017} with the FR scheme in space which leads to a general method that can be applied to any PDE system. No Jacobian or other computations are required and only the flux function is used in this method. The time derivatives of the flux in the Taylor expansion are approximated using a finite difference method. This allows us to write a generic solver that can easily applied to solve any PDE system. The single step method is achieved by using an approximation to the time average flux. A numerical flux is required for the time average flux and this needs to be computed carefully in order to obtain good CFL numbers. In some previous works, the solution at the current time level has been used to estimate the dissipative part of the numerical flux; however, it reduces the CFL number based on Fourier stability analysis and does not lead to an upwind flux, even for the linear advection equation. Here we propose to use the time average solution to compute the numerical flux, which leads to an upwind scheme for linear problems, and also increases the CFL numbers, which are comparable to other single step methods like ADER-DG scheme. An interesting observation we make is that the method at fifth order has a mild instability even though we use the CFL number determined from Fourier stability analysis. This mild instability seems to be present even in some RKDG schemes and in ADER-DG schemes. The central part of the numerical flux can be computed either by extrapolating it from the solution points to the faces or by directly estimating them at the faces by applying the approximate Lax-Wendroff procedure. These two methods perform differently for non-linear problems, with the extrapolation method leading to loss of convergence rate at odd polynomial degrees and also having larger errors compared to RK scheme. The alternate method proposed in this work performs uniformly well at all polynomial degrees and shows comparable accuracy to RK schemes. The LW method is also developed for hyperbolic systems like Euler equations, where many commonly used numerical fluxes based on approximate Riemann solvers like Roe, HLL, HLLC, can be used, along with the modifications that enhance the CFL number.  The method is described up to fifth order accuracy and it is cast in terms of matrix-vector operations. The new LW method is applied to linear and non-linear problems in one and two dimensions and representative results are shown to exhibit its accuracy and ability to compute smooth and discontinuous solutions.

The rest of the paper is organized as follows. The discretization of the domain and function approximation by polynomials are presented in Section~\ref{sec:scl}. Section~\ref{sec:rk} discusses the one dimensional FR method in the context of RK schemes and Section~\ref{sec:lw} introduces the LW method. The computation of the numerical flux and improvements over existing methods is presented in Section~\ref{sec:numflux}. The Fourier stability analysis in 1-D is performed in Section~\ref{sec:four1d} and the TVD limiter is discussed in Section~\ref{sec:lim}. Section~\ref{sec:res1d} and \ref{sec:res1dsys} present some numerical results in 1-D for scalar and system problems, to demonstrate the convergence rates and effect of correction functions, solution points and numerical flux schemes. Section~\ref{sec:2d} presents the LW scheme in two dimensions and Section~\ref{sec:res2d} presents some numerical results. Section~\ref{sec:sum} presents a summary of the new scheme. The Appendix discusses extension of some numerical fluxes for systems of PDEs to the LW case where time average fluxes and solutions are used, presents Fourier stability analysis in 2-D and the equivalence of two flux reconstruction approaches.

\section{Scalar conservation law}\label{sec:scl}
\begin{figure}
\begin{center}
\begin{tabular}{cc}
\includegraphics[width=0.45\textwidth]{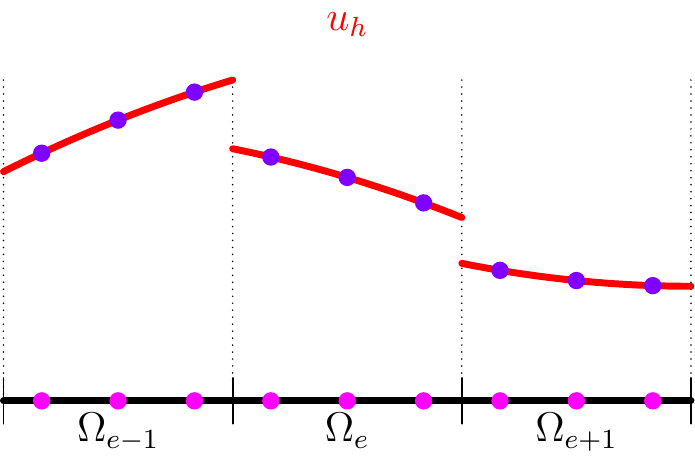} &
\includegraphics[width=0.45\textwidth]{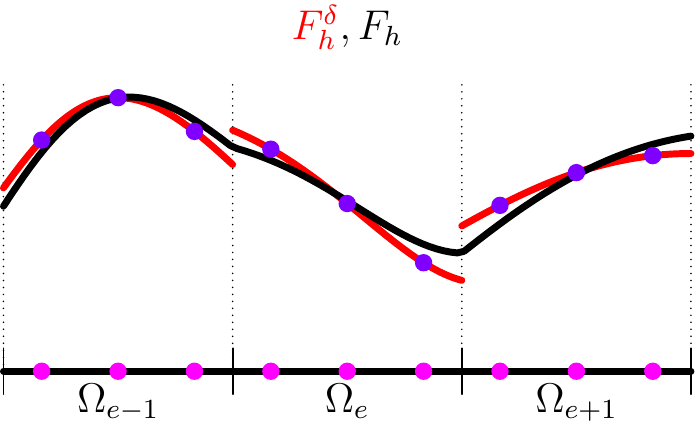} \\
(a) & (b)
\end{tabular}
\end{center}
\caption{(a) Piecewise polynomial solution at time $t_n$, and (b) discontinuous and continuous flux.}
\label{fig:solflux1}
\end{figure}

Let us consider a scalar conservation law of the form
\[
u_t + f(u)_x = 0
\]
where $u$ is some conserved quantity, $f(u)$ is the corresponding flux, together with some initial and boundary conditions. We will divide the computational domain $\Omega$ into disjoint elements $\Omega_e$, with
\[
\Omega_e = [x_\emh, x_\eph] \qquad\textrm{and}\qquad \Delta x_e = x_\eph - x_\emh
\]
Let us map each element to a reference element, $\Omega_e \to [0,1]$, by
\[
x \to \xi = \frac{x - x_\emh}{\Delta x_e}
\]
Inside each element, we approximate the solution by degree $N \ge 0$ polynomials belonging to the set $\poly_N$. For this, choose $N+1$ distinct nodes
\[
0 \le \xi_0 < \xi_1 < \cdots < \xi_N \le 1
\]
which will be taken to be Gauss-Legendre (GL) or Gauss-Lobatto-Legendre (GLL) nodes, and will also be referred to as {\em solution points}. There are associated quadrature weights $w_j$ such that the quadrature rule is exact for polynomials of degree up to $2N+1$ for GL points and upto degree $2N-1$ for GLL points. Note that the nodes and weights we use are with respect to the interval $[0,1]$ whereas they are usually defined for the interval $[-1,+1]$. The solution inside an element is given by
\[
x \in \Omega_e: \qquad u_h(\xi,t) = \sum_{j=0}^N u_j^e(t) \ell_j(\xi)
\]
where each $\ell_j$ is a Lagrange polynomial of degree $N$ given by
\[
\ell_j(\xi) = \prod_{i=0, i\ne j}^N \frac{\xi - \xi_i}{\xi_j - \xi_i} \in \poly_N, \qquad \ell_j(\xi_i) = \delta_{ij}
\]
Figure~(\ref{fig:solflux1}a) illustrates a piecewise polynomial solution at some time $t_n$ with discontinuities at the element boundaries. Note that the coefficients $u_j^e$ which are the basic unknowns or {\em degrees of freedom} (dof), are the solution values at the solution points.
The numerical method will require spatial derivatives of certain quantities.  We can compute the spatial derivatives on the reference interval using a differentiation matrix $\vD = [D_{ij}]$ whose entries are given by
\[
D_{ij} = \ell_j'(\xi_i), \qquad 0 \le i,j \le N
\]
For example, we can obtain the spatial derivatives of the solution $u_h$ at all the solution points by a matrix-vector product as follows
\[
\begin{bmatrix}
\partial_x u_h(\xi_0,t) \\
\vdots \\
\partial_x u_h(\xi_N,t)
\end{bmatrix} = \frac{1}{\Delta x_e} \vD \vu(t), \qquad
\vu = \begin{bmatrix}
u_0^e \\ \vdots \\ u_N^e \end{bmatrix}
\]
We will use symbols in sans serif font like $\vD,\vu$, etc. to denote matrices or vectors defined with respect to the solution points. The entries of the differentiation matrix are given by
\[
D_{ij} = \frac{W_j}{W_i} \frac{1}{(\xi_i - \xi_j)}, \qquad i \ne j \qquad\textrm{and}\qquad D_{ii} = -\sum_{j=0, j \ne i}^N D_{ij}
\]
where the $W_i$ are barycentric weights given by
\[
W_i = \frac{1}{\prod_{j=0,j\ne i}^N (\xi_i - \xi_j)}, \qquad 0 \le i \le N
\]
Define the Vandermonde matrices corresponding to the left and right boundaries of a cell by
\begin{equation}
\vV_L = [\ell_0(0), \ell_1(0), \ldots, \ell_N(0)]^\top, \qquad \vV_R = [\ell_0(1), \ell_1(1), \ldots, \ell_N(1)]^\top
\label{eq:VlVr}
\end{equation}
which is used to extrapolate the solution and/or flux to the cell faces for the computation of inter-cell fluxes.
\section{Runge-Kutta FR scheme}\label{sec:rk}
The RKFR scheme is based on a FR spatial discretization leading to a system of ODE followed by application of an RK scheme to march forward in time. The key idea is to construct a continuous polynomial approximation of the flux which is then used in a collocation scheme to update the nodal solution values. At some time $t$, we have the piecewise polynomial solution defined inside each cell; the FR scheme can be described by the following steps.

\paragraph{Step 1.} In each element, we construct the flux approximation by interpolating the flux at the solution points leading to a polynomial of degree $N$, given by
\[
f_h^\delta(\xi,t) = \sum_{j=0}^N f(u_j^e(t)) \ell_j(\xi)
\]
The above flux is in general discontinuous across the elements similar to the red curve in Figure~(\ref{fig:solflux1}b).

\paragraph{Step 2.} We build a continuous flux approximation by adding some correction terms at the element boundaries
\[
f_h(\xi,t) = \left[f_\emh(t) - f_h^\delta(0,t) \right] g_L(\xi) + f_h^\delta(\xi,t) + \left[f_\eph(t) - f_h^\delta(1,t) \right] g_R(\xi)
\]
where
\[
f_\eph(t) = f(u_h(x_\eph^-,t), u_h(x_\eph^+,t))
\]
is a numerical flux function that makes the flux unique across the cells. The continuous flux approximation is illustrated by the black curve in Figure~(\ref{fig:solflux1}b). The functions $g_L, g_R$ are the correction functions which must be chosen to obtain a stable scheme.

\paragraph{Step 3.} We obtain the system of ODE by collocating the PDE at the solution points
\[
\od{u_j^e}{t}(t) = - \frac{1}{\Delta x_e} \pd{f_h}{\xi} (\xi_j, t), \qquad 0 \le j \le N
\]
which is solved in time by a Runge-Kutta scheme.

\paragraph{Correction functions.} The correction functions $g_L,g_R$ should satisfy the end point conditions
\begin{align*}
g_L(0) = 1, \qquad & g_R(0) = 0 \\
g_L(1) = 0, \qquad & g_R(1) = 1
\end{align*}
which ensures the continuity of the flux, i.e., $f_h(x_\eph^-,t) = f_h(x_\eph^+,t) = f_\eph(t)$. Moreover, we want them to be close to zero inside the element. There is a wide family of correction functions available in the literature~\cite{Huynh2007,Vincent2011a}. A family of correction functions depending on a parameter $c$ were developed in~\cite{Vincent2011a} based on stability in a Sobolev-type norm. Two of these functions, the Radau and g2 correction functions, are of major interest since they correspond to commonly used DG formulations. The Radau correction function is a polynomial of degree $N+1$ which belongs to the family of~\cite{Vincent2011a} corresponding to the parameter $c=0$ and given by
\[
g_L(\xi) = \frac{(-1)^N}{2}[L_N(2\xi-1) - L_{N+1}(2\xi-1)], \qquad g_R(\xi) = \frac{1}{2}[L_N(2\xi-1) + L_{N+1}(2\xi-1)]
\]
where $L_N : [-1,1] \to \re$ is the Legendre polynomial of degree $N$. The resulting RKFR scheme can be shown to be identical to the nodal RKDG scheme using Gauss-Legendre nodes for solution points and quadrature. In the general class of~\cite{Vincent2011a}, $g_2$ correction function of degree $N+1$ corresponds to $c = \frac{2(N+1)}{(2N+1)N(a_N N!)^2}$ where $a_N$ is the leading coefficient of $L_N$; they are given by
\begin{align*}
g_L(\xi) &= \frac{(-1)^N}{2}\left[L_N(2\xi-1) - \frac{(N+1)L_{N-1}(2\xi-1)+NL_{N+1}(2\xi-1)}{2N+1}  \right] \\
g_R(\xi) &= \frac{1}{2}\left[L_N(2\xi-1) + \frac{(N+1)L_{N-1}(2\xi-1)+NL_{N+1}(2\xi-1)}{2N+1} \right]
\end{align*}
The resulting RKFR scheme can be shown to be identical to the nodal RKDG scheme using Gauss-Lobatto-Legendre points as solution points and for quadrature. We will perform Fourier stability analysis of the Lax-Wendroff scheme based on these correction functions in a later section. Note that the correction functions are usually defined in the interval $[-1,1]$ but here we have written them for our reference interval which is $[0,1]$.
\section{Lax-Wendroff FR scheme}\label{sec:lw}
The Lax-Wendroff scheme combines the spatial and temporal discretization into a single step. The starting point is a Taylor expansion in time following  the Cauchy-Kowalewski procedure where the PDE is used to rewrite some of the time derivatives in the Taylor expansion as spatial derivatives. Using Taylor expansion in time around $t = t_n$, we can write the solution at the next time level as
\[
u^{n+1} = u^n + \sum_{m=1}^{N+1} \frac{\Delta t^m}{m!} \partial_t^m u^n + O(\Delta t^{N+2})
\]
Since the spatial error is expected to be of $O(\Delta x^{N+1})$, we retain terms up to $O(\Delta t^{N+1})$ in the Taylor expansion, so that the over all accuracy is of order $N+1$ both space and time.  Using the PDE, $\partial_t u = - \partial_x f$, we re-write time derivatives of the solution in terms of spatial derivatives of the flux
\[
\partial_t^m u = - \partial_t^{m-1} \partial_x f = - (\partial_t^{m-1} f)_x, \qquad m=1,2,\ldots
\]
so that
\begin{eqnarray}
\nonumber
u^{n+1} &=& u^n - \sum_{m=1}^{N+1} \frac{\Delta t^m}{m!} (\partial_t^{m-1} f)_x + O(\Delta t^{N+2}) \\
\nonumber
&=& u^n - \Delta t \left[ \sum_{m=0}^N \frac{\Delta t^m}{(m+1)!} \partial_t^m f \right]_x + O(\Delta t^{N+2}) \\
\label{eq:lwtay}
&=& u^n - \Delta t \pd{F}{x}(u^n) + O(\Delta t^{N+2})
\end{eqnarray}
where
\begin{equation} \label{eq:tavgflux}
F(u) = \sum_{m=0}^N \frac{\Delta t^m}{(m+1)!} \partial_t^m f(u) = f(u) + \frac{\Delta t}{2} \partial_t f(u) + \ldots + \frac{\Delta t^N}{(N+1)!} \partial_t^N f(u)
\end{equation}
Note that $F(u^n)$ is an approximation to the time average flux in the interval $[t_n, t_{n+1}]$ since it can be written as
\begin{equation}
F(u^n) = \frac{1}{\Delta t} \int_{t_n}^{t_{n+1}} \left[ f(u^n) + (t-t_n) \partial_t f(u^n) + \ldots + \frac{(t - t_n)^N}{N!} \partial_t^N f(u^n) \right] \ud t \label{eq:tvgproperty}
\end{equation}
where the quantity inside the square brackets is the truncated Taylor expansion of the flux $f$ in time.  Equation~\eqref{eq:lwtay} is the basis for the construction of the Lax-Wendroff method. Following the ideas in the RKFR scheme, we will first reconstruct the time average flux $F$ inside each element by a continuous polynomial $F_h(\xi)$. Then truncating equation~\eqref{eq:lwtay}, the solution at the nodes is updated by a collocation scheme as follows
\begin{equation}
(u_j^e)^{n+1} = (u_j^e)^n - \frac{\Delta t}{\Delta x_e} \od{F_h}{\xi}(\xi_j), \qquad 0 \le j \le N
\label{eq:uplwfr}
\end{equation}
This is the single step Lax-Wendroff update scheme for any order of accuracy. The major work in the above scheme is involved in the construction of the time average flux approximation $F_h$ which is explained in subsequent sections.
\subsection{Conservation property}
The computation of correct weak solutions for non-linear conservation laws in the presence of discontinuous solutions requires the use of conservative numerical schemes.  The Lax-Wendroff theorem shows that if a consistent, conservative method converges, then the limit is a weak solution. The method~\eqref{eq:uplwfr} is also conservative though it is not directly apparent; to see this multiply~\eqref{eq:uplwfr} by the quadrature weights associated with the solution points and sum over all the points in the $e$'th element,
\[
\sum_{j=0}^N w_j (u_j^e)^{n+1} = \sum_{j=0}^N w_j (u_j^e)^n - \frac{\Delta t}{\Delta x_e} \sum_{j=0}^N w_j \pd{F_h}{\xi}(\xi_j)
\]
The correction functions are of degree $N+1$ and the flux $F_h$ is a polynomial of degree $\le N+1$. If the quadrature is exact for polynomials of degree atleast $N$, which is true for both GLL and GL points, then the quadrature is exact for the flux derivative term and we can write it as an integral, which leads to
\begin{equation}
\int_{\Omega_e} u_h^{n+1} \ud x = \int_{\Omega_e} u_h^{n} \ud x - \Delta t [F_\eph - F_\emh]
\label{eq:upmean}
\end{equation}
This shows that the total mass inside the cell changes only due to the boundary fluxes and the scheme is hence conservative.
\subsection{Reconstruction of the time average flux}\label{sec:reconstruction}
To complete the description of the LW method~\eqref{eq:uplwfr}, we must explain the method for the computation of the time average flux $F_h$. The flux reconstruction $F_h(\xi)$ for a time interval $[t_n, t_{n+1}]$ is performed in three steps.

\paragraph{Step 1.} Use the approximate Lax-Wendroff procedure to compute the time average flux $F$ at all the solution points
\[
F_j^e \approx F(\xi_j), \qquad 0 \le j \le N
\]
The approximate LW procedure is explained in a subsequent section.

\paragraph{Step 2.} Build a local approximation of the time average flux inside each element by interpolating at the solution points
\[
F_h^\delta(\xi) = \sum_{j=0}^N F_j^e \ell_j(\xi)
\]
which however may not be continuous across the elements. This is illustrated in Figure~(\ref{fig:solflux1}b).

\paragraph{Step 3.} Modify the flux approximation $F_h^\delta(\xi)$ so that it becomes continuous across the elements. Let $F_\eph$ be some numerical flux function that approximates the flux $F$ at $x=x_\eph$. Then the continuous flux approximation is given by
\[
F_h(\xi) = \left[F_\emh - F_h^\delta(0) \right] g_L(\xi) + F_h^\delta(\xi) + \left[F_\eph - F_h^\delta(1) \right] g_R(\xi)
\]
which is illustrated in Figure~(\ref{fig:solflux1}b). The correction functions $g_L,g_R$ are chosen from the FR literature.

\paragraph{Step 4.}
The derivatives of the continuous flux approximation at the solution points can be obtained as
\[
\myvector{\partial_\xi F_h} = \left[F_\emh - \vV_L^\top \vF \right] \vb_L +
\vD \vF +
\left[F_\eph - \vV_R^\top \vF \right] \vb_R, \qquad \vb_L = \begin{bmatrix}
g_L'(\xi_0) \\
\vdots \\
g_L'(\xi_N) \end{bmatrix}, \qquad \vb_R = \begin{bmatrix}
g_R'(\xi_0) \\
\vdots \\
g_R'(\xi_N) \end{bmatrix}
\]
which can also be written as
\begin{equation}\label{eq:fder}
\myvector{\partial_\xi F_h} = F_\emh \vb_L + \vD_1 \vF +
F_\eph \vb_R, \qquad \vD_1 =  \vD - \vb_L \vV_L^\top - \vb_R \vV_R^\top
\end{equation}
where $\vV_L, \vV_R$ are Vandermonde matrices which are defined in~\eqref{eq:VlVr}.  The quantities $\vb_L, \vb_R, \vV_L, \vV_R, \vD, \vD_1$ can be computed once and re-used in all subsequent computations. They do not depend on the element and are computed on the reference element. Equation~\eqref{eq:fder} contains terms that can be computed inside a single cell (middle term) and those computed at the faces (first and third terms) where it is required to use the data from two adjacent cells. The computation of the flux derivatives can thus be performed by looping over cells and then the faces.
\subsection{Direct flux reconstruction (DFR) scheme}\label{sec:DFR}
An alternate approach to flux reconstruction which does not require the choice of a correction function is based on the idea of direct flux reconstruction~\cite{Romero2016}, which we adopt in the Lax-Wendroff scheme as follows. Let us take the solution points to be the $N+1$ Gauss-Legendre nodes, and define
\[
\xi_{-1} = 0, \qquad \xi_{N+1} = 1
\]
The Lagrange polynomials corresponding to the $N+3$ points $\{\xi_i, i=-1,0,\ldots,N+1\}$ are given by
\[
\tilde\ell_j(\xi) = \prod_{i=-1,i\ne j}^{N+1} \frac{\xi - \xi_i}{\xi_j - \xi_i} \in \poly_{N+2}
\]
We approximate the continuous flux in terms of these polynomials
\[
F_h(\xi) = F_\emh \tilde\ell_{-1}(\xi) + \sum_{j=0}^N F_j^e \tilde\ell_j(\xi) + F_\eph \tilde\ell_{N+1}(\xi)
\]
We can compute the spatial derivatives using a differentiation matrix $\tilde \vD \in \re^{(N+1) \times (N+3)}$
\[
\tilde D_{ij} = \tilde\ell_j'(\xi_i), \qquad 0 \le i \le N, \quad -1 \le j \le N+1
\]
Define $\vb_L, \vb_R$ to be the first and last column of the matrix $\tilde \vD$ and $\vD_1$ to be the remaining columns
\[
\vb_L = \tilde \vD{\tt(:,-1)}, \qquad \vb_R = \tilde \vD{\tt(:,N+1)}, \qquad \vD_1 = \tilde \vD{\tt(0:N, 0:N)}
\]
The flux derivatives at all the solution points can be computed as follows
\[
\myvector{\partial_\xi F_h} = F_\emh \vb_L + \vD_1 \vF + F_\eph \vb_R
\]
Note that the above equation has the same structure as~\eqref{eq:fder} from the FR procedure but $\vb_L,\vb_R,\vD_1$ are obtained using different idea. In this DFR approach, we cannot use GLL points since then the boundary points $\xi_{-1} = \xi_0$, $\xi_{N+1} = \xi_N$ would be repeated and the Lagrange interpolation is not well-defined; if we use GL points, the resulting scheme is identical to the LWFR approach using Radau correction  function in combination with GL points as solution points, as shown in Appendix~\ref{sec:frdfr}.
\subsection{Approximate Lax-Wendroff procedure}\label{sec:alw}
The time average flux at the solution points $F_j^e$ must be computed using~\eqref{eq:tavgflux}. The usual approach is to use the PDE and replace time derivatives with spatial derivatives, but this leads to large amount of algebraic computations since we need to evaluate the flux Jacobian and its higher tensor versions. To avoid this process, we follow the ideas in~\cite{Zorio2017,Burger2017} and adopt an approximate Lax-Wendroff procedure. To present this idea in a concise and efficient form, we introduce the notation
\[
u^{(m)} = \Delta t^m \partial_t^m u, \qquad f^{(m)} = \Delta t^m \partial_t^m f, \qquad m=1,2,\ldots
\]
The time derivatives of the solution are computed using the PDE
\[
u^{(m)} = - \Delta t \ \partial_x f^{(m-1)}, \qquad m=1,2,\ldots
\]
The approximate Lax-Wendroff procedure uses a finite difference approximation applied at the solution points to compute the time derivatives of the fluxes. For example, a second order approximation is given by
\[
f_t(\xi,t) \approx \frac{f(u(\xi,t+\Delta t)) - f(u(\xi,t-\Delta t))}{2\Delta t}
\]
The arguments to the flux are in turn approximated by a Taylor expansion in time
\[
u(\xi,t\pm\Delta t) \approx u(\xi,t) \pm u_t(\xi,t) \Delta t
\]
Using this approximation at the $j$'th solution point in an element, we get
\[
f_j^{(1)} = f_t(\xi_j,t) \Delta t \approx  \half [ f(u_j + u_j^{(1)}) - f(u_j - u_j^{(1)})], \qquad u_j^{(1)} = u_t(\xi_j,t) \Delta t = - \frac{\Delta t}{\Delta x_e} f_\xi(\xi_j,t) \approx - \frac{\Delta t}{\Delta x_e} (\vD \vf)_j
\]
It can be shown that the above approximation to $f_t$ is second order accurate in $\Delta t$. Such approximations can be written for higher accuracy and for higher time derivatives~\cite{Zorio2017,Burger2017}, and we summarize them below at different orders of accuracy which are used in the paper. The neglected term in the Taylor expansion~\eqref{eq:tavgflux}  is of $O(\Delta t^{N+1})$, and hence the derivative approximation $\partial_t^m f$ must be computed to at least $O(\Delta t^{N+1-m})$ accuracy. The Lax-Wendroff procedure is applied in each element and so for simplicity of notation, we do not show the element index in the following sub-sections.
\subsubsection{Second order scheme, $N=1$}
The time average flux at the solution points is given by
\[
\vF = \vf + \frac{1}{2} \vf^{(1)}
\]
where
\begin{eqnarray*}
\vu^{(1)} &=& - \frac{\Delta t}{\Delta x_e} \vD \vf \\
\vf^{(1)} &=& \frac{1}{2}\left[ f\left(\vu +  \vu^{(1)}\right) - f\left(\vu - \vu^{(1)}\right) \right]
\end{eqnarray*}
\subsubsection{Third order scheme, $N=2$}
The time average flux at the solution points is given by
\[
\vF = \vf + \frac{1}{2} \vf^{(1)} + \frac{1}{6} \vf^{(2)}
\]
where
\begin{eqnarray*}
\vu^{(1)} &=& - \frac{\Delta t}{\Delta x_e} \vD \vf \\
\vf^{(1)} &=& \frac{1}{2} \left[ f\left(\vu + \vu^{(1)}\right) - f\left(\vu -  \vu^{(1)}\right) \right] \\
\vu^{(2)} &=& - \frac{\Delta t}{\Delta x_e} \vD \vf^{(1)} \\
\vf^{(2)} &=& f\left(\vu + \vu^{(1)} + \half  \vu^{(2)}\right) - 2 f(\vu) + f\left(\vu -  \vu^{(1)} + \half  \vu^{(2)} \right)
\end{eqnarray*}
\subsubsection{Fourth order scheme, $N=3$}
For the fourth order scheme, the time average flux  at the solution points reads as
\[
\vF = \vf + \frac{1}{2} \vf^{(1)} + \frac{1}{6} \vf^{(2)} + \frac{1}{24} \vf^{(3)}
\]
where
\begin{eqnarray*}
\vu^{(1)} &=& - \frac{\Delta t}{\Delta x_e} \vD \vf \\
\vf^{(1)}&=&\frac{1}{12} \left[  - f \left(\vu + 2 \vu^{(1)}\right) +8f \left(\vu + \vu^{(1)} \right)
 -8f \left(\vu - \vu^{(1)}\right) + f \left(\vu - 2 \vu^{(1)}\right) \right] \\
\vu^{(2)} &=& - \frac{\Delta t}{ \Delta x_e} \vD \vf^{(1)} \\
\vf^{(2)} &=&  f\left(\vu + \vu^{(1)} + \half \vu^{(2)}\right) - 2 f(\vu) + f\left(\vu - \vu^{(1)} + \half  \vu^{(2)} \right) \\
\vu^{(3)} &=& - \frac{\Delta t}{\Delta x_e} \vD \vf^{(2)} \\
\vf^{(3)}&=& \frac{1}{2} \left[f \left(\vu + 2 \vu^{(1)} + \frac{2^2}{2!}\vu^{(2)} + \frac{2^3}{3!}\vu^{(3)}
 \right)
 -2 f \left(\vu + \vu^{(1)} + \frac{1}{2!}\vu^{(2)} + \frac{1}{3!}\vu^{(3)} \right) \right. \\
 && \quad \left. + 2 f \left( \vu - \vu^{(1)} + \frac{1}{2!}\vu^{(2)} - \frac{1}{3!}\vu^{(3)} \right)
 - f \left(\vu - 2 \vu^{(1)} + \frac{2^2}{2!}\vu^{(2)} - \frac{2^3}{3!} \vu^{(3)} \right) \right]
\end{eqnarray*}
\subsubsection{Fifth order scheme, $N=4$}
The time average flux at the solution points  for the fifth order scheme  takes the form
\[
\vF = \vf + \frac{1}{2} \vf^{(1)} + \frac{1}{6} \vf^{(2)} + \frac{1}{24} \vf^{(3)} + \frac{1}{120} \vf^{(4)}
\]
where

\begin{eqnarray*}
\vu^{(1)} &=& - \frac{\Delta t}{\Delta x_e} \vD \vf \\
\vf^{(1)} &=& \frac{1}{12} \left[ -f \left(\vu + 2 \vu^{(1)} \right) + 8 f \left(\vu + \vu^{(1)} \right)
 - 8 f \left( \vu - \vu^{(1)} \right) + f \left( \vu - 2 \vu^{(1)} \right) \right] \\
\vu^{(2)} &=& - \frac{\Delta t}{\Delta x_e} \vD \vf^{(1)} \\
\vf^{(2)} &=&
  \frac{1}{12}\left[-f\left(\vu + 2 \vu^{(1)} + \frac{2^2}{2!} \vu^{(2)} \right)
  + 16 f \left(\vu + \vu^{(1)} + \frac{1}{2!} \vu^{(2)} \right)
  - 30 f(\vu) \right.\\
  && \qquad \left.
  +16 f\left(\vu - \vu^{(1)} + \frac{1}{2!} \vu^{(2)} \right)
  -f \left(\vu - 2 \vu^{(1)} + \frac{2^2}{2!} \vu^{(2)} \right)
  \right]\\
\vu^{(3)} &=& - \frac{\Delta t}{\Delta x_e} \vD \vf^{(2)} \\
\vf^{(3)}&=& \frac{1}{2} \left[ f \left(\vu + 2 \vu^{(1)} + \frac{2^2}{2!} \vu^{(2)} + \frac{2^3}{3!} \vu^{(3)} \right)
 - 2 f \left(\vu + \vu^{(1)} + \frac{1}{2!} \vu^{(2)} + \frac{1}{3!} \vu^{(3)} \right) \right. \\
 && \quad \left. + 2 f \left( \vu - \vu^{(1)} + \frac{1}{2!} \vu^{(2)} - \frac{1}{3!} \vu^{(3)} \right)
 - f \left( \vu - 2 \vu^{(1)} + \frac{2^2}{2!} \vu^{(2)} - \frac{2^3}{3!} \vu^{(3)} \right) \right]\\
\vu^{(4)} &=& - \frac{\Delta t}{\Delta x_e} \vD \vf^{(3)} \\
\vf^{(4)} &=& \left[  f \left(\vu + 2 \vu^{(1)}
+ \frac{2^2}{2!} \vu^{(2)} + \frac{2^3}{3!} \vu^{(3)} + \frac{2^4}{4!} \vu^{(4)} \right)\right.\\
&& \left. - 4 f \left( \vu + \vu^{(1)}
+ \frac{1}{2!} \vu^{(2)} + \frac{1}{3!} \vu^{(3)} + \frac{1}{4!} \vu^{(4)} \right)
+ 6 f(\vu) \right.\\
&& \left. - 4 f \left(\vu - \vu^{(1)} + \frac{1}{2!} \vu^{(2)} - \frac{1}{3!} \vu^{(3)} + \frac{1}{4!} \vu^{(4)} \right) \right.\\
&& + \left.
f \left(\vu - 2 \vu^{(1)}
+ \frac{2^2}{2!} \vu^{(2)} - \frac{2^3}{3!} \vu^{(3)} + \frac{2^4}{4!} \vu^{(4)} \right)
\right]
\end{eqnarray*}
The above set of formulae shows the sequence of steps that have to be performed to compute the time average flux at various orders. The arguments of the fluxes used on the right hand side in these steps are built in a sequential manner. Note that all the equations are vectorial equations and are applied at each solution point.
\section{Numerical flux}\label{sec:numflux}
The numerical flux couples the solution between two neighbouring cells in a discontinuous Galerkin type method. In RK methods, the numerical flux is a function of the trace values of the solution at the faces. In the Lax-Wendroff scheme, we have constructed the time average flux at all the solution points inside the element and we want to use this information to compute the time averaged numerical flux at the element faces. The simplest numerical flux is based on Lax-Friedrich type approximation and is given by~\cite{Qiu2005b}
\begin{equation} \label{eq:nfdiss1}
F_\eph = \half[ F_\eph^- + F_\eph^+] - \half \lambda_\eph [ u_h(x_\eph^+, t_n) - u_h(x_\eph^-, t_n)]
\tag{D1}
\end{equation}
which consists of a central flux and a dissipative part. For linear advection equation $u_t + a u_x = 0$, the coefficient in the dissipative part of the flux is taken as $\lambda_\eph = |a|$, while for a non-linear PDE like Burger's equation, we take it to be
\[
\lambda_\eph = \max\{ |f'(\au_e^n)|, |f'(\au_{e+1}^n)| \}
\]
where $\au_e^n$ is the cell average solution in element $\Omega_e$ at time $t_n$, and will be referred to as Rusanov or local Lax-Friedrich~\cite{Rusanov1962} approximation. Note that the dissipation term in the above numerical flux is evaluated at time $t_n$ whereas the central part of the flux uses the time average flux. Since the dissipation term contains solution difference at faces, we still expect to obtain optimal convergence rates, which is verified in numerical experiments. This numerical flux depends on the following quantities: $\{ \au_e^n, \au_{e+1}^n, u_h(x_\eph^-,t_n), u_h(x_\eph^+,t_n), F_\eph^-, F_\eph^+ \}$.

The numerical flux of the form~\eqref{eq:nfdiss1} leads to somewhat reduced CFL numbers as shown by Fourier stability analysis in a later section, and also does not have upwind property even for linear advection equation. An alternate form of the numerical flux is obtained by evaluating the dissipation term using the time average solution, leading to the formula
\begin{equation} \label{eq:nfdiss2}
F_\eph = \half[ F_\eph^- + F_\eph^+] - \half \lambda_\eph [ U_\eph^+ - U_\eph^-]
\tag{D2}
\end{equation}
where
\begin{equation} \label{eq:tavgsol}
U = \sum_{m=0}^N \frac{\Delta t^m}{(m+1)!} \partial_t^m u = u + \frac{\Delta t}{2} \partial_t u + \ldots + \frac{\Delta t^N}{(N+1)!} \partial_t^N u
\end{equation}
is the time average solution. In this case, the numerical flux depends on the following quantities: $\{ \au_e^n, \au_{e+1}^n, U_\eph^-, U_\eph^+, F_\eph^-, F_\eph^+ \}$. We will refer to the above two forms of dissipation as D1 and D2, respectively. The dissipation model D2 is not computationally expensive compared to the D1 model since all the quantities required to compute the time average solution $U$ are available during the Lax-Wendroff procedure. Numerical fluxes for the case of systems of hyperbolic equations are described in the Appendix. It remains to explain how to compute $F_\eph^\pm$ appearing in the central part of the numerical flux, which can be accomplished in two different ways, which we term AE and EA in the next two sub-sections.
\begin{remark}
In case of constant linear advection equation, $u_t + a u_x = 0$, $f^{(m)} = a u^{(m)}$ so that $F_j^e = a U_j^e$. Then, since $\lambda_\eph=|a|$, the numerical flux~\eqref{eq:nfdiss2} becomes the upwind flux
\[
F_\eph = \begin{cases}
F_h^\delta(x_\eph^-), & a \ge 0 \\
F_h^\delta(x_\eph^+), & a < 0
\end{cases}
\]
but the flux~\eqref{eq:nfdiss1} does not have this upwind property. For a variable coefficient advection problem with flux $f = a(x)u$, we get $F_j^e = a(x_j) U_j^e$, the numerical flux~\eqref{eq:nfdiss2} is
\begin{equation}
F_\eph = \half[ F_\eph^- + F_\eph^+] - \half |a(x_\eph)| [ U_\eph^+ - U_\eph^-]
\label{eq:nflin1}
\end{equation}
which does not reduce to an upwind flux due to interpolation errors, though it will be close to it in the well resolved cases. In this case, we can define the upwind numerical flux as
\begin{equation}
F_\eph = \begin{cases}
F_\eph^-, & a(x_\eph) \ge 0 \\
F_\eph^+, & a(x_\eph) < 0
\end{cases}
\label{eq:nflin2}
\end{equation}
which is defined in terms of the time average flux only and does not make use of the solution.
\end{remark}

\begin{remark}
For non-linear problems, we can also consider the global Lax-Friedrich and Roe type dissipation models which are given by
\[
\lambda_\eph = \lambda = \max_e |f'(\au_e)|, \qquad \lambda_\eph = \left| f'\left( \frac{\au_e + \au_{e+1}}{2} \right) \right|
\]
respectively. In the global Lax-Friedrich flux, the maximum is taken over the whole grid. For Burger's equation, we can consider an Osher type flux~\cite{Engquist1981} which is given by
\[
F_\eph = \begin{cases}
F_\eph^- & \au_e, \au_{e+1} > 0 \\
F_\eph^+ & \au_e, \au_{e+1} < 0 \\
F_\eph^- + F_\eph^+ & \au_e \ge 0 \ge \au_{e+1}\\
0 & \textrm{otherwise}
\end{cases}
\]
\end{remark}
\subsection{Numerical flux -- average and extrapolate to face (AE)}
In each element, the time average flux $F_h^\delta$ has been constructed using the Lax-Wendroff procedure. The simplest approximation that can be used for $F_\eph^\pm$ in the central part of the numerical flux is to extrapolate the flux $F_h^\delta$ to the faces,
\[
F_\eph^\pm = F_h^\delta(x_\eph^\pm)
\]
We will refer to this approach with the abbreviation {\bf AE}. However, as shown in the numerical results, this approximation can lead to sub-optimal convergence rates for some non-linear problems. Hence we propose another method for the computation of the inter-cell flux which overcomes this problem as explained next.
\subsection{Numerical flux -- extrapolate to face and average (EA)}
Instead of extrapolating the time average flux from the solution points to the faces, we can instead build the time average flux at the faces directly using the approximate Lax-Wendroff procedure that is used at the solution points. The flux at the faces is constructed after the solution is evolved at all the solution points. In the following equations, $\alpha$ denotes either the left face ($L$) or the right face ($R$) of a cell. For $\alpha \in \{L,R\}$, we compute the time average flux at the faces of the $e$'th element by the following steps, where we suppress the element index since all the operations are performed inside one element.
\begin{equation*}
\begin{aligned}[t]
\textrm{\bf Degree $N=1$} \quad\quad\quad & \\
u_\alpha &= \vV_\alpha^\top \vu \\
u_\alpha^{\pm} &= \vV_\alpha^\top (\vu \pm \vu^{(1)}) \\
f^{(1)}_\alpha &= \half [ f(u_\alpha^+) - f(u_\alpha^-) ] \\
F_\alpha &= f(u_\alpha) + \half f_\alpha^{(1)}
\end{aligned}
\qquad\qquad\qquad\qquad
\begin{aligned}[t]
\textrm{\bf Degree $N=2$} \quad\quad\quad & \\
u_\alpha &= \vV_\alpha^\top \vu \\
u_\alpha^{\pm} &= \vV_\alpha^\top \left( \vu \pm \vu^{(1)} + \half  \vu^{(2)} \right) \\
f^{(1)}_\alpha &= \frac{1}{2} \left[ f(u_\alpha^+) - f(u_\alpha^-) \right] \\
f^{(2)}_\alpha &= f(u_\alpha^-) - 2 f(u_\alpha) + f(u_\alpha^+) \\
F_\alpha &= f(u_\alpha) + \frac{1}{2} f^{(1)}_\alpha + \frac{1}{6} f_\alpha^{(2)}
\end{aligned}
\end{equation*}
\paragraph{Degree $N=3$}
\begin{eqnarray*}
u_\alpha &=& \vV_\alpha^\top \vu \\
u_\alpha^{\pm} &=& \vV_\alpha^\top \left( \vu \pm \vu^{(1)} + \frac{1}{2!}\vu^{(2)} \pm \frac{1}{3!}\vu^{(3)} \right) \\
u_\alpha^{\pm2} &=& \vV_\alpha^\top \left( \vu \pm 2 \vu^{(1)} + \frac{2^2}{2!}\vu^{(2)} \pm \frac{2^3}{3!} \vu^{(3)} \right) \\
f^{(1)}_\alpha &=& \frac{1}{12} \left[  - f \left( u^{+2}_\alpha \right) + 8 f \left( u^+_\alpha \right) - 8 f \left( u^-_\alpha \right) + f \left( u^{-2}_\alpha \right) \right] \\
f^{(2)}_\alpha &=& f(u_\alpha^-) - 2 f(u_\alpha) + f(u_\alpha^+) \\
f^{(3)}_\alpha &=& \frac{1}{2} \left[f \left( u^{+2}_\alpha \right)
 -2 f \left( u^+_\alpha \right) + 2 f \left( u^-_\alpha \right)
 - f \left( u^{-2}_\alpha \right) \right] \\
F_\alpha &=& f(u_\alpha) + \frac{1}{2} f^{(1)}_\alpha + \frac{1}{6} f_\alpha^{(2)} + \frac{1}{24} f_\alpha^{(3)}
\end{eqnarray*}
\paragraph{Degree $N=4$}
\begin{eqnarray*}
u_\alpha &=& \vV_\alpha^\top \vu \\
u_\alpha^{\pm} &=& \vV_\alpha^\top \left( \vu \pm \vu^{(1)} + \frac{1}{2!}\vu^{(2)} \pm \frac{1}{3!}\vu^{(3)}+ \frac{1}{4!}\vu^{(4)} \right) \\
u_\alpha^{\pm2} &=& \vV_\alpha^\top \left( \vu \pm 2 \vu^{(1)} + \frac{2^2}{2!}\vu^{(2)} \pm \frac{2^3}{3!} \vu^{(3)}+ \frac{2^3}{3!} \vu^{(3)} \right) \\
f^{(1)}_\alpha &=& \frac{1}{12} \left[  - f \left( u^{+2}_\alpha \right) + 8 f \left( u^+_\alpha \right) - 8 f \left( u^-_\alpha \right) + f \left( u^{-2}_\alpha \right) \right] \\
f^{(2)}_\alpha &=& \frac{1}{12} \left[ -f \left(u^{+2}_\alpha \right) +16 f\left(u^+_\alpha \right)
-30 f\left(u_\alpha \right)+16 f\left( u^-_\alpha \right) - f \left( u^{-2}_\alpha \right)
\right]\\
f^{(3)}_\alpha &=& \frac{1}{2} \left[f \left( u^{+2}_\alpha \right)
 -2 f \left( u^+_\alpha \right) + 2 f \left( u^-_\alpha \right)
 - f \left( u^{-2}_\alpha \right) \right] \\
 f^{(4)}_\alpha &=&  \left[f \left( u^{+2}_\alpha \right)
 -4 f \left( u^+_\alpha \right) +6 f \left ( u_\alpha \right)- 4 f \left( u^-_\alpha \right)
 + f \left( u^{-2}_\alpha \right) \right] \\
F_\alpha &=& f(u_\alpha) + \frac{1}{2} f^{(1)}_\alpha + \frac{1}{6} f_\alpha^{(2)} + \frac{1}{24} f_\alpha^{(3)}+ \frac{1}{120} f_\alpha^{(4)}
\end{eqnarray*}
We see that the solution is first extrapolated to the cell faces and the same finite difference formulae for the time derivatives of the flux which are used at the solution points, are also used at the faces.  The numerical flux is computed using the time average flux built as above at the faces; the central part of the flux $F_\eph^\pm$ in equations~\eqref{eq:nfdiss1},~\eqref{eq:nfdiss2} are computed as
\[
F_\eph^- = (F_R)_e, \qquad F_\eph^+ = (F_L)_{e+1}
\]
We will refer to this method with the abbreviation {\bf EA}. The dissipative part of the flux is computed using either the solution at time $t_n$ or the time average solution, which are extrapolated to the faces, leading to the D1 and D2 models, respectively. 
\begin{remark}
The two methods AE and EA are different only when there are no solution points at the faces. E.g., if we use GLL solution points, then the two methods yield the same result since there is no interpolation error. For the constant coefficient advection equation, the above two schemes for the numerical flux lead to the same approximation but they differ in case of variable coefficient advection problems and when the flux is non-linear with respect to $u$. The effect of this non-linearity and the performance of the two methods are shown later using some numerical experiments.
\end{remark}
\section{Fourier stability analysis in 1-D}\label{sec:four1d}
We now perform Fourier stability analysis of the LW schemes applied to the linear advection equation $u_t + a u_x = 0$ where $a$ is the constant advection speed. We will assume that the advection speed $a$ is positive and denote the CFL number by
\[
\sigma = \frac{a \Delta t}{\Delta x}
\]
Since $f^{(m)} = a u^{(m)}$, the time average flux at all the solution points is given by
\[
\vF_e = a \vU_e \quad  \textrm{where} \quad
\vU_e = \vT \vu_e  \quad\textrm{and}\quad \vT = \sum_{m=0}^N \frac{1}{(m+1)!} (-\sigma \vD)^m
\]
Then the discontinuous flux at the cell boundaries are given by
\[
F_h^\delta(x_\emh^+) = \vV_L^\top \vF_e , \qquad F_h^\delta(x_\eph^-) = \vV_R^\top \vF_e
\]
We can write the update equation in the form
\begin{equation}\label{eq:laup}
\vu_e^{n+1}= -\sigma \vA_{-1} \vu_{e-1}^n + (\vI - \sigma \vA_0) \vu_e^n - \sigma \vA_{+1} \vu_{e+1}^n
\end{equation}
where the matrices $\vA_{-1}, \vA_0, \vA_{+1}$ depend on the choice of the dissipation model in the numerical flux. The EA and AE schemes for the flux are identical for this linear problem, and hence we do not make any distinction between them for Fourier stability analysis.
\paragraph{Dissipation model D1.}
The numerical flux is given by
\[
F_\eph = \half [ \vV_R^\top \vF_e + \vV_L^\top \vF_{e+1}] - \half a (\vV_L^\top \vu_{e+1} - \vV_R^\top \vu_e)
\]
Since the flux difference at the faces is
\[
F_\emh - F_h^\delta(x_\emh^+) = \half a \vV_R^\top (\vT + \vI) \vu_{e-1} - \half a \vV_L^\top (\vT + \vI) \vu_e
\]
\[
F_\eph - F_h^\delta(x_\eph^-) = \half a \vV_L^\top(\vT - \vI)\vu_{e+1} - \half a \vV_R^\top (\vT - \vI) \vu_e
\]
the flux derivative at the solution points is given by
\[
\partial_\xi \vF_h  = \half a \vb_L \vV_R^\top (\vT + \vI) \vu_{e-1} + \left[ a \vD \vT - \half a \vb_L \vV_L^\top (\vT + \vI) - \half a \vb_R \vV_R^\top (\vT-\vI) \right] \vu_e + \half a \vb_R \vV_L^\top(\vT - \vI) \vu_{e+1}
\]
Thus the matrices in~\eqref{eq:laup} are given by
\[
\vA_{-1} = \half \vb_L \vV_R^\top(\vT+\vI), \qquad \vA_{+1} = \half \vb_R \vV_L^\top(\vT-\vI), \qquad
\vA_0 = \vD \vT - \half \vb_L \vV_L^\top (\vT+\vI) - \half \vb_R \vV_R^\top (\vT-\vI)
\]
\paragraph{Dissipation model D2.}
The numerical flux is given by
\[
F_\eph = \vV_R^\top \vF_e = a \vV_R^\top \vT \vu_e
\]
and the flux differences at the face are
\[
F_\emh - F_h^\delta(x_\emh^+) = a \vV_R^\top \vT \vu_{e-1} - a \vV_L^\top \vT \vu_e, \qquad F_\eph - F_h^\delta(x_\eph^-) = 0
\]
so that the flux derivative at the solution points is given by
\[
\partial_\xi \vF_h = (a \vV_R^\top \vT \vu_{e-1} - a \vV_L^\top \vT \vu_e)  \vb_L + a \vD \vT \vu_e = a  \vb_L \vV_R^\top \vT \vu_{e-1} + a (\vD \vT - \vb_L \vV_L^\top \vT) \vu_e
\]
Thus the matrices in~\eqref{eq:laup} are given by
\[
\vA_{-1} = \vb_L \vV_R^\top \vT, \qquad \vA_0 = \vD \vT - \vb_L \vV_L^\top \vT, \qquad \vA_{+1} = 0
\]
The upwind character of the flux leads to $\vA_{+}=0$ and the right cell does not appear in the update equation.
\paragraph{Stability analysis.}
We assume a solution of the form $ \vu_e^n =  { \hat {{\vu}}_k^n} \exp(\im k x_e)$ where $\im = \sqrt{-1}$,  $k$ is the wave number which is an integer and $\hat \vu_k^n \in \re^{N+1}$ are the Fourier amplitudes; substituting this ansatz in~\eqref{eq:laup}, we find that the amplitudes evolve according to the equation
\[
\hat \vu_k^{n+1} = \vH(\sigma, \kappa) \hat \vu_k^n, \qquad \vH = \vI - \sigma \vA_0 - \sigma \vA_{-1} \exp(-\im \kappa) - \sigma \vA_{+1} \exp(\im \kappa), \qquad \kappa = k\Delta x
\]
where $\kappa$ is the non-dimensional wave number. The eigenvalues of $\vH$ depend on the CFL number $\sigma$ and wave number $\kappa$, i.e., $\lambda = \lambda(\sigma,\kappa)$; for stability, all the eigenvalues of $\vH$ must have magnitude less than or equal to one for all $\kappa \in [0, 2\pi]$, i.e.,
\[
\lambda(\sigma) = \max_\kappa |\lambda(\sigma, \kappa)| \le 1
\]
The CFL number is the maximum value of $\sigma$ for which above stability condition is satisfied. This CFL number is determined approximately by sampling in the wave number space; we partition $\kappa \in [0,2\pi]$ into a large number of uniformly spaced points $\kappa_j$ and determine
\[
\lambda(\sigma) = \max_j |\lambda(\sigma, \kappa_j)|
\]
The values of $\sigma$ are also sampled in some interval $[\sigma_{\min}, \sigma_{\max}]$ and the largest value of $\sigma_l \in [\sigma_{\min}, \sigma_{\max}]$ for which $\lambda(\sigma_l) \le 1$ is determined in a Python code. We start with a large interval $[\sigma_{\min}, \sigma_{\max}]$ and then progressively reduce the size of this interval so that the CFL number is determined to about three decimal places. The results of this numerical investigation of stability are shown in Table~\eqref{tab:cfl} for two correction functions and different polynomial degrees.
\begin{table}
    \centering
    \begin{tabular}{|c|c|c|c|c|c|c|}
    \hline
    $N$ & \multicolumn{3}{|c|}{Radau} & \multicolumn{3}{|c|}{$g_2$} \\
    \hline
        &  D1 & D2 & Ratio & D1 & D2 & Ratio \\
    \hline
    1 & 0.226 & 0.333 &1.47 & 0.465 & 1.000 &2.15 \\
    2 & 0.117 & 0.170 &1.45& 0.204  & 0.333 &1.63\\
    3 & 0.072 & 0.103 &1.43 & 0.116 & 0.170 & 1.47\\
    4 & 0.049 & 0.069 &1.40 & 0.060 & 0.103 & 1.72\\
    \hline
    \end{tabular}
    \caption{CFL numbers for 1-D LWFR using the two dissipation models and correction functions}
    \label{tab:cfl}
\end{table}
We see that dissipation model D2 has a higher CFL number compared to dissipation model D1. The CFL numbers for the $g_2$ correction function are also significantly higher than those for the Radau correction function. The LWFR scheme with Radau correction function is identical to DG scheme and the CFL numbers found here agree with those from the ADER-DG scheme~\cite{Dumbser2008,Gassner2011a}. The optimality of these CFL numbers has been verified by experiment on the linear advection test case~\eqref{sec:cla}, i.e., the solution eventually blows up if the time step is slightly higher than what is allowed by the CFL condition.
\begin{remark}{(Dirichlet boundary condition)}
The boundary conditions for hyperbolic problems are usually implemented in a weak manner through the fluxes. The boundary condition on the solution can be specified only at inflow boundaries, i.e., where the characteristics are entering the domain. For example, at the left boundary of the domain, say $x=0$, the boundary condition can be specified if $f' > 0$. Assuming this is the case for our problem, let the boundary condition be given as $u(0,t) = g(t)$. It will be enforced by defining an upwind numerical flux at the boundary face, which is given by
\[
F_{lb} = F_{lb}^- \approx \frac{1}{\Delta t} \int_{t_n}^{t_{n+1}} f(g(t)) \ud t
\]
If the integral cannot be computed analytically, then it is approximated by quadrature in time. From~\eqref{eq:tvgproperty}, we see that integral must be at least accurate to $O(\Delta t^{N+1})$ which is of the same order as the neglected terms in~\eqref{eq:tvgproperty}. In the numerical tests, we use $(N+1)$-point Gauss-Legendre quadrature which ensures the required accuracy. If the right boundary is an outflow boundary, i.e., $f' > 0$, then the flux across the right boundary is computed in an upwind manner using the interior solution, i.e., $F_{rb} = F_{rb}^-$ where $F_{rb}^-$ is obtained from the Lax-Wendroff procedure.
\end{remark}
\section{TVD limiter}\label{sec:lim}
The computation of discontinuous solutions with high order methods can give oscillatory solutions which must be limited by some non-linear process. The a posteriori limiters developed in the context of RKDG schemes~\cite{Cockburn1991,Cockburn1989a} can be adopted in the framework of LWFR schemes. The limiter is applied in a post-processing step after the solution is updated to the new time level. The limiter is thus applied only once for each time step unlike in RKDG scheme where it has to be applied after each RK stage update. Let $u_h(x)$ denote the solution at time $t_{n+1}$. In element $\Omega_e$, let the average solution be $\au_e$; define the backward and forward differences of the solution and cell averages by
\[
\Delta^- u_e = \au_e - u_h(x_\emh^+), \qquad \Delta^+ u_e = u_h(x_\eph^-) - \au_e
\]
\[
\Delta^- \au_e = \au_e - \au_{e-1}, \qquad \Delta^+ \au_e = \au_{e+1} - \au_e
\]
We limit the solution by comparing its variation within the cell with the difference of the neighbouring cell averages through a limiter function,
\[
\Delta^- u_e^{m} = \minmod(\Delta^- u_e, \Delta^- \au_e, \Delta^+ \au_e), \qquad
\Delta^+ u_e^{m} = \minmod(\Delta^+ u_e, \Delta^- \au_e, \Delta^+ \au_e)
\]
where
\[
\minmod(a,b,c) = \begin{cases}
s \min(|a|, |b|, |c|), & \textrm{if } s = \sign(a) = \sign(b) = \sign(c) \\
0, & \textrm{otherwise}
\end{cases}
\]
If $\Delta^- u_e^m \ne \Delta^- u_e$ or $\Delta^+ u_e^m \ne \Delta^+ u_e$, then the solution is deemed to be locally oscillatory and we modify the solution inside the cell by replacing it as a linear polynomial with a limited slope, which is taken to be the average limited slope. The limited solution polynomial in cell $\Omega_e$ is given by
\[
u_h|_{\Omega_e} = \au_e + \frac{\Delta^- u_e^m + \Delta^+ u_e^m}{2} (2 \xi - 1), \qquad \xi \in [0,1]
\]
This limiter is known to clip smooth extrema since it cannot distinguish them from jump discontinuities. A small modification based on the idea of TVB limiters~\cite{Cockburn1991} can be used to relax the amount of limiting that is performed which leads to improved resolution of smooth extrema. The $\minmod$ function is replaced by
\[
\overline{\minmod}(a,b,c) = \begin{cases}
a, & |a| \le M \Delta x^2 \\
\minmod(a,b,c), & \textrm{otherwise}
\end{cases}
\]
which requires the choice of a parameter $M$, which is an estimate of the second derivative of the solution at smooth extrema. In the case of systems of equations, the limiter is applied to the characteristic variables, which is known to yield better control on the spurious numerical oscillations~\cite{Cockburn1989}. The limiters used in this work are not able to provide high order accuracy and the development of better limiters for LW schemes is deferred to future work. Here we would like to show that with the use of TVD-type limiters, the LW and RK schemes provides similar resolution.
\section{Numerical results in 1-D: scalar problems}\label{sec:res1d}
In this section, we present some numerical results to show the convergence rates and the comparison of different scheme parameters like correction function, solution points and dissipation model. For each problem in this section, the corresponding CFL number is chosen from Table~\ref{tab:cfl}. Here after, when we use the CFL numbers obtained using the Fourier stability analysis, we multiply it with a safety factor of 0.95. When D1 and D2 schemes are compared together, the CFL numbers of D1 schemes are used as these are lower; the same CFL numbers are used for the RKFR schemes. Up to degree $N=3$, RKFR schemes use Runge-Kutta time integration of order $N+1$ with $N+1$ stages. In the $N=4$ case, for non-linear problems there is no five stage Runge-Kutta method of order 5, see Chapter 32 of~\cite{Butcher2016}. However, for linear, autonomous problems, the five stage SSPRK method in~\cite{Gottlieb2001} is fifth order accurate,  and we make use of it for the constant advection test cases with periodic boundary conditions and refer to it as SSPRK55. For non-linear or non-periodic problems, to make a fair comparison of LW and RK, we make use of the six stage, fifth order Runge-Kutta (RK65) time integration introduced in~\cite{Tsitouras2011}.

The RKFR and LWFR schemes are illustrated at a high level in Algorithm 1 and Algorithm 2, respectively, for solving a hyperbolic conservation law in a time interval $[0,T]$. Here we assume that an a posteriori limiter like a TVD/TVB limiter and a positivity limiter are applied in a post-processing step after the solution is updated. The LWFR scheme requires the application of the limiter only once per time step while the RKFR scheme requires multiple applications of the limiter depending on the number of RK stages. The limiter can be costly to apply for systems of equations where a characteristic approach and/or WENO-type limiter is used. In the present work, we use a simple TVD/TVB limiter but use characteristic limiting for systems.
\begin{algorithm}
$t=0$\;
\While{$t < T$}
{
\For{each RK stage}
{
Loop over cells and assemble rhs\;
Loop over faces and assemble rhs\;
Update solution to next RK stage\;
Apply a posteriori limiter\;
Apply positivity limiter\;
}
$t = t + \Delta t$\;
}
\caption{RKFR scheme}
\end{algorithm}

\begin{algorithm}
$t=0$\;
\While{$t < T$}
{
Loop over cells and assemble rhs\;
Loop over faces and assemble rhs\;
Update solution to next time level\;
Apply a posteriori limiter\;
Apply positivity limiter\;
$t = t + \Delta t$\;
}
\caption{LWFR scheme}
\end{algorithm}

\subsection{Linear advection equation: constant speed \label{sec:cla}}
We first consider the 1-D linear advection equation $u_t+au_x=0$ with speed $a=1$ and periodic boundary condition.
\subsubsection{Smooth solutions}
\begin{figure}
\centering
\begin{tabular}{cc}
\includegraphics[width=0.40\textwidth]{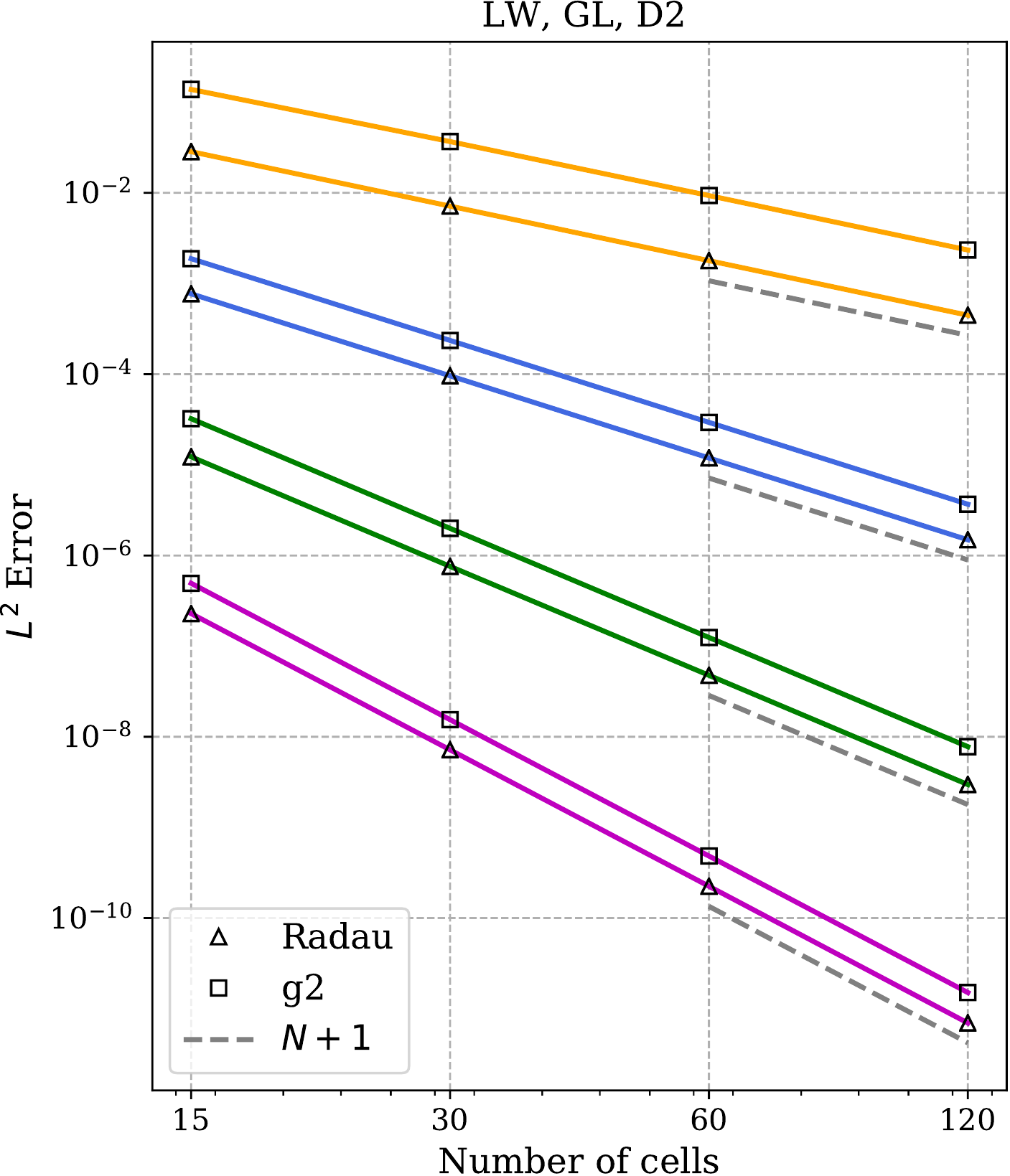} &
\includegraphics[width=0.40\textwidth]{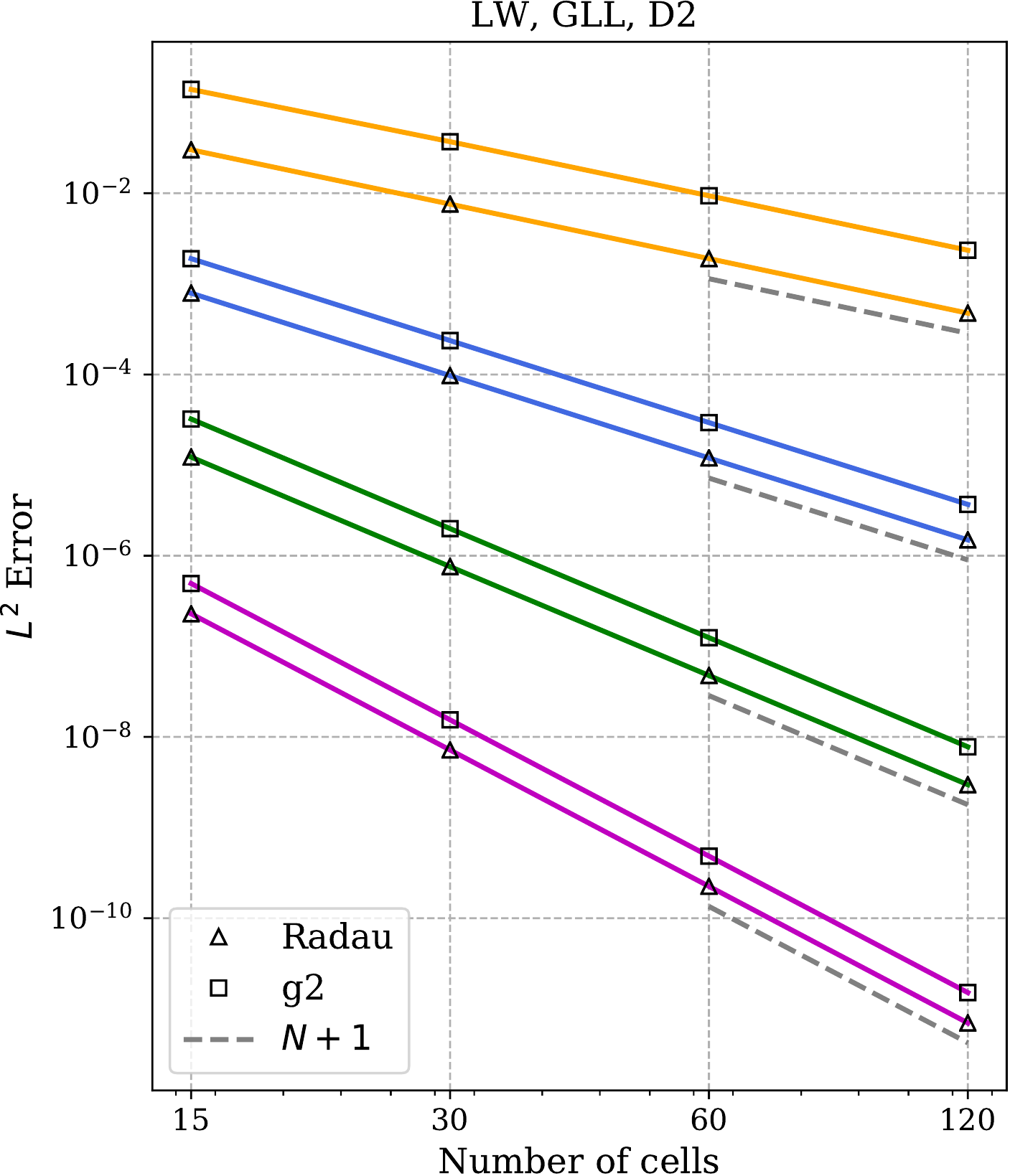} \\
(a) & (b)
\end{tabular}
\caption{Error convergence for constant linear advection; (a) GL points, (b) GLL points. The different colors correspond to degrees $N=1,2,3,4$ from top to bottom.}
\label{fig:cla1}
\end{figure}

For the initial condition $u(x,0)=\sin(2\pi x)$ with periodic boundaries  on $[0,1]$, we perform grid convergence studies using various parameters like correction functions and solution points. The error norms are computed at time $t=2$ units. In Figure~(\ref{fig:cla1}) we compare the convergence behaviour for Radau and g2 correction functions and for both choices of solution points using the D2 dissipation model. It is clear that the errors due to Radau are consistently smaller than those with g2 correction function. The choice of the solution points does not significantly affect the error in the solution.

\begin{figure}
\begin{center}
\begin{tabular}{cc}
\includegraphics[width=0.40\textwidth]{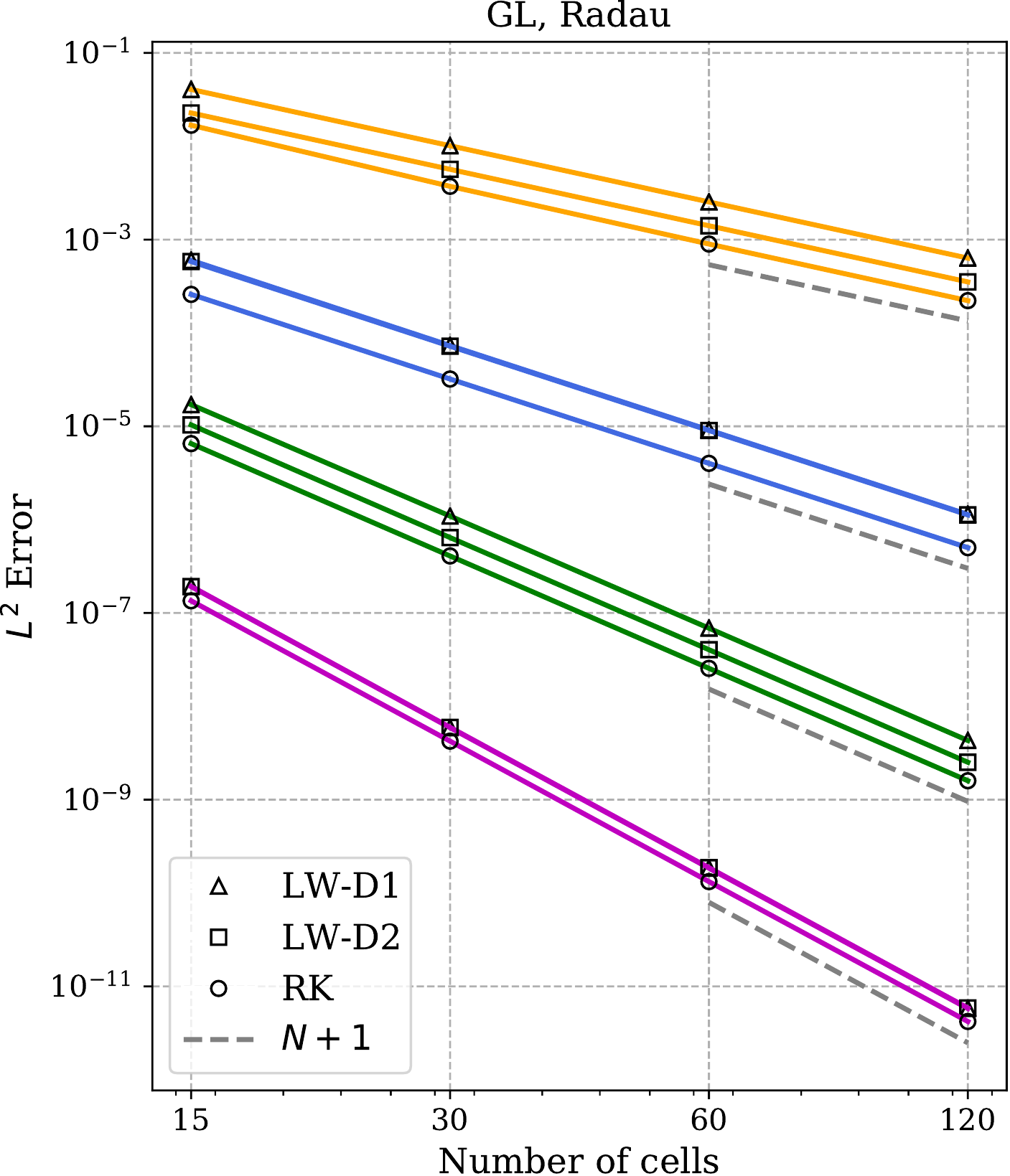} &
\includegraphics[width=0.40\textwidth]{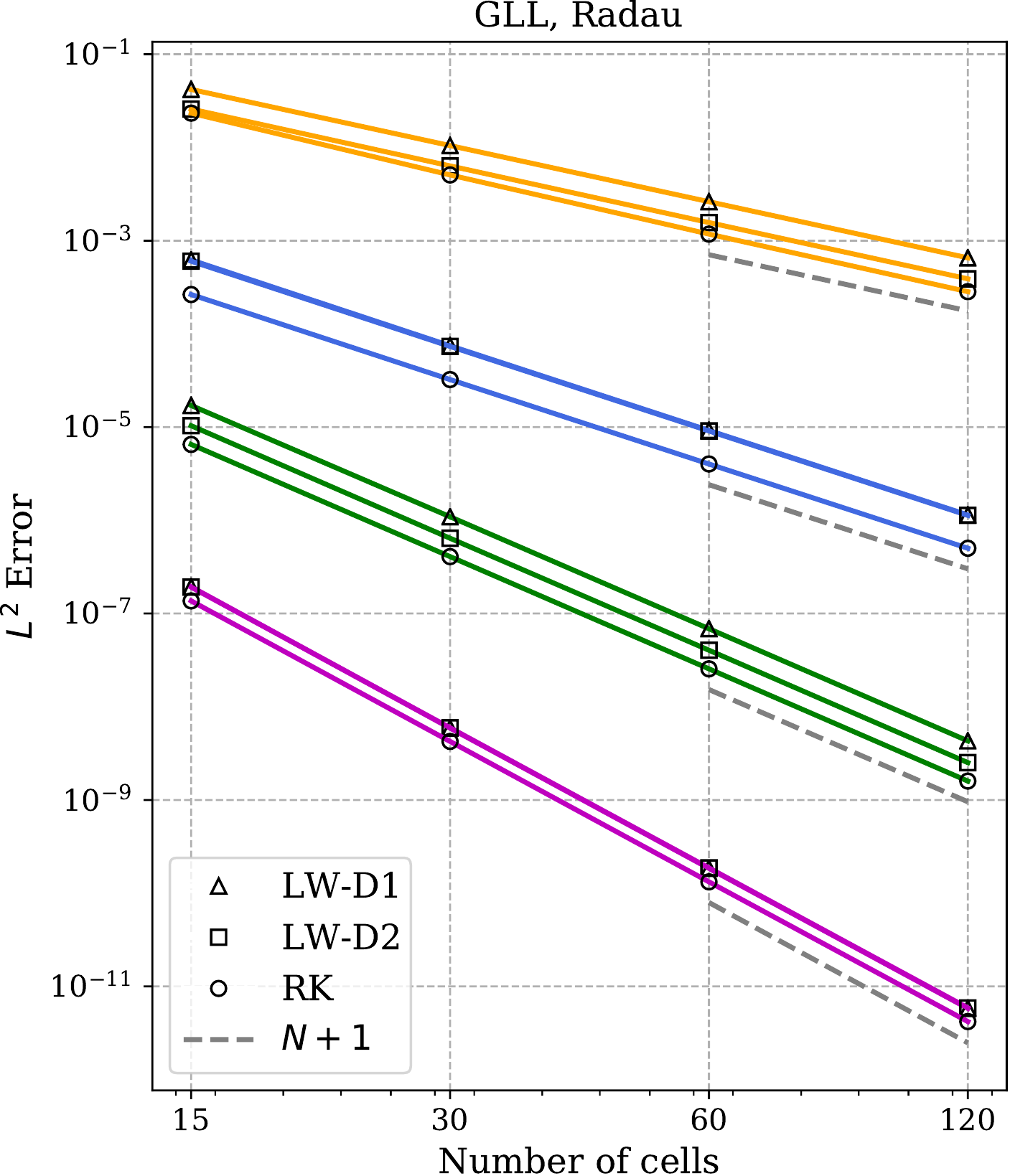} \\
(a) & (b)
\end{tabular}
\end{center}
\caption{Error convergence for constant linear advection; (a) GL points , (b) GLL points. The different colors correspond to degrees $N=1,2,3,4$ from top to bottom.}
\label{fig:cla2}
\end{figure}

\begin{figure}
\begin{center}
\begin{tabular}{cc}
\includegraphics[width=0.40\textwidth]{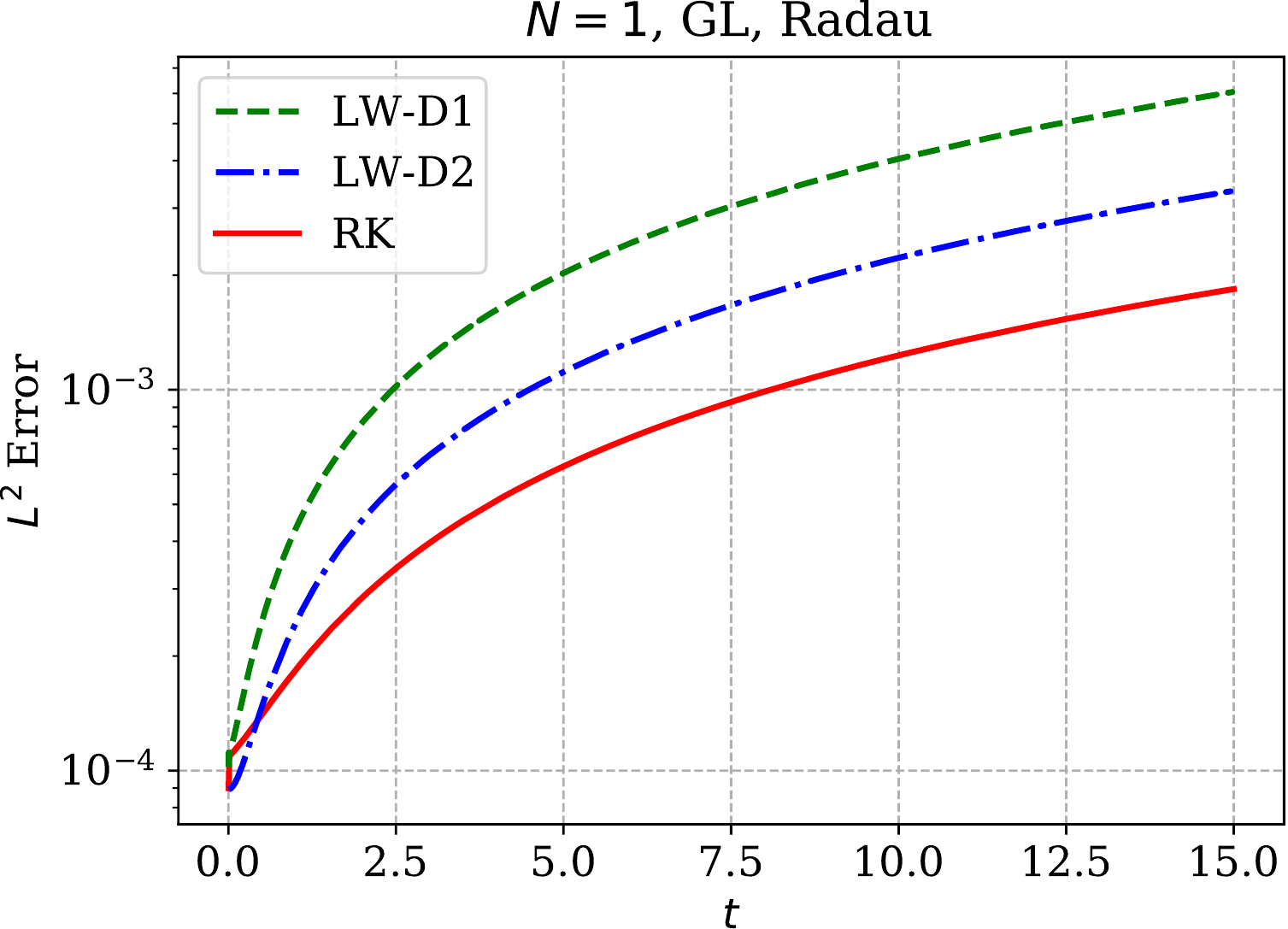} &
\includegraphics[width=0.40\textwidth]{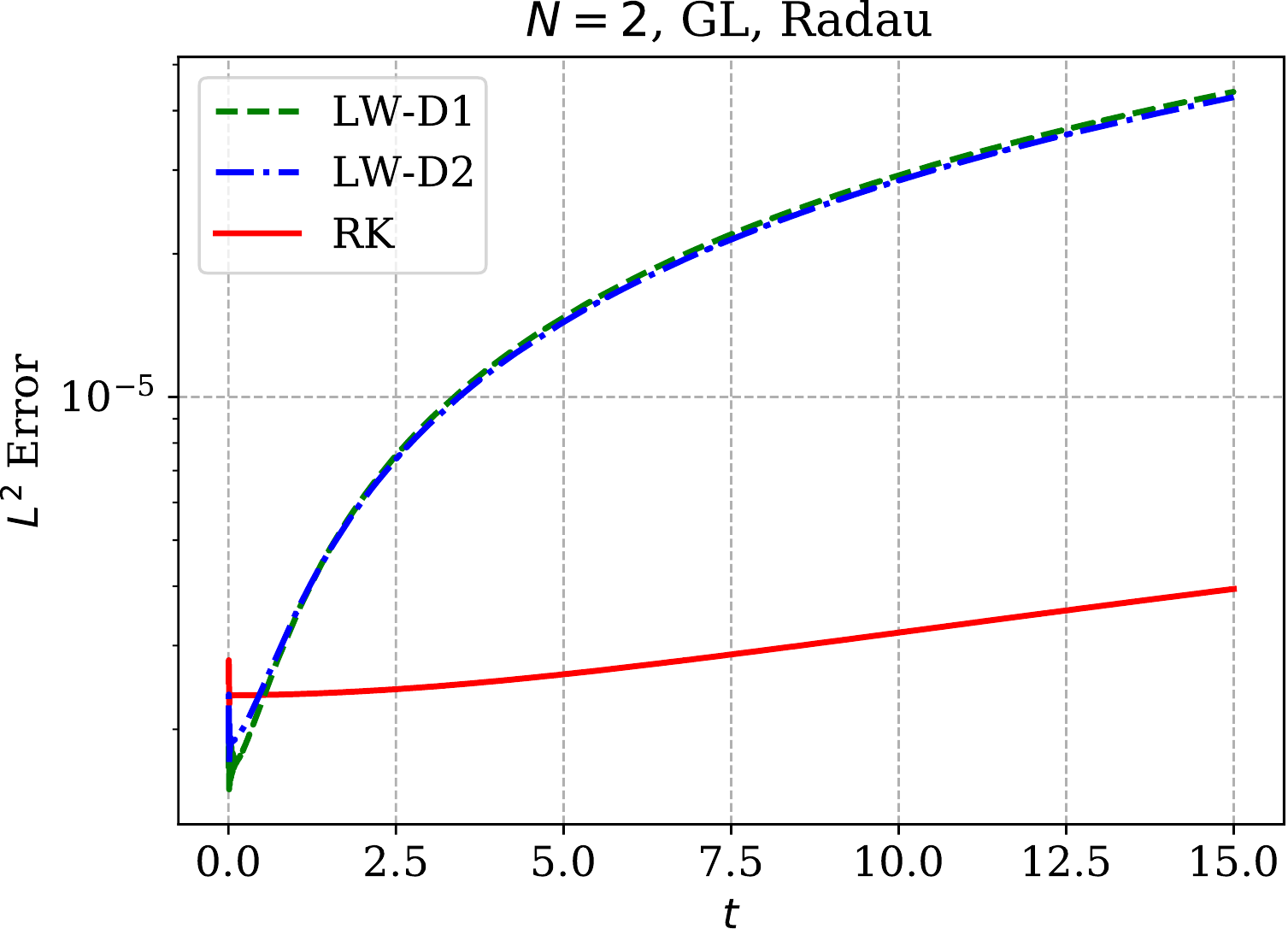} \\
(a) & (b) \\
\includegraphics[width=0.40\textwidth]{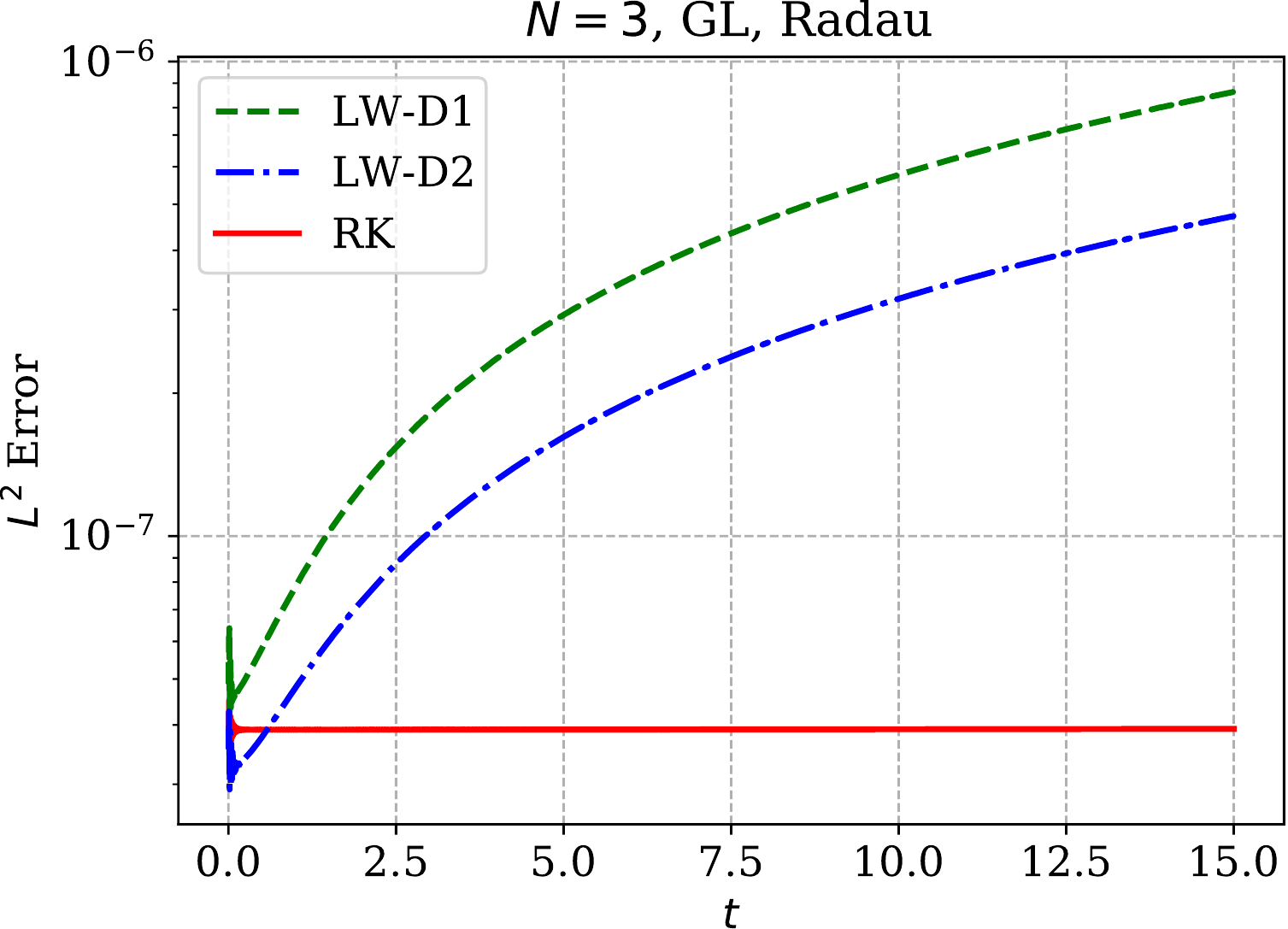} &
\includegraphics[width=0.40\textwidth]{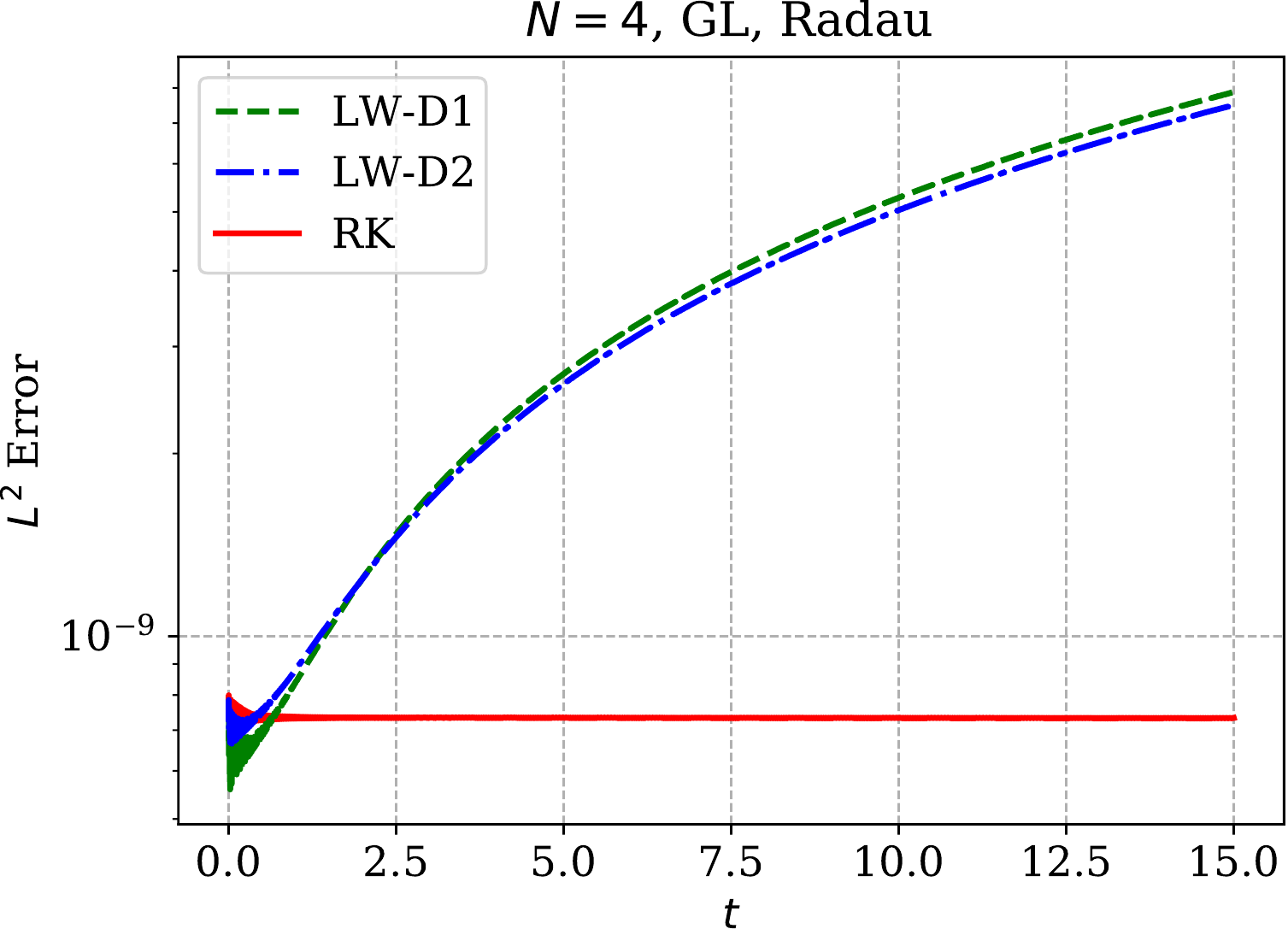} \\
(c) & (d)
\end{tabular}
\end{center}
\caption{Error versus time for constant linear advection, initial condition $u(x, 0) = \sin(2 \pi x)$, $x\in [0,1]$, for different polynomial degrees, each with 200 dofs; GL solution points and Radau correction. (a) $N=1$, (b) $N=2$, (c) $N=3$, (d) $N=4$.}
\label{fig:cla3}
\end{figure}

Figure~(\ref{fig:cla2}) shows the comparison of LW and RK schemes using Radau correction and two types of solution points. There is a small difference in the error levels between the two dissipation models, with the D2 model performing better for odd $N$. The RK scheme has slightly smaller errors than the LW scheme. We can see this more clearly by plotting the error norm versus time as shown in Figure~(\ref{fig:cla3}), where all the four degrees consist of the same number of total dofs which is 200. We see that the error growth with time is higher for the LW scheme than for the RK scheme. The superior performance of the RK scheme is already known in the literature~\cite{Guo2015} and is due to its super-convergence properties. It is possible to construct LW schemes which are also super-convergent as done in~\cite{Guo2015} but the resulting schemes are computationally more expensive as they involve a stronger coupling with the neighbouring cell solutions, than what is used in the standard LW schemes. Hence we do not pursue that approach for our LW schemes. Note that this super-convergence occurs for constant coefficient linear problems on uniform grids and these advantages of RK schemes are not present when we consider non-linear problems, as shown in later results.

\begin{figure}
\begin{center}
\begin{tabular}{cc}
\includegraphics[width=0.40\textwidth]{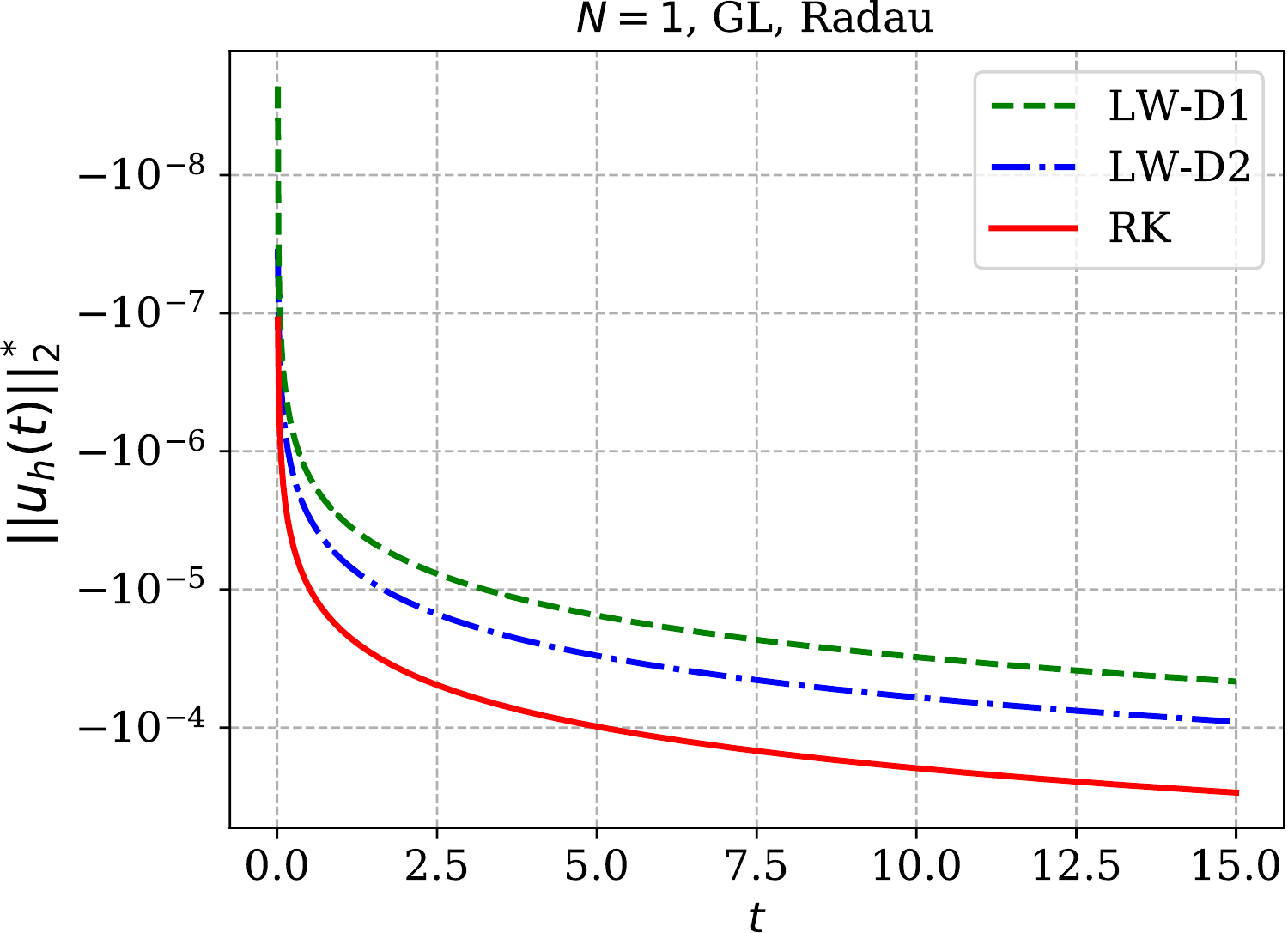} &
\includegraphics[width=0.40\textwidth]{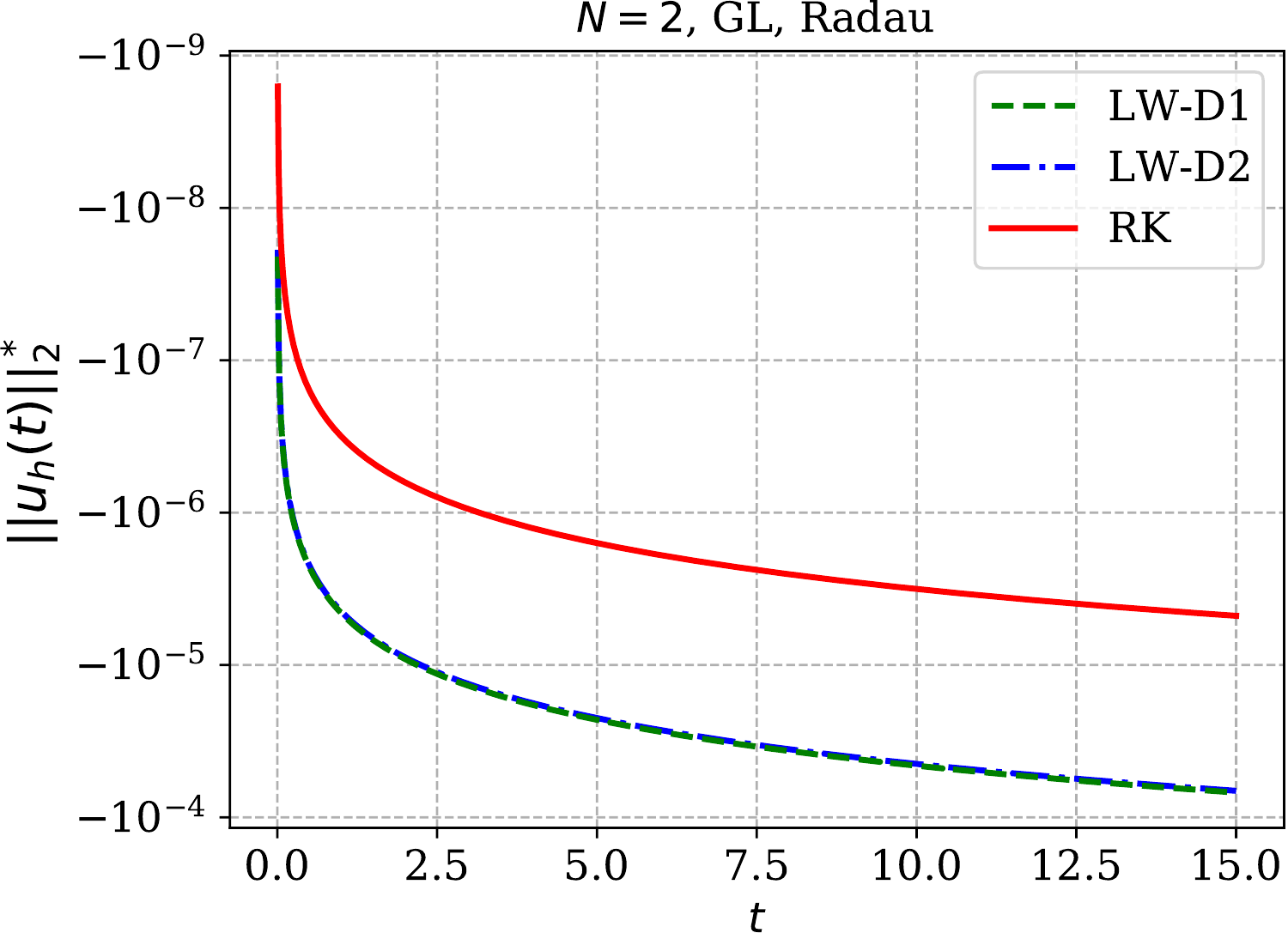} \\
(a) & (b) \\
\includegraphics[width=0.40\textwidth]{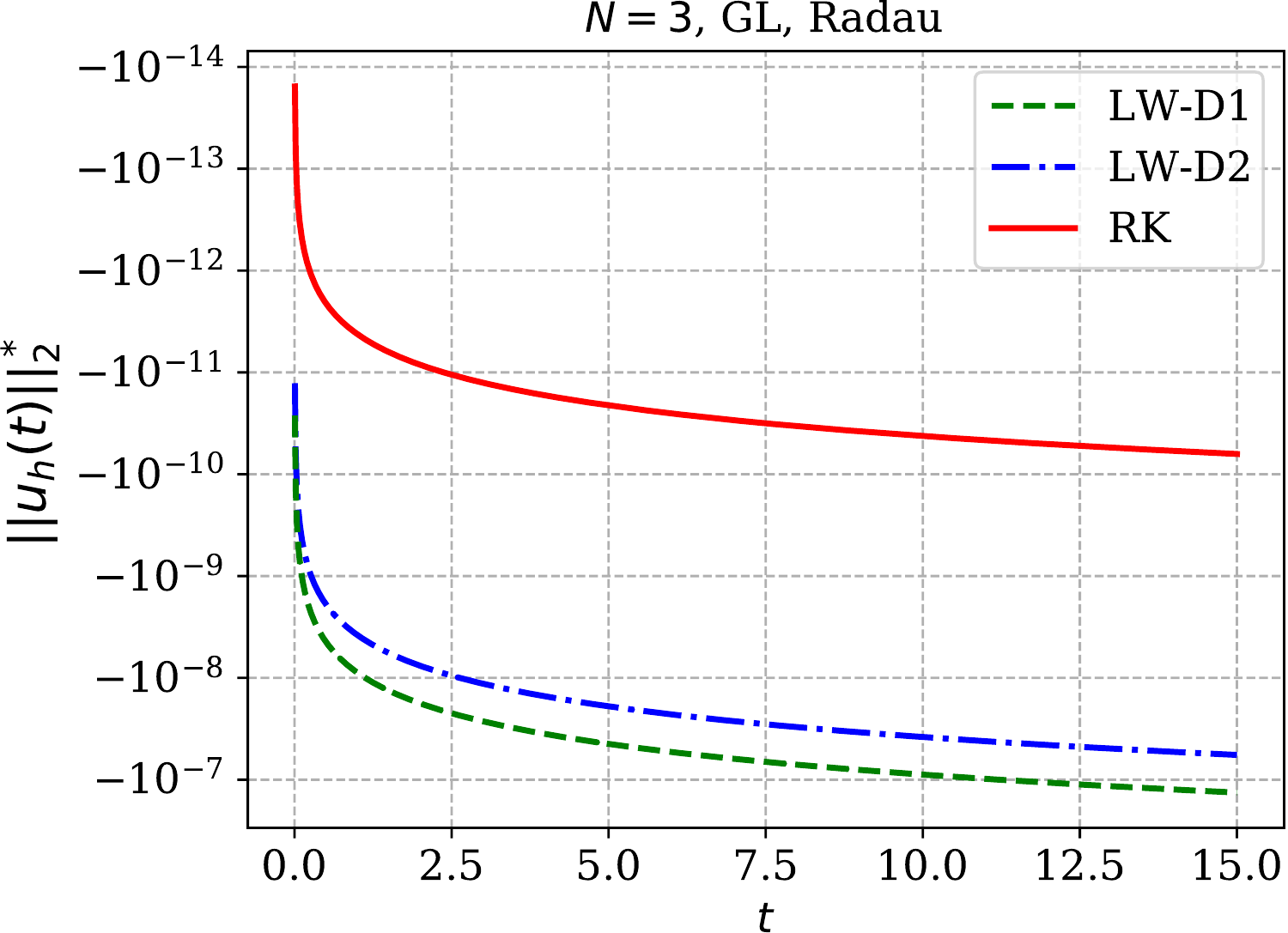} &
\includegraphics[width=0.40\textwidth]{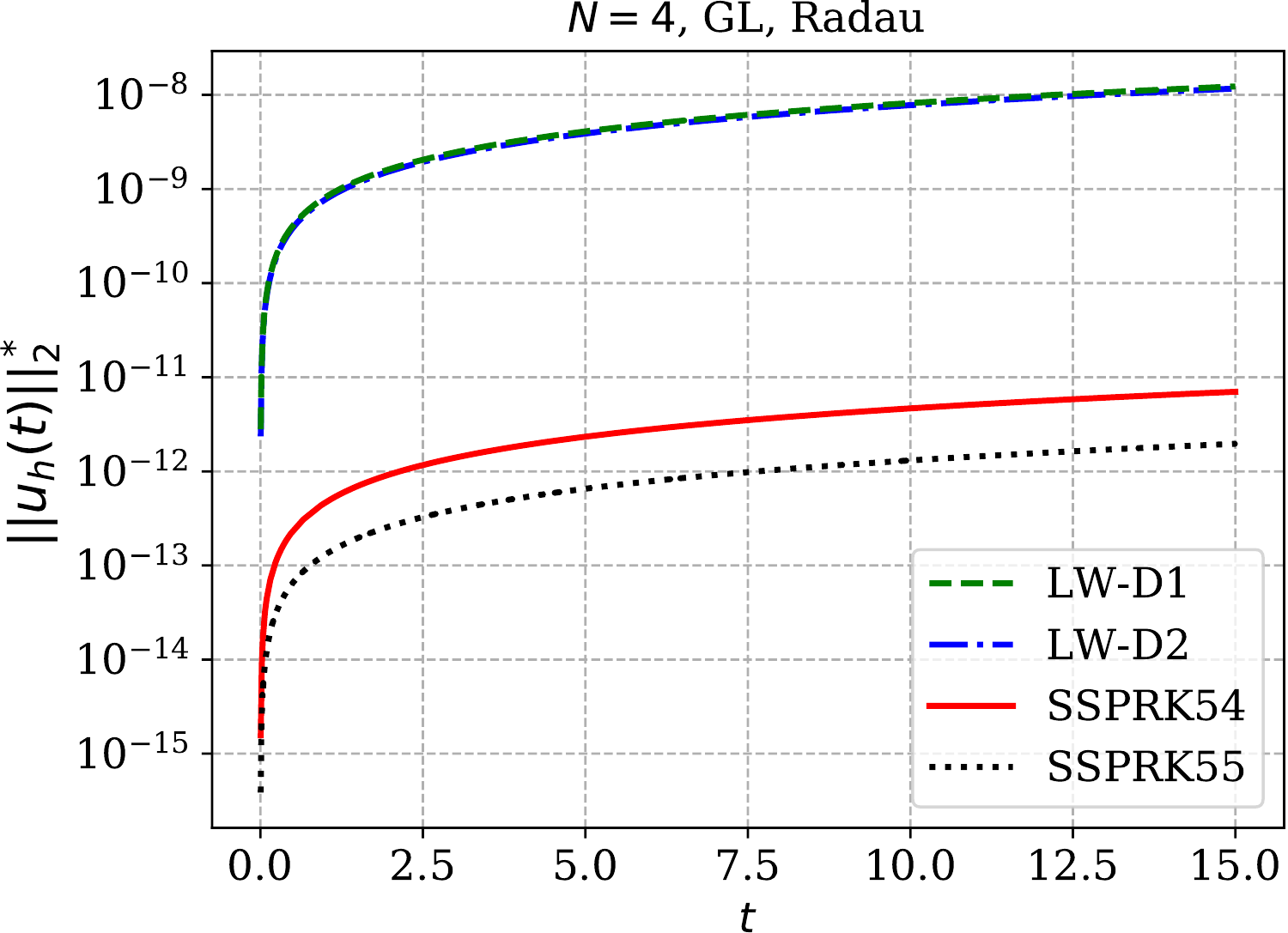} \\
(c) & (d)
\end{tabular}
\end{center}
\caption{Semi-log plot of relative change in $L^2$ norm versus time for constant linear advection with the initial condition $u(x, 0) = \sin(2 \pi x)$, $x\in [0,1]$ for different polynomial degrees, each with 200 dofs; GL solution points and Radau correction.}
\label{fig:cla4}
\end{figure}

Figure~(\ref{fig:cla4}) analyzes the behaviour of $L^2$ norm of the solution where we plot the relative change in the $L^2$ norm with respect to the initial value, defined as  $\|u_h(t)\|_2^*=\frac{\|u_h(t)\|_2-\|u_h(0)\|_2}{\|u_h(0)\|_2}$. For $N=1$, we see that LW is less dissipative than RK and thus better at conserving energy while for $N=2,3$, RK schemes perform better. For $N=4$, we see a mild instability for both LWFR and RKFR schemes. For $N=4$, we compare LW with SSPRK55 which is fifth order only for linear problems and SSPRK54~\cite{Spiteri2002}, which is more relevant for non-linear problems and is fourth order accurate. Choosing time step size by CFL numbers of the LW-D1 scheme, we observe the mild instability for both the time integration schemes. The instability of RKDG scheme has been studied in~\cite{Xu2019} which can be remedied by using an RK scheme with different number of stages; however the use of six stage RK65 method with a limiter may dampen the solution too much, as we discuss later in Figure~(\ref{fig:hat3}). The solution norm grows linearly, with a very small slope for both LW-D1 and LW-D2 (approximately 6.177e-10 and 5.415e-10) schemes, and also for SSPRK54 and SSPRK55 (approximately 2.862e-13, 1.908e-13) schemes. This mild instability does not affect the error convergence rates for finite time problems, but it might lead to problems for very long time simulations with periodic boundary conditions, especially for non-linear problems and needs to be investigated further. This type of mild instability for $N=4$ seems to be present in other single step methods like ADER-DG schemes also.

\begin{figure}
\centering
\begin{tabular}{cc}
\includegraphics[width=0.40\textwidth]{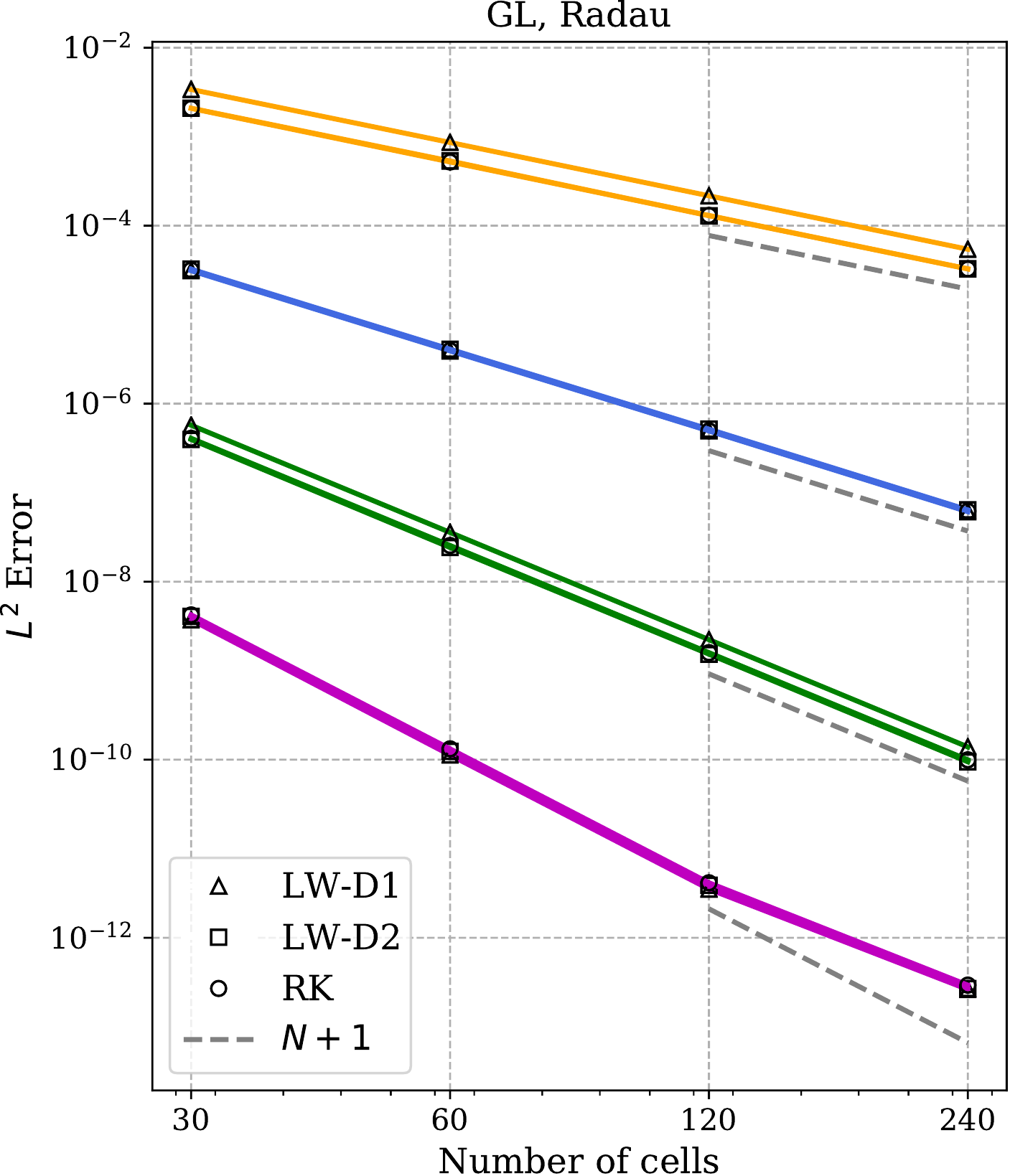} &
\includegraphics[width=0.40\textwidth]{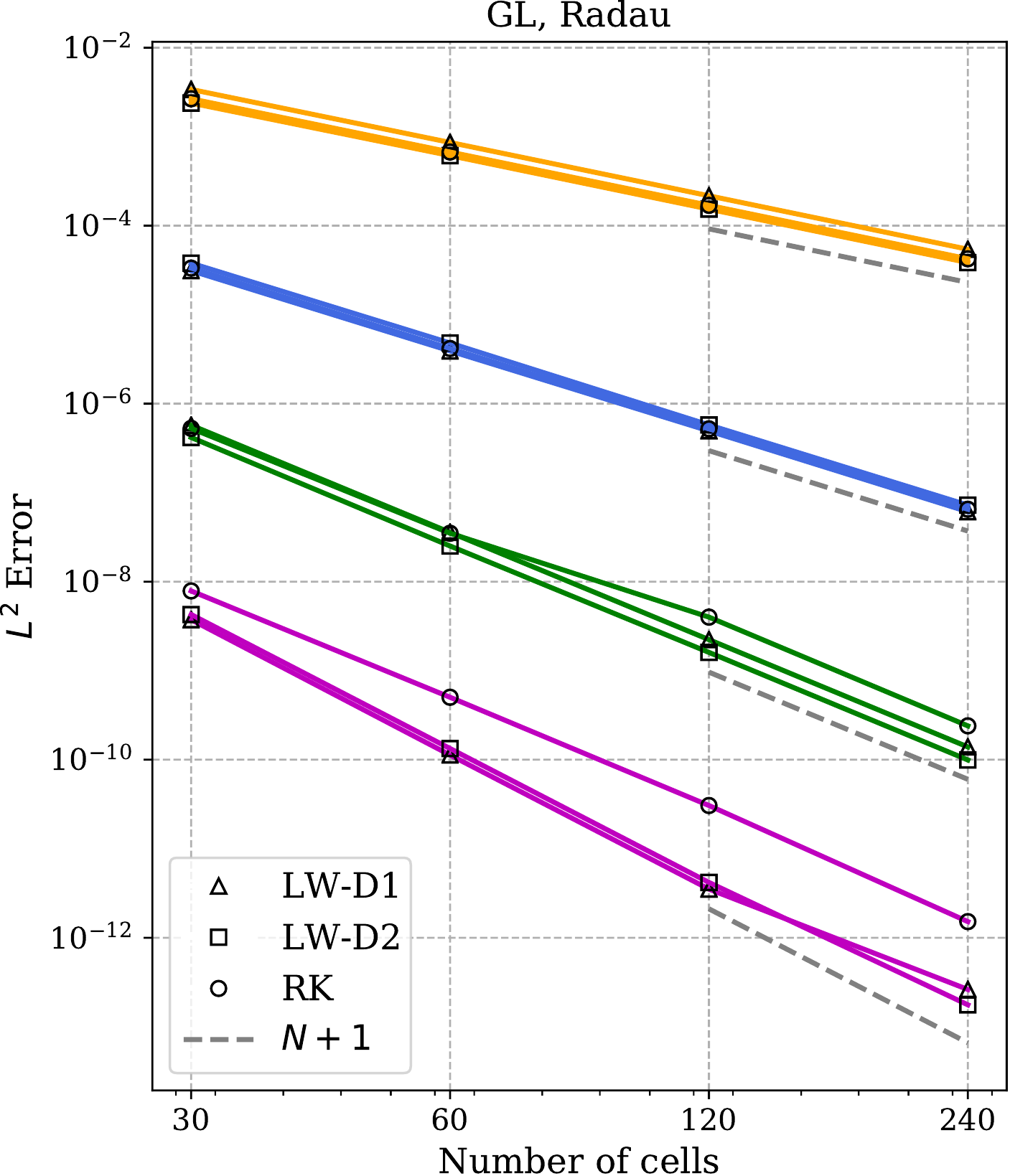} \\
(a) & (b)
\end{tabular}
\caption{Convergence for constant linear advection with Dirichlet boundary conditions; (a) using  CFL numbers of LW-D1 for all schemes, (b) using corresponding CFL number for each scheme. The different colors correspond to degrees $N=1,2,3,4$ from top to bottom.}
\label{fig:cla2_dirichlet}
\end{figure}

\begin{figure}
\centering
\begin{tabular}{cc}
\includegraphics[width=0.46\textwidth]{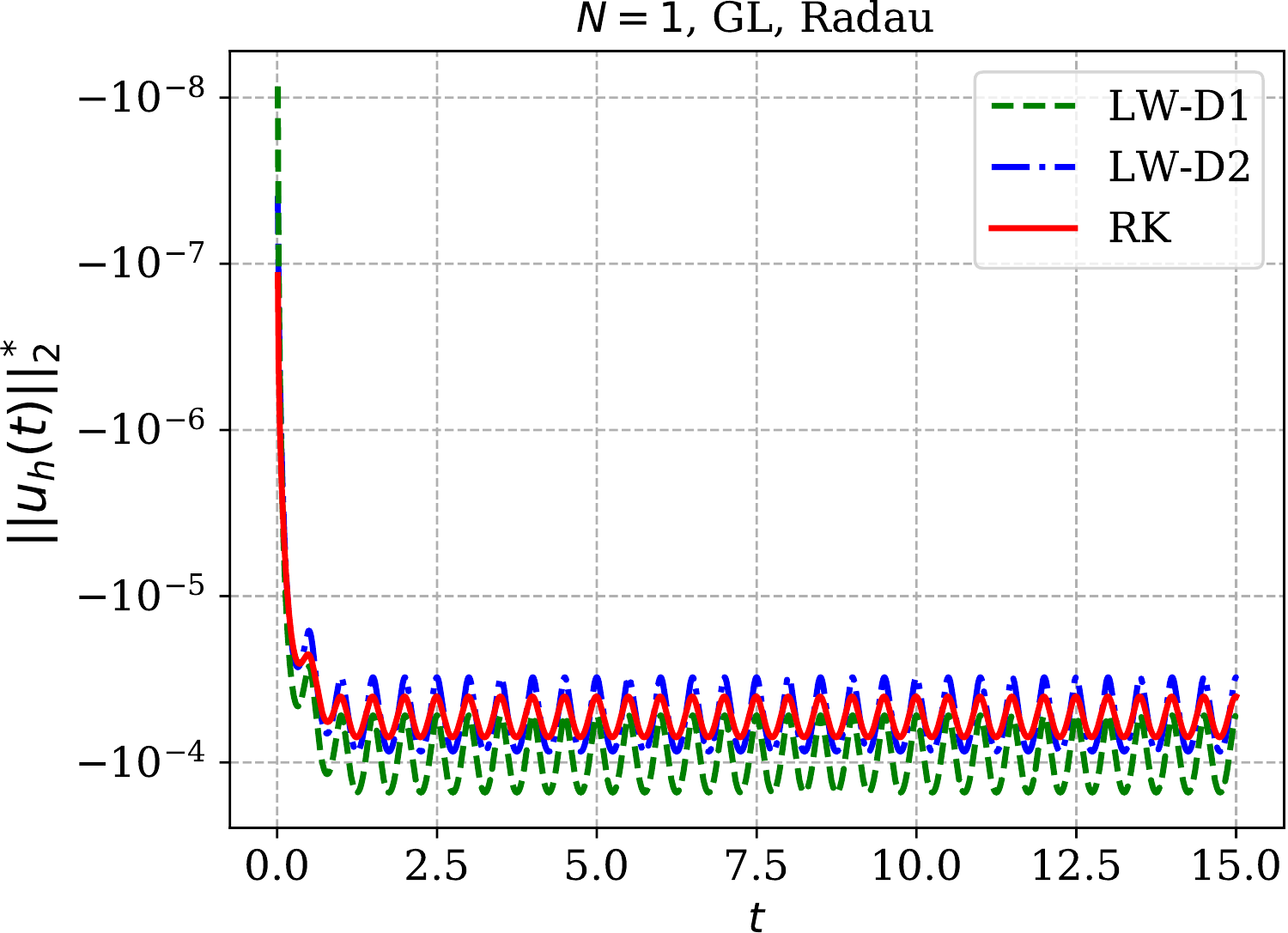} &
\includegraphics[width=0.46\textwidth]{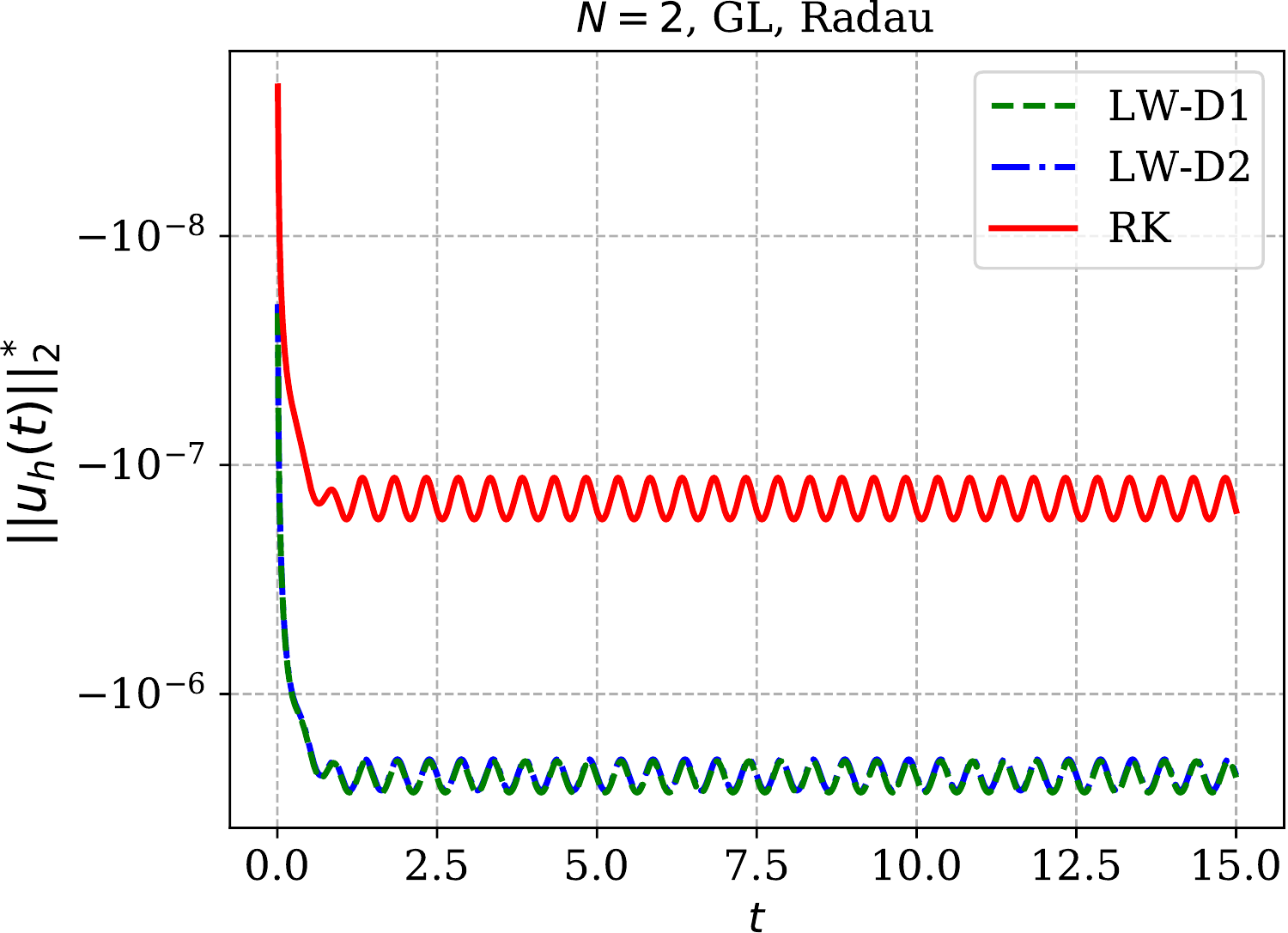} \\
(a) & (b) \\
\includegraphics[width=0.46\textwidth]{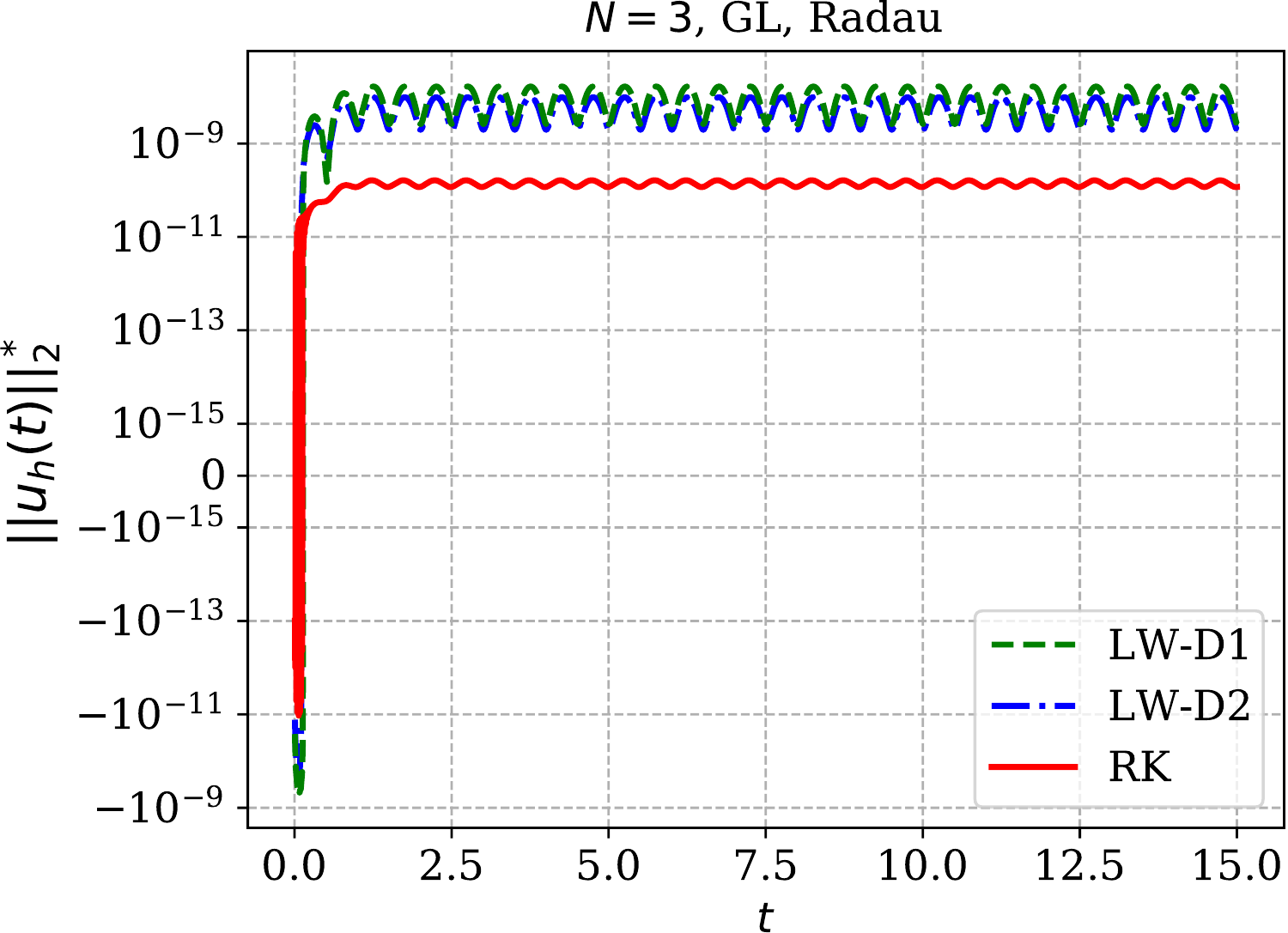} &
\includegraphics[width=0.46\textwidth]{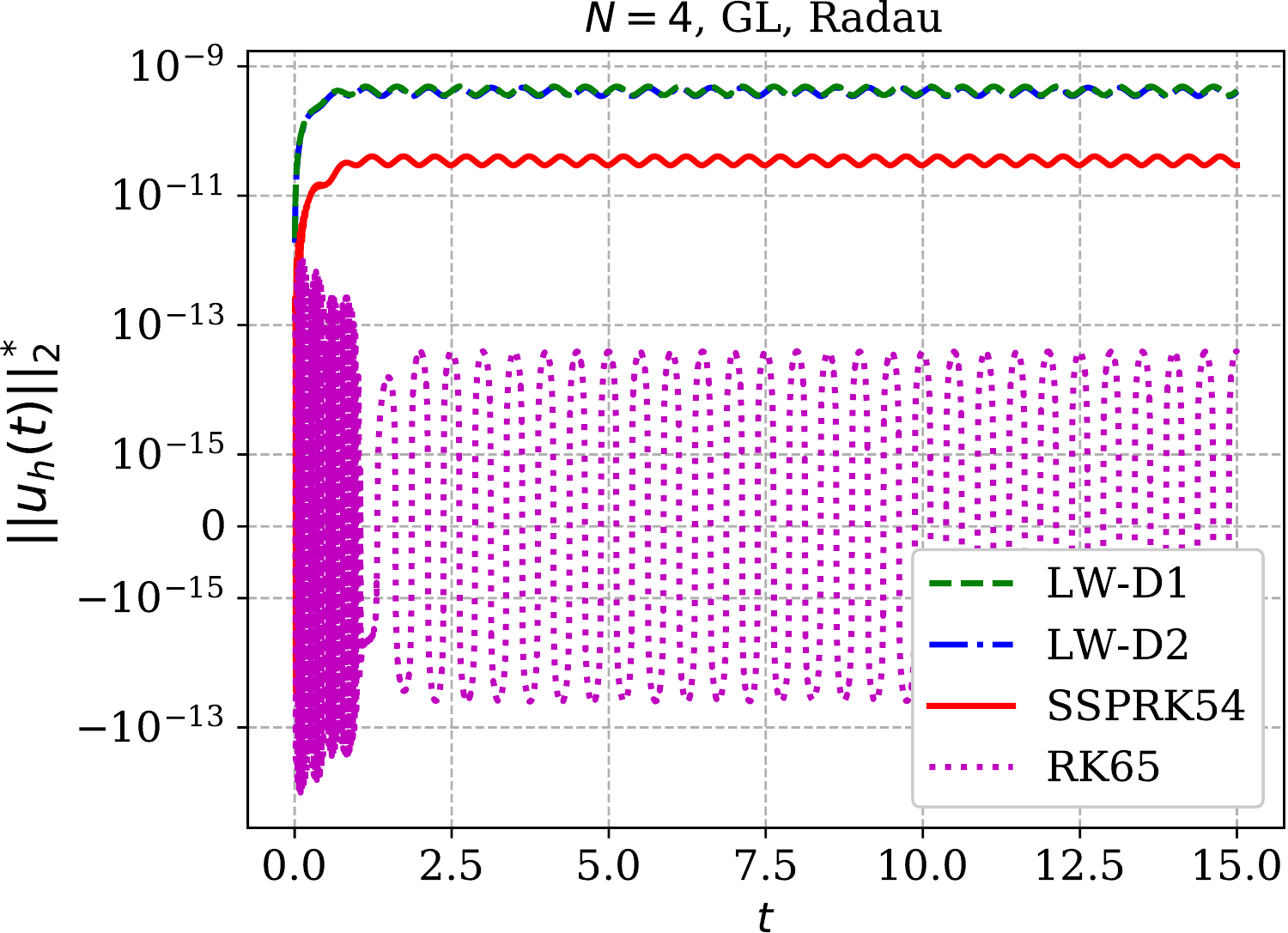} \\
(c) & (d) \\
\end{tabular}
\caption{Semi-log plot of relative change in $L^2$ norm versus time for constant linear advection, with initial condition $u(x, 0) = \sin(2 \pi x)$, $x\in [0,1]$ together with Dirichlet boundary conditions, for different polynomial degrees, each with 200 dofs; GL solution points and Radau correction.}
\label{fig:cla5}
\end{figure}

We now solve the problem with the same initial condition but using Dirichlet boundary condition at the left side of the domain; the exact solution remains same as before.  The fifth order SSPRK scheme of~\cite{Gottlieb2001} is only for autonomous systems, so here we use RK65~\cite{Tsitouras2011} for $N=4$. Figure~(\ref{fig:cla2_dirichlet}a) shows the error convergence when we use the CFL of LW-D1 scheme for all schemes, since this is the smallest. All the schemes show optimal convergence rates with the LW-D2 and RK schemes showing very similar errors. In Figure~(\ref{fig:cla2_dirichlet}b), we perform the same error convergence study but using the stable CFL number of each scheme; we still observe optimal convergence rates in each scheme, but the RK scheme shows slightly larger errors at degrees $N=3,4$, when the error level has become small. The issue with RK schemes may be related with the way Dirichlet conditions are implemented inside the RK stages as studied in~\cite{Carpenter1995}. Figure~(\ref{fig:cla5}) shows the time history of the relative change in $L^2$ norm of the numerical solution; for degrees $N=1,2$ the norm does not increase relative to the initial value but for degrees $N=3,4$, there is some increase in the norm at initial times for some schemes. In Figure~(\ref{fig:cla5}), we also make comparison of LW and RK for $N=4$, with RK time integration performed with SSPRK54~\cite{Spiteri2002}. However, in all cases, the norm does not grow monotonically but reaches a periodic oscillatory behaviour. Unlike the case of periodic boundary conditions, the inflow and outflow boundary conditions may lead to increase and decrease in energy respectively; if the two mechanisms aren't exactly balanced, we can observe the oscillatory behaviour in the solution norm.

\begin{figure}
\begin{center}
\begin{tabular}{cc}
\includegraphics[width=0.47\textwidth]{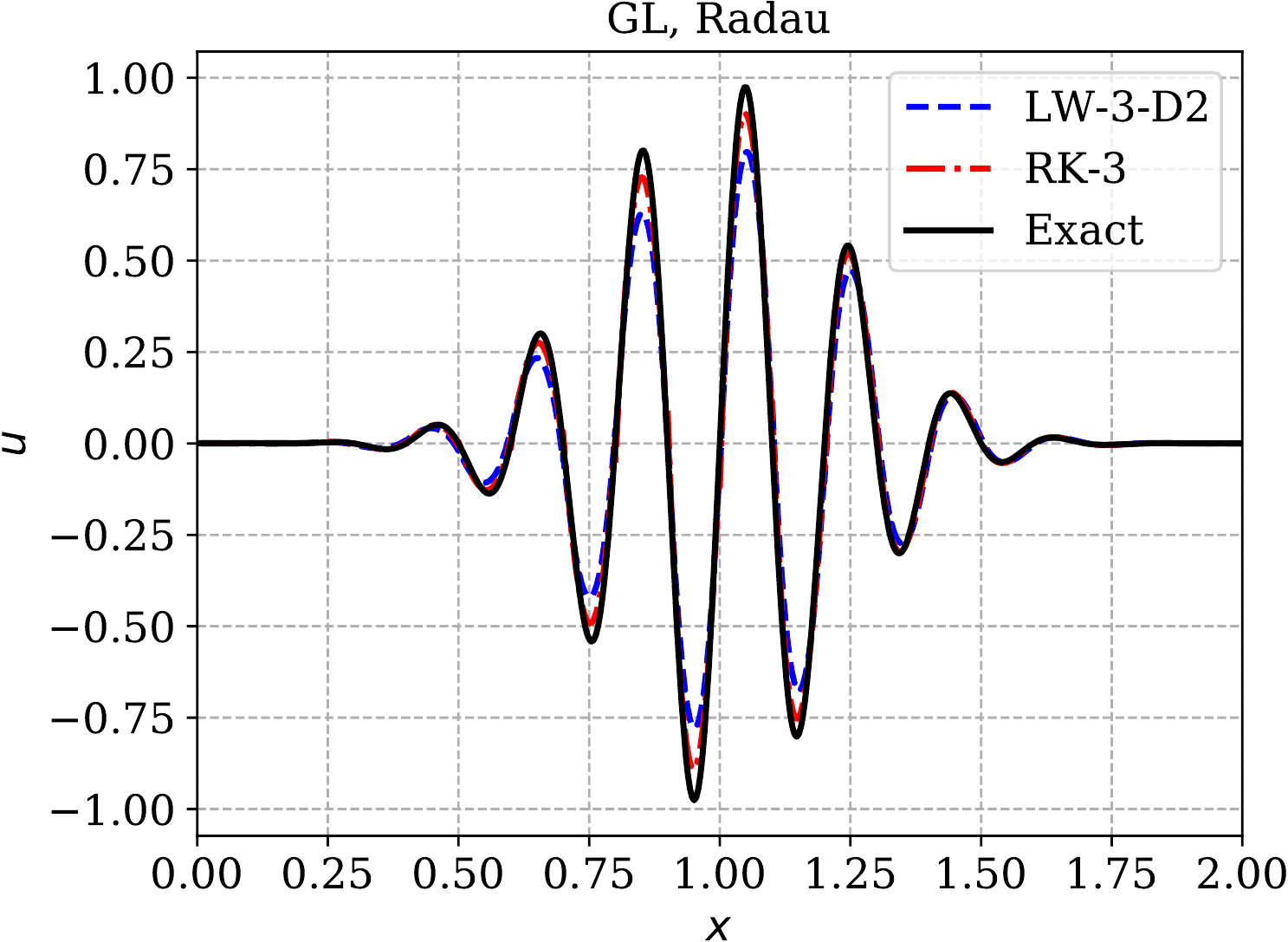} &
\includegraphics[width=0.47\textwidth]{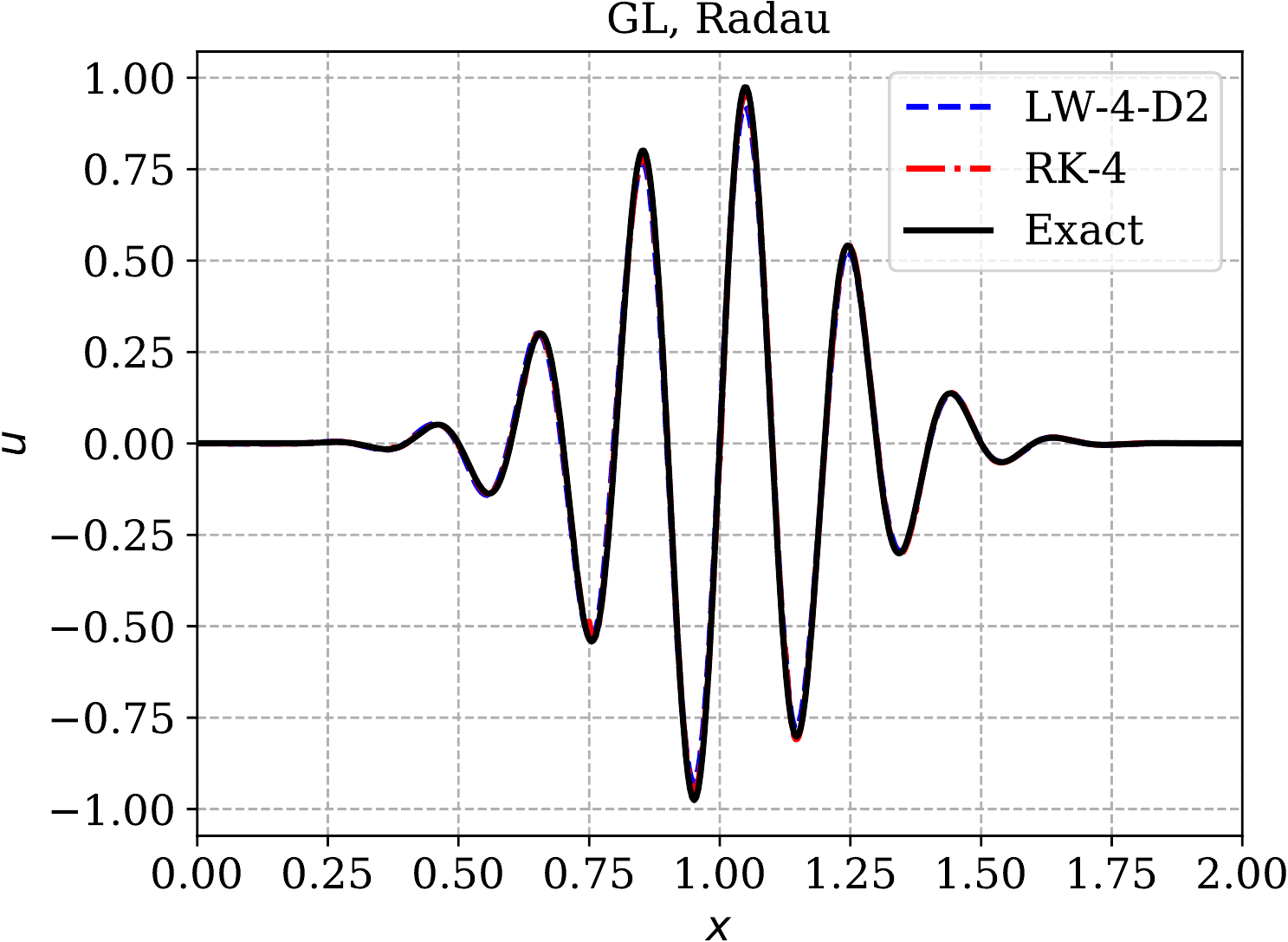} \\
(a) $N=3$ & (b) $N=4$
\end{tabular}
\end{center}
\caption{Constant linear advection of a wave packet; solution at time $t=1.0$ with 160 dofs using polynomial degree (a) $N=3$, (b) $N=4$.}
\label{fig:wp1}
\end{figure}

\begin{figure}
\begin{center}
\begin{tabular}{cc}
\includegraphics[width=0.40\textwidth]{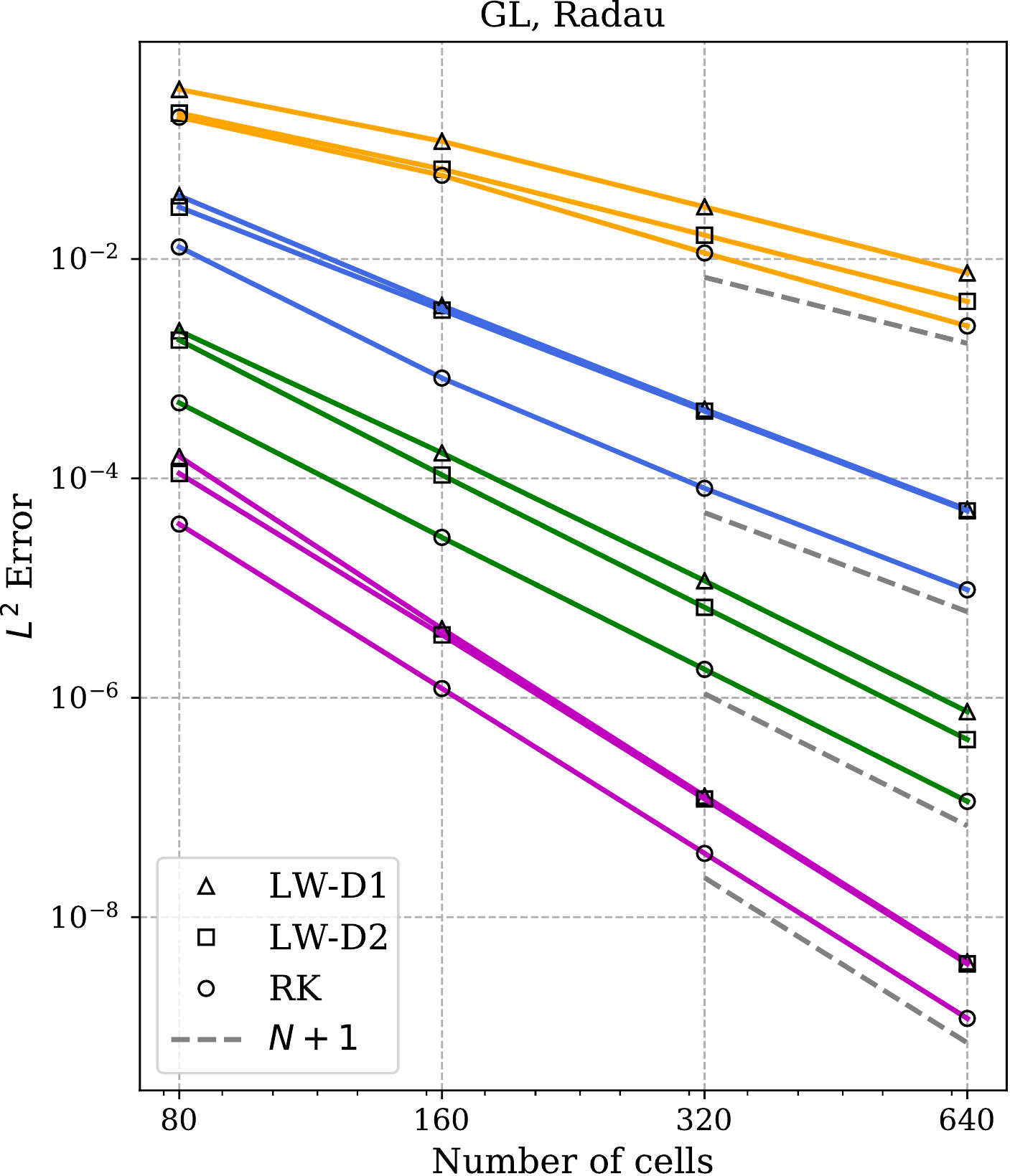} &
\includegraphics[width=0.40\textwidth]{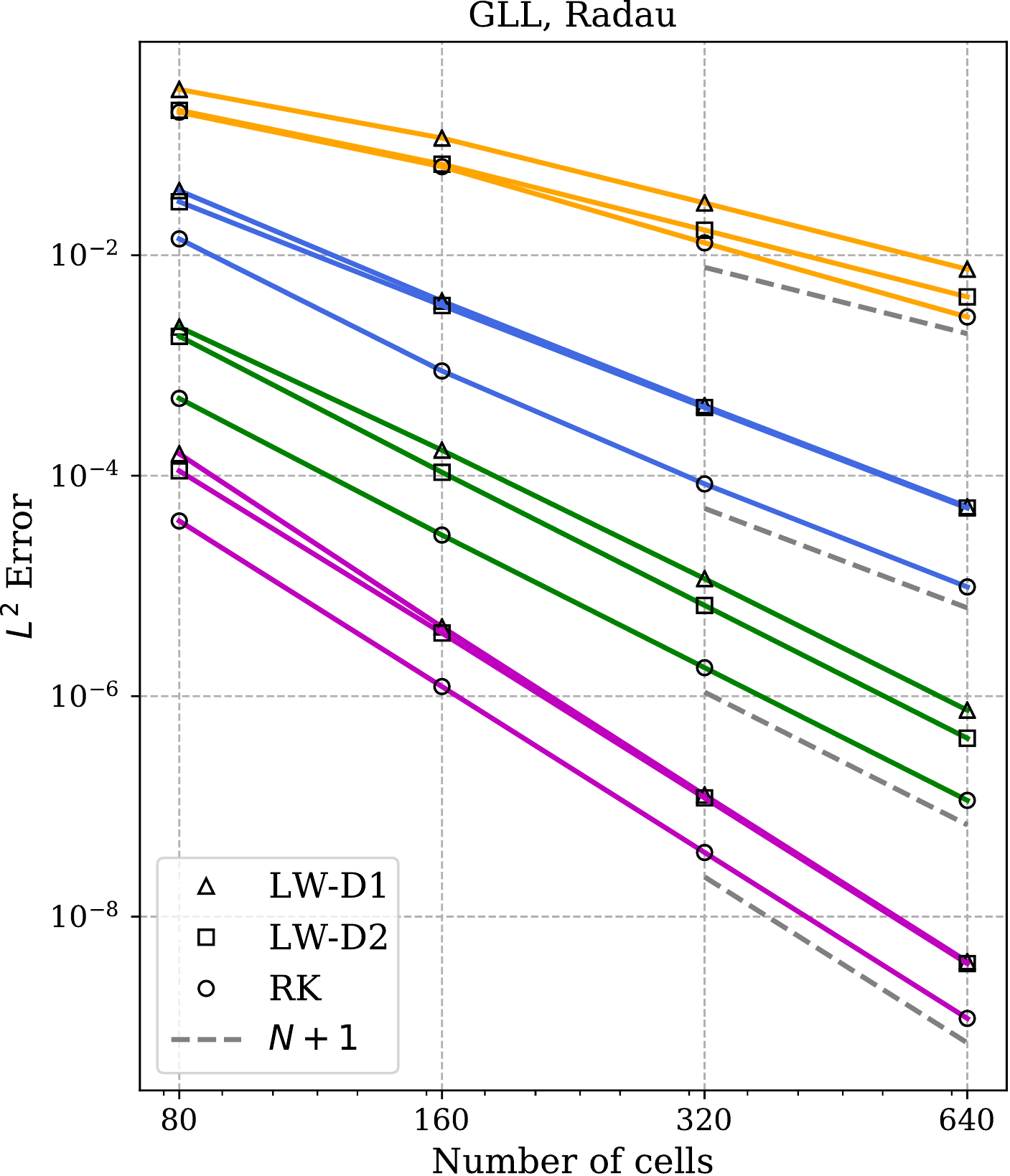} \\
(a) GL points & (b) GLL points
\end{tabular}
\end{center}
\caption{Error convergence for constant linear advection of a wave packet; (a) GL points, (b) GLL points. The different colors correspond to degrees $N=1,2,3,4$ from top to bottom.}
\label{fig:wp2}
\end{figure}

Next we perform error convergence studies for an initial condition of a wave packet given by $u(x,0) = e^{-10x^2}\sin(10\pi x)$ with periodic boundary conditions. This initial condition has a more broadband Fourier spectrum than the previous case which had only one Fourier mode. Figure~(\ref{fig:wp1}) shows the solutions obtained for $N=3,4$ and using 200 dofs in each case. The solutions are more accurate in case of $N=4$ compared to $N=3$ showing the benefits of a higher order method. We  see that RK schemes are able to capture the peak solution more accurately than LW schemes, especially in case of $N=3$, but the difference between the two schemes reduces for $N=4$ case. Figure~(\ref{fig:wp2}) shows the error convergence plot with GL points; as before, we see that RK schemes show smaller errors than the LW schemes due to their super-convergence property. For odd degrees, the D2 dissipation has slightly smaller errors than the D1 model, while for even degrees, the difference between the two models is negligible.

\subsubsection{Non-smooth solutions}\label{sec:non-smooth1d}
\begin{figure}
\centering
\begin{tabular}{cc}
\includegraphics[width=0.45\textwidth]{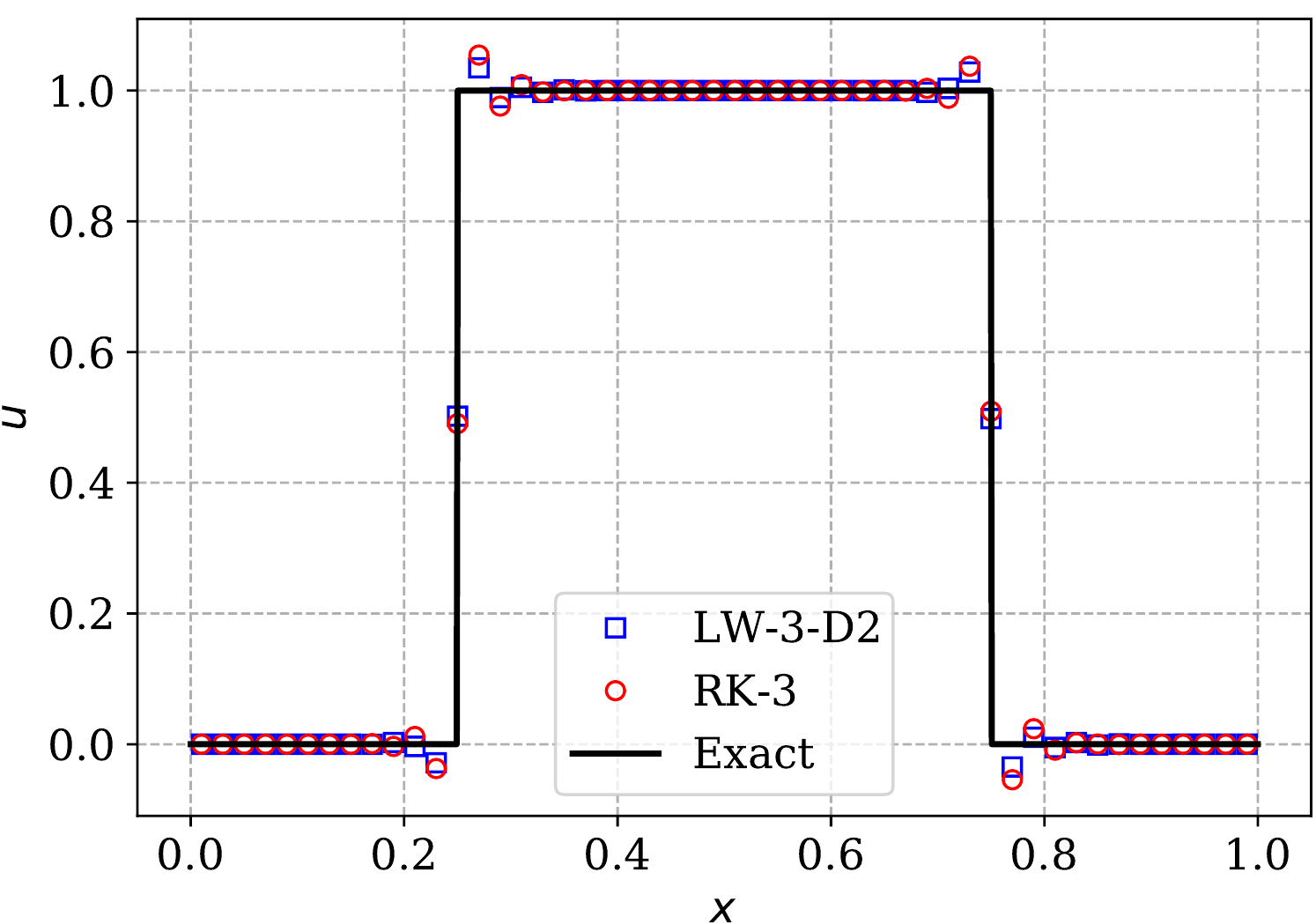} &
\includegraphics[width=0.45\textwidth]{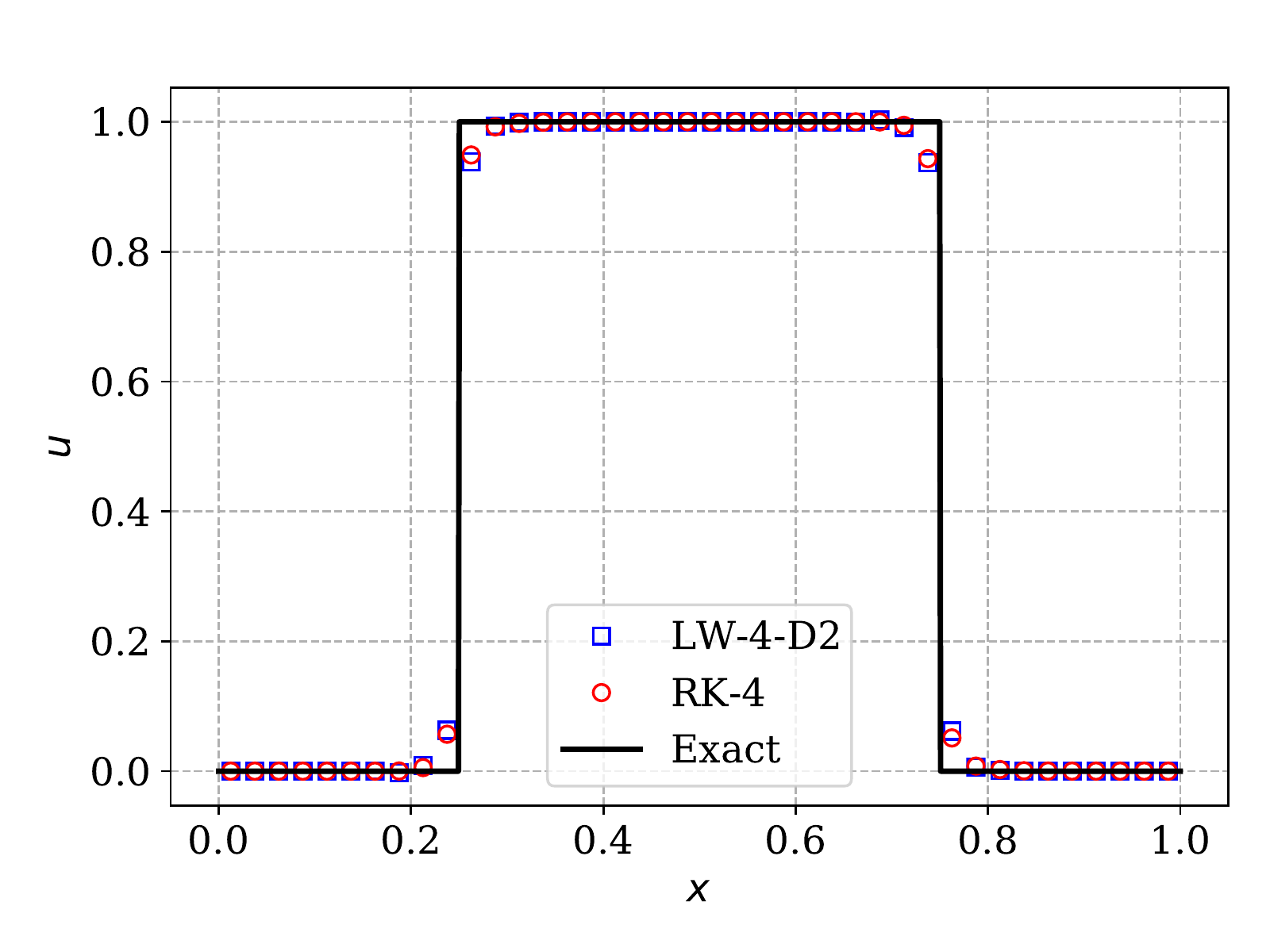} \\
(a) $N=3$ & (b) $N=4$
\end{tabular}
\caption{Constant linear advection of hat profile without limiter. The solution is shown at time $t=1$ with 200 dofs using polynomial degree (a) $N=3$, (b) $N=4$.}
\label{fig:hat1}
\end{figure}

\begin{figure}
\centering
\begin{tabular}{cc}
\includegraphics[width=0.45\textwidth]{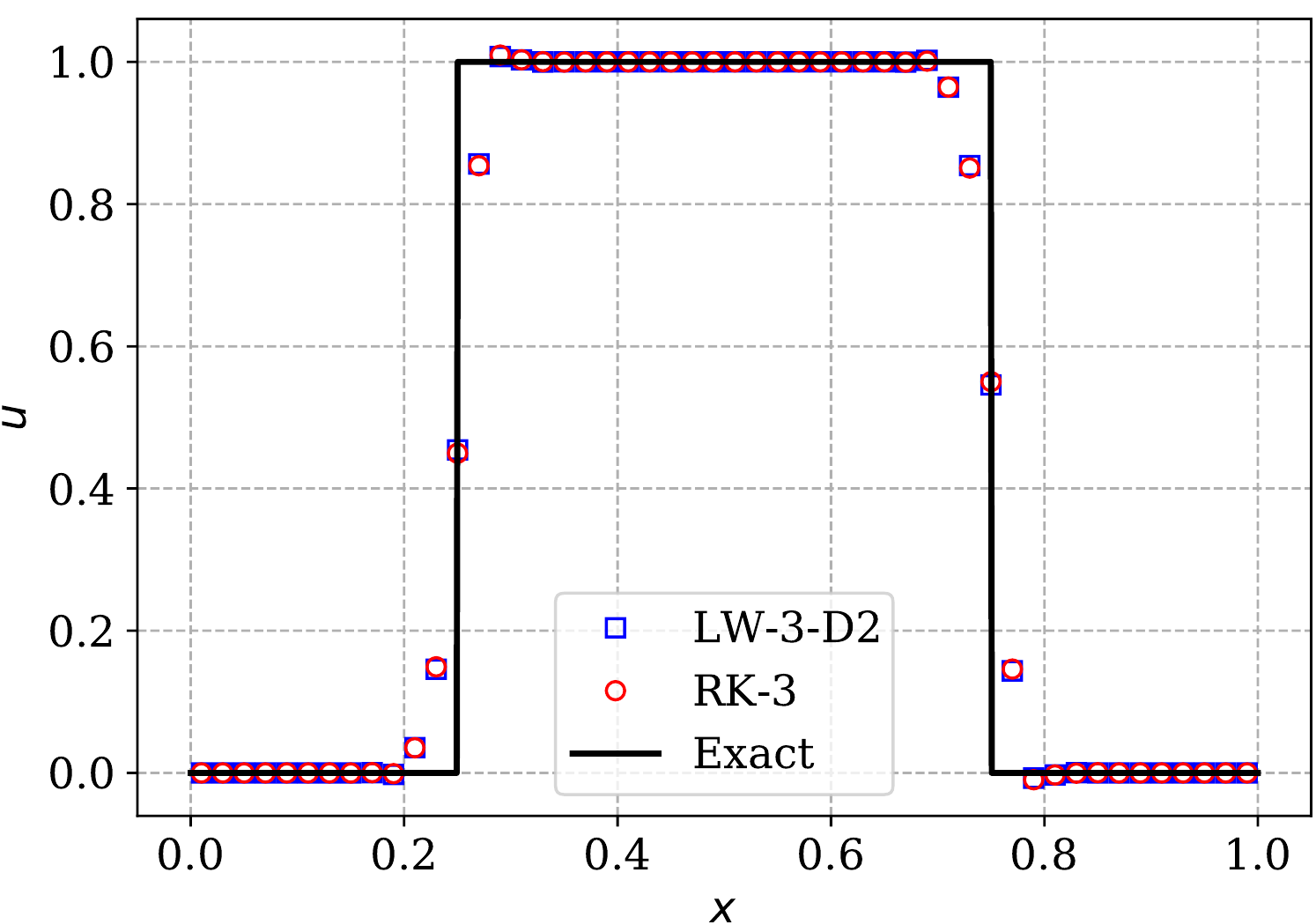} &
\includegraphics[width=0.45\textwidth]{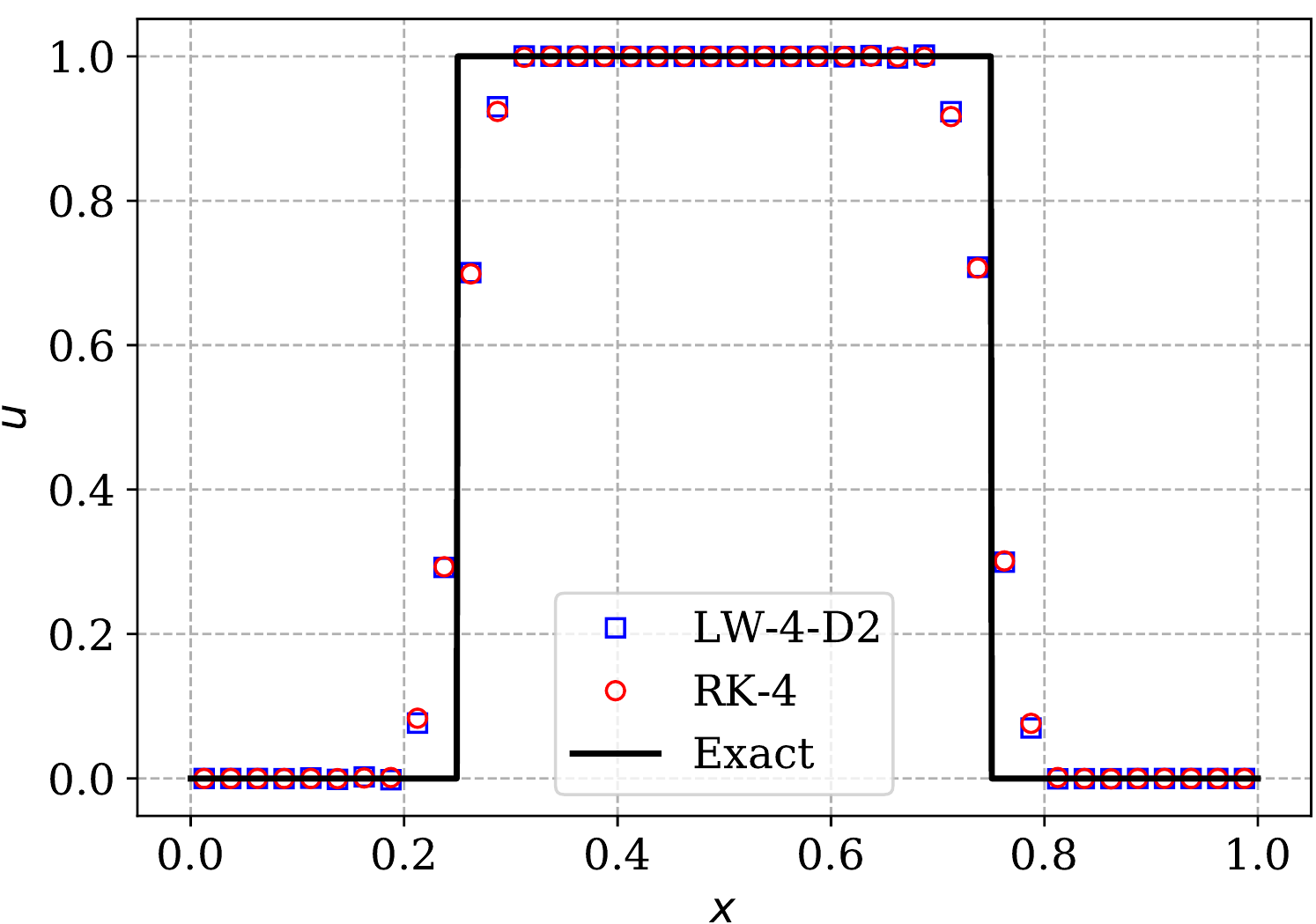} \\
(a) $N=3$ & (b) $N=4$
\end{tabular}
\caption{Constant linear advection of hat profile with TVB limiter ($M=100$). The solution is shown at time $t=1$ and 200 dofs using polynomial degree (a) $N=3$, (b) $N=4$.}
\label{fig:hat2}
\end{figure}

\begin{figure}
\begin{center}
\includegraphics[width=0.45\textwidth]{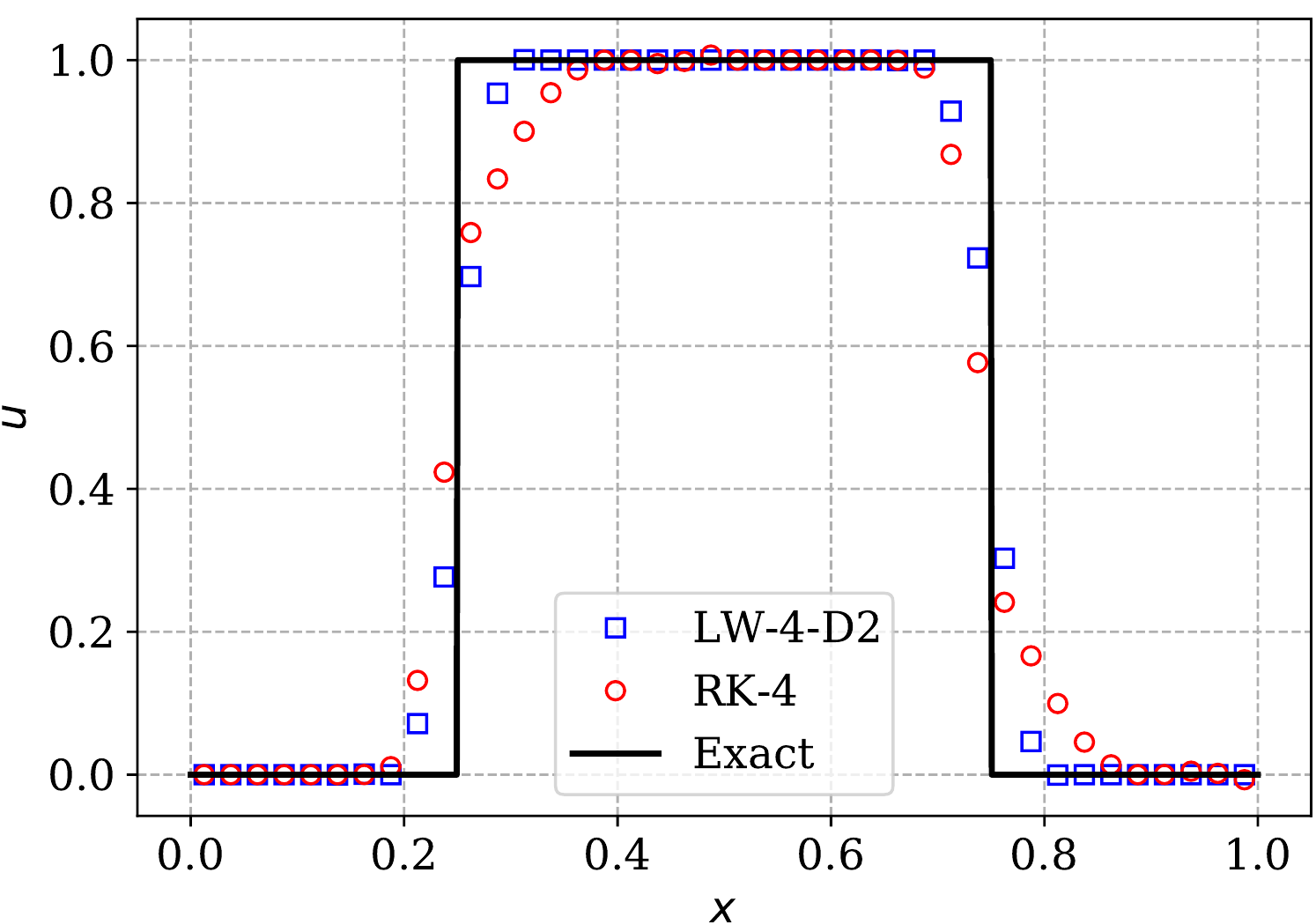}
\end{center}
\caption{Constant linear advection of hat profile with TVB limiter ($M=100$) where Runge-Kutta time integration is performed using RK65~\cite{Tsitouras2011}. The solution is shown at time $t=1$ and 200 dofs using polynomial degree $N=4$.}
\label{fig:hat3}
\end{figure}

If the initial condition is not smooth and has a jump discontinuity, then high order methods will generate oscillatory solutions due to the non-monotone property of the schemes as shown by Godunov. For such problems, we need some form of limiter to control the oscillations and we use the TVB-type limiters which are applied in a posteriori manner as explained in Section~\ref{sec:lim}. Consider the initial condition consisting of a square hat function,
\[
u(x,0) = \begin{cases}
1, & x \in (0.25, 0.75) \\
0, & x \in [0,0.25) \cup (0.75, 1]
\end{cases}
\]
and which is extended by periodicity. We compute the solution upto the time $t=1$ unit when the solution returns to its initial position. Figure~(\ref{fig:hat1}) shows the solutions obtained with degree $N=3,4$ and without applying any limiter. We observe oscillations in case of $N=3$ but no significant oscillations are seen for the $N=4$ case. The oscillations are however localized around the discontinuity and do not corrupt the rest of the solution. When TVB limiter is applied, these oscillations disappear as seen in Figure~(\ref{fig:hat2}) but the jumps are smeared over more cells. If we use the RK65 scheme which is fifth order accurate but has six stages, then the results are shown in Figure~(\ref{fig:hat3}) where we observe increased smearing of the jump in the RK scheme. Overall, we see that the limiter smears the discontinuity over a few cells in case of both  LW and RK schemes; but we also observe that the solutions obtained with the LW schemes are very similar to the RK schemes.

\begin{figure}
\begin{center}
\begin{tabular}{cc}
\includegraphics[width=0.45\textwidth]{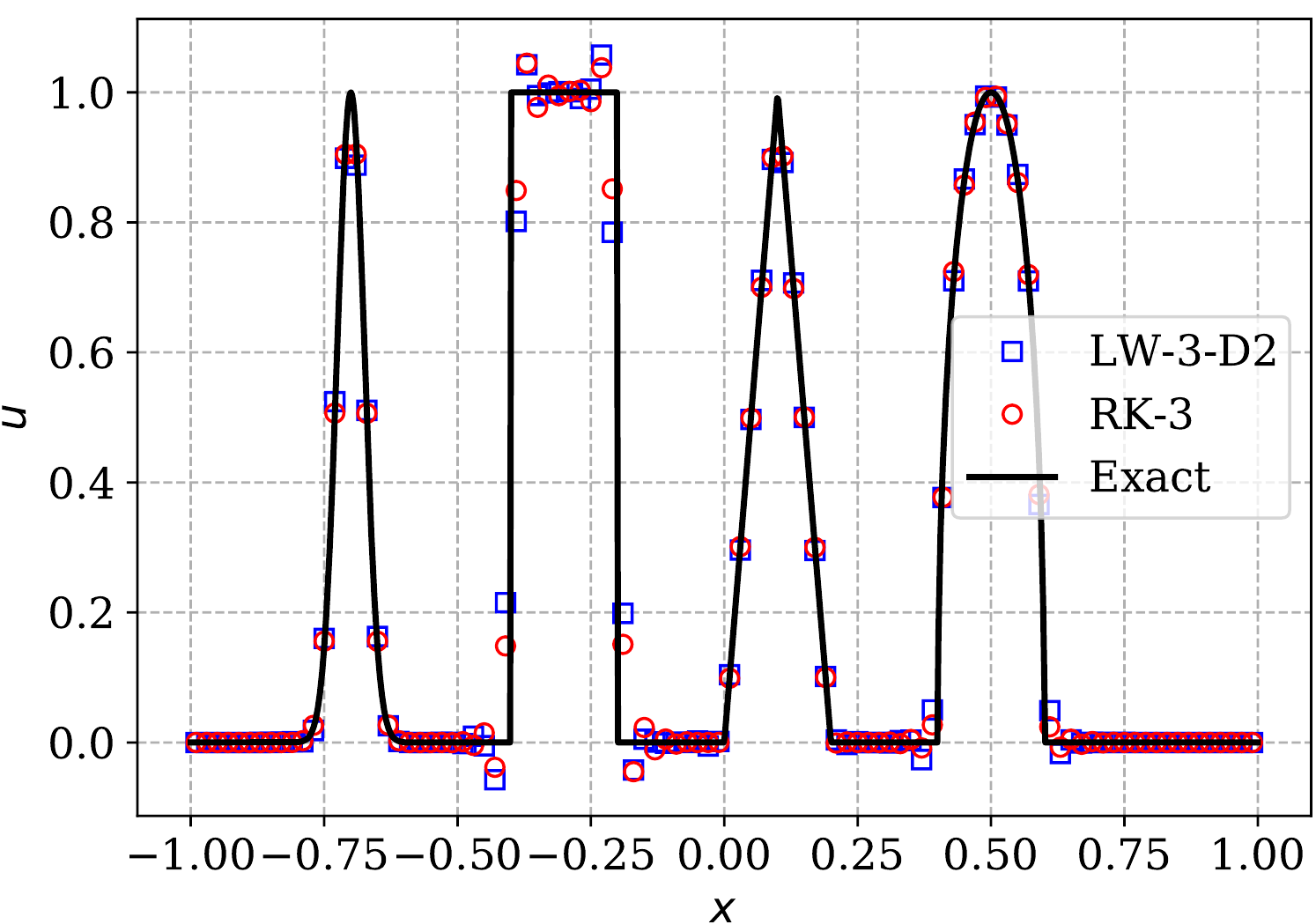} &
\includegraphics[width=0.45\textwidth]{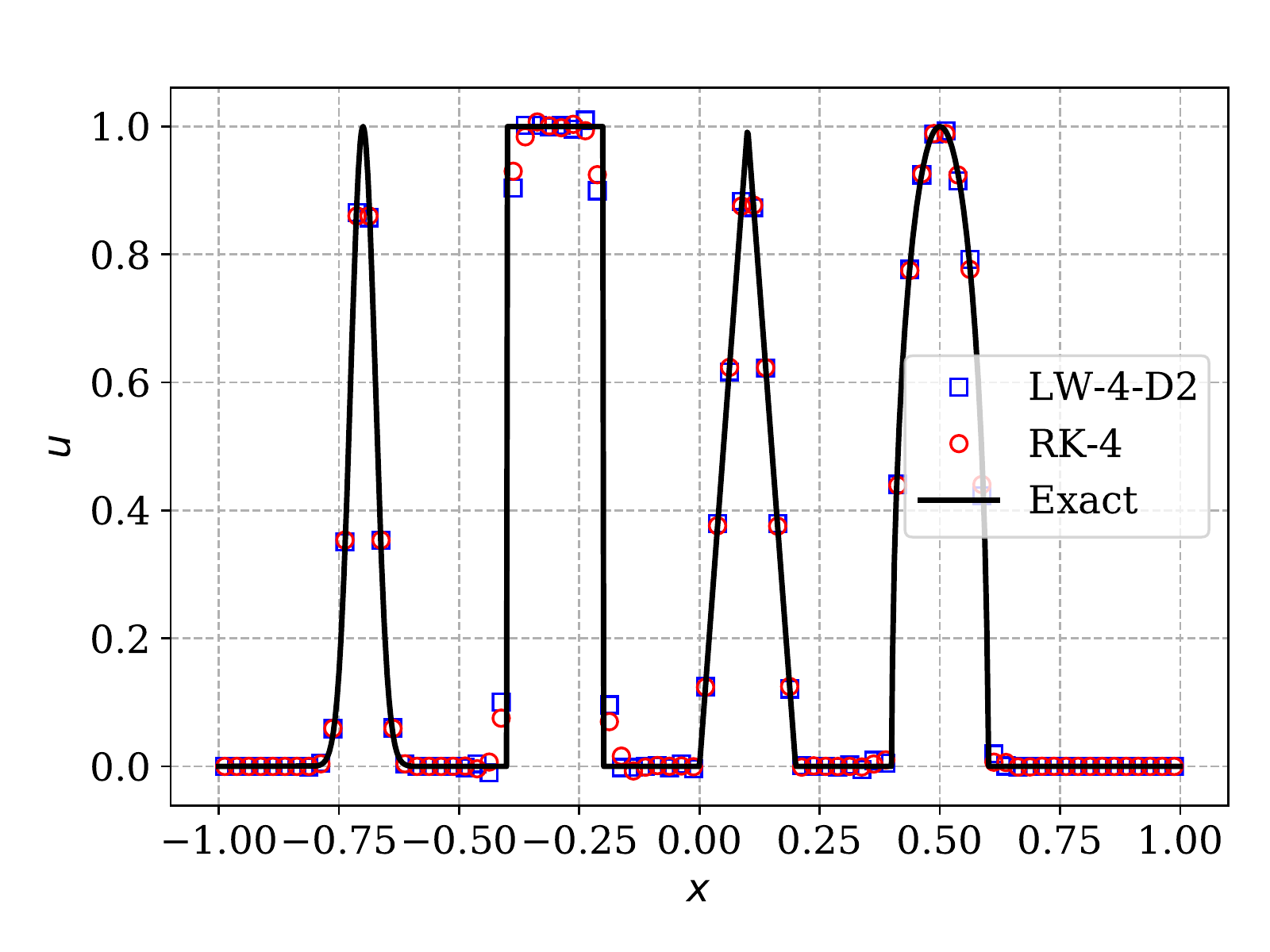} \\
(a) $N=3$ & (b) $N=4$
\end{tabular}
\end{center}
\caption{Constant linear advection of a composite profile without limiter. The solution is shown at time $t=8$ using 400 dofs in each case and polynomial degree (a) $N=3$, (b) $N=4$.}
\label{fig:mult1}
\end{figure}

\begin{figure}
\begin{center}
\begin{tabular}{cc}
\includegraphics[width=0.45\textwidth]{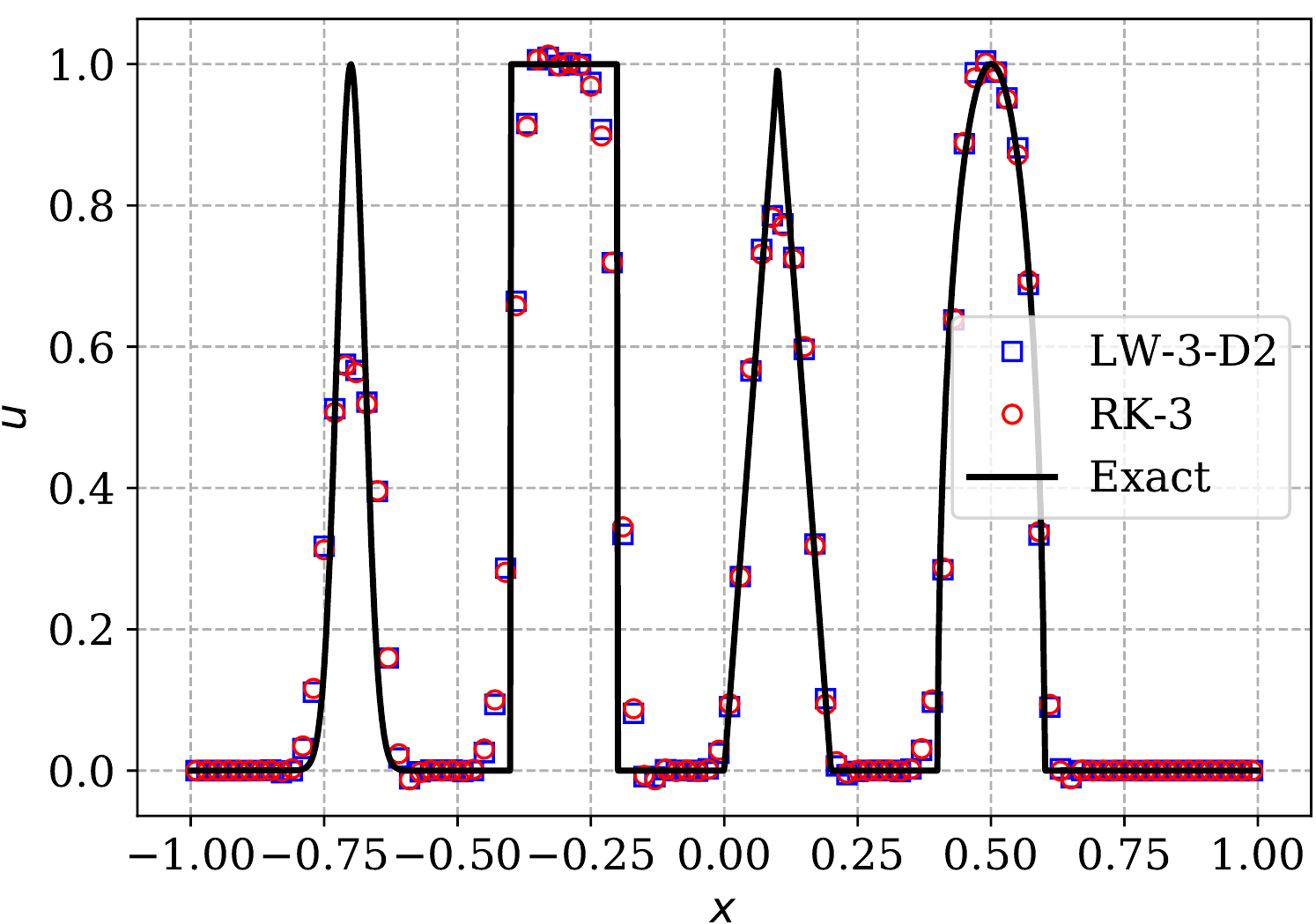} &
\includegraphics[width=0.45\textwidth]{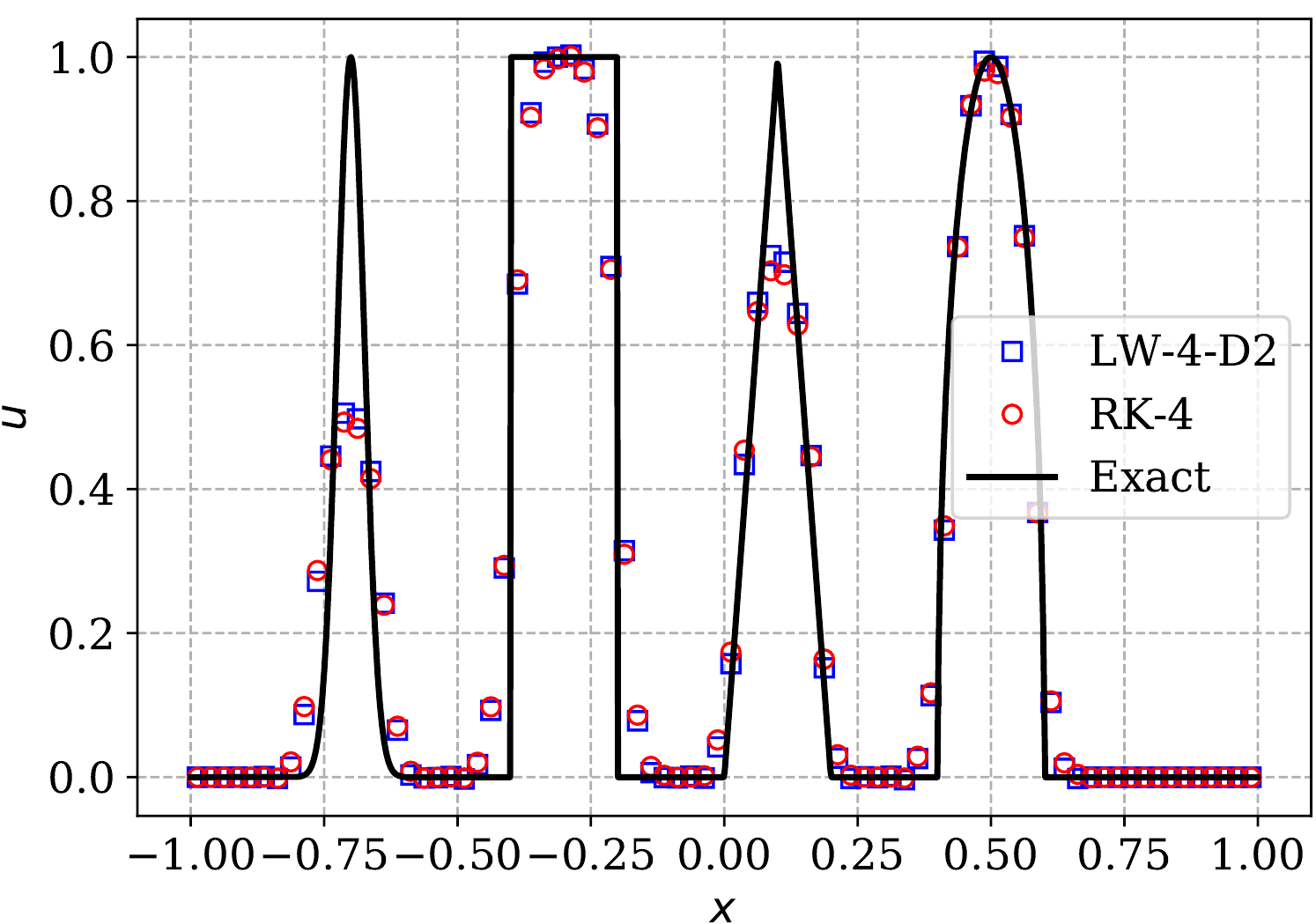} \\
(a) $N=3$ & (b) $N=4$
\end{tabular}
\end{center}
\caption{Constant linear advection of a composite profile with TVB limiter ($M=50$). The solution is shown at time $t=8$ using 400 dofs in each case and polynomial degree (a) $N=3$, (b) $N=4$.}
\label{fig:mult2}
\end{figure}

We next consider a composite signal consisting of profiles with different levels of smoothness whose initial condition is a slightly different version of~\cite{Jiang1996} given by
\[
u(x,0)=\begin{cases}
G(x,\beta,z), & -0.8\le x \le -0.6,\\
1, & -0.4\le x \le -0.2\\
1 - |10(x-0.1)|, & 0\le x\le 0.2\\
F(x,\alpha,a),& 0.4\le x \le 0.6\\
0,&  \textrm{elsewhere}
\end{cases}
\]
where $G(x,\beta,z)=\mathrm e^{-\beta(x-z)^2}$, $F(x,\alpha,a)=\sqrt{1-\alpha^2(x-a)^2}$ with the constants $a=0.5$, $z=-0.7$, $\delta=0.005$,  $\alpha=10$ and $\beta=\frac{\log 2}{36\delta^2}$.  This initial condition is composed of the succession of a Gaussian, rectangular, triangular and parabolic signals. We compute the numerical solutions at $t=8$ (after $4$ periods) and for degrees $N=3,4$ but with 400 dofs in total for each case. The results without any limiter are shown in Figure~(\ref{fig:mult1}); the profiles which are more regular are captured accurately by the numerical schemes, while the hat profile shows some oscillations. These oscillations are larger for $N=3$ than for $N=4$ case. When the TVB limiter is used, the corresponding solutions are shown in Figure~(\ref{fig:mult2}). Now the oscillations in the hat profile are controlled but there is more numerical dissipation as is evident in the reduced amplitude of the smooth profiles. We observe that the results from the LW scheme are very similar to those of the RK scheme.

\subsection{Linear advection equation: variable speed \label{sec: vla}}
Now we consider the linear advection equation with spatially varying speed which is given by
\[
u_t + f(x,u)_x = 0, \qquad f(x,u) = a(x) u
\]
This problem is non-linear in the spatial variable, i.e., if $I_h$ is the interpolation operator, then $I_h(a u_h) \ne I_h(a)I_h(u_h)$. This can lead to different behaviour of the numerical schemes compared to the linear case, depending on AE and EA methods for the numerical flux. To study the effect of non-linearity, we consider different types of speeds with different degree of non-linearity from \cite{Offner2019}.

\begin{figure}
\begin{center}
\begin{tabular}{ccc}
\includegraphics[width=0.30\textwidth]{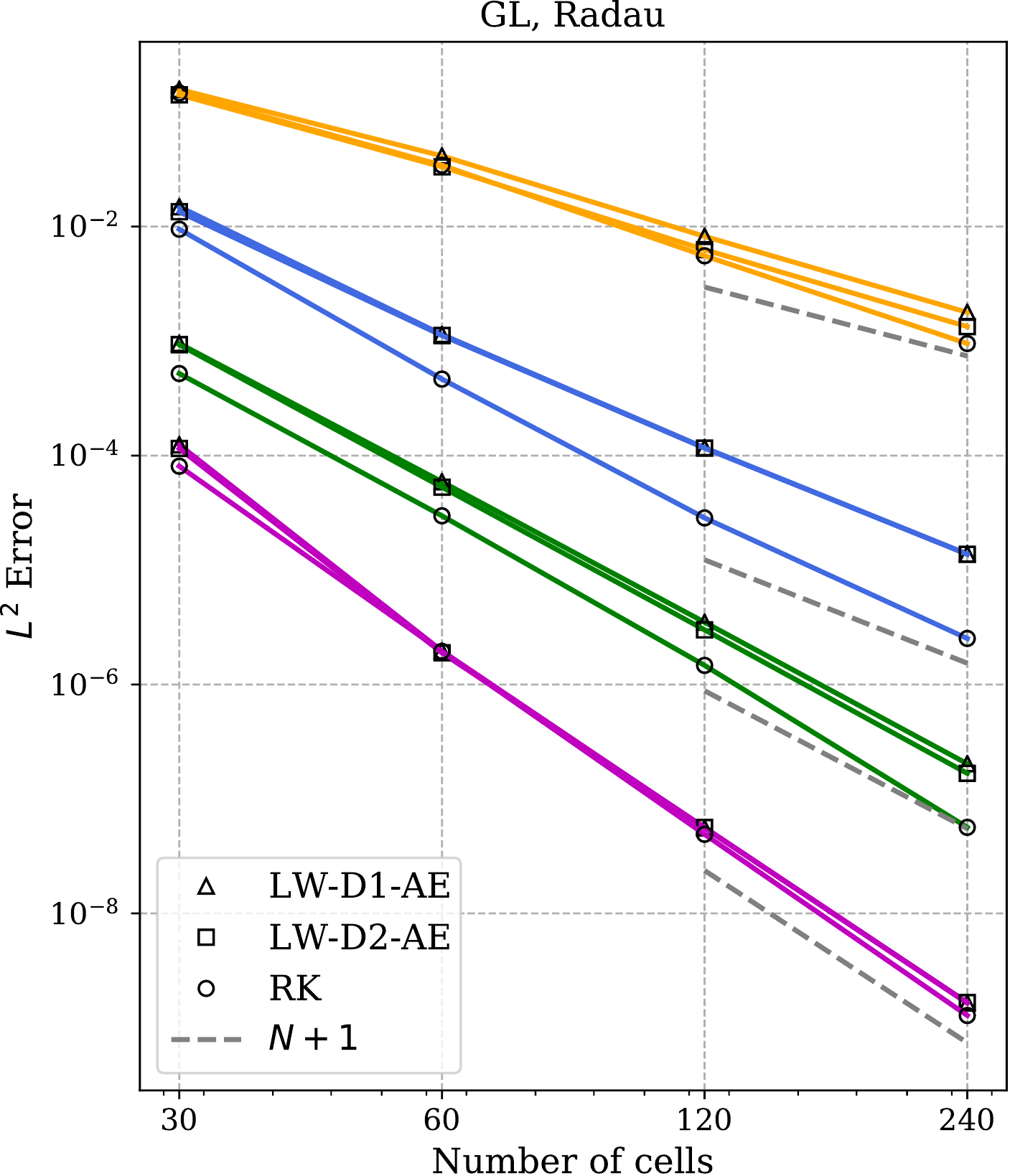} &
\includegraphics[width=0.30\textwidth]{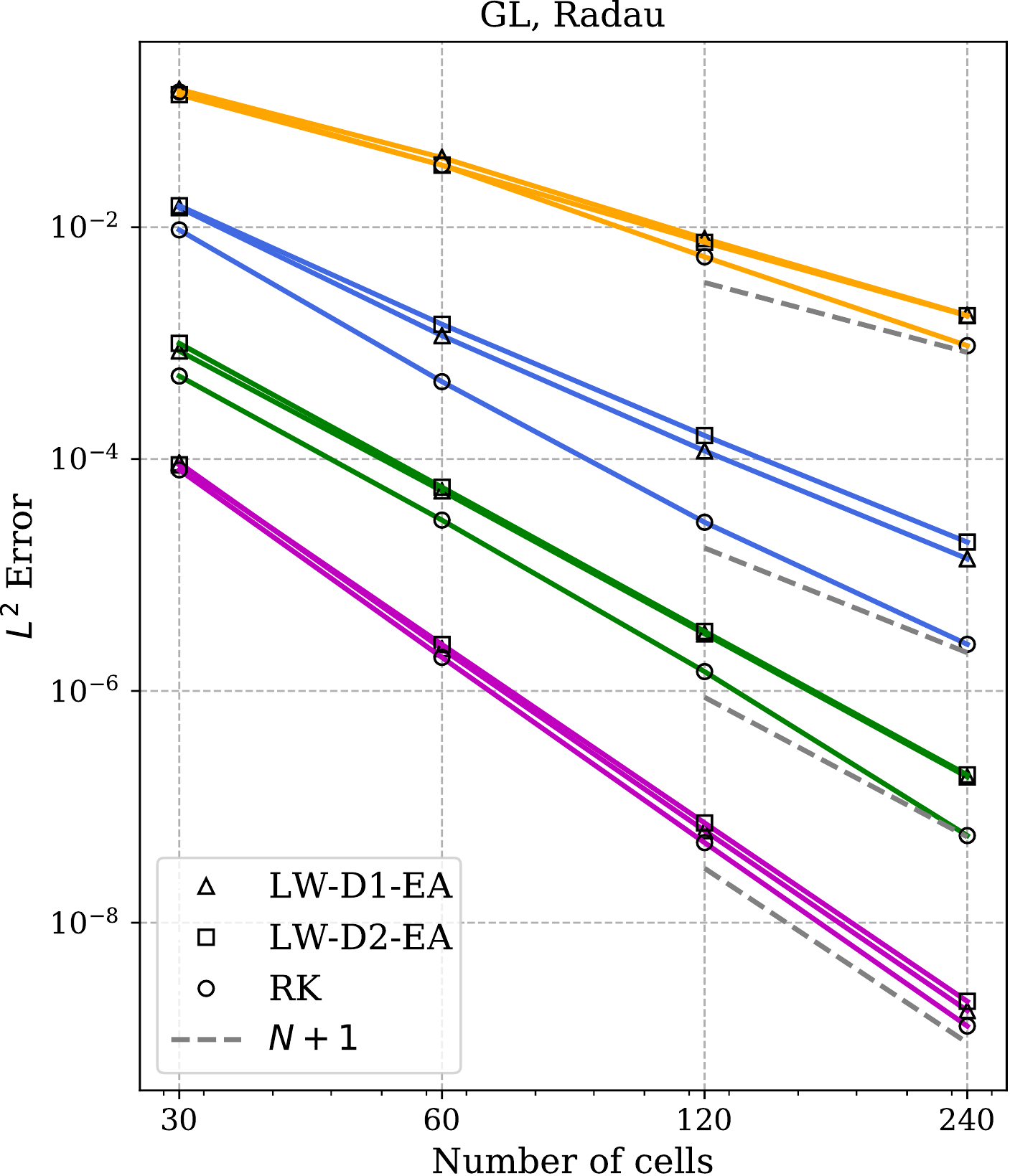} &
\includegraphics[width=0.30\textwidth]{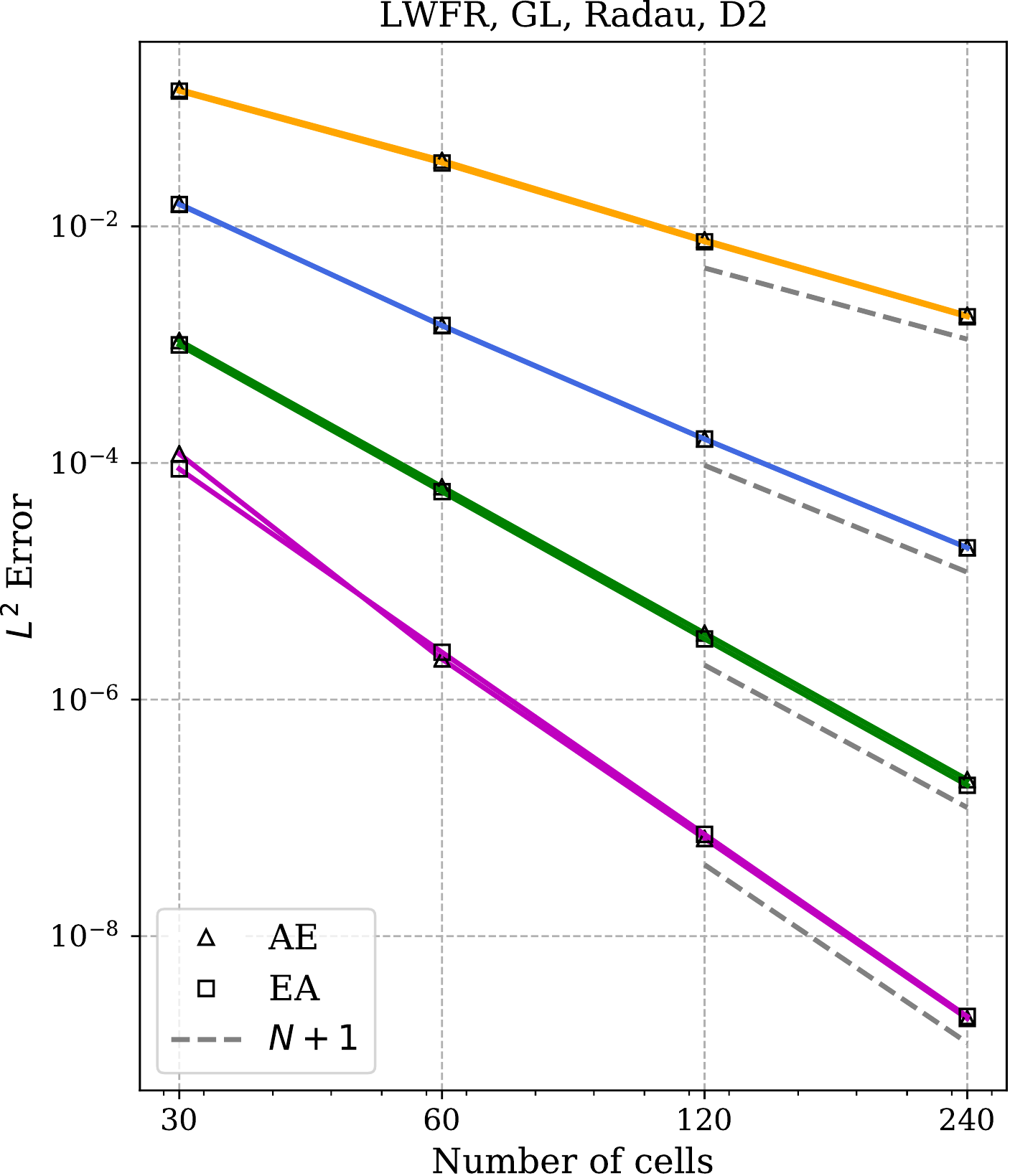} \\
(a) & (b) & (c)
\end{tabular}
\end{center}
\caption{Error convergence for variable linear advection with $a(x)=x$: (a) AE scheme, (b) EA scheme, (c) AE vs EA.}
\label{fig:vla1}
\end{figure}

Figure~(\ref{fig:vla1}) shows the error convergence for the AE and EA schemes, and for the speed $a(x) = x$ with initial condition $u_0(x)=\sin(12(x-0.1))$. The domain is $[0.1,2\pi]$ and we use Dirichlet boundary conditions at $x=0.1$ and outflow condition at $x=2\pi$ so that the exact solution is given by $u(x,t)=e^{-t}u_0(xe^{-t})$. As mentioned earlier, upwind flux is used to enforce the boundary condition at inflow boundaries. The LW scheme with either AE or EA method yields correct convergence rates, while the RK scheme exhibits a small super-convergence. The right most figure shows that the error levels with AE and EA are nearly same. The non-linearity in this problem is small enough that it does not spoil the error and convergence behaviour of the LW schemes, for both AE and EA methods.

\begin{figure}
\begin{center}
\begin{tabular}{ccc}
\includegraphics[width=0.30\textwidth]{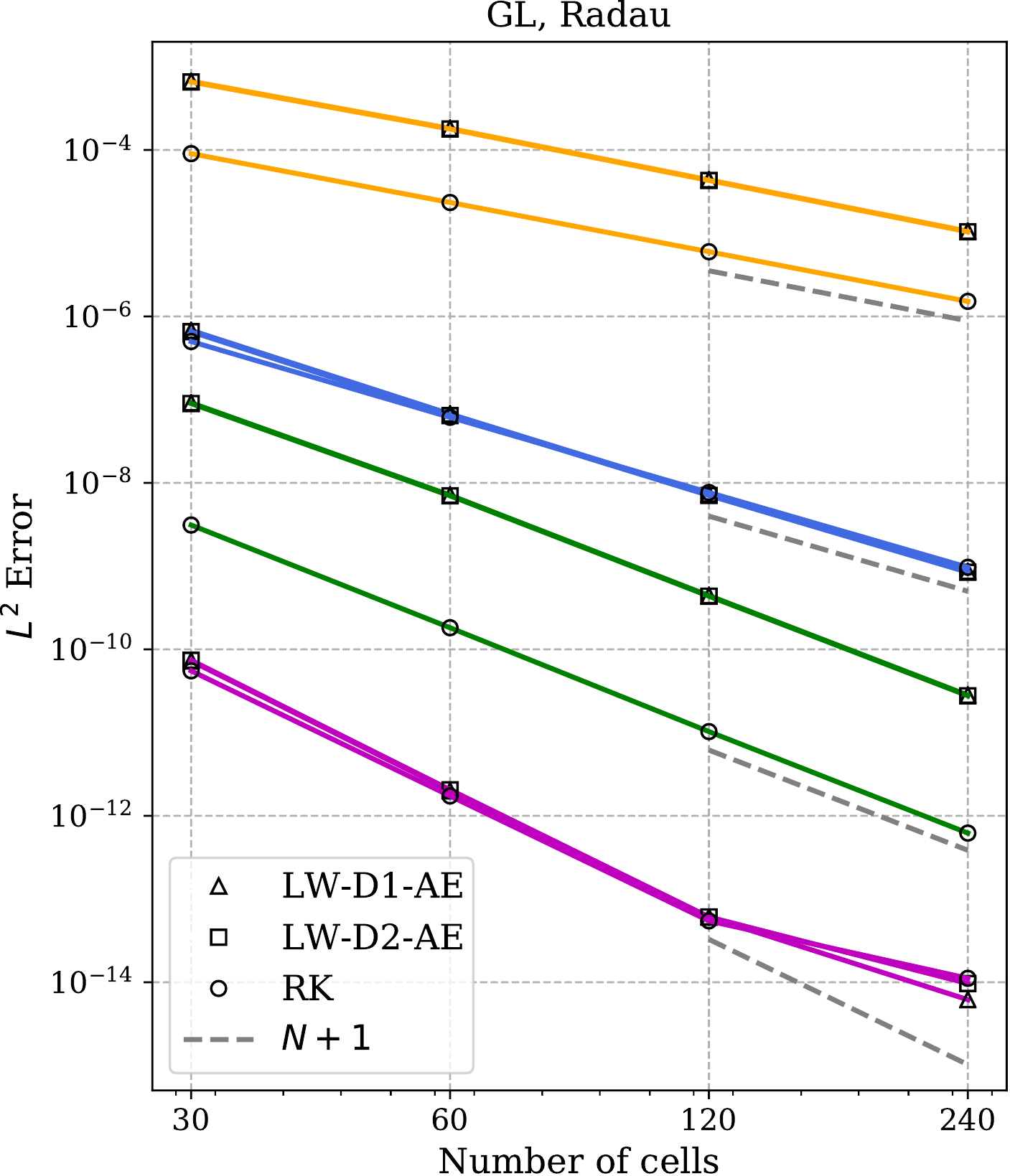} &
\includegraphics[width=0.30\textwidth]{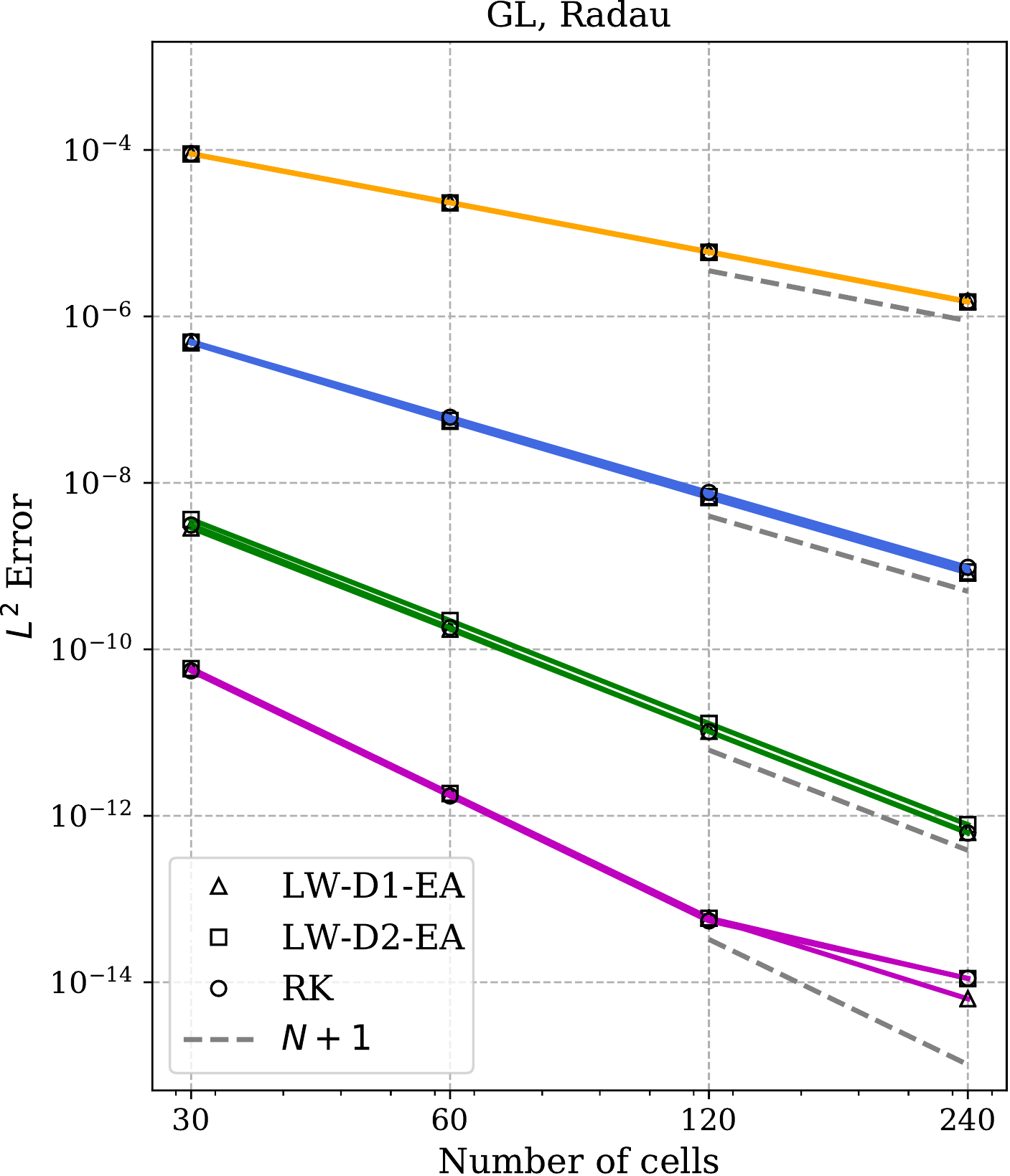} &
\includegraphics[width=0.30\textwidth]{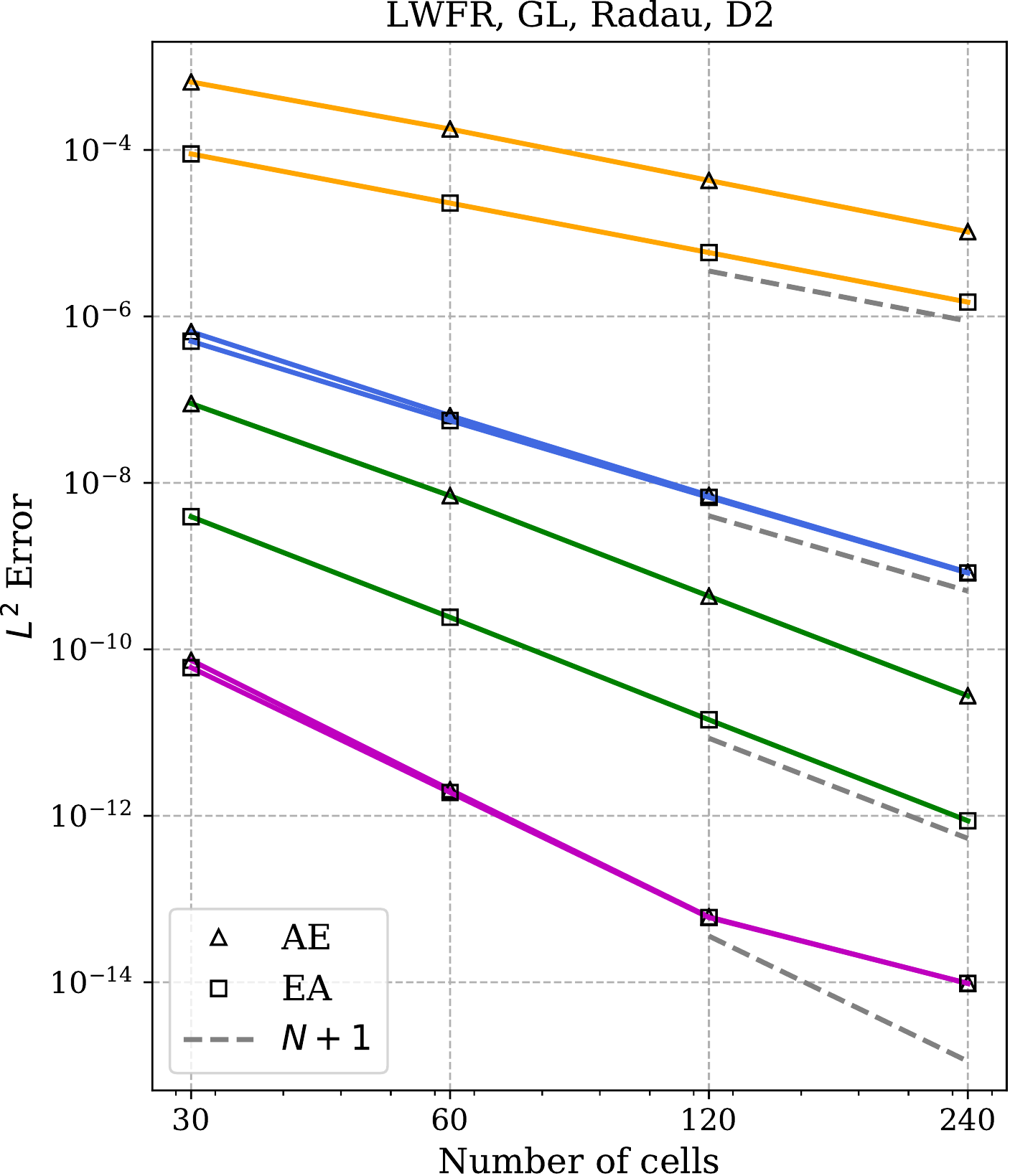} \\
(a) & (b) & (c)
\end{tabular}
\end{center}
\caption{Error convergence for variable linear advection with $a(x)=x^2$: (a) AE scheme, (b) EA scheme, (c) AE vs EA.}
\label{fig:vla2}
\end{figure}

\begin{figure}
\begin{center}
\begin{tabular}{cc}
\includegraphics[width=0.46\textwidth]{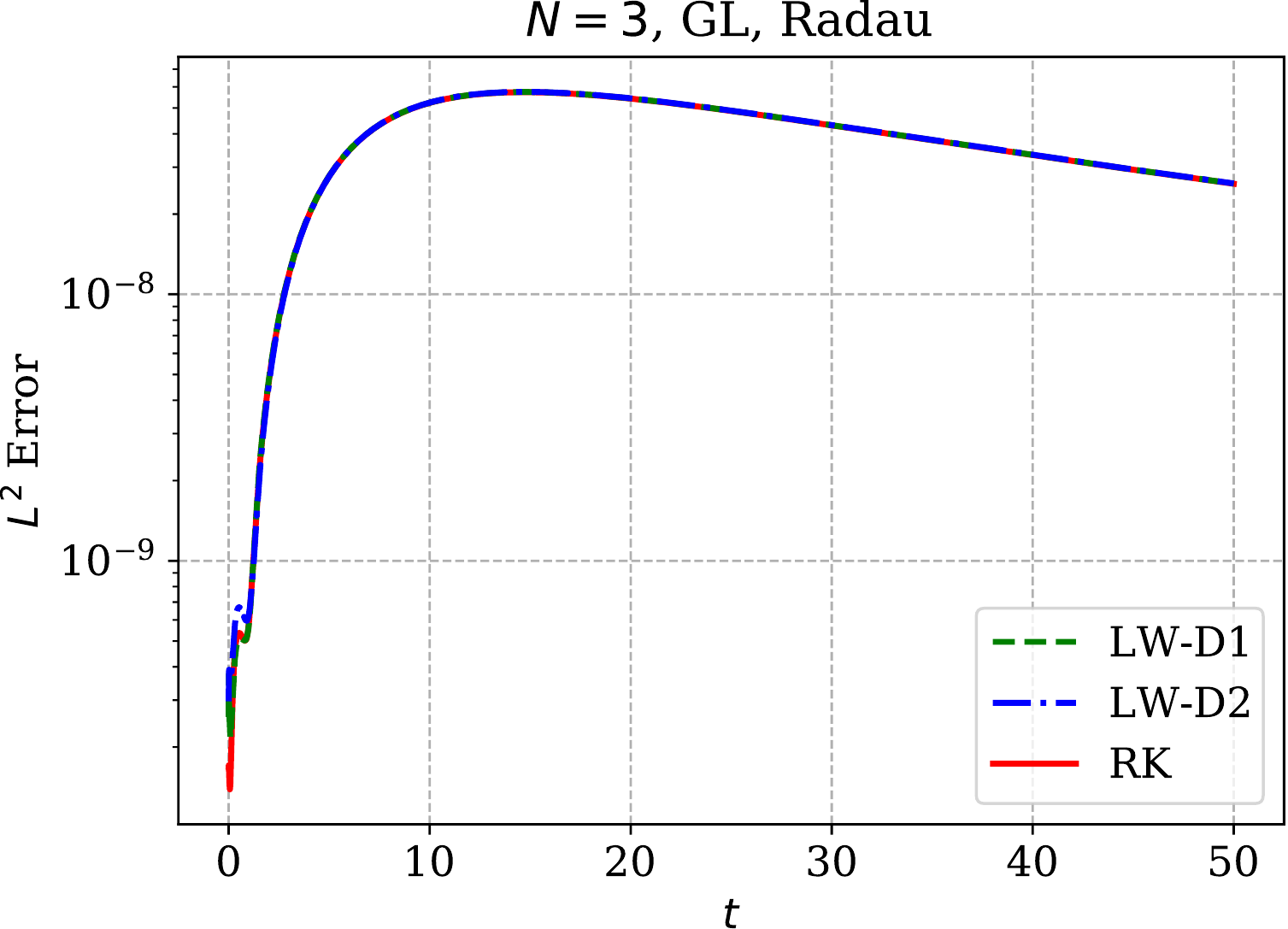} &
\includegraphics[width=0.46\textwidth]{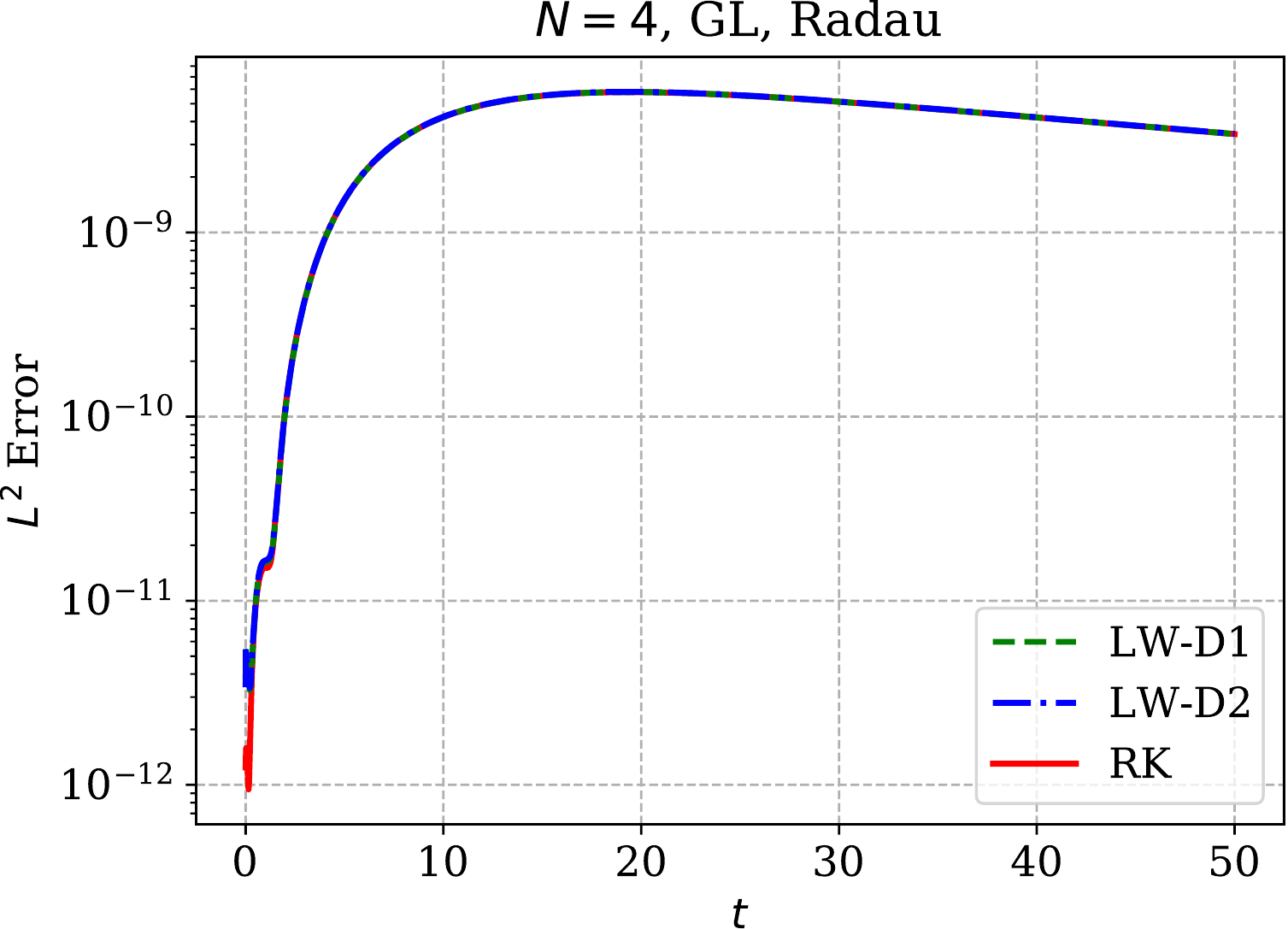} \\
(a) $N=3$ & (b) $N=4$
\end{tabular}
\end{center}
\caption{Error versus time for linear advection with wave speed $a(x)=x^2$ for different polynomial degrees; GL solution points, Radau correction and polynomial degree (a) $N=3$, (b) $N=4$.}
\label{fig:vla3}
\end{figure}

Figure~(\ref{fig:vla2}) shows the error convergence for the AE and EA schemes, and for the non-linear speed $a(x) = x^2$, with initial condition $u_0(x)=\cos(\pi x/2)$. The domain is $[0.1,1]$, and we use Dirichlet boundary conditions at $x=0.1$ which is an inflow boundary, and outflow condition at $x=1$ so that the exact solution is given by $u(x,t)=u_0(x/(1+tx))/{(1+tx)}^2$. For odd degrees, the LW scheme with AE shows larger errors compared to the RK scheme though the convergence rate is optimal. The LW scheme with EA shown in the middle figure, is as accurate as the RK scheme at all degrees. The last figure compares AE and EA schemes using GL solution points, Radau correction function and D2 dissipation; we clearly see that EA scheme has smaller errors than AE scheme at odd degrees, while they are very similar for even degrees.  Figure~(\ref{fig:vla3}) shows the error versus time plots for degrees $N=3,4$; we see that the LW and RK schemes have very similar error levels and the superior performance of RK schemes observed for constant linear advection is not realized in this non-linear case.

We have observed the same behaviour in all other non-linear test cases given in~\cite{Offner2019} but the results are not shown here, i.e., the LW schemes with EA perform at par with RK schemes for non-linear problems.

\subsection{Inviscid Burger's equation}
\begin{figure}
\begin{center}
\begin{tabular}{ccc}
\includegraphics[width=0.30\textwidth]{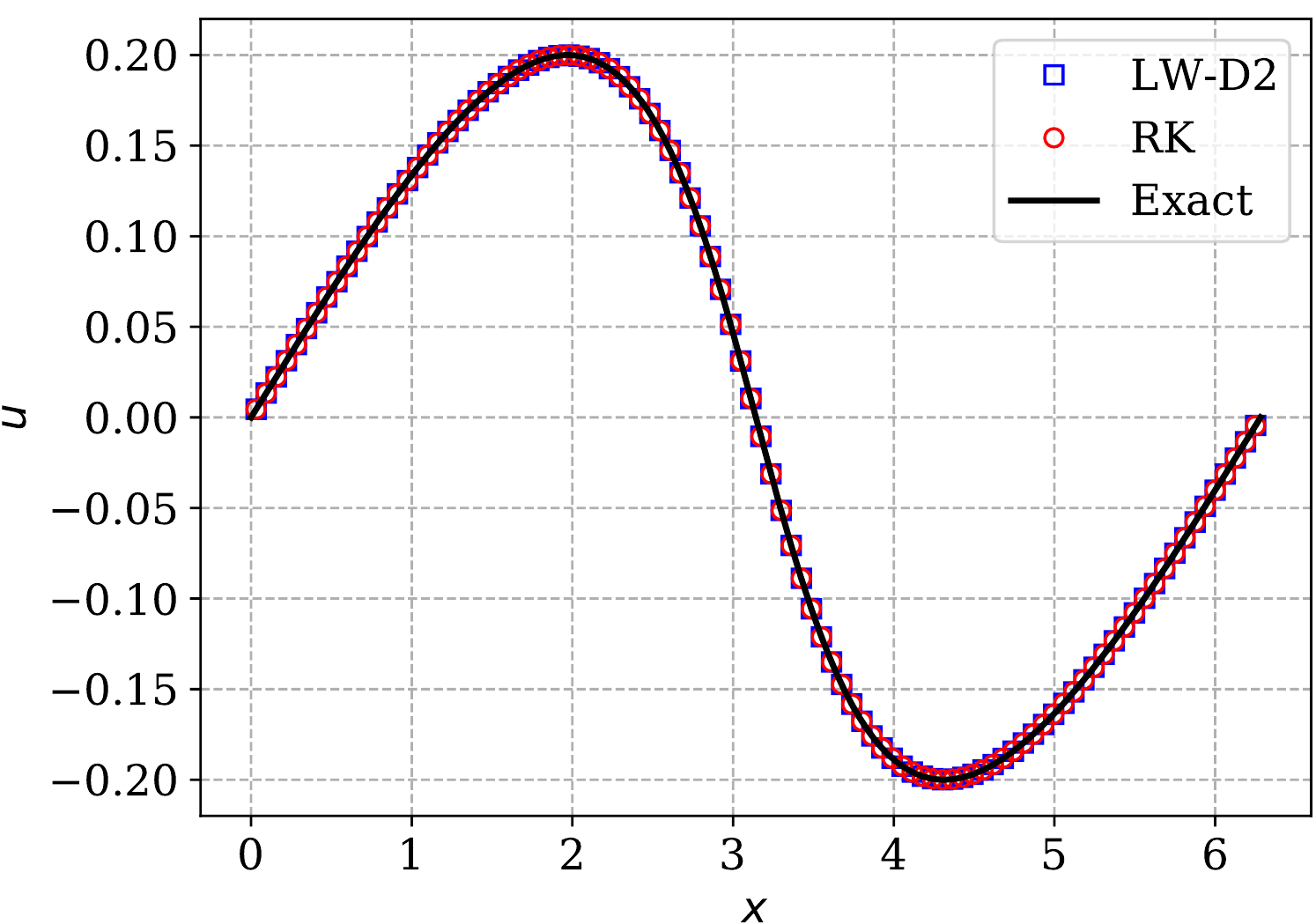} &
\includegraphics[width=0.30\textwidth]{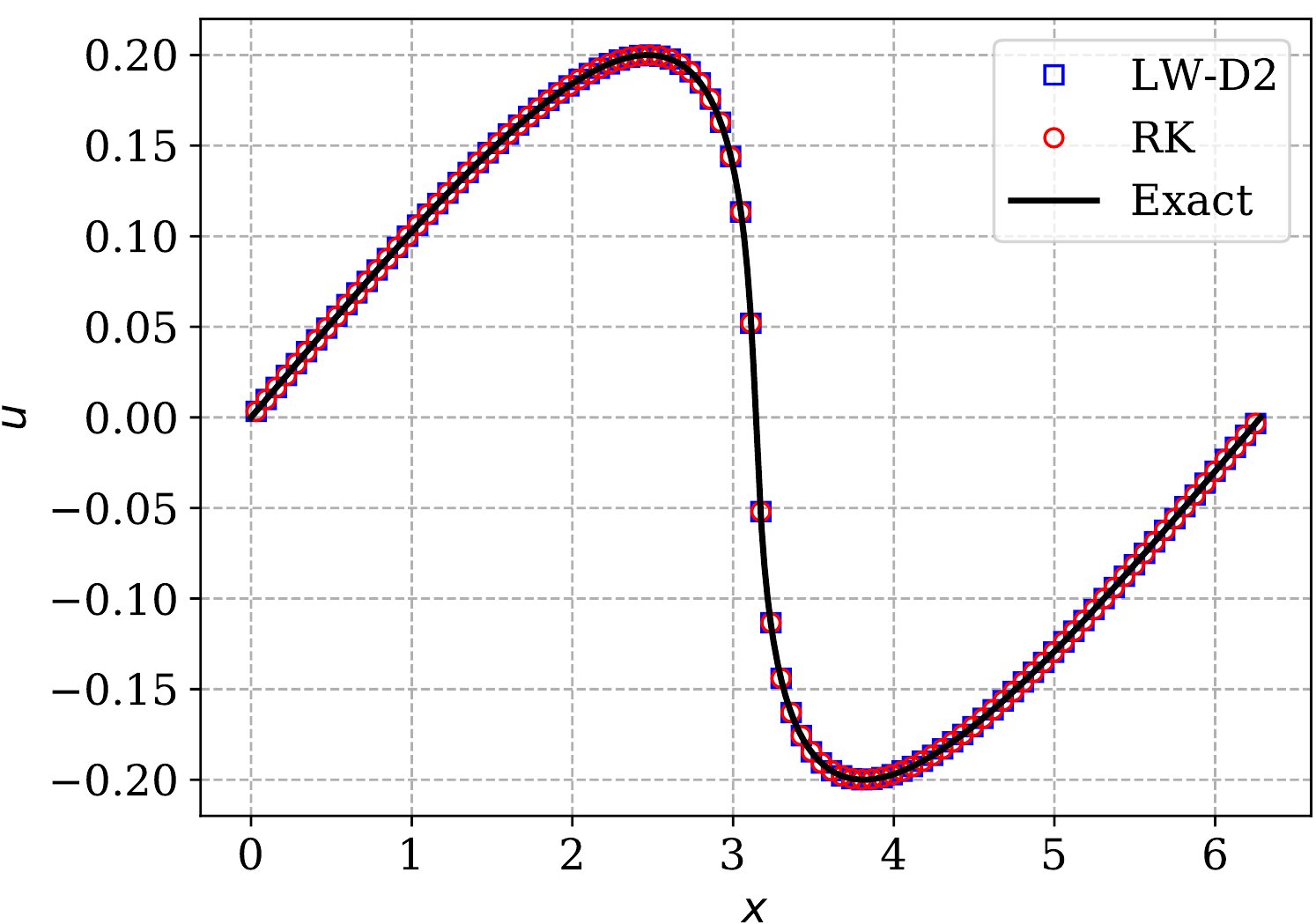} &
\includegraphics[width=0.30\textwidth]{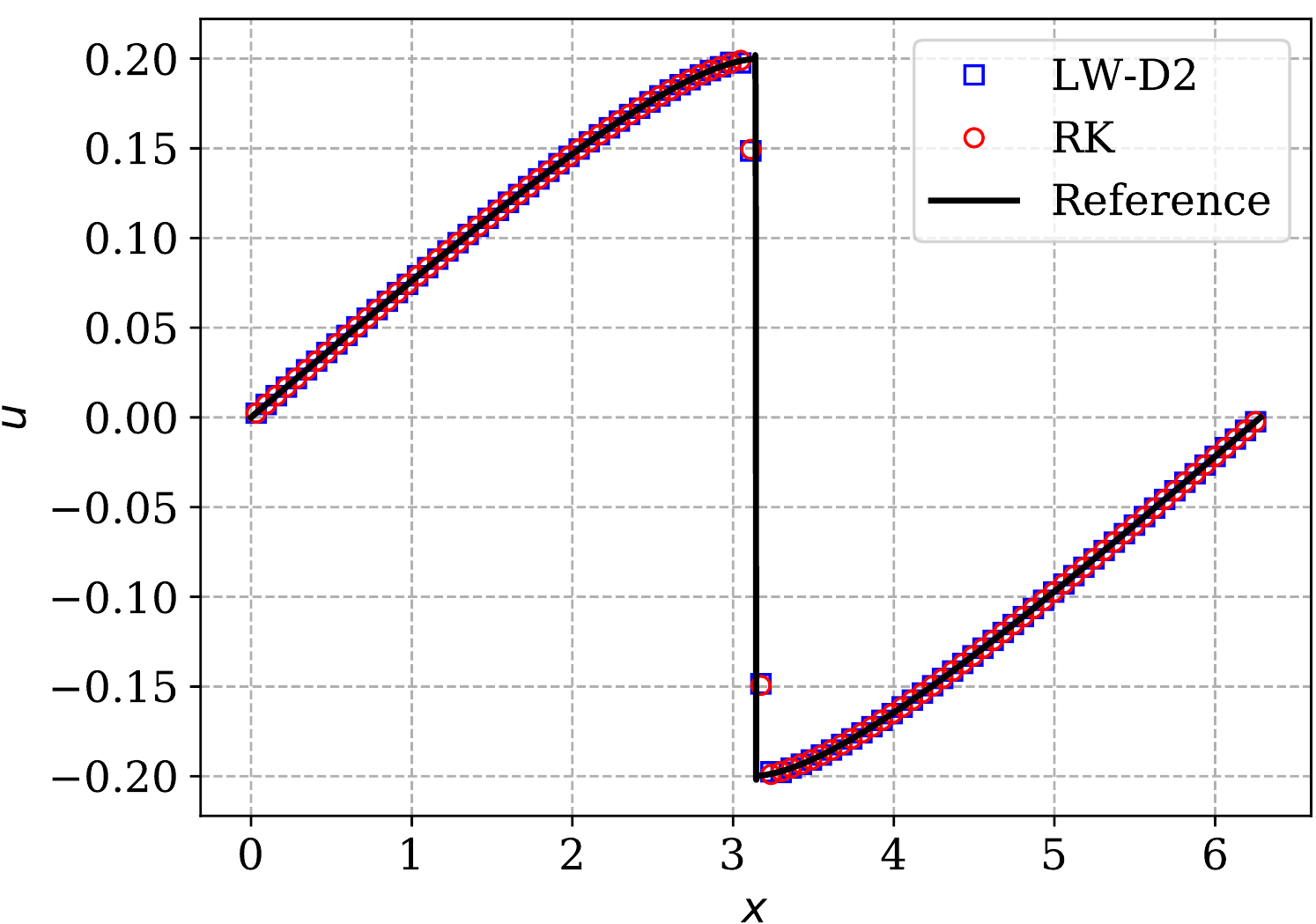} \\
(a) $t=2$ & (b) $t=4.5$ & (c) $t=8$
\end{tabular}
\end{center}
\caption{Solution of 1-D Burger's equation with $N=3$ and 100 cells at different time instants. TVB limiter $(M=1)$ is used. The reference solution is computed using RKFR, degree $N=1$, on a mesh of 3500 cells.}
\label{fig:burg1}
\end{figure}

\begin{figure}
\begin{center}
\begin{tabular}{ccc}
\includegraphics[width=0.30\textwidth]{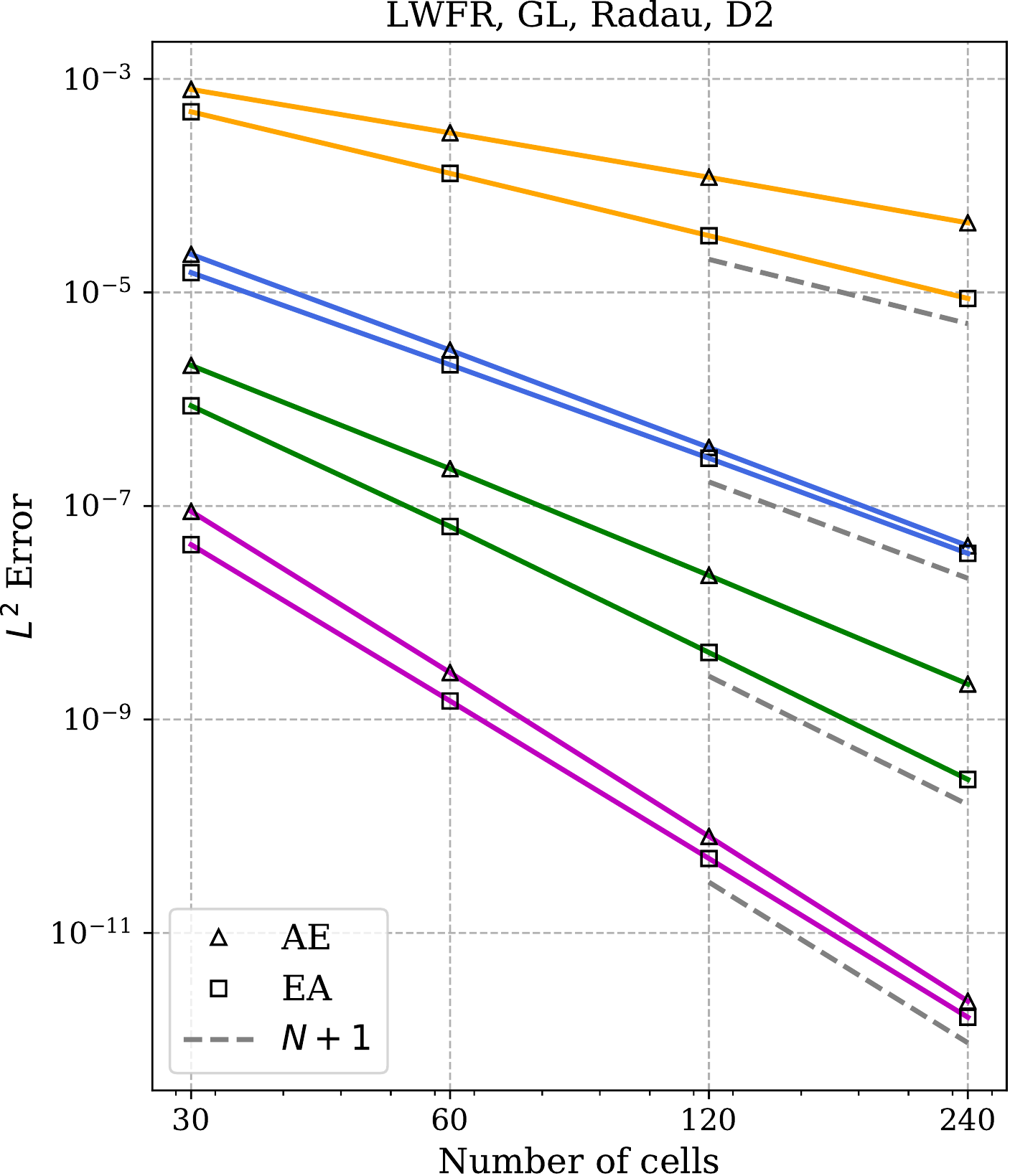} &
\includegraphics[width=0.30\textwidth]{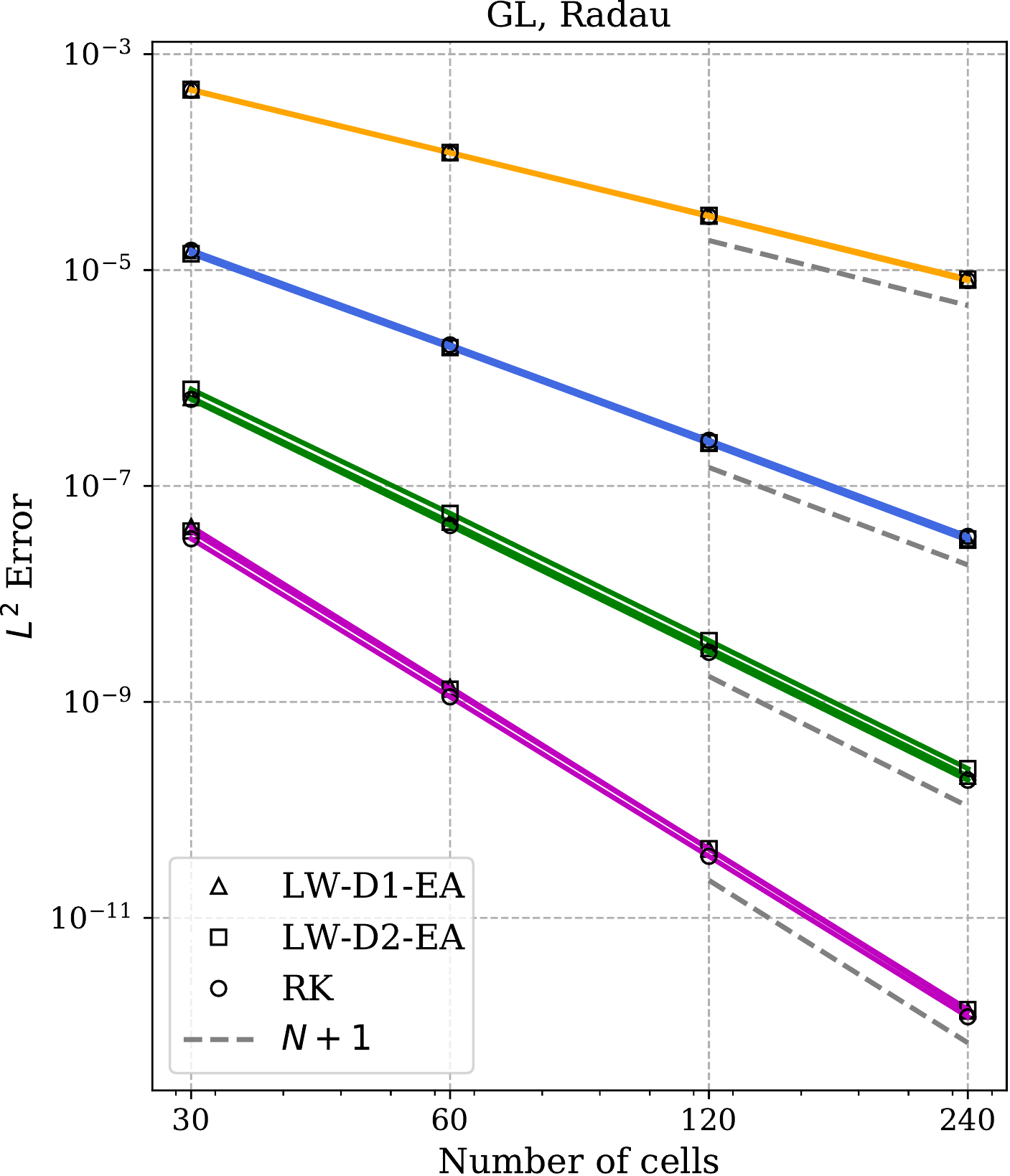} &
\includegraphics[width=0.30\textwidth]{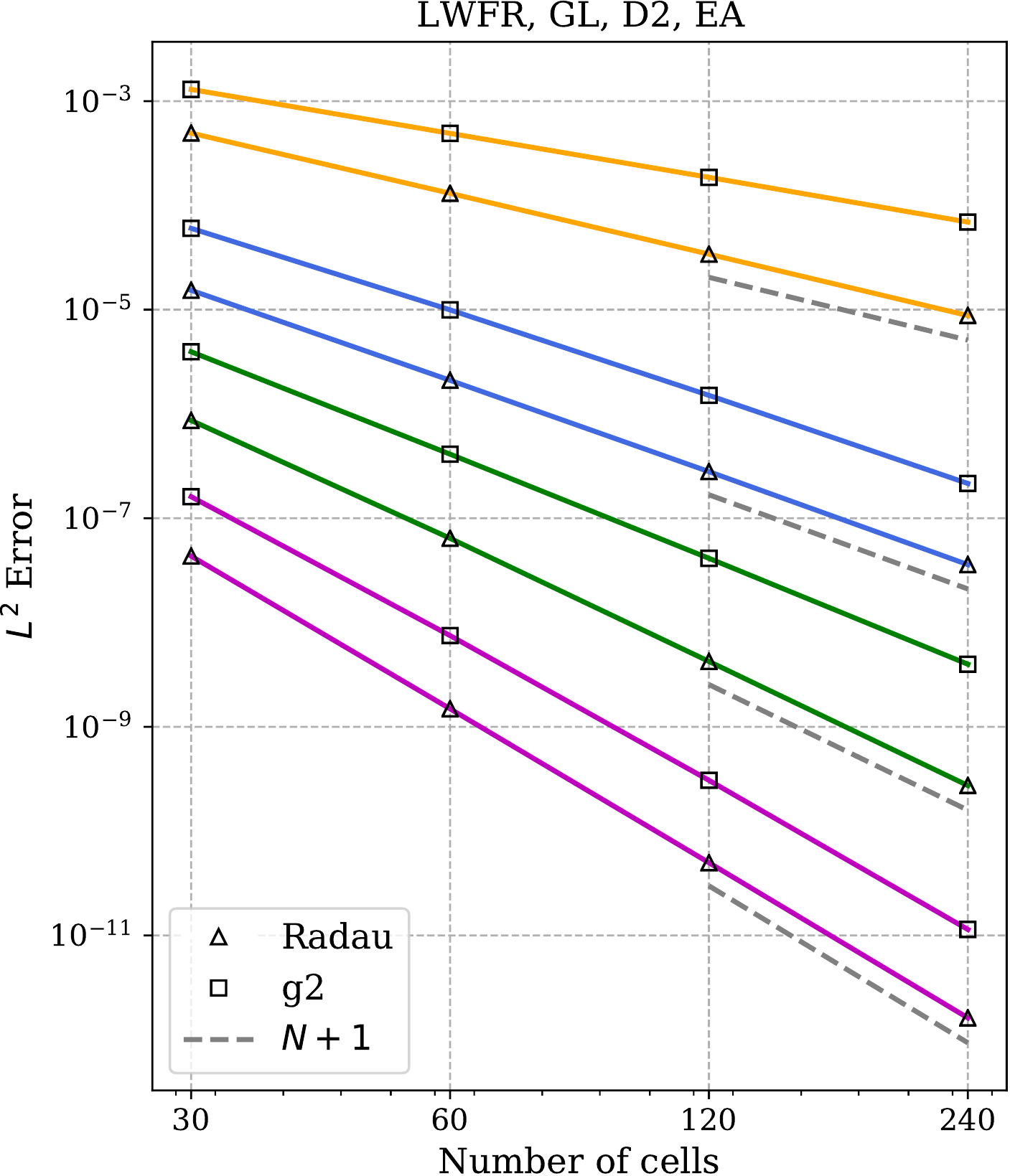} \\
(a) & (b) & (c)
\end{tabular}
\end{center}
\caption{Error convergence for 1-D Burger's equation at time $t=2$. (a) AE vs EA, (b) Radau vs g2, EA scheme, (c) LW-EA vs RK.}
\label{fig:burg2}
\end{figure}

The one dimensional Burger's equation is a conservation law of the form $u_t + f(u)_x=0$ with the quadratic flux $f(u) = u^2/2$.  For the smooth initial condition $u(x,0)=0.2\sin (x)$, we compute the numerical solution at different times $t\in\{2.0,4.5,8.0\}$ with periodic boundary condition in the domain $[0,2\pi]$. The TVB limiter with parameter $M=1$ is used. A stationary discontinuity is formed at $x=\pi$ and time $t_c=5.0$. The solutions are shown in Figure~(\ref{fig:burg1}) for degree $N=3$ and compared with the results from the RK method. We see that the discontinuity is captured accurately and without any oscillations, and the LW results compare very well with the RK results. At time $t=2$, the solution is still smooth and we can obtain the exact solution, using which, error norms and convergence rates can be estimated, see Figure~(\ref{fig:burg2}).

Figure~(\ref{fig:burg2}a) compares the error norms for the AE and EA methods for the Rusanov numerical flux, and using GL solution points, Radau correction and D2 dissipation; at odd degrees, the convergence rate of AE is less than optimal and close to $O(h^{N+1/2})$, while at even degrees, we obtain the optimal $O(h^{N+1})$ rate.  In Figure~(\ref{fig:burg2}b), we see that error norms of LW-EA and RK schemes are very close.  In Figure~(\ref{fig:burg2}c), we compare the two correction functions using the EA scheme and Rusanov flux; it is clear that the errors with Radau correction are significantly smaller than those with g2 correction. 

\begin{figure}
\begin{center}
\begin{tabular}{cc}
\includegraphics[width=0.40\textwidth]{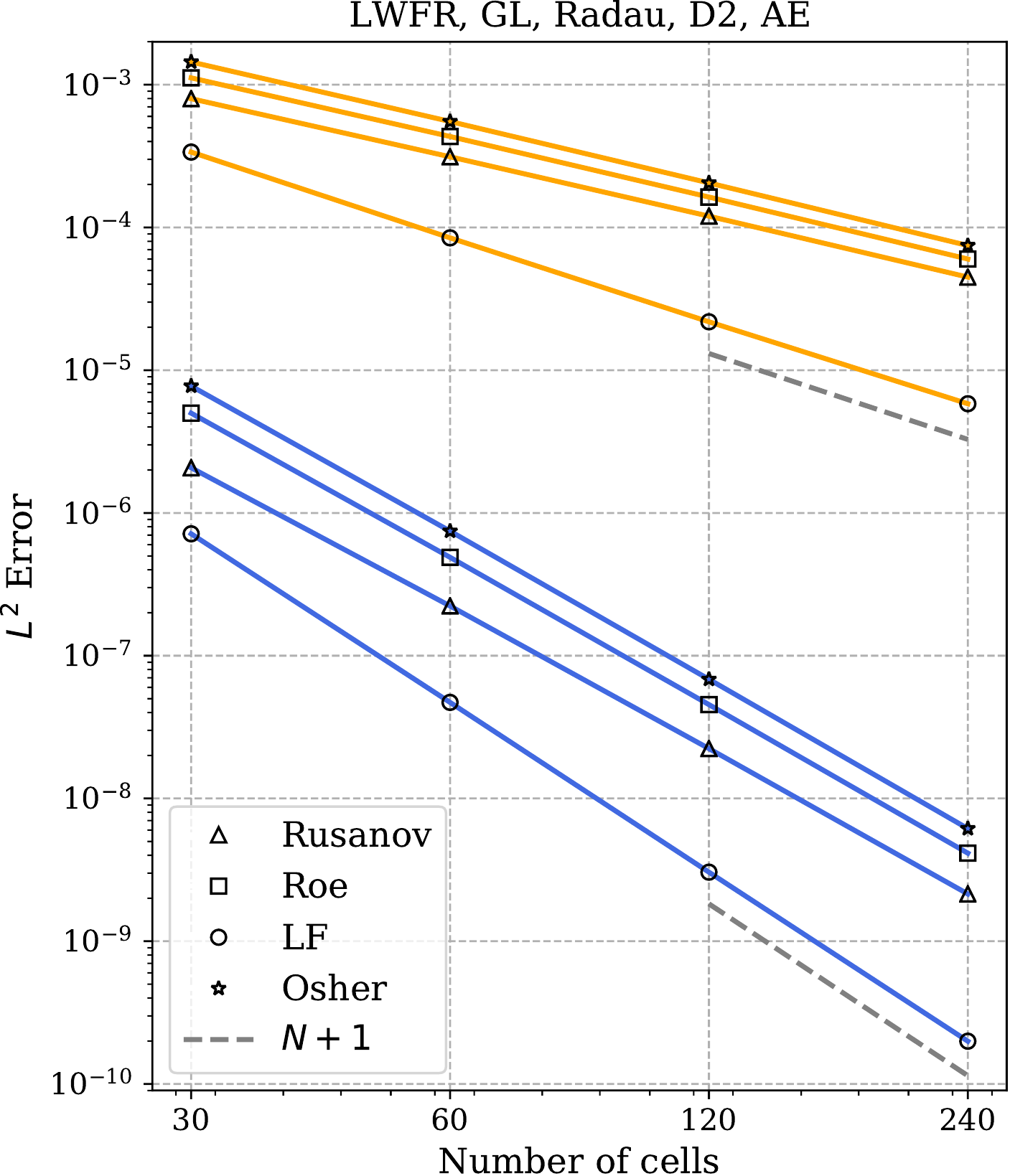} &
\includegraphics[width=0.40\textwidth]{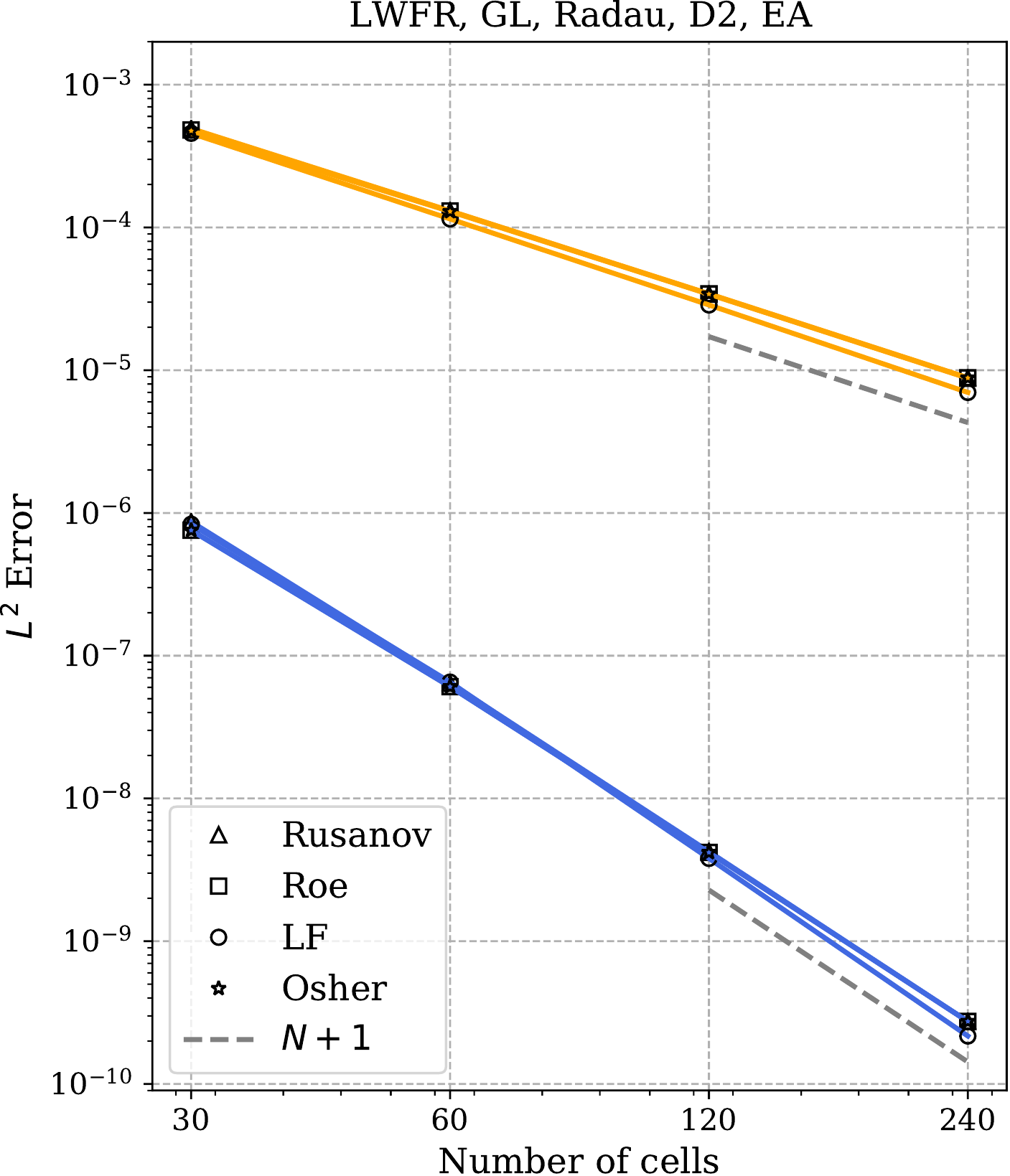} \\
(a) & (b)
\end{tabular}
\end{center}
\caption{Error convergence for 1-D Burger's equation at time $t=2$; effect of numerical fluxes for $N=1,3$. (a) AE scheme, (b) EA scheme.}
\label{fig:burg4}
\end{figure}

Next, we study the effect of different numerical fluxes in Figure~(\ref{fig:burg4}) for odd degrees $N=1,3$. With the AE scheme, only the global Lax-Friedrich flux is able to achieve the correct convergence rates and has the smallest errors compared to other fluxes which is a surprising result since it is a very dissipative flux. When the EA scheme is used as shown in the right of Figure~(\ref{fig:burg4}), all the numerical fluxes perform very similarly and achieve the optimal convergence rate. An examination of the error distribution in space shows that the AE scheme in combination with any numerical flux other than global Lax-Friedrich, produces large errors around the region of sonic points where $f'(u)=0$; however this happens only for odd degrees and the reason for this behaviour is not known at present. For initial data where the solution does not have a sonic point as in Example 2 of~\cite{Lou2020}, we get optimal convergence rates for all degrees even with the AE scheme.

\subsection{Non-convex problem: Buckley-Leverett equation}
\begin{figure}
\begin{center}
\begin{tabular}{cc}
\includegraphics[width=0.45\textwidth]{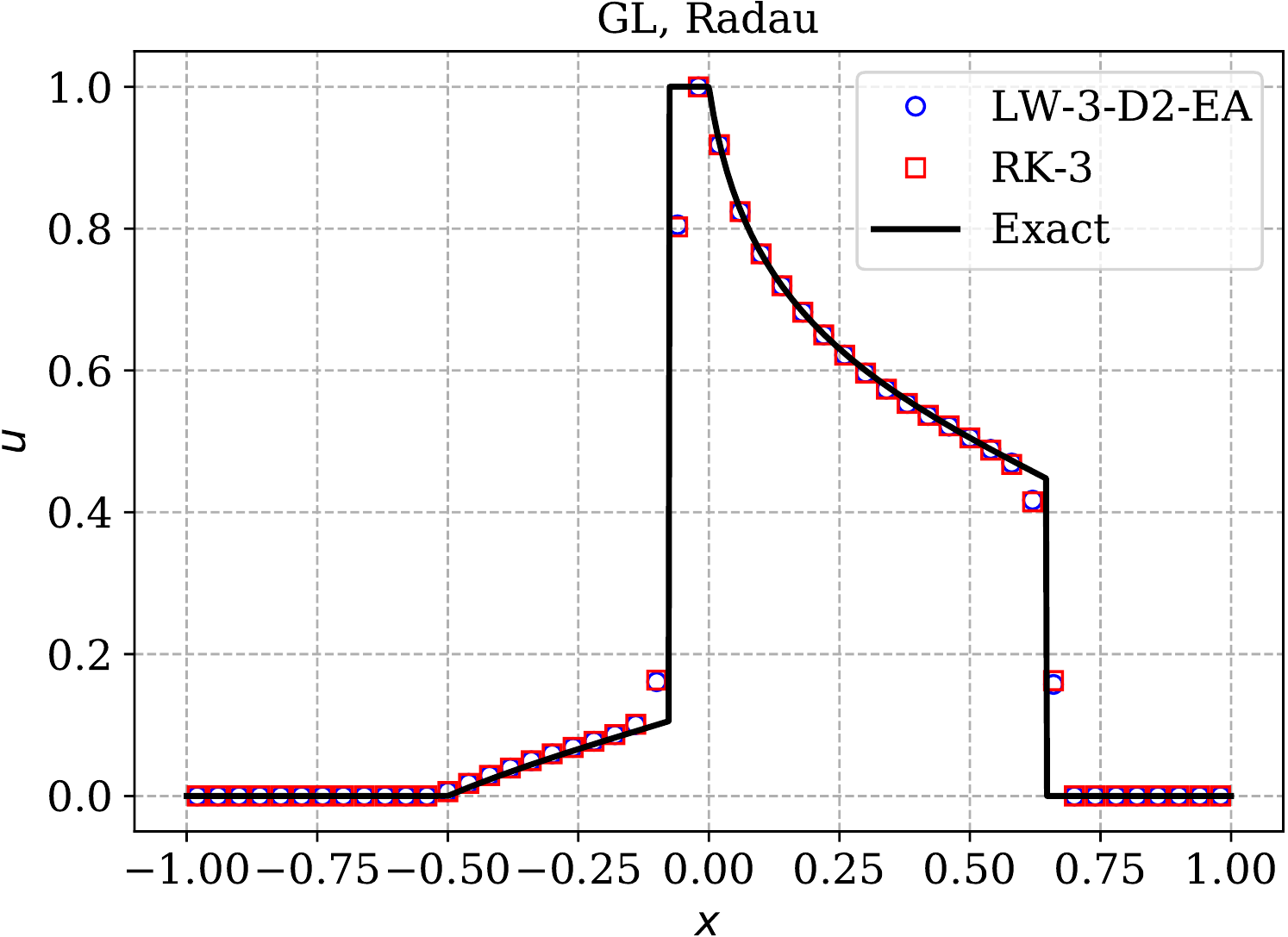} &
\includegraphics[width=0.45\textwidth]{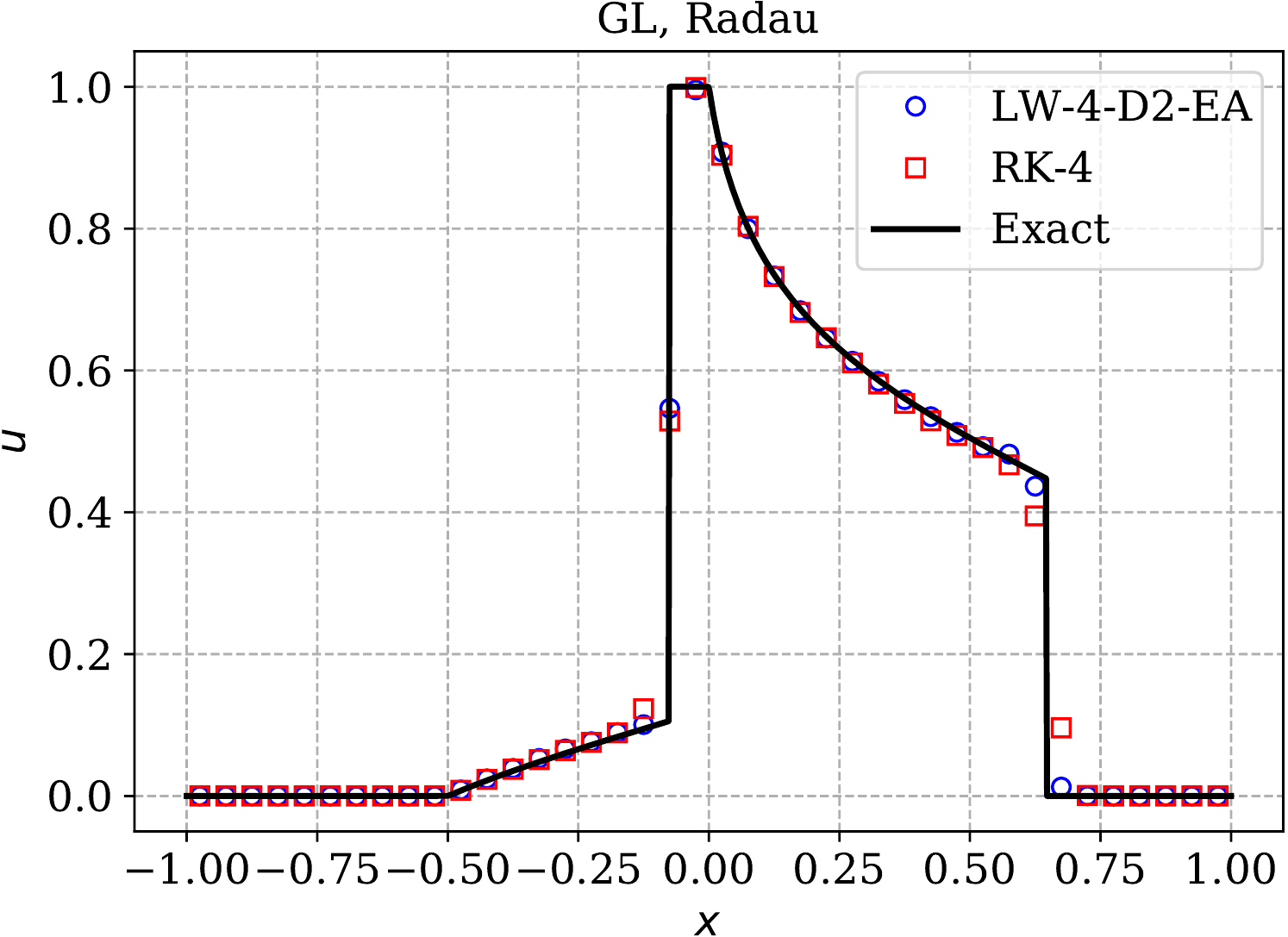} \\
(a) $N=3$ & (b) $N=4$
\end{tabular}
\end{center}
\caption{Solution of Buckley-Leverett model with TVD limiter using polynomial degrees $N=3,4$ with 200 dofs in each case. A positivity preserving scaling limiter~\cite{Zhang2010b} has been used to keep the solution in $[0,1]$.}
\label{fig:bucklev1}
\end{figure}

We consider the Buckley-Leverett equation $u_t+f(u)_x=0,$ where the flux $f(u)=\dfrac{4u^2}{4u^2+(1-u)^2}$ is convex and concave in different regions of the solution space. The numerical solution is computed up to the time $t=0.4$ with the initial condition
\[
u(x,0) = \begin{cases}
1, & x \in [-1/2,0] \\
0, & \textrm{elsewhere}
\end{cases}
\]
At $t=0.4$, the characteristics that originate from  the two discontinuities do not intersect, and thus we only have to deal with the two uncoupled Riemann problems. Solutions to Riemann problems for piecewise strictly convex-concave fluxes can be computed explicitly. In case of the Buckley-Leverett model, we can split the state-space $[0,1]$ into $[0,u_\text{buck}]$ and $[u_\text{buck},1]$ so that $f$ is strictly convex in $(0,u_\text{buck})$ and strictly concave in $(u_\text{buck},1)$. Thus, the solution to a Riemann problem with states $0,1$ would compose of a rarefaction and shock, and the exact solution corresponding to the above defined initial condition is composed of two rarefaction-shock combinations. Since the solution measures saturation of displacing fluid, it should lie in the interval $[0,1]$ and we try to enforce this by applying a positivity preserving scaling limiter~\cite{Zhang2010b}. For the LW schemes, we cannot strictly prove the positivity of the resulting scheme but numerical results show that this holds in practice, with slightly reduced CFL number compared to the Fourier CFL number. For $N=4$, the optimal CFL conditions preserve the bounds, while a slightly reduced CFL of 0.079 was needed for $N=3$. Figure~(\ref{fig:bucklev1}) shows the results at the final time obtained using degree $N=3,4$,  respectively. Since the flux is monotone in solution space, an upwind numerical flux is used at cell interfaces, i.e., $F_\eph = F_\eph^-$. The numerical solutions are able to resolve all the waves well including correct shock location, and they compare well with the results from the RK scheme.

\section{Numerical results in 1-D: Euler equations}\label{sec:res1dsys}
As an example of system of non-linear hyperbolic equations, we consider the one-dimensional Euler equations of gas dynamics given by
\begin{equation}\label{eq:1deuler}
\pd{}{t} \begin{pmatrix}
\rho \\
\rho v \\
E
\end{pmatrix} +
\pd{}{x} \begin{pmatrix}
\rho v \\
p + \rho v^2 \\
(E+p)v
\end{pmatrix} = 0
\end{equation}
where $\rho, v, p$ and $E$ denote the  density, velocity, pressure and total energy of the gas, respectively. For a polytropic gas, an equation of sate $E=E(\rho, v, p)$ which leads to a closed system is given by
\begin{equation}\label{eq:state}
E = E(\rho, v, p) = \frac{p}{\gamma -1}+\frac{1}{2} \rho v^2
\end{equation}
where $\gamma > 1$ is the adiabatic constant, that will be taken as $1.4$ which is the value for air. The numerical fluxes for Euler equations are explained in the Appendix~\ref{apendix:numfluxes}. In the following results, wherever it is not mentioned, we use the Rusanov flux.

In the scalar results, we see that the LW scheme with Radau correction function is superior to that with $g_2$ correction function in terms of error reduction. In light of this, for the 1-D Euler case we compare the performance of LW scheme with RK scheme using the Radau correction function. It is also observed that the EA scheme is more accurate than AE in the scalar case. So, for the 1-D Euler case we present only those results obtained using EA schemes.

Note that wherever it is not specified we use the CFL numbers of Table~(\ref{tab:cfl}) and whenever we compare the numerical solutions of LW scheme with that of  RK scheme, both are run with the CFL numbers of LW scheme. In the time integration of the RK schemes, for degree $N=1$ and $2$, we use $(N+1)$-stage SSPRK method  of order $N+1$. For $N=3,$ we use a five stage, SSPRK method of order four~\cite{Spiteri2002} as there is no four stage SSPRK method. In smooth test cases and for $N=4$, we use a six stage Runge-Kutta method of order five~\cite{Tsitouras2011}. In those test cases where the solution is  not sufficiently smooth, the SSP property of the RK time integration is useful in obtaining non-oscillatory solutions. As we do not have SSPRK method of order five with positive coefficients~\cite{Ruuth2002}, we use the SSPRK method ~\cite{Spiteri2002} of order four when $N=4$.

\subsection{Smooth solution}
\begin{figure}
\centering
\begin{tabular}{cc}
\includegraphics[width=0.40\textwidth]{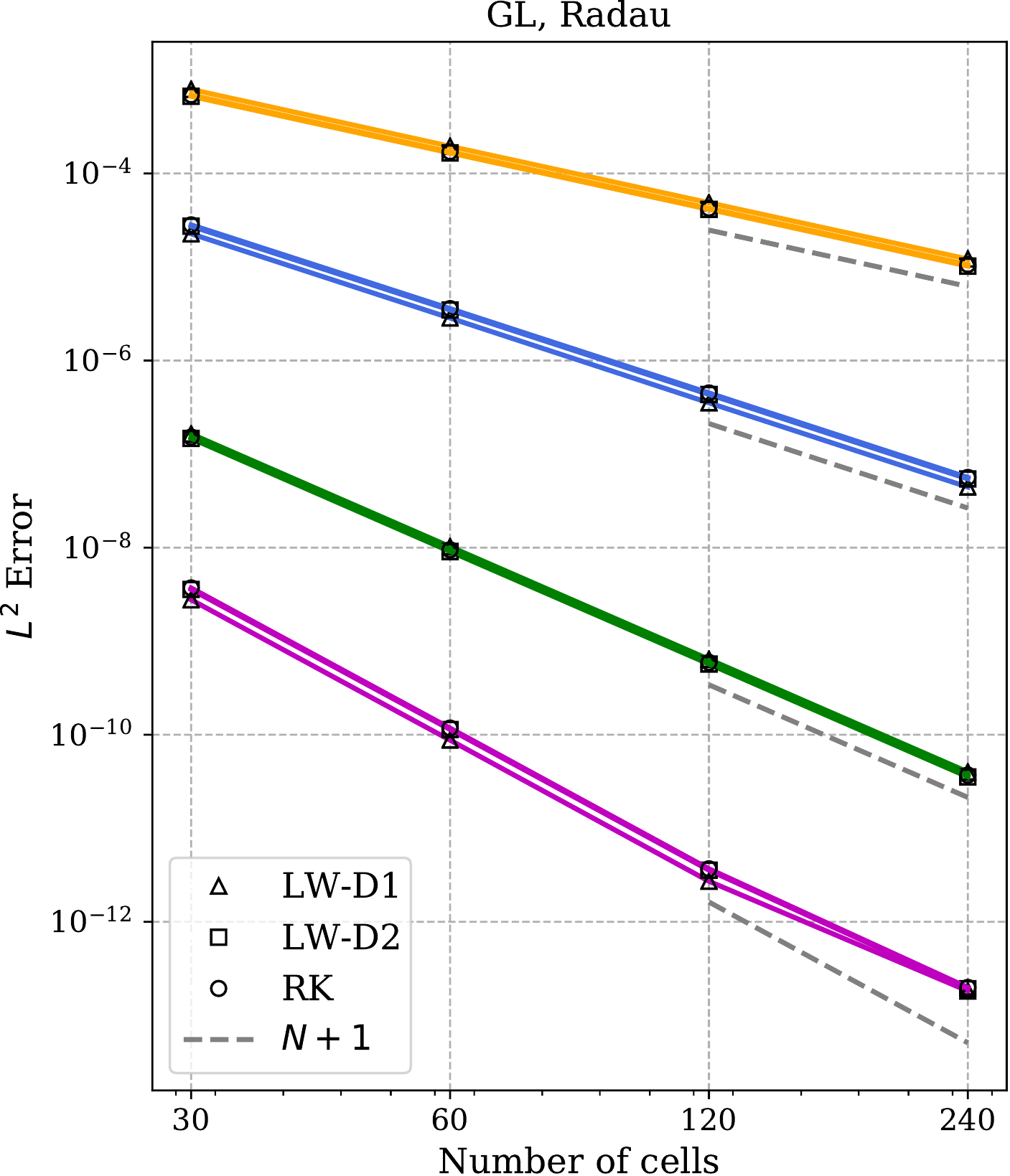} &
\includegraphics[width=0.40\textwidth]{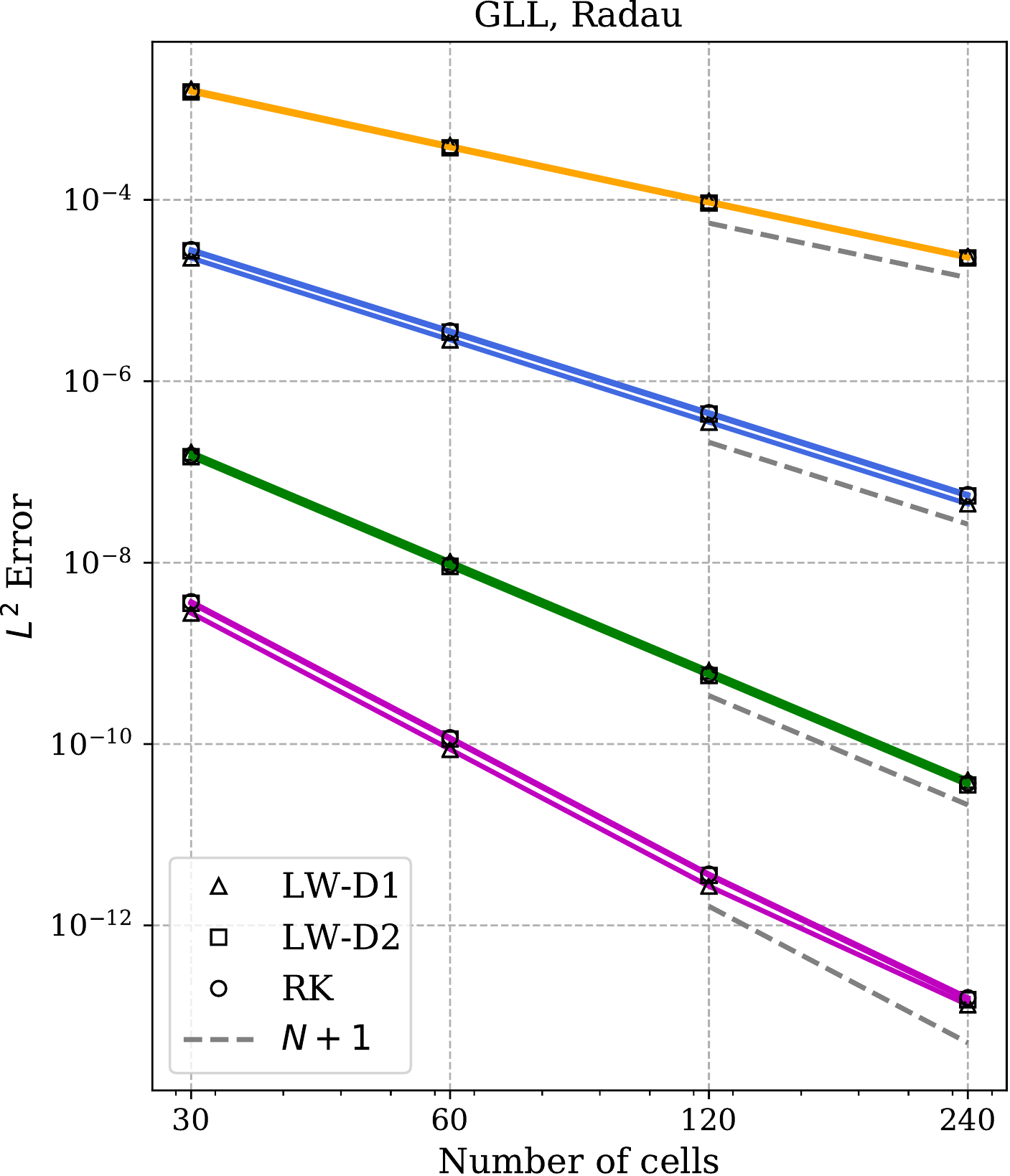} \\
(a) GL points & (b) GLL points
\end{tabular}
\caption{Density error convergence for 1-D Euler's equation at time $t=1$. The different colors correspond to degrees $N=1,2,3,4$ from top to bottom. (a) GL points, (b) GLL points.}
\label{fig:dwave}
\end{figure}

To verify the accuracy of the proposed scheme, we solve the Euler equations \eqref{eq:1deuler},~\eqref{eq:state} with a smooth initial condition
\[
\rho(x,0)=1+0.5\sin(2\pi x), \qquad v(x,0)=1, \qquad p(x,0)=1
\]
in the domain $[0,1]$ with periodic boundary condition for all the variables. The corresponding exact solution  is a density wave, i.e., it consists of a translation of the initial density at constant speed of one, and is given by
\[
\rho(x,t)= 1+0.5\sin(2\pi (x-t)), \qquad v(x,t) = 1, \qquad p(x,t) = 1
\]
We compute the solution up to the time $t=1$ and estimate the error norms. The linear nature of this test case makes EA and AE schemes equivalent, and we show only the EA results. We plot the error in the density obtained using the LW and RK schemes and the corresponding results are given in Figure~(\ref{fig:dwave}). In each case we observe the expected order of accuracy and the  error reduction is close to that of RK scheme.
\begin{remark}
Based on the scalar test cases and the above smooth test case for Euler equations, in  all the remaining 1-D Euler test cases, we present only the results obtained using D2 scheme since it is advantageous in terms of having higher CFL number and performs as well as or better than the D1 scheme.
\end{remark}
\subsection{Sod's shock tube problem}
\begin{figure}
\centering
\begin{tabular}{cc}
\includegraphics[width=0.45\textwidth]{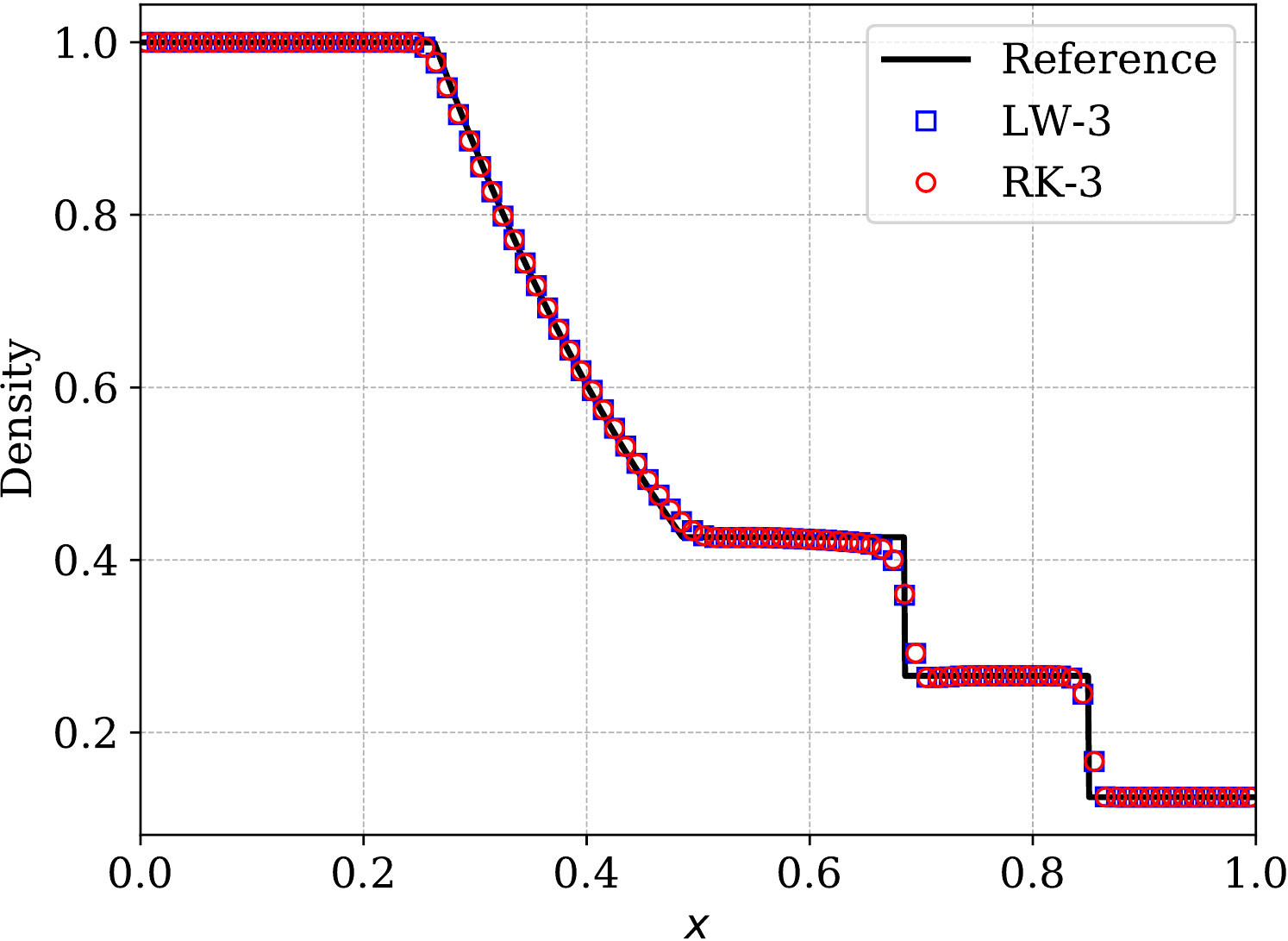} &
\includegraphics[width=0.45\textwidth]{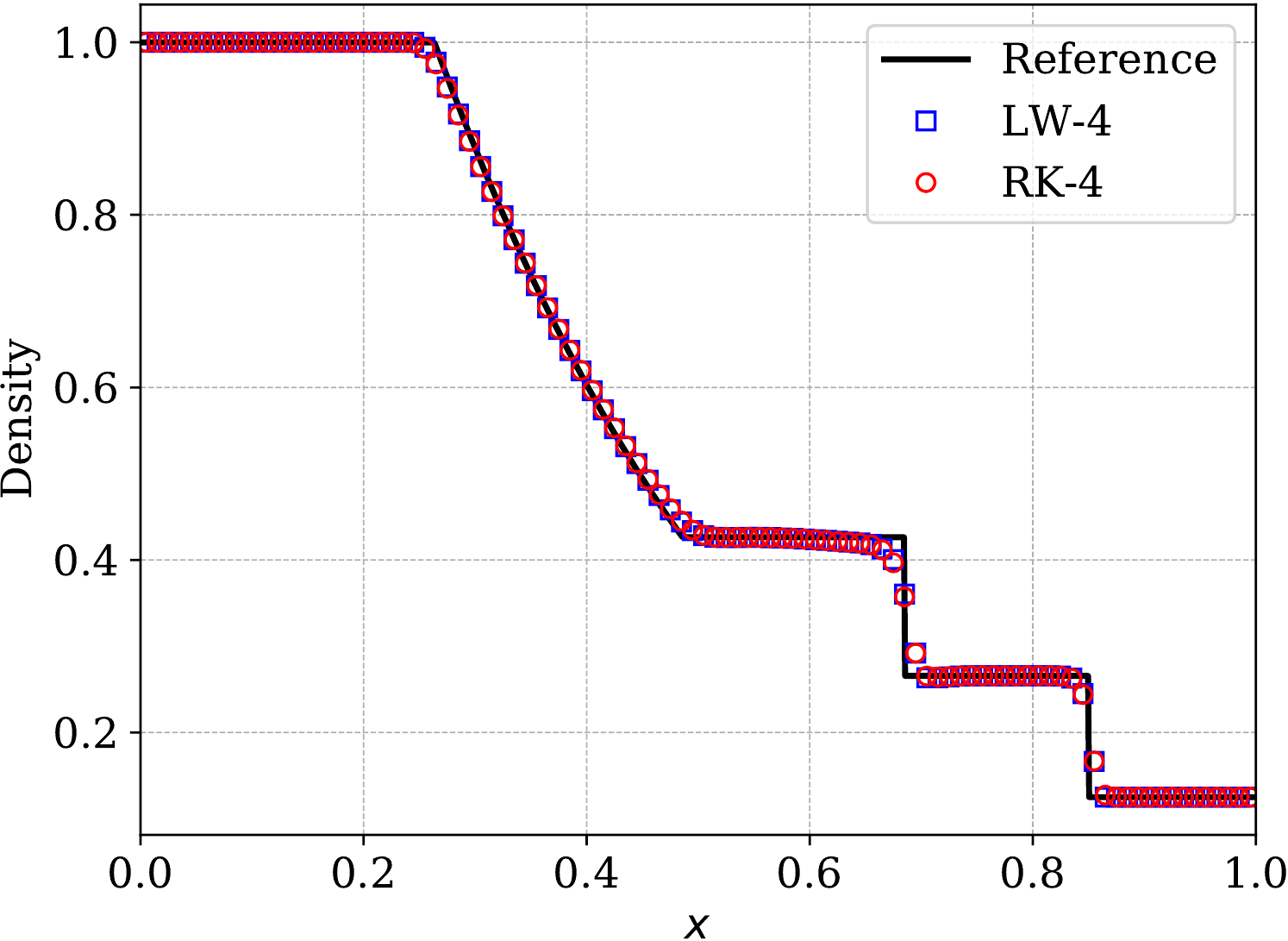} \\
(a) & (b)
\end{tabular}
\caption{Numerical solutions of 1-D Euler equations (Sod's test case) obtained by LW and RK schemes for polynomial degree (a) $N=3$ and (b) $N=4$ using Radau correction function and GL solution points. The solutions are shown at time $t=0.2$ on a mesh of $100$ cells with Rusanov flux and D2 dissipation. The TVB limiter is used with the parameter $M=10$.}
\label{fig:sod}
\end{figure}

The Sod's shock tube problem~\cite{Sod1978} is a Riemann problem which models a shock tube where gas at two different conditions initially is allowed to interact, with the formation of many waves including shock, contact and rarefaction waves. The Euler equations~\eqref{eq:1deuler} are solved with the initial condition
\begin{equation}\label{eq:sod}
(\rho,v,p) = \begin{cases}
(1.0,~0,~1), & x < 0.5 \\
(0.125,~0,~0.1), & x > 0.5
\end{cases}
\end{equation}
for which the exact solution is composed of a left rarefaction, a contact discontinuity and a right shock wave. The approximated solution is computed in the domain $[0,1]$ with the outflow boundary conditions on both the ends $x=0$ and $x=1$. We run the numerical scheme up to time  $t=0.2$  with 100 cells using the TVB limiter with parameter $M=10$. The  density profile  obtained using the  LW and RK schemes for $N=3$ and $N=4$ are shown in Figure~(\ref{fig:sod}) together with the exact solution. From the plots it is visible that  the results obtained using the LW scheme agree very well with that of RK scheme.

\subsection{Lax problem}
\begin{figure}
\centering
\begin{tabular}{cc}
\includegraphics[width=0.45\textwidth]{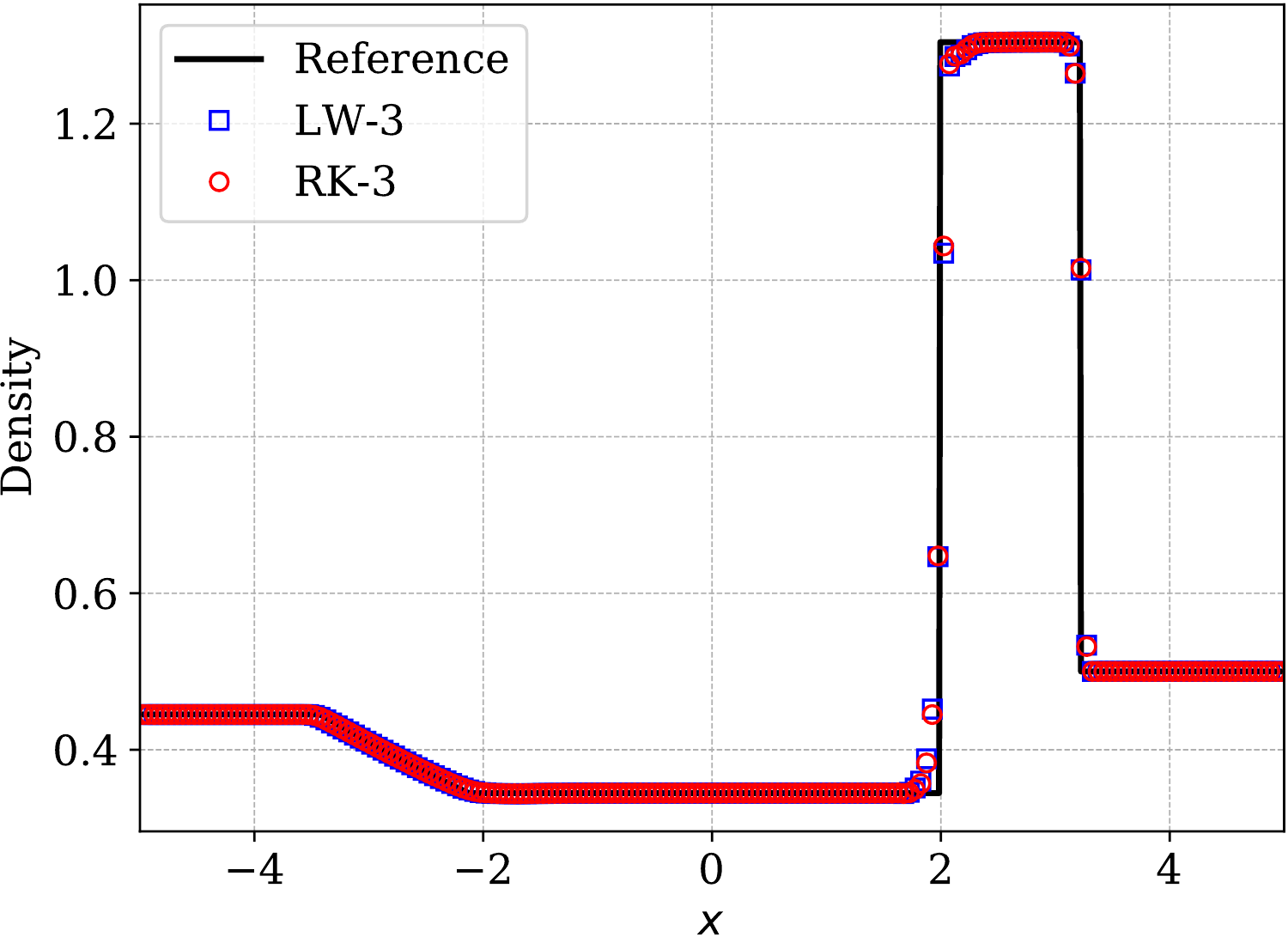} &
\includegraphics[width=0.45\textwidth]{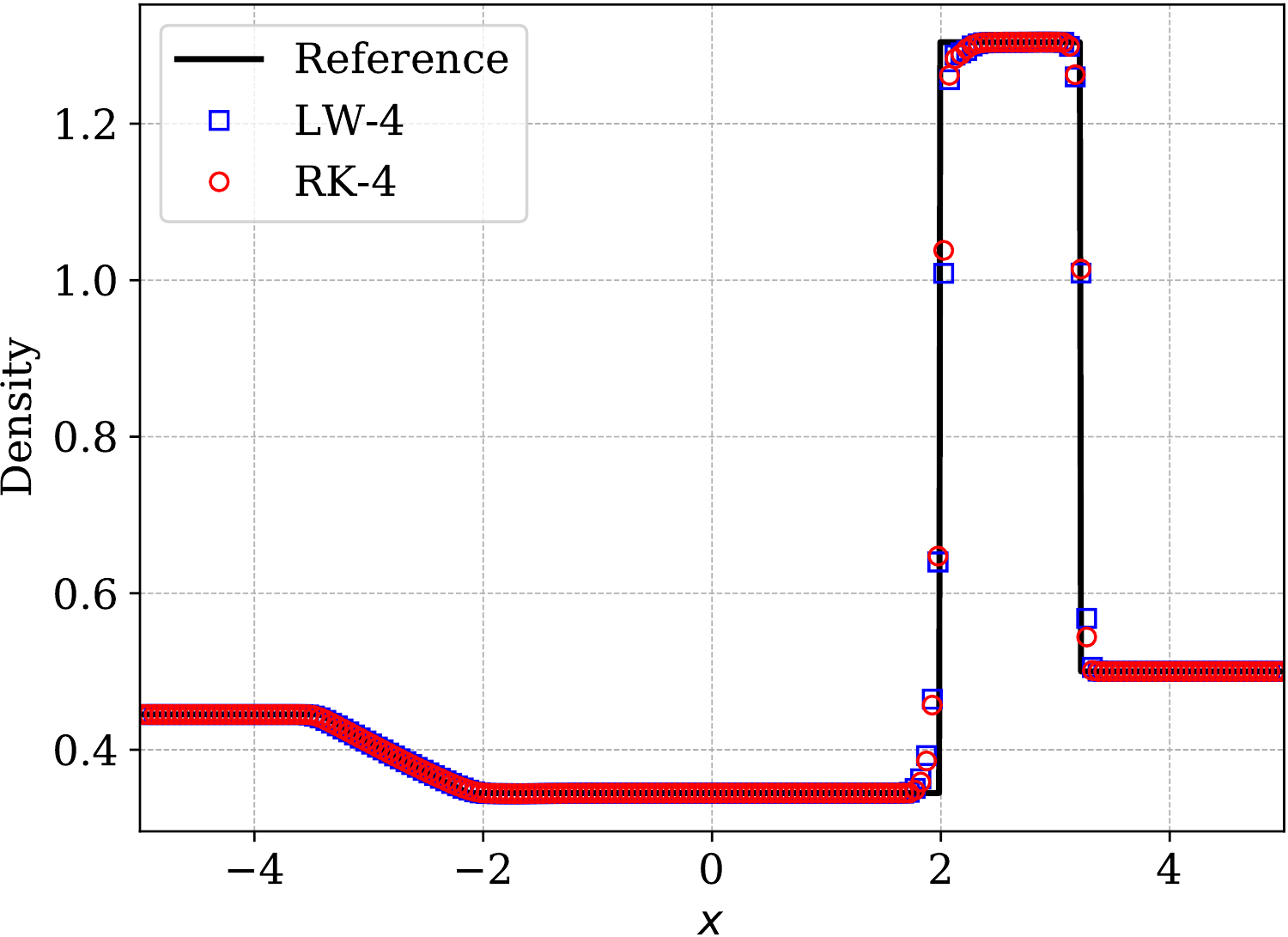} \\
(a) & (b)
\end{tabular}
\caption{Numerical solutions of 1-D Euler equations (Lax's test case) obtained by LW and RK schemes for  polynomial degree (a) $N=3$ and (b) $N=4$ with Radau correction function and GL solution points. The solutions are computed on a mesh of $200$ cells with dissipation model D2 and Rusanov numerical flux and are shown at time $t=1.3$. The TVB limiter is used with parameter $M=1$.}
\label{fig:lax}
\end{figure}

We consider the Lax problem given in~\cite{Lax1954,Hirsch1990} which solves a Riemann problem for the system of equations \eqref{eq:1deuler}  with initial condition
\begin{equation}
(\rho,v,p) = \begin{cases}(0.445,0.698,3.528), & x < 0 \\
(0.5,0,0.571), & x > 0
\end{cases}
\end{equation}
where, unlike the Sod's shock tube problem, the initial velocity is not zero. The exact solution of this  Riemann problem is known and it consists of  a rarefaction, a right moving contact discontinuity and shock. For a detailed description of this problem, see~\cite{Hirsch1990}.  This is a demanding test case in the sense that,  high order schemes are prone to produce oscillations near the contact discontinuity. The numerical solution is computed up to time  $t=1.3$ in the domain $[-5,5]$ using 100 cells and using TVB limiter with parameter $M=1$. The approximate solutions are computed for polynomial degrees $N=3$ and $N=4$, and the corresponding density profiles are shown in Figure~(\ref{fig:lax}) along with the exact solution. We observe that  the LW scheme captures the wave structures accurately without oscillations and the numerical solutions are very close to that of RK scheme.
\subsection{Shu-Osher problem}
\begin{figure}
\centering
\begin{tabular}{cc}
\includegraphics[width=0.45\textwidth]{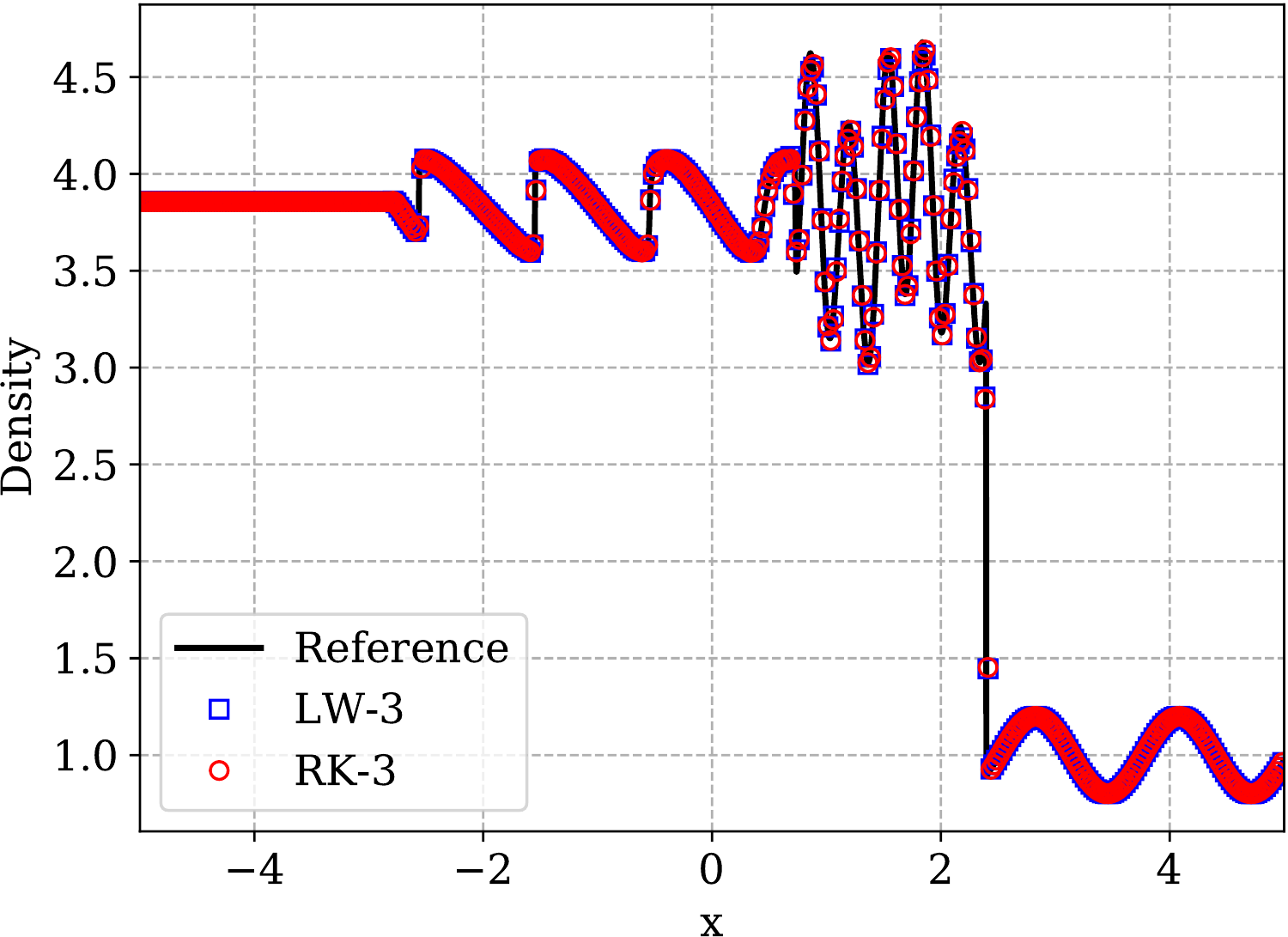} &
\includegraphics[width=0.45\textwidth]{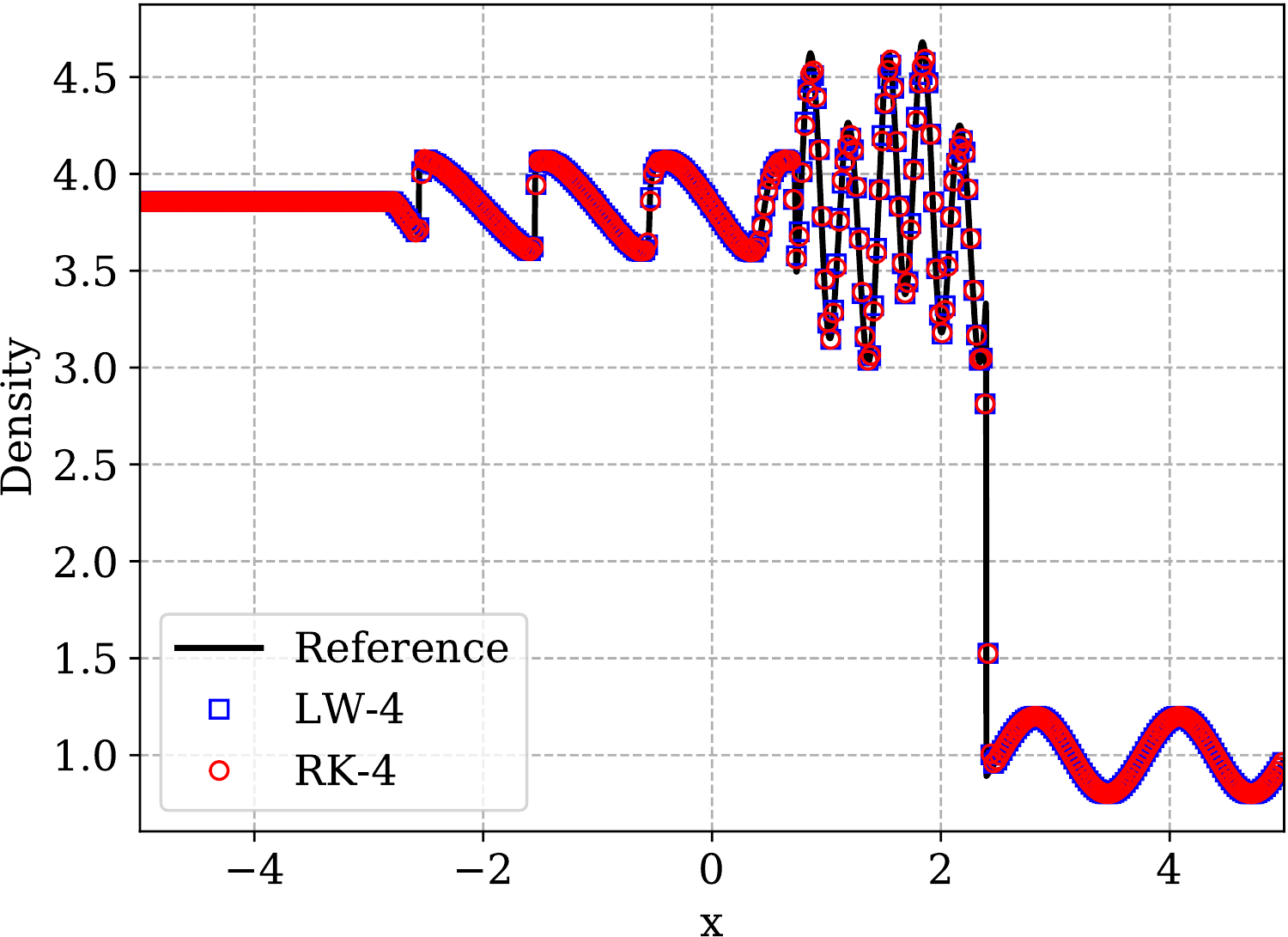} \\
(a) & (b) \\
\includegraphics[width=0.45\textwidth]{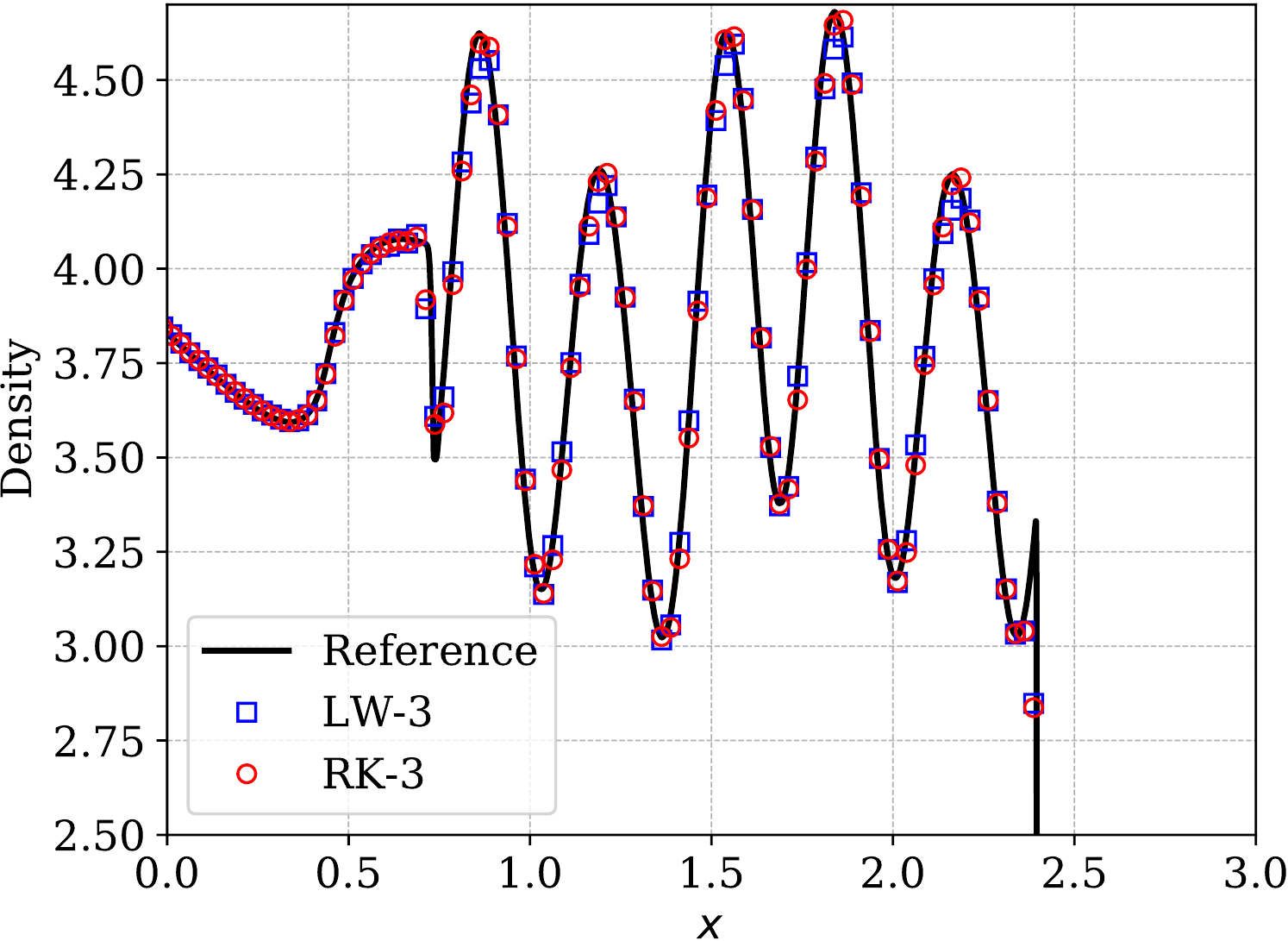} &
\includegraphics[width=0.45\textwidth]{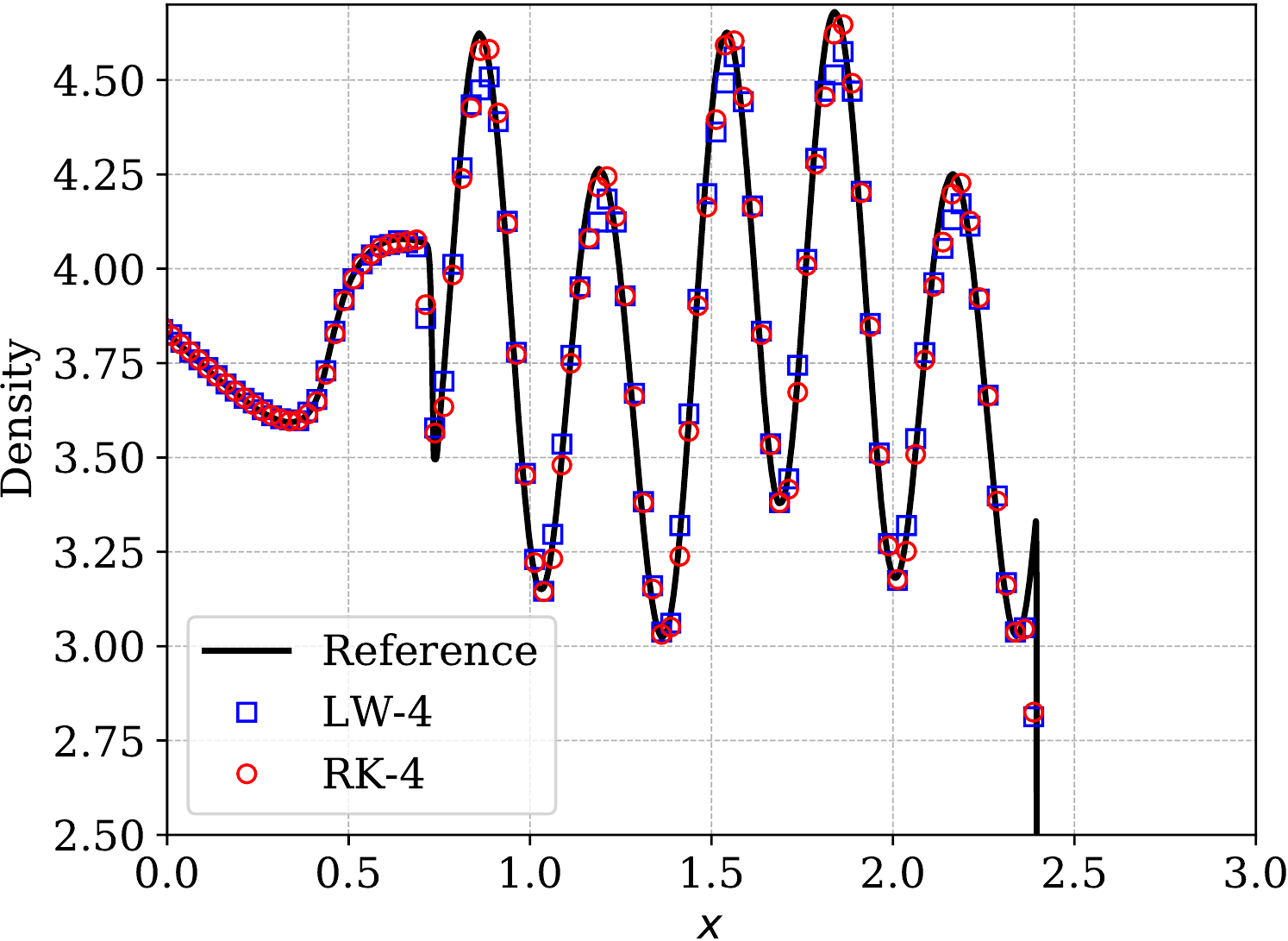} \\
(c) & (d)
\end{tabular}
\caption{Numerical solutions of 1-D Euler equations (Shu-Osher problem) obtained by LW and RK schemes for  (a,c) $N=3$ and (b,d) $N=4$ with Radau correction function and GL solution points. The enlarged plot of the oscillatory portion is given in the bottom row. The solutions are shown at time $t=1.8$ on a mesh of 400 cells with dissipation model D2 and Rusanov numerical flux.  The TVB limiter is used with parameter $M=300$.}
\label{fig:ShuOsher}
\end{figure}

This problem was developed in~\cite{Shu1989}  to test the numerical scheme's capability to accurately capture a shock wave and its interaction with a smooth density field, which then  propagates  downstream of the shock. Here, we compute the numerical solution of  \eqref{eq:1deuler} with initial condition
\begin{equation}\label{eq:shuosher}
(\rho,v,p)=\begin{cases}
(3.857143, 2.629369, 10.333333) & \mbox{ if } x<-4\\
(1+0.2\sin(5x),0,1) & \mbox{ if }x\geq -4
\end{cases}
\end{equation}
prescribed in the domain $[-5,5]$ at time $t=1.8$.  The smooth density profile passes through the shock and appears on the other side, and its accurate computation is challenging due to numerical dissipation. We discretize the spatial domain with 400 cells  and to control the spurious oscillations we use the TVB limiter with parameter $M=300$~\cite{Qiu2005b}. The density component of  the approximate solutions computed using  LW and RK schemes for $N=3$ and $N=4$ are  plotted against a reference solution obtained using  a very fine mesh, and are given in Figure~(\ref{fig:ShuOsher}).  We observe that the complicated wave patterns are accurately captured by the proposed LW scheme. Furthermore, the enlarged plots of the oscillatory portion  indicate that the numerical solutions corresponding to LW scheme are comparable with that of RK schemes.

\subsection{Blast wave}
\begin{figure}
\centering
\begin{tabular}{cc}
\includegraphics[width=0.45\textwidth]{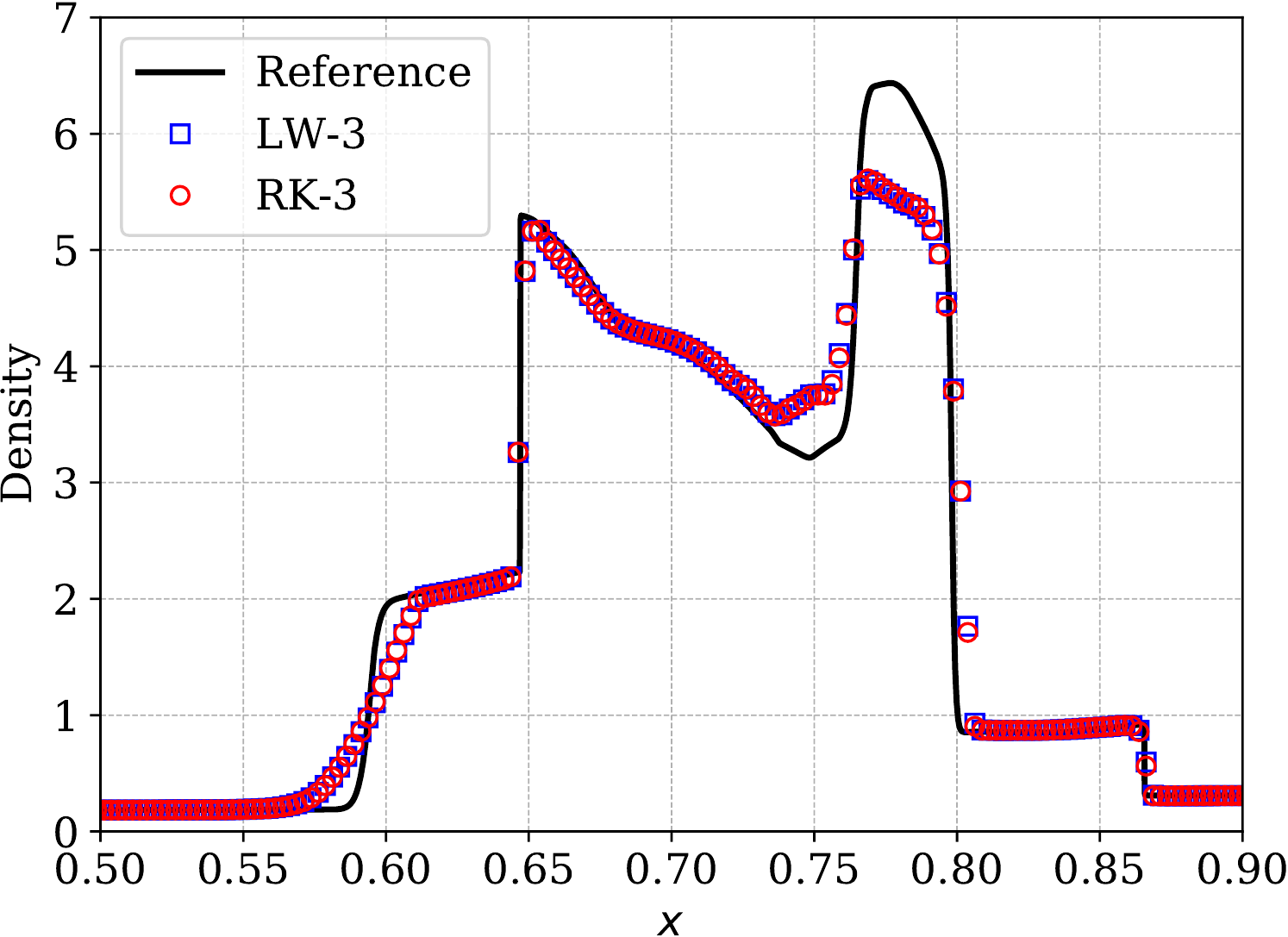} &
\includegraphics[width=0.45\textwidth]{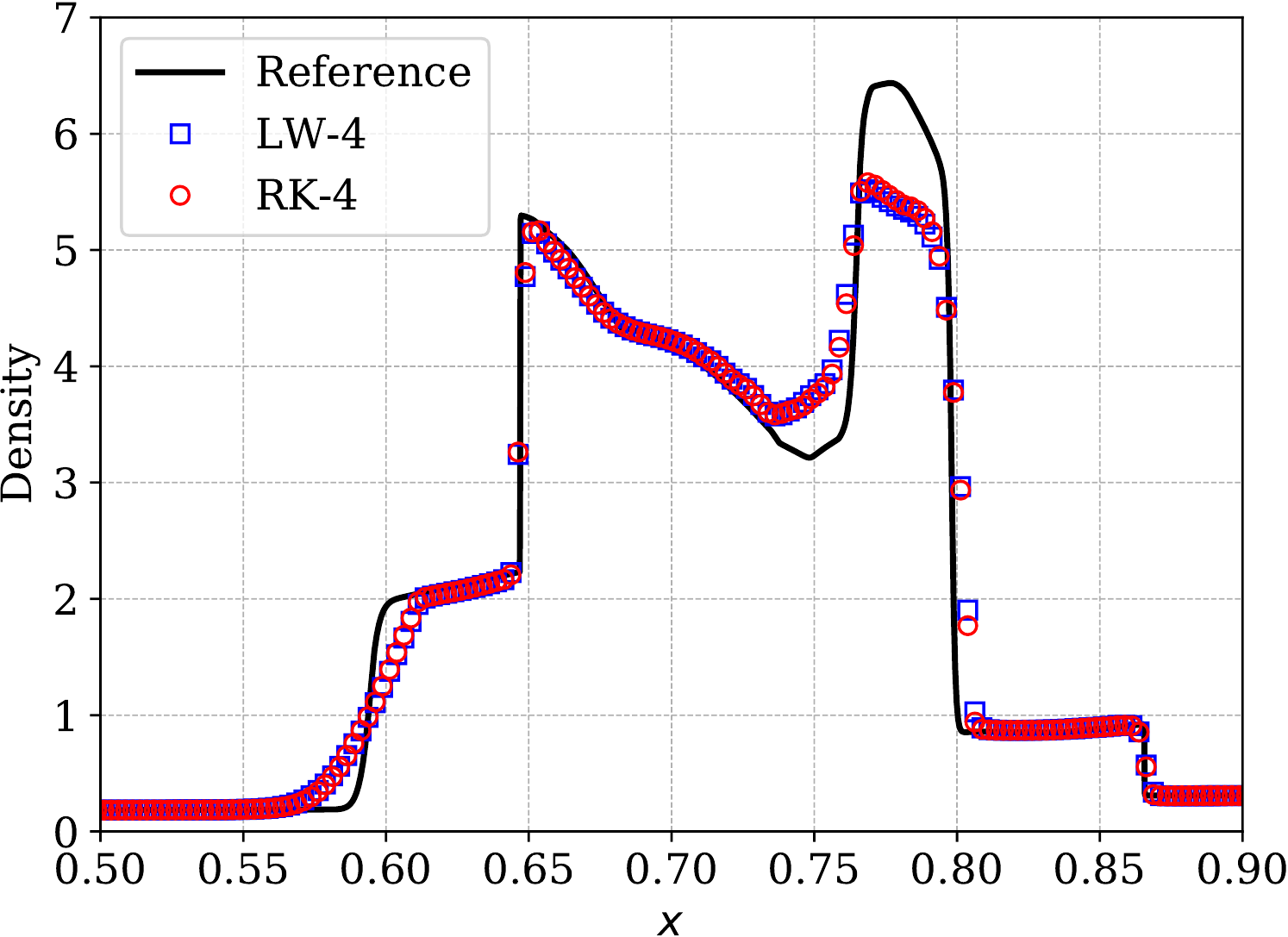} \\
(a) & (b)
\end{tabular}
\caption{Numerical solutions of 1-D Euler equations (Blast wave) obtained by LWFR and RKFR schemes for  (a) $N=3$ and (b) $N=4$ using Radau correction function and GL solution points. The solutions are plotted at time $t=0.038$   on a mesh  of $400$ cells with dissipation model D2 and Rusanov numerical flux.  The corresponding CFL numbers of LWFR scheme are chosen from the Table~\ref{tab:cfl} and kept same for the RKFR schemes. TVB limiter is used with the parameter $M=300$ and the EA scheme is used for numerical flux.}
\label{fig:blast}
\end{figure}

In this test case the Euler equations \eqref{eq:1deuler} are solved with the initial condition
\begin{equation*}
(\rho,v,p)=\begin{cases}
(1,0, 1000) & \mbox{ if } x<0.1\\
(1,0,0.01) & \mbox{ if }  0.1 < x <0.9\\
(1,0, 100) & \mbox{ if } x> 0.9
\end{cases}
\end{equation*}
in the domain $[0,1]$. It is originally introduced in~\cite{Woodward1984} to check the capability of the numerical scheme to accurately capture the shock-shock interaction scenario. The boundaries are set as  solid walls by imposing the  reflecting boundary conditions at $x=0$ and $x=1$. The solution consists of reflection of shocks and expansion waves off the boundary wall and several wave interactions inside the domain. With a grid of 400 cells, we run the simulation till the time $t=0.038$ where a high density peak profile is produced. The TVB limiter as in~\cite{Qiu2005b} with parameter $M=300$ is used along with a positivity preserving limiter~\cite{Zhang2010b}. We compare the performance of the LW scheme with the RK scheme and analyse how well they predict the density profile and its peak amplitude. For $N=3$ and $N=4$ cases, the results are given in Figure~(\ref{fig:blast}) where the approximated density profiles are compared with a reference solution computed using a very fine mesh. From the plots it is evident that the computed density profile obtained using LW scheme are  close to that of RK scheme.
\subsection{Numerical Fluxes: LF, Roe, HLL and HLLC}
\begin{figure}
\centering
\begin{tabular}{cc}
\includegraphics[width=0.45\textwidth]{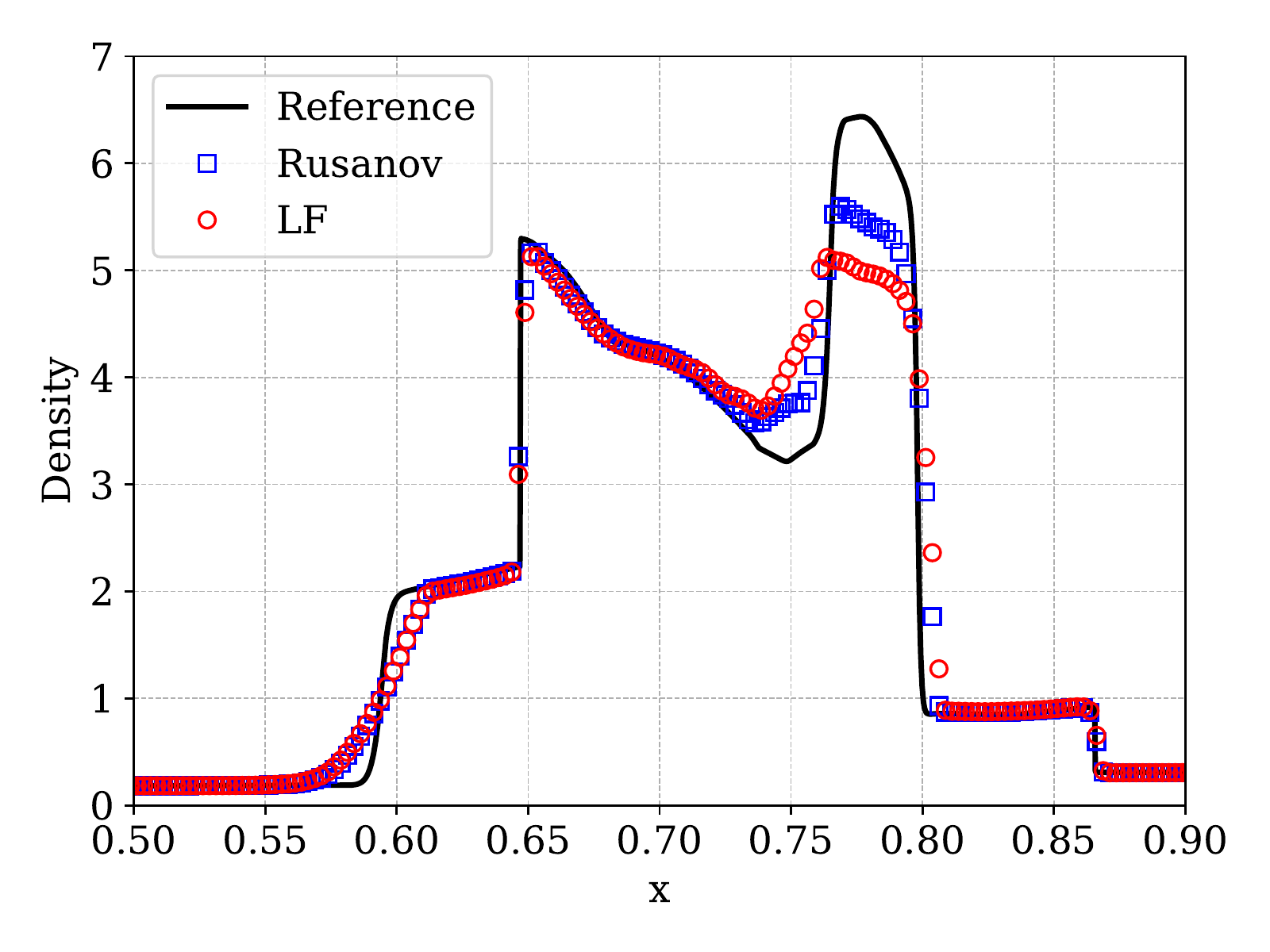} &
\includegraphics[width=0.45\textwidth]{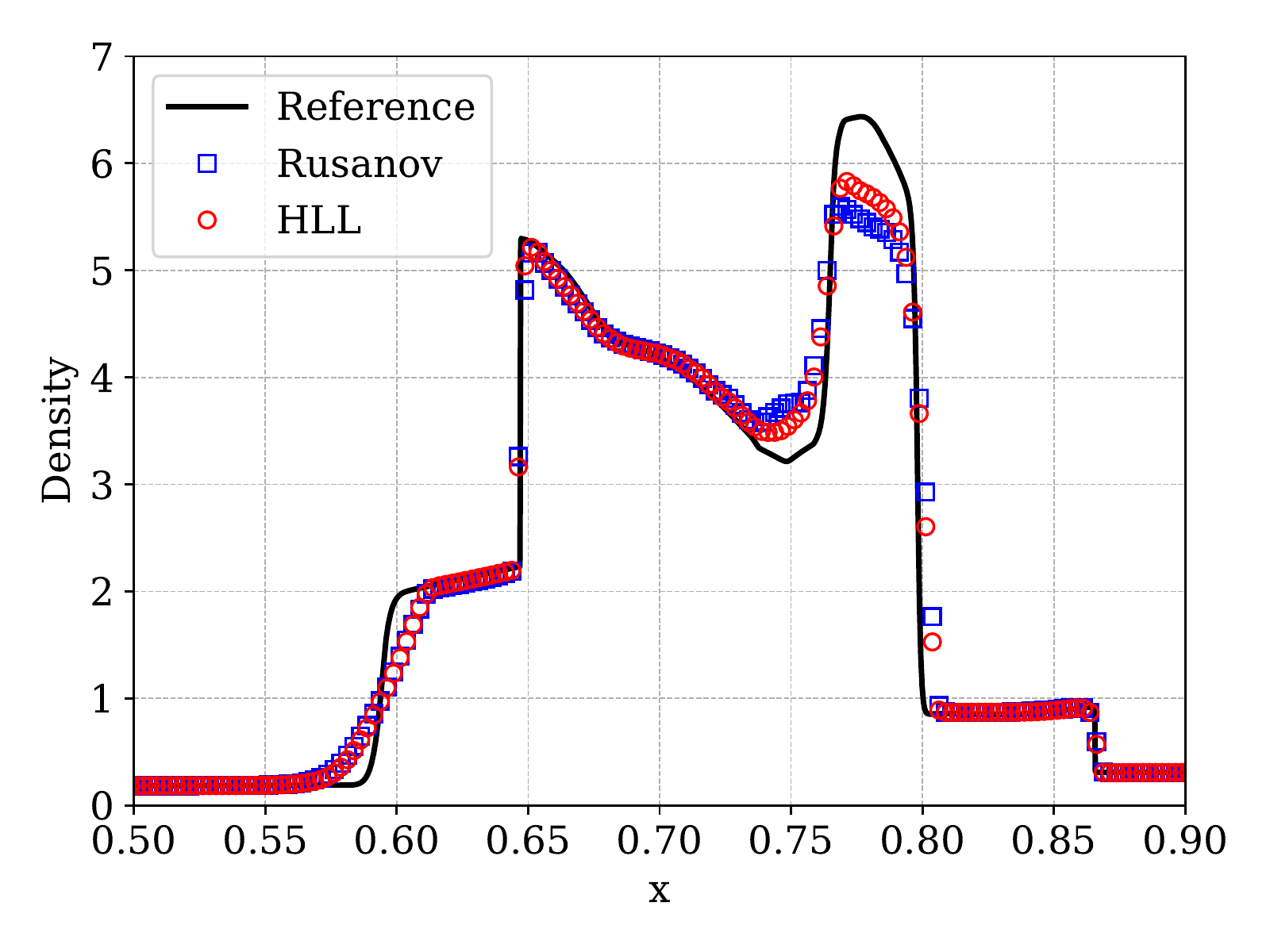} \\
(a) & (b) \\
\includegraphics[width=0.45\textwidth]{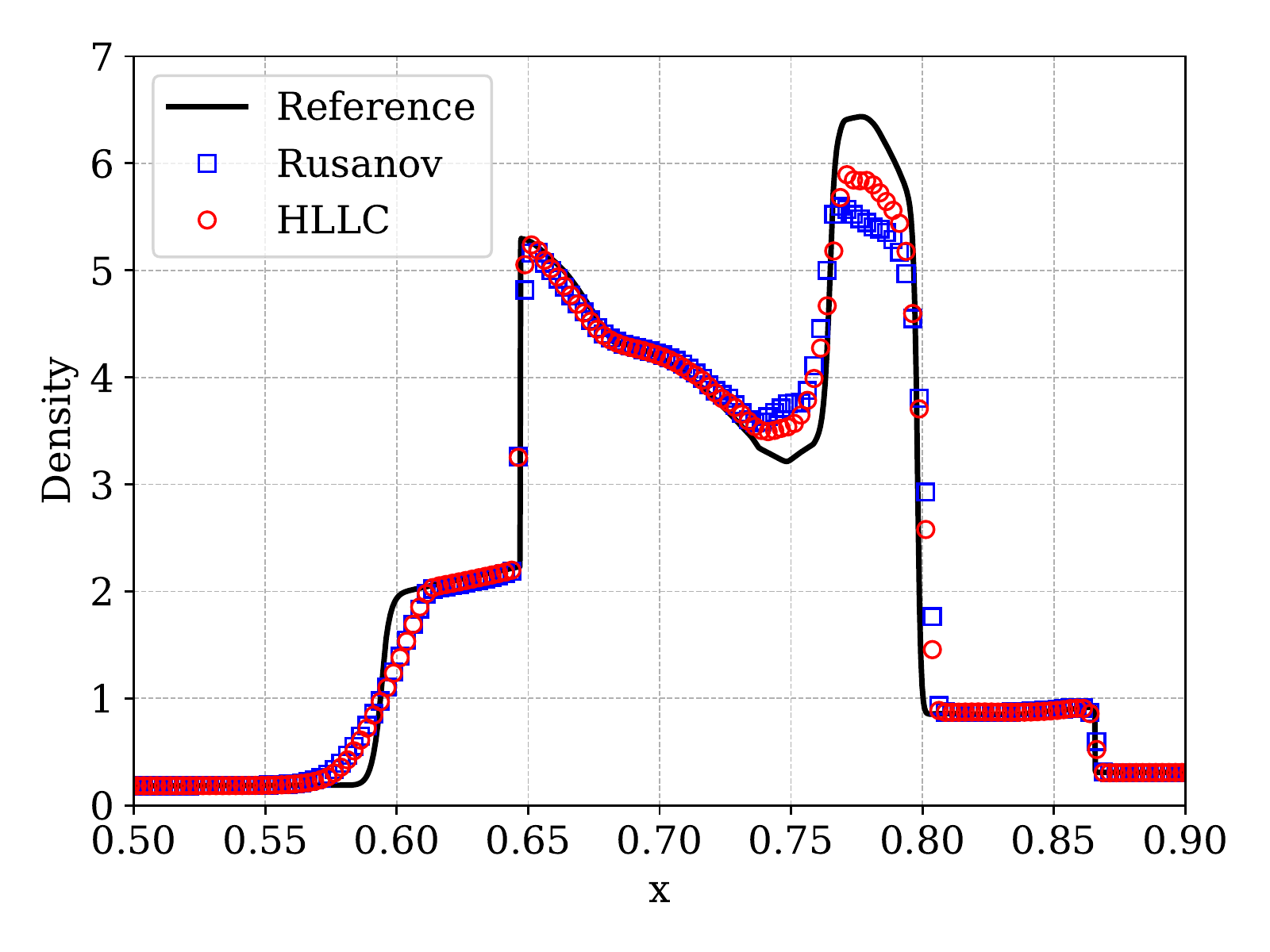} &
\includegraphics[width=0.45\textwidth]{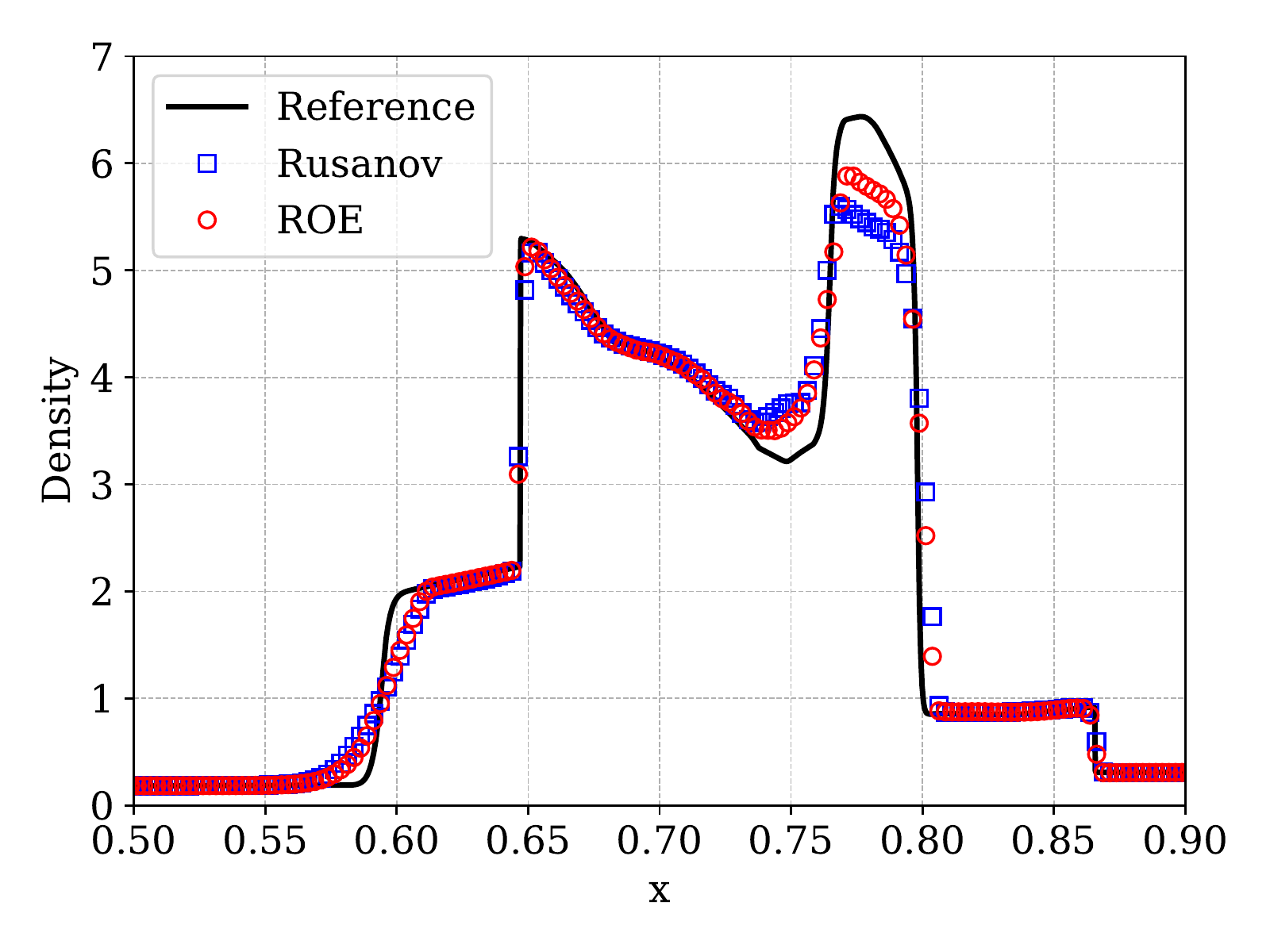} \\
(c) & (d)
\end{tabular}
\caption{Numerical solutions of 1-D Euler equations (Blast wave) obtained by LW schemes with different numerical fluxes (a) LF, (b) HLL, (c) HLLC and (d) ROE compared with Rusanov flux, for  $N=3$ using Radau correction function and GL solution points. All other parameters remain the same as in Figure~\ref{fig:blast}.}
\label{fig:numflx}
\end{figure}

The previous Euler results used Rusanov flux. In the LW scheme, we can also use other numerical fluxes  like HLL, HLLC, Roe and  global Lax-Friedrichs and compare its performance with the standard Rusanov flux. A description of these fluxes is given in the Appendix~\ref{apendix:numfluxes}. Fluxes like HLL, HLLC and Roe may be desirable in some problems due to their upwind character, unlike Lax-Friedrich/Rusanov type fluxes. Moreover HLLC and Roe fluxes also model the linear contact and shear waves which can lead to better approximations of these waves. We have tested the numerical fluxes in all the test cases, however to save space we present only the blast test case for $N=3$. The results are  given in Figure~(\ref{fig:numflx}) which compare these fluxes with Rusanov flux. The high density peak region is better approximated by the LW schemes using HLL, HLLC and Roe fluxes, as compared to the Rusanov flux. The global LF flux is found to be less accurate in this respect when compared to the Rusanov flux, which is expected due to the larger amount of numerical dissipation in the global Lax-Friedrich flux.

\begin{figure}
\centering
\begin{tabular}{cc}
\includegraphics[width=0.47\textwidth]{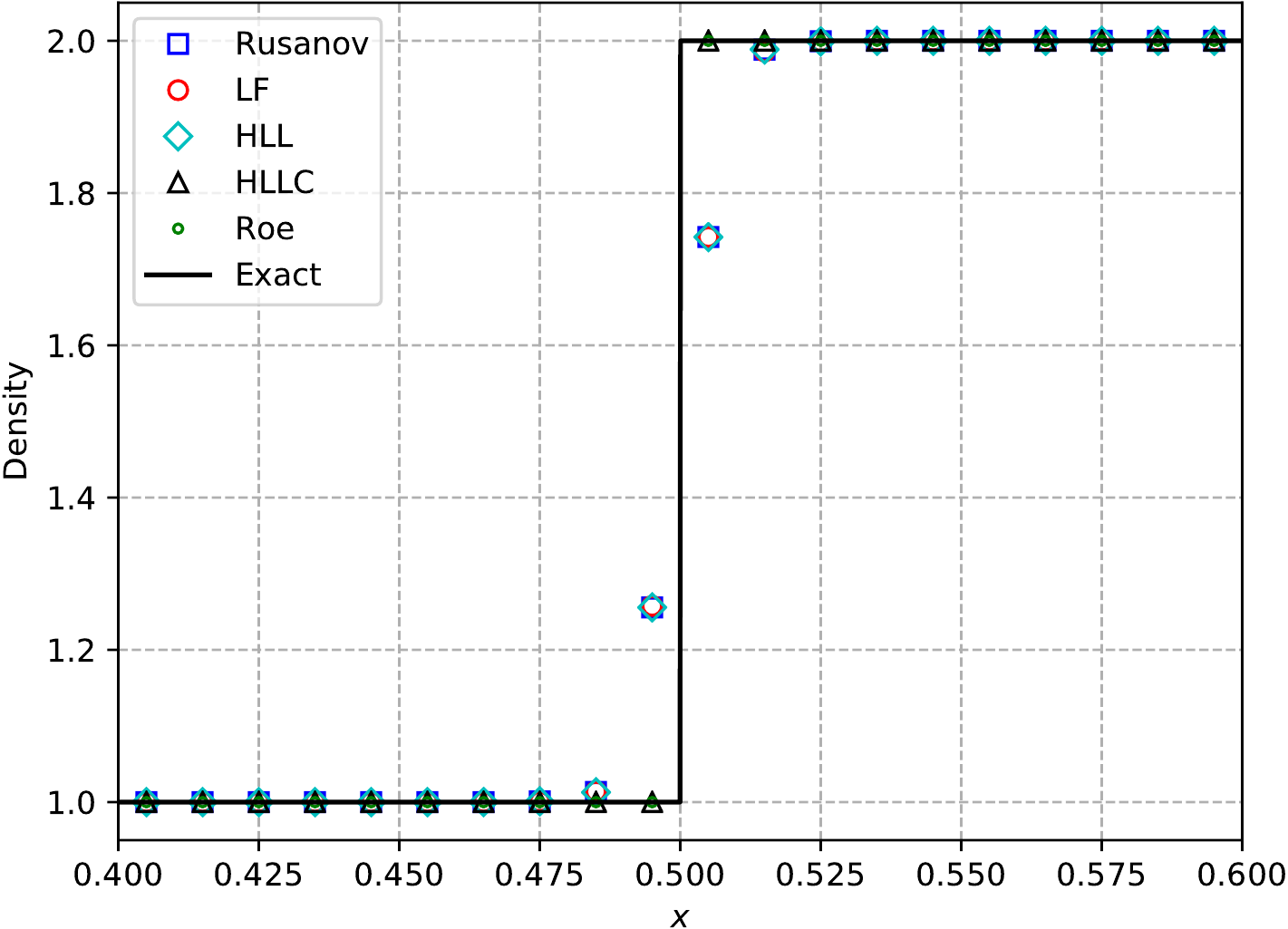} &
\includegraphics[width=0.47\textwidth]{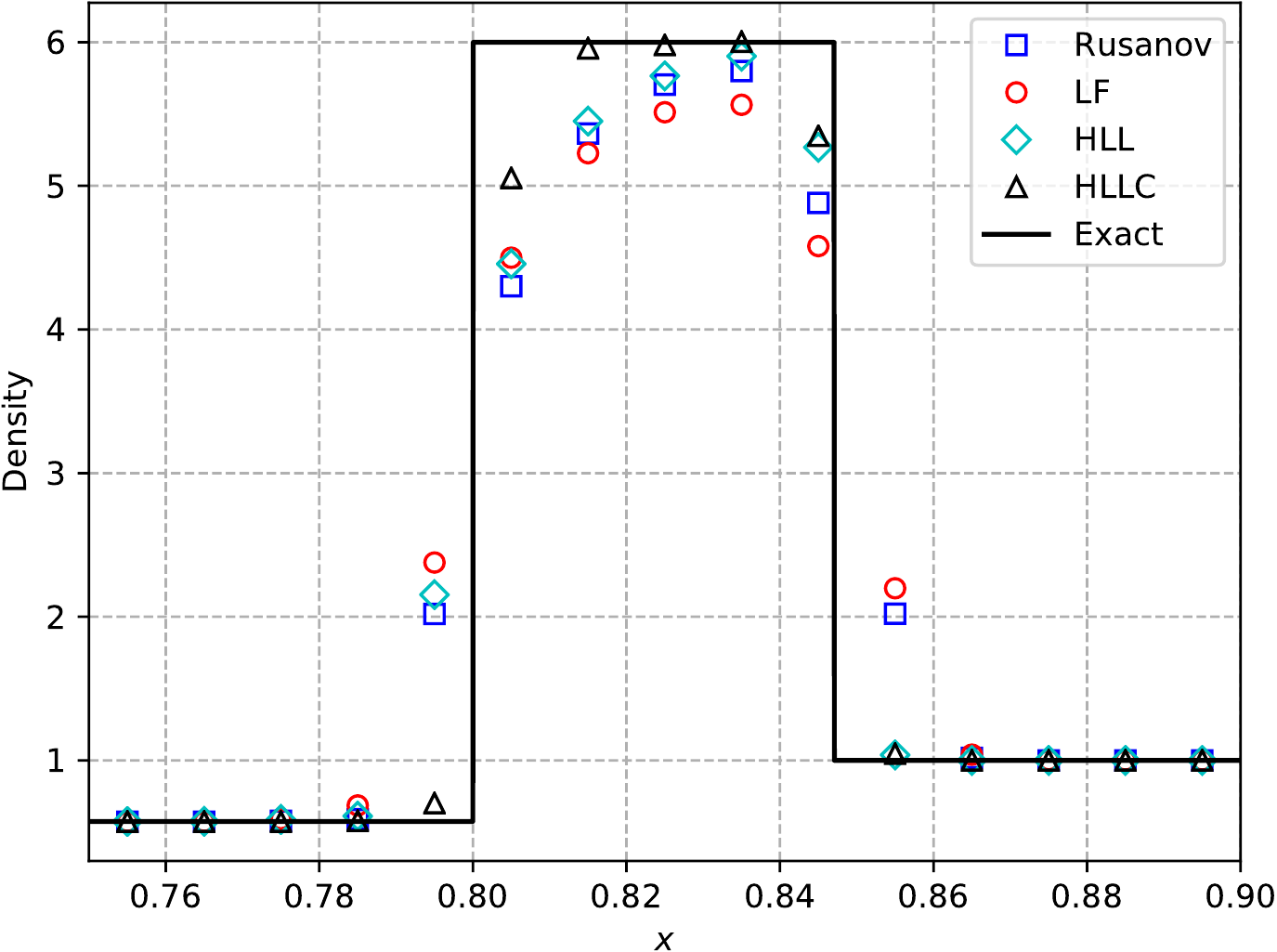} \\
(a) & (b)
\end{tabular}
\caption{Numerical solutions of 1-D Euler's equations for (a) stationary contact test, (b) Toro's Test 5 obtained by LW schemes with different numerical fluxes for polynomial degree $N=4$ using Radau correction function and GL solution points,  TVB $(M=1)$ limiter on a grid of 100 cells.}
\label{fig:numflx_toro5_contact}
\end{figure}

Since HLLC and Roe schemes contain more information about the wave structure, they are better at resolving contact discontinuities which are linearly degenerate waves that can be severely affected by numerical dissipation. We illustrate this through two Riemann problems containing stationary contact waves. The first one consists of an initial density jump that leads to a stationary contact wave, with initial condition given by,
\begin{equation*}
(\rho,v,p)=\begin{cases}
(1,0, 1), & \mbox{ if } x < 0.5\\
(2,0, 1), & \mbox{ if } x > 0.5
\end{cases}
\end{equation*}
In Figure~(\ref{fig:numflx_toro5_contact}a), we show the comparison of numerical fluxes for this stationary solution test case, zoomed near the discontinuity, at $t=1.0$ using LW schemes with D2 dissipation model for degree $N=4$ on a grid of 100 cells together with TVB ($M=1$) limiter. As expected, we see that Roe and HLLC fluxes are able to resolve the contact discontinuity exactly, while the other fluxes smear the jump over two cells.

The second Riemann problem is a tough test case with respect to maintaining positivity of pressure and is taken from~\cite{Toro2009}. The initial condition is given by
\begin{equation*}
(\rho,v,p)=\begin{cases}
(1,-19.59745, 1000), & \mbox{ if } x < 0.8\\
(1,-19.59745, 0.01), & \mbox{ if } x > 0.8
\end{cases}
\end{equation*}
The solution develops a stationary contact at the location of the initial discontinuity $x=0.8$ and a right moving shock wave. In Figure~(\ref{fig:numflx_toro5_contact}b), we show the comparison of numerical fluxes, zoomed near the contact discontinuity, at $t=0.012$ obtained using LW scheme with D2 dissipation model for polynomial degree $N=4$ on a grid of 100 cells and TVB ($M=1$) limiter. As in the previous case, the HLLC flux captures the contact discontinuity more accurately than the other fluxes.

\subsection{Comparison of correction functions}
\begin{figure}
\centering
\begin{tabular}{cc}
\includegraphics[width=0.45\textwidth]{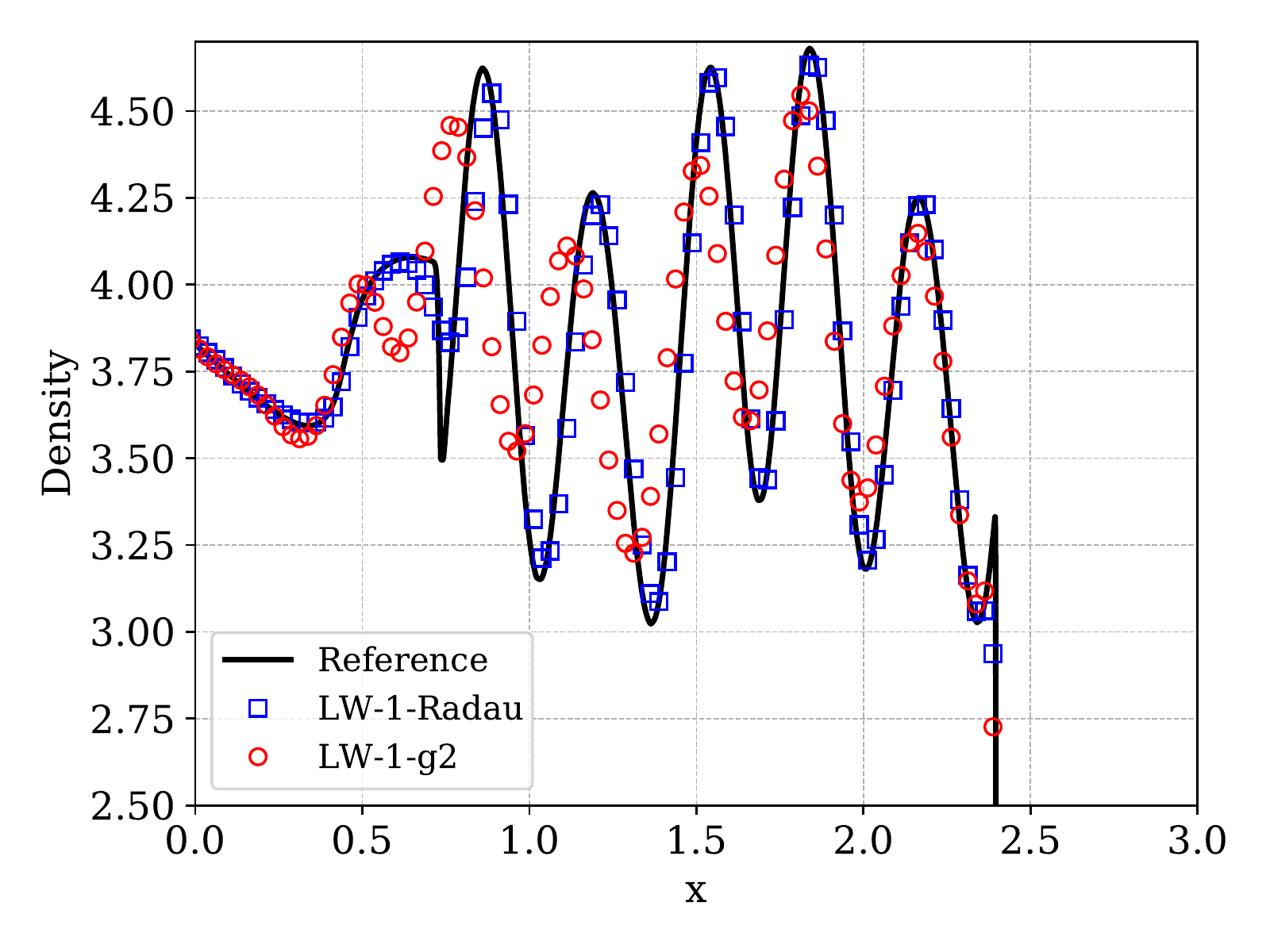} &
\includegraphics[width=0.45\textwidth]{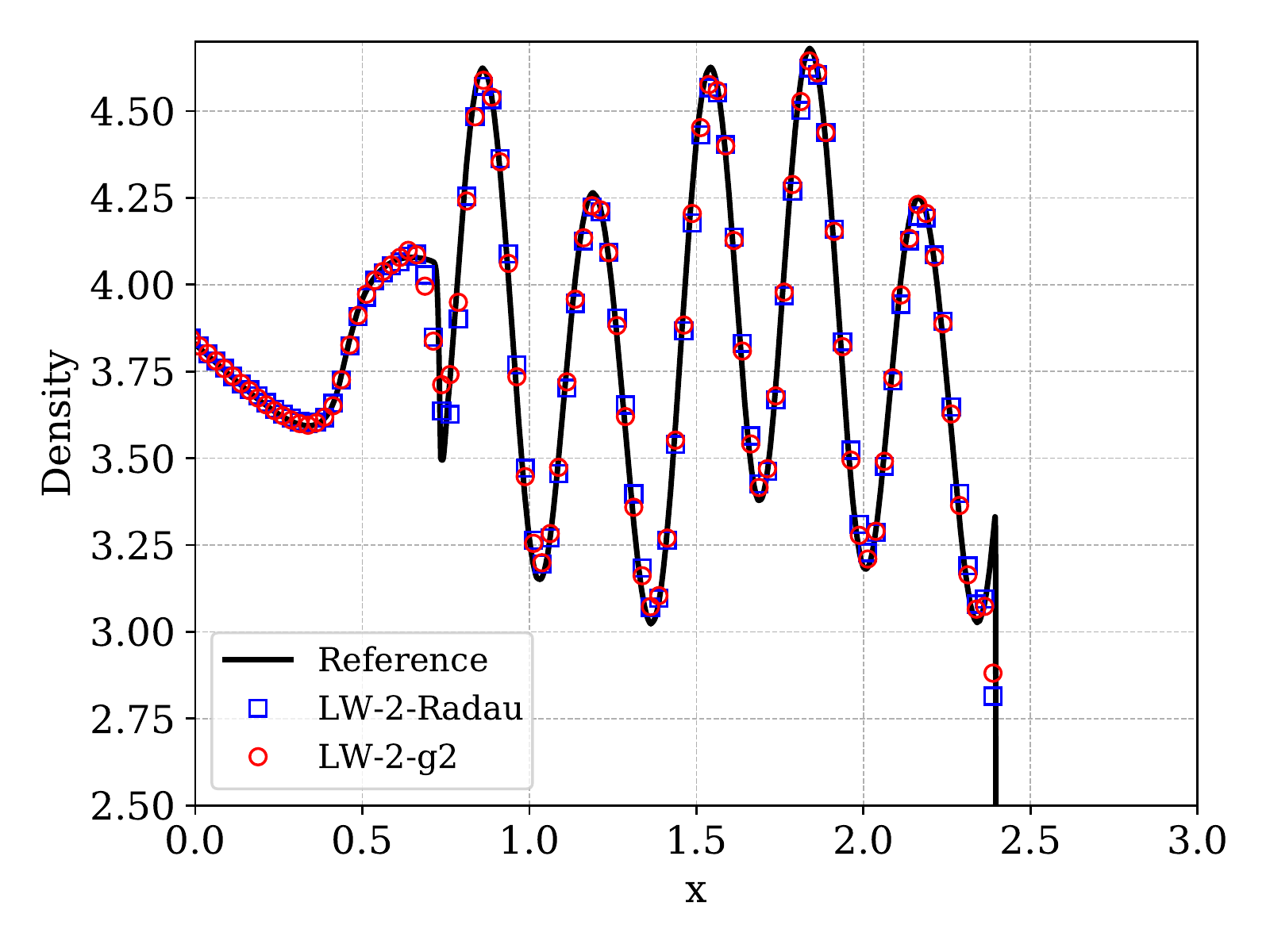} \\
(a) & (b) \\
\includegraphics[width=0.45\textwidth]{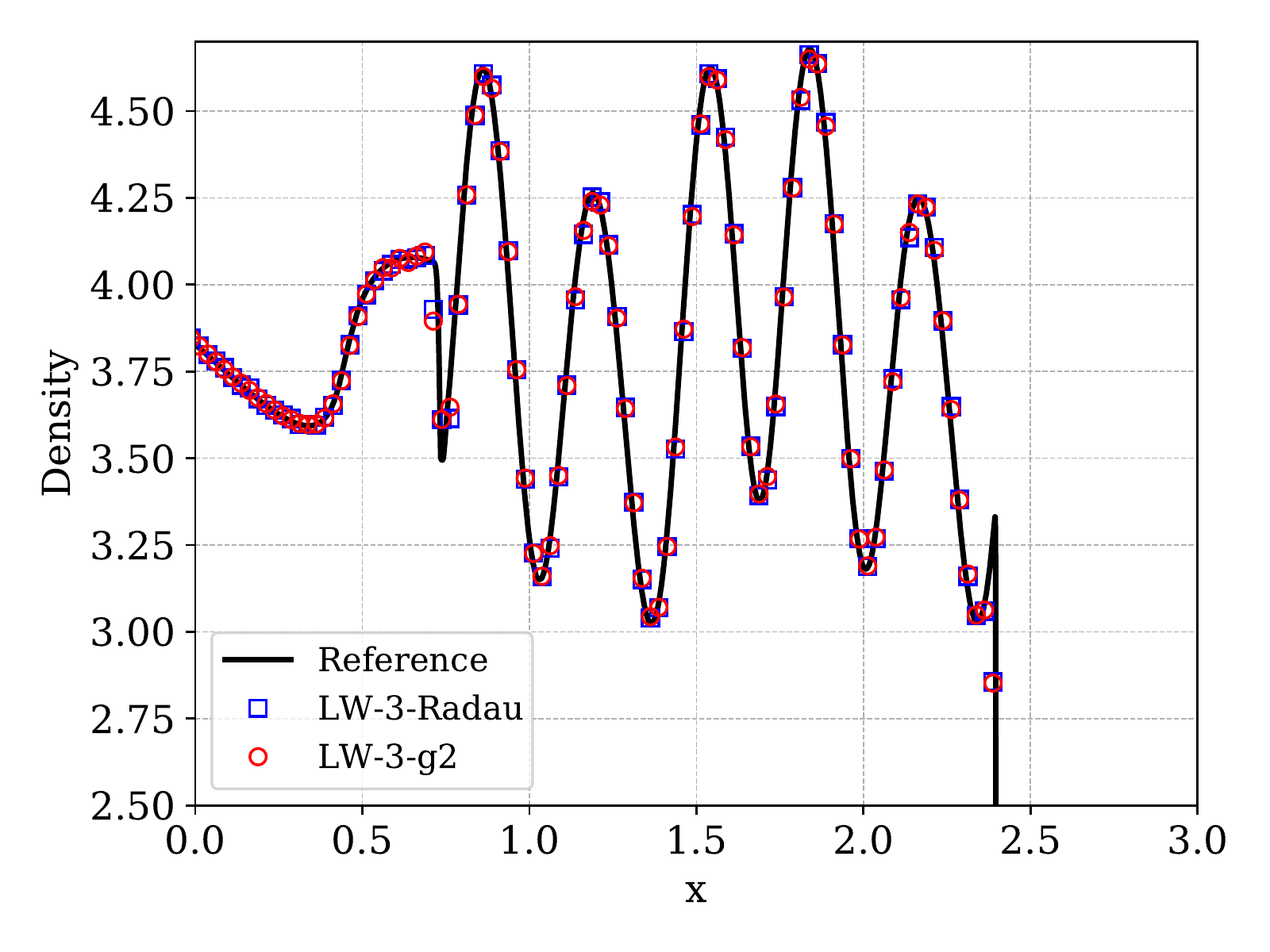} &
\includegraphics[width=0.45\textwidth]{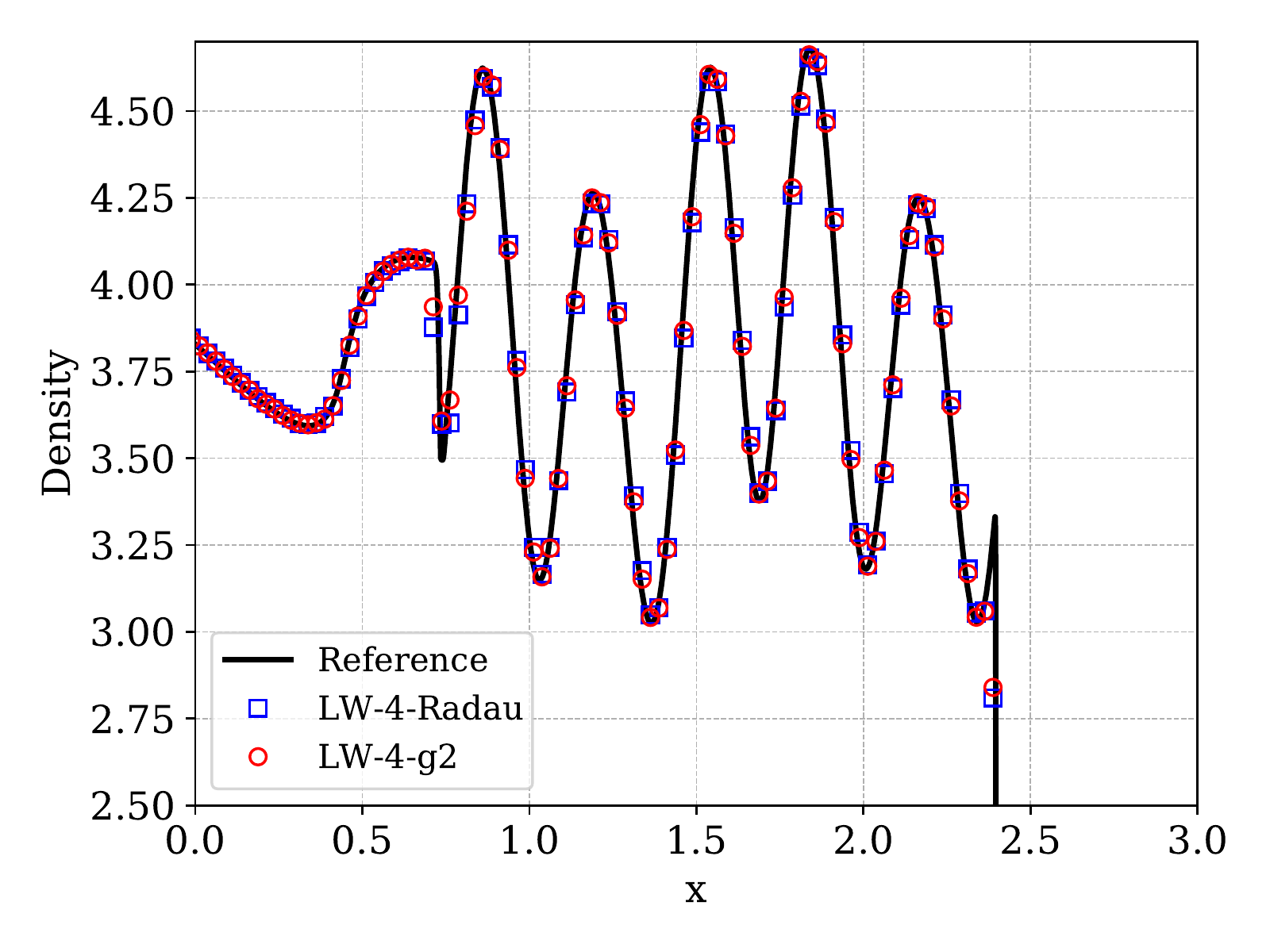} \\
(c) & (d)
\end{tabular}
\caption{Numerical solutions of 1-D Euler equations (Shu-Osher problem) for (a) $N=1$, (b) $N=2$, (c) $N=3$, (d) $N=4$. Comparison of  LW scheme with GL solution points for two correction functions, Radau and g2, with their own CFL numbers chosen from Table~\ref{tab:cfl}, except for g2 correction function with $N=1$, where we choose CFL=0.44. The enlarged  oscillatory portion of the solutions are shown. The solutions are computed at time $t=1.8$ on a mesh of 400 cells with dissipation model D2 and HLLC numerical flux. The TVB limiter is used with parameter $M=300$.}
\label{fig:ShuOsherCorr}
\end{figure}

We compare the robustness and accuracy of the two correction functions, Radau and g2, in the LW scheme when applied to the Euler equations~\eqref{eq:1deuler} with GL solution points. The numerical experiments are conducted for the Shu-Osher test case and the corresponding  results are obtained with the HLLC numerical flux for polynomial degrees $N=1,2,3,4$, see  Figure~(\ref{fig:ShuOsherCorr}). For this test case,  it is observed that the LW scheme with g2 correction function fails to work for $N=1$ with the optimal CFL of Table (\ref{tab:cfl}) due to loss of positivity of pressure. So, we use a smaller  CFL number of $0.44$  to compute the $N=1$ case in Figure~(\ref{fig:ShuOsherCorr}). For $N \ge 2$, the solutions computed using the g2 correction function are found to be close to that of Radau correction function. However, it fails to perform consistently, as we see in the $N=1$ case. With this observation and also the behaviour for other problems, we see that it is desirable to use the Radau correction function in the LW scheme.
\section{Two dimensional scheme}\label{sec:2d}
\begin{figure}
\centering
\includegraphics[width=0.7\textwidth]{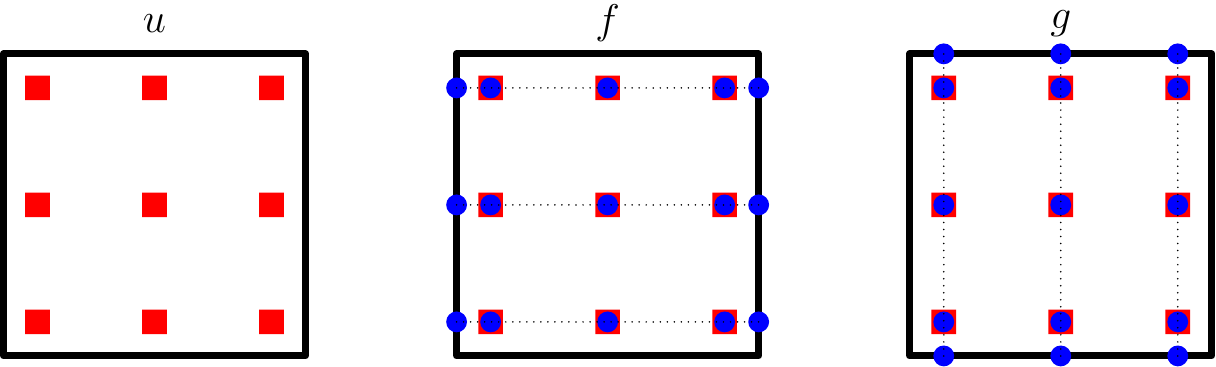}
\caption{Location of solution points and flux points. Numerical flux is required at the blue points on the faces.}
\label{fig:dofs2d}
\end{figure}
The extension of the 1-D scheme to two dimensions is performed by applying the 1-D ideas along each coordinate direction. Consider a 2-D conservation law of the form
\begin{equation}
u_t + f(u)_x + g(u)_y = 0
\end{equation}
where $(f,g)$ are Cartesian components of the flux vector. Using Taylor expansion in time, we can write the solution at $t=t_{n+1}$ as
\[
u^{n+1} = u^n - \Delta t \left[ \pd{F}{x}(u^n) + \pd{G}{y}(u^n) \right] + O(\Delta t^{N+2})
\]
where $F,G$ are time average fluxes given by
\begin{equation*}
F(u) = \sum_{m=0}^N \frac{\Delta t^m}{(m+1)!} \partial_t^m f(u), \qquad
G(u) = \sum_{m=0}^N \frac{\Delta t^m}{(m+1)!} \partial_t^m g(u) 
\end{equation*}
We will consider a Cartesian mesh and map each element $\Omega_e$ to the reference element $\Kref = [0,1] \times [0,1]$. Inside the reference element, the solution points are chosen to be tensor product of 1-D solution points, which may be either GL or GLL points. Figure~(\ref{fig:dofs2d}) shows an example of 2-D solution points based on tensor product of 1-D GL points. The solution inside an element $\Omega_e$ is approximated by a tensor product polynomial of degree $N$,
\begin{equation*}
(x,y) \in \Omega_e : \qquad u_h = \sum_{i=0}^N \sum_{j=0}^N u_{ij}^e \ell_i(\xi) \ell_j(\eta)
\end{equation*}
where $(\xi,\eta)$ are coordinates in the reference element, and $\ell_i(\xi), \ell_j(\eta)$ are the 1-D Lagrange polynomials based on the solution points. The discontinuous fluxes are approximated by interpolating at the solution points,
\[
F_h^\delta(\xi,\eta) = \sum_{i=0}^N \sum_{j=0}^N F_{ij}^e \ell_i(\xi) \ell_j(\eta), \qquad
G_h^\delta(\xi,\eta) = \sum_{i=0}^N \sum_{j=0}^N G_{ij}^e \ell_i(\xi) \ell_j(\eta)
\]
where $F_{ij}^e, G_{ij}^e$ are time average fluxes obtained from the Lax-Wendroff procedure applied at each solution point.  The continuous flux along the $x$ and $y$ axes are constructed using the one dimensional algorithm along the $\eta=\xi_j =$ constant and $\xi=\xi_i=$ constant lines, respectively, see Figure~(\ref{fig:dofs2d}),
\[
F_h(\xi,\xi_j) = [F_{\emh,j} - F_h^\delta(0,\xi_j)] g_L(\xi) + F_h^\delta(\xi,\xi_j) + [F_{\eph,j} - F_h^\delta(1,\xi_j)] g_R(\xi), \qquad 0 \le j \le N
\]
\[
G_h(\xi_i,\eta) = [G_{\emh,i} - G_h^\delta(\xi_i,0)] g_L(\eta) + G_h^\delta(\xi_i,\eta) + [G_{\eph,i} - G_h^\delta(\xi_i,1)] g_R(\eta), \qquad 0 \le i \le N
\]
Note that the above equations are obtained by applying the FR idea along the horizontal and vertical lines in Figure~(\ref{fig:dofs2d}). The quantities $F_\emh, F_\eph$ are $x$-directional numerical fluxes on the left and right faces, while $G_\emh, G_\eph$ are the $y$-directional numerical fluxes across the bottom and top faces, respectively.  The update equation is given by a collocation procedure at each solution point,
\begin{equation}
(u_{ij}^e)^{n+1} = (u_{ij}^e)^n - \Delta t \left[ \frac{1}{\Delta x_e} \pd{F_h}{\xi}(\xi_i,\xi_j) + \frac{1}{\Delta y_e} \pd{G_h}{\eta}(\xi_i,\xi_j) \right], \qquad 0 \le i,j \le N
\end{equation}
where the flux derivatives can be computed from
\[
\partial_\xi F_h(:,\xi_j) = \left[F_{\emh,j} - F_h^\delta(0,\xi_j) \right] \vb_L + \partial_\xi F_h^\delta(:,\xi_j) + \left[F_{\eph,j} - F_h^\delta(1,\xi_j) \right] \vb_R, \qquad 0 \le j \le N
\]
\[
\partial_\eta G_h(\xi_i,:) = \left[G_{\emh,i} - G_h^\delta(\xi_i,0) \right] \vb_L + \partial_\eta G_h^\delta(\xi_i,:) + \left[G_{\eph,i} - G_h^\delta(\xi_i,1) \right] \vb_R, \qquad 0 \le i \le N
\]
We can cast the update equation in matrix form. For this, define the flux matrices
\[
\vF_e(i,j) = F_{ij}^e, \qquad \vG_e(i,j) = G_{ij}^e, \qquad 0 \le i,j \le N
\]
Then we can compute the derivatives of the discontinuous flux at all the solution points by a matrix-matrix product
\[
\partial_\xi F_h^\delta(:,:) = \vD \vF_e, \qquad \partial_\eta G_h^\delta(:,:) = \vG_e \vD^\top
\]
where $\vD$ is the 1-D differentiation matrix. The update equation can be written in matrix form,
\begin{equation}\label{eq:up2d}
\vu_e^{n+1} = \vu_e^n - \left[ \frac{\Delta t}{\Delta x_e} \vD_1 \vF_e + \frac{\Delta t}{\Delta y_e} \vG_e \vD_1^\top \right] - \frac{\Delta t}{\Delta x_e} \left[ \vb_L \vF_\emh^\top + \vb_R \vF_\eph^\top \right] - \frac{\Delta t}{\Delta y_e} \left[ \vG_\emh \vb_L^\top + \vG_\eph \vb_R^\top \right]
\end{equation}
where the quantities $\vD_1, \vb_L, \vb_R$ have been defined before in the description of the 1-D scheme. 

The time average fluxes are computed by the approximate Lax-Wendroff procedure. To describe this, let us define the flux matrices
\[
\vf_e(i,j) = f(u_{ij}^e), \qquad \vg_e(i,j) = g(u_{ij}^e), \qquad 0 \le i,j \le N
\]
The time derivatives of the solution at all solution points are obtained from the PDE by the following matrix equation,
\[
\vu^{(m)}_e = - \frac{\Delta t}{\Delta x_e} \vD \vf^{(m-1)}_e - \frac{\Delta t}{\Delta y_e} \vg^{(m-1)}_e \vD^\top, \quad m=1,2,\ldots,N
\]
and the time average solution and fluxes are given by
\[
\vU_e = \sum_{m=0}^N \frac{\vu_e^{(m)}}{(m+1)!}, \qquad
\vF_e = \sum_{m=0}^N \frac{\vf_e^{(m)}}{(m+1)!}, \qquad
\vG_e = \sum_{m=0}^N \frac{\vg_e^{(m)}}{(m+1)!}
\]
The time derivatives of the fluxes $\vf_e^{(m)}$, $\vg_e^{(m)}$ are approximated using finite differences in time as in the 1-D case given in Section~\ref{sec:alw}; those formulae are applied to both components of the flux.  The stable time step is determined by considering the linear advection equation in 2-D and applying Fourier stability analysis to the LW scheme, see Appendix~\ref{sec:fourier2d}.
\section{Numerical results in 2-D : scalar problems}\label{sec:res2d}
We present results to test the error convergence properties of the LW schemes for some 2-D problems and compare them to RK scheme. For each problem in this section, the corresponding CFL numbers are chosen based on Fourier stability analysis which are given in Table~\ref{tab:2Dcfl}. We compare  Lax-Wendroff scheme with D2 dissipation model and Runge-Kutta schemes in this section, and the CFL numbers of the former are used for both schemes. For the RKFR scheme, we use SSP Runge-Kutta time integration~\cite{Gottlieb2001} for $N=1$ and 2,  the classical four stage Runge-Kutta method of order four for  $N=3$, and six-stage, fifth order Runge-Kutta (RK65) time integration for $N=4$~\cite{Tsitouras2011} implemented in {\tt DifferentialEquations.jl}~\cite{Rackauckas2017}. All the results in this section are produced using code written in Julia~\cite{Bezanson2017}; the design and optimization of the code was inspired by {\tt Trixi.jl}~\cite{Ranocha2022}.
\subsection{Advection of a smooth signal}
\begin{figure}
\begin{center}
\begin{tabular}{cc}
\includegraphics[width=0.40\textwidth]{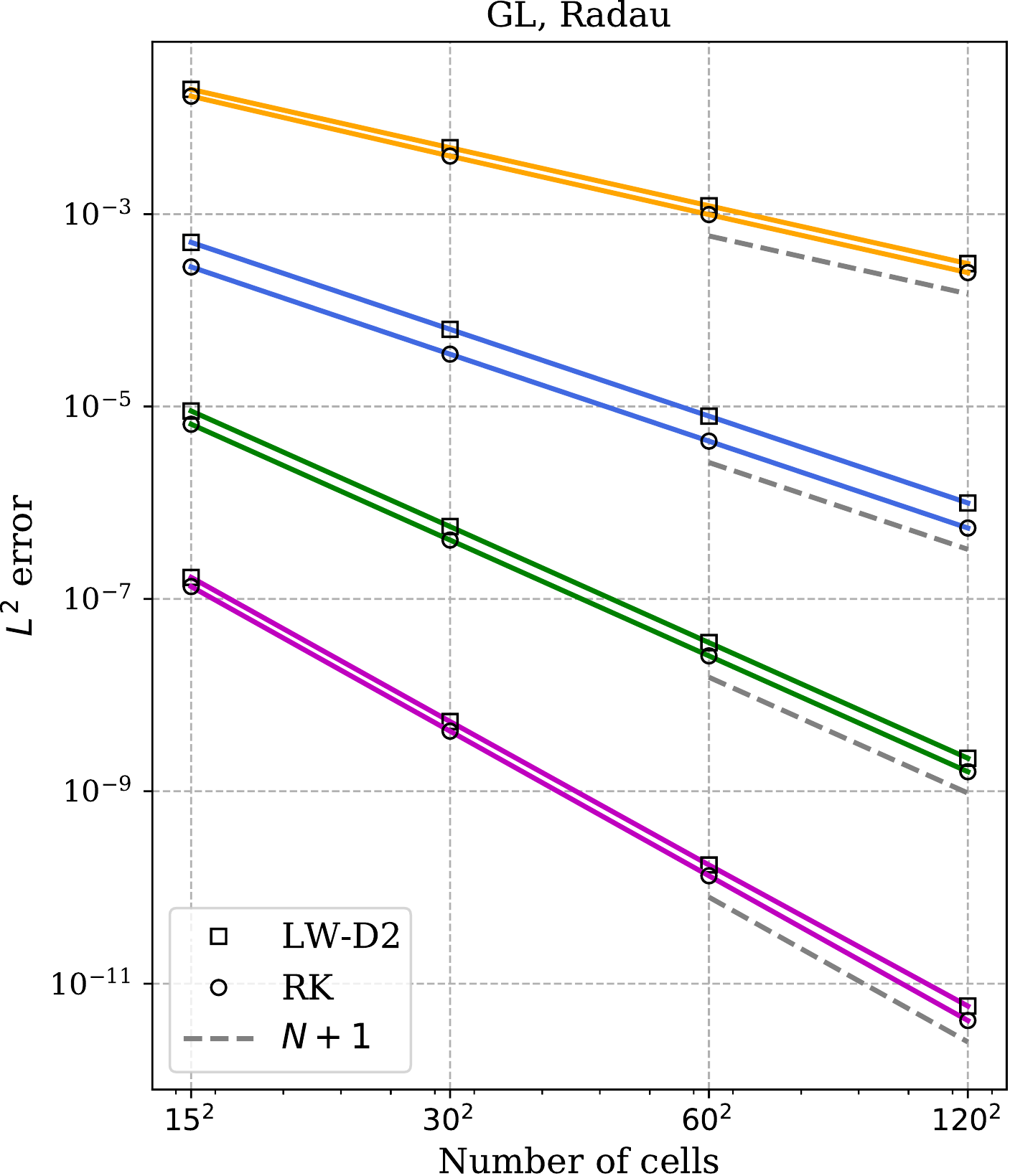} &
\includegraphics[width=0.40\textwidth]{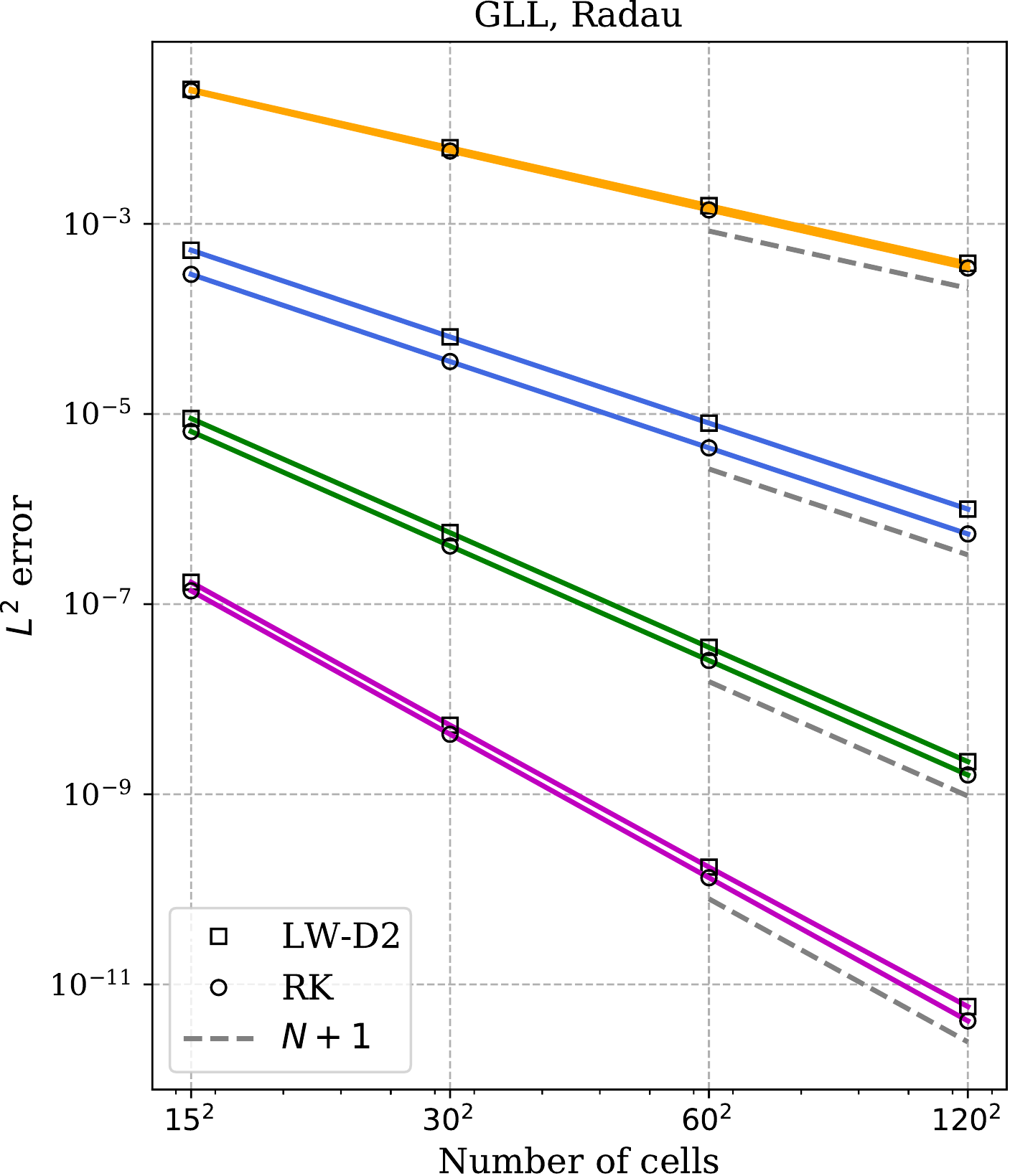} \\
(a)  & (b)
\end{tabular}
\end{center}
\caption{Error convergence test for 2-D linear advection equation with velocity $\bm{a}=(1,1)$ at $t=1.0$, initial data $u_0(x,y)=\sin(2\pi x)\sin(2\pi y)$; (a) GL points, (b) GLL points. The different colors correspond to degrees $N=1,2,3,4$ from top to bottom.}
\label{fig:conv linear adv 2d}
\end{figure}

We consider  the advection equation in two dimensions
\begin{equation}\label{eq:2dvaradv}
u_t+\nabla \cdot [\bm{a}(x,y)u] =0,
\end{equation}
with two types of divergence-free advection velocity, namely  a constant velocity $\bm a=(1,1)$ and a variable velocity $\bm a=(-y,x)$. For the second velocity, the flux components are $(f,g) = (-y u, x u)$ so that $F_h^\delta$ is of degree $N$ in $x$ and $G_h^\delta$ is of degree $N$ in $y$ variable. Since $F_h^\delta$ is interpolated along $x$ direction and $G_h^\delta$ is interpolated along $y$ direction, there is no interpolation error due to non-linearity of the flux and the EA and AE schemes are equivalent. Due to this reason, we only show the results with EA scheme. In order to verify the accuracy of the LWFR scheme we consider the equation~\eqref{eq:2dvaradv}  with a smooth initial condition and perform the simulation for both the advection velocities. For the velocity $\bm a=(1,1)$, the characteristic curves are straight lines and we use periodic boundary conditions on the domain $[0,1]\times[0,1]$ with initial condition $u_0(x,y)=\sin(2\pi x)\sin(2\pi y)$. The error convergence plots are shown in Figure~(\ref{fig:conv linear adv 2d}) using Radau correction function. The optimal convergence rate is attained by both LW and RK schemes and there is no significant difference between GL and GLL points. The errors of RK scheme are slightly smaller than those of the LW scheme, similar to the 1-D case.

\begin{figure}
\begin{center}
\begin{tabular}{ccc}
\includegraphics[width=0.32\textwidth]{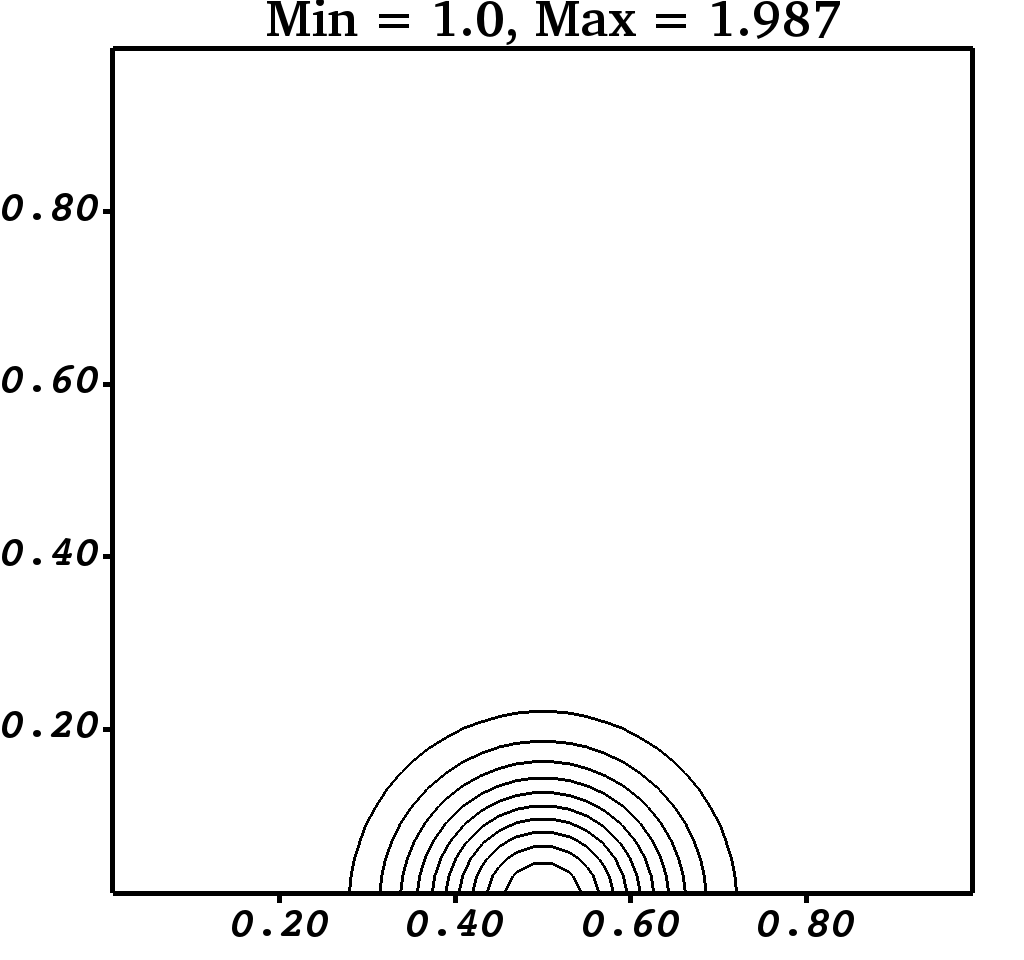} &
\includegraphics[width=0.32\textwidth]{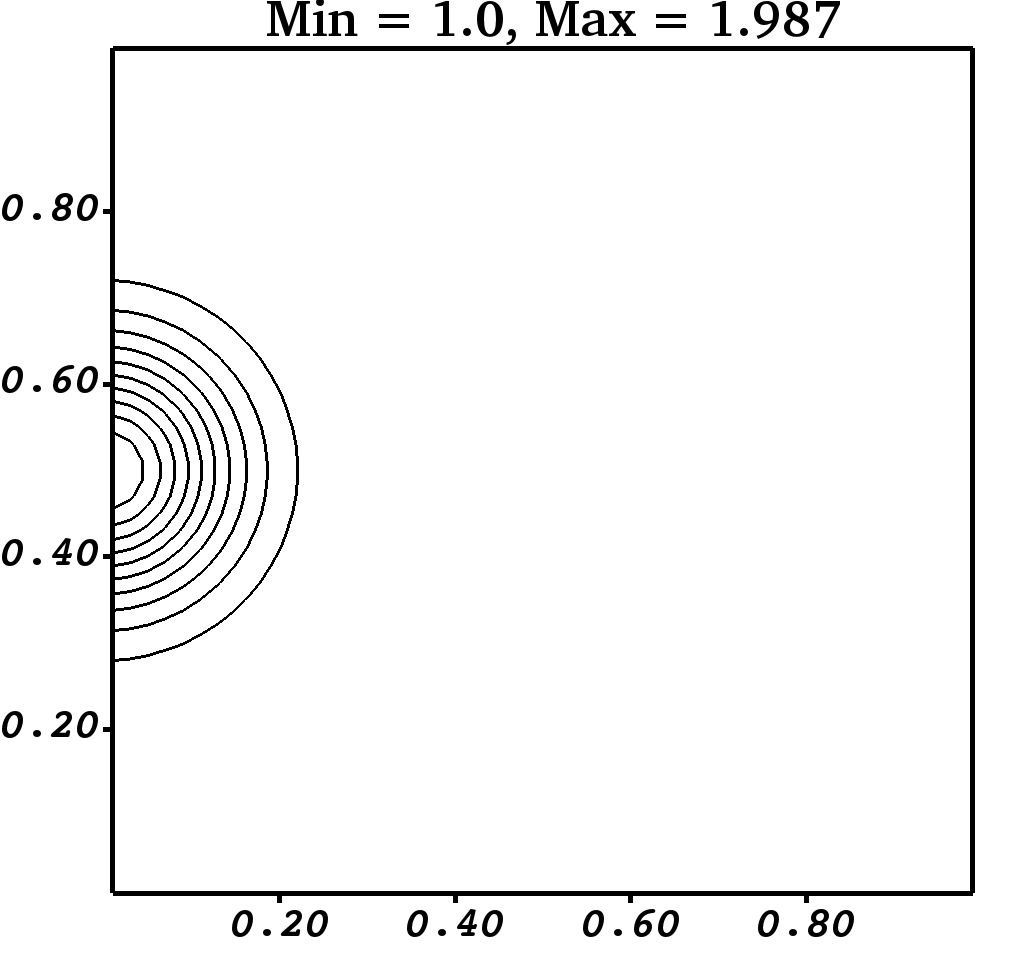} &
\includegraphics[width=0.32\textwidth]{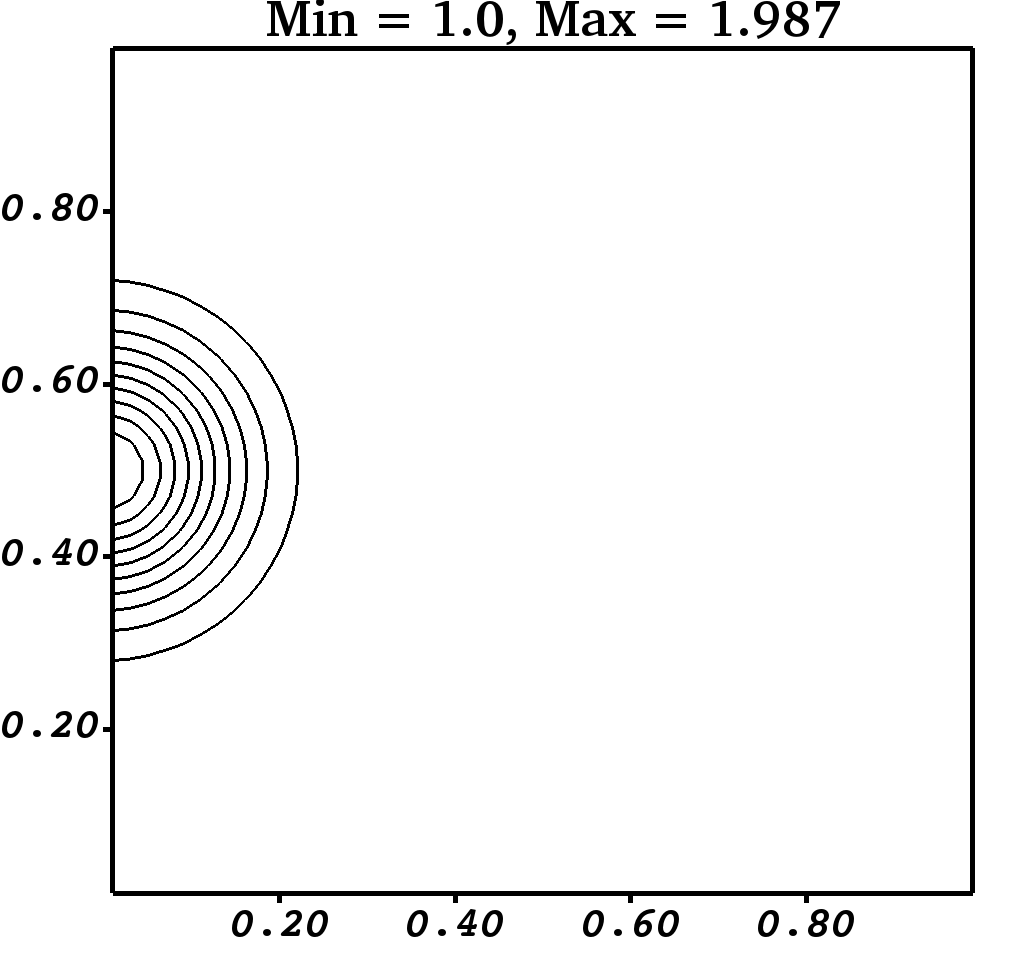} \\
(a) & (b) & (c)
\end{tabular}
\end{center}
\caption{Linear advection with velocity $\bm{a}=(-y,x)$ on $[0,1]\times[0,1]$ with inflow/outflow boundary condition. The solutions are shown on a mesh of $50 \times 50$ cells with polynomial degree $N=3$; (a) initial solution, (b) LWFR, $t= \frac{\pi}{2}$ (c) RKFR, $t= \frac{\pi}{2}$.}
\label{fig:lin2d_rotate_soln}
\end{figure}

\begin{figure}
\begin{center}
\begin{tabular}{cc}
\includegraphics[width=0.40\textwidth]{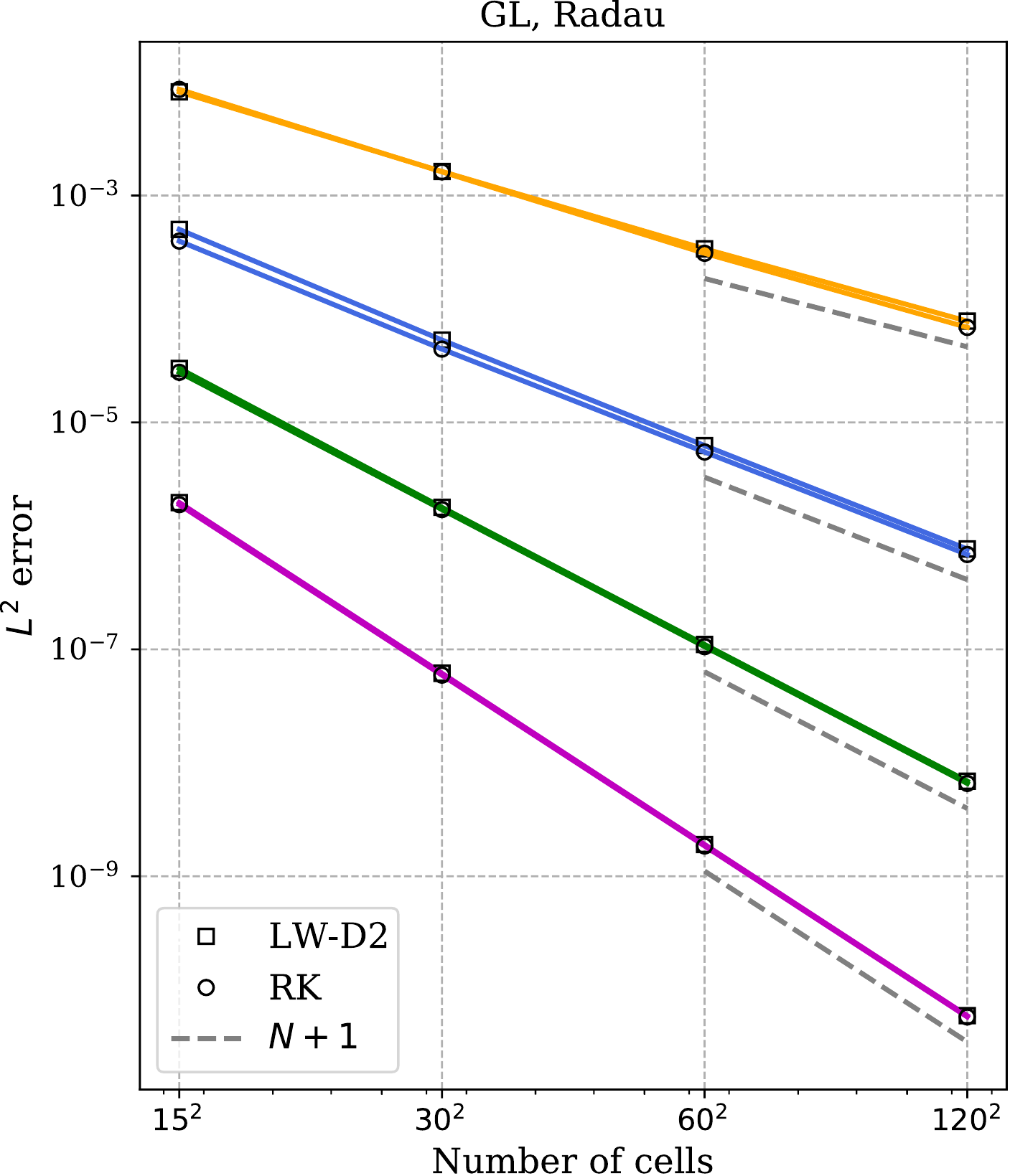} &
\includegraphics[width=0.40\textwidth]{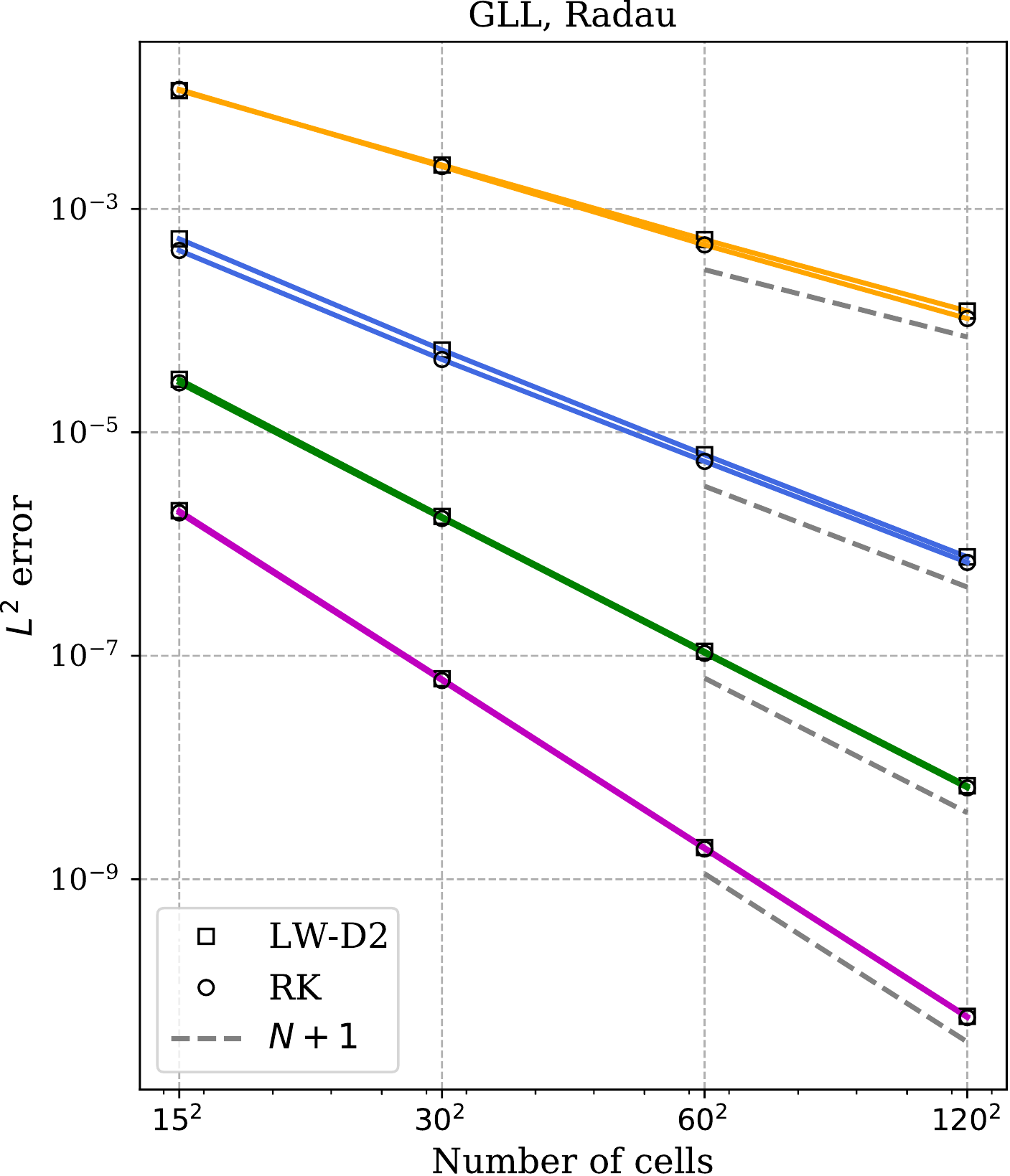} \\
(a) & (b)
\end{tabular}
\end{center}
\caption{Error convergence test for 2-D linear advection equation with velocity $\bm{a}=(-y,x)$ at $t=\frac{\pi}{2}$, initial data $u_0(x,y)=1+\exp(-50((x-1/2)^2+y^2))$ using (a) GL points, (b) GLL points.}
\label{fig:conv linear rotate 2d}
\end{figure}

For the variable velocity $\bm a=(-y,x)$, the characteristic curves are circles whose center is at the origin and we take the domain $\Omega = [0,1] \times [0,1]$. The exact solution is given by $u(x,y,t)=u_0(x\cos(t)+y\sin(t),-x\sin(t)+y\cos(t))$ with  the initial condition $u_0(x,y)=1+\exp(-50((x-1/2)^2+y^2))$.  At the bottom and right boundaries, we use inflow conditions while on top and left side of the boundary, we use outflow conditions. The initial condition advects along the circular characteristic curves in the counter clock-wise direction.   A contour plot of the numerical solution is visualized in Figure~(\ref{fig:lin2d_rotate_soln}), and the error convergence analysis is made in Figure~(\ref{fig:conv linear rotate 2d}). The error convergence agrees with the optimal convergence rates and the error values of the LW scheme are comparable to those from the RK scheme at all orders shown in the figures.

\subsection{Rotation of a composite signal}
\begin{figure}
\begin{subfigure}{0.33\textwidth}
\centering
\includegraphics[width=\textwidth]{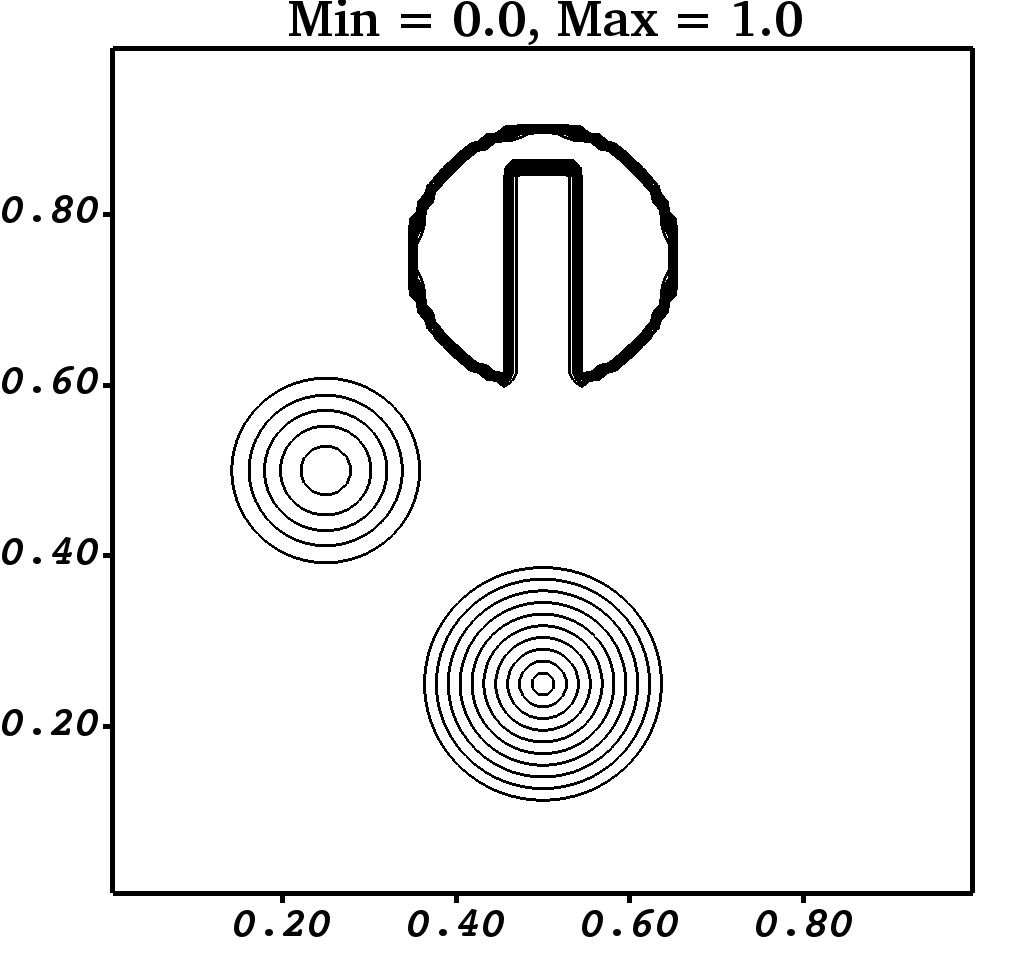}
\caption*{(a) Exact}
\end{subfigure}%
\begin{subfigure}{0.33\textwidth}
\begin{subfigure}{\textwidth}
\renewcommand\thesubfigure{\alph{subfigure}1}
\centering
\includegraphics[width=\textwidth]{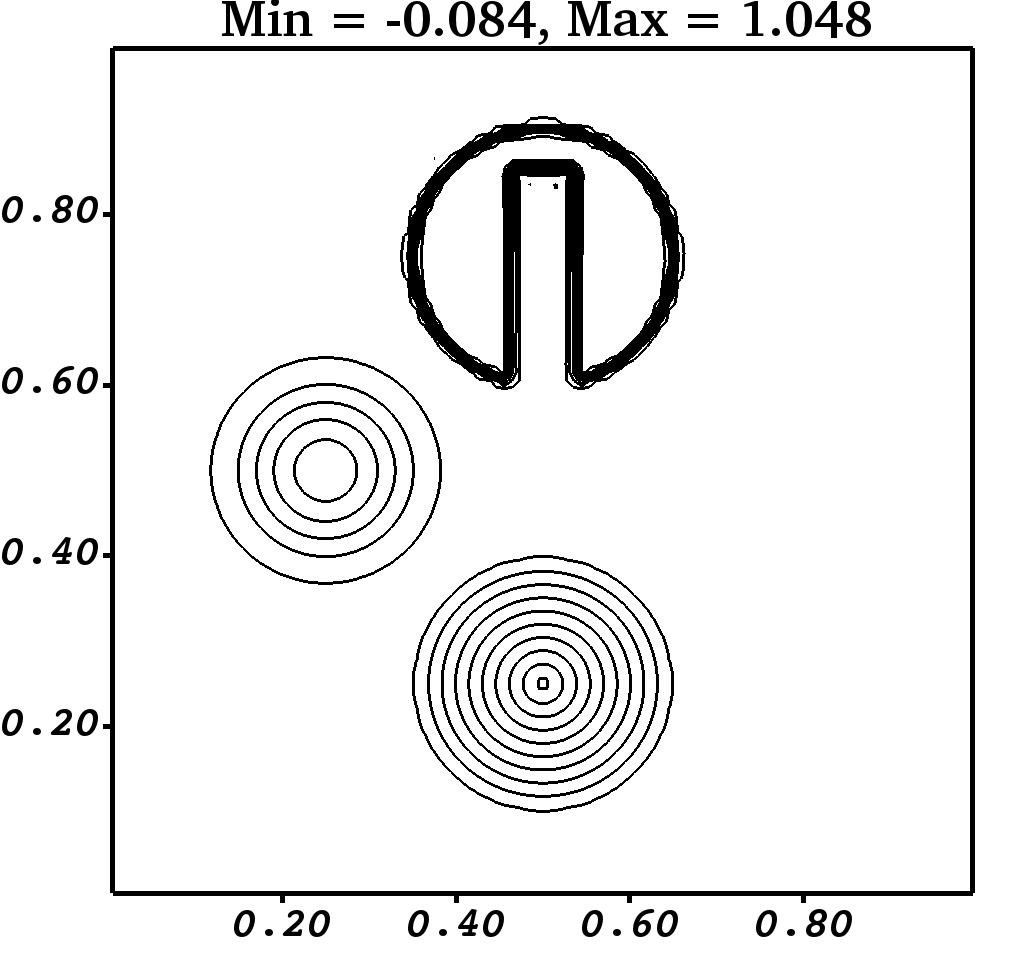}
\caption*{(b) LWFR}
\end{subfigure}
\begin{subfigure}{\textwidth}
\renewcommand\thesubfigure{\alph{subfigure}2}
\centering
\includegraphics[width=\textwidth]{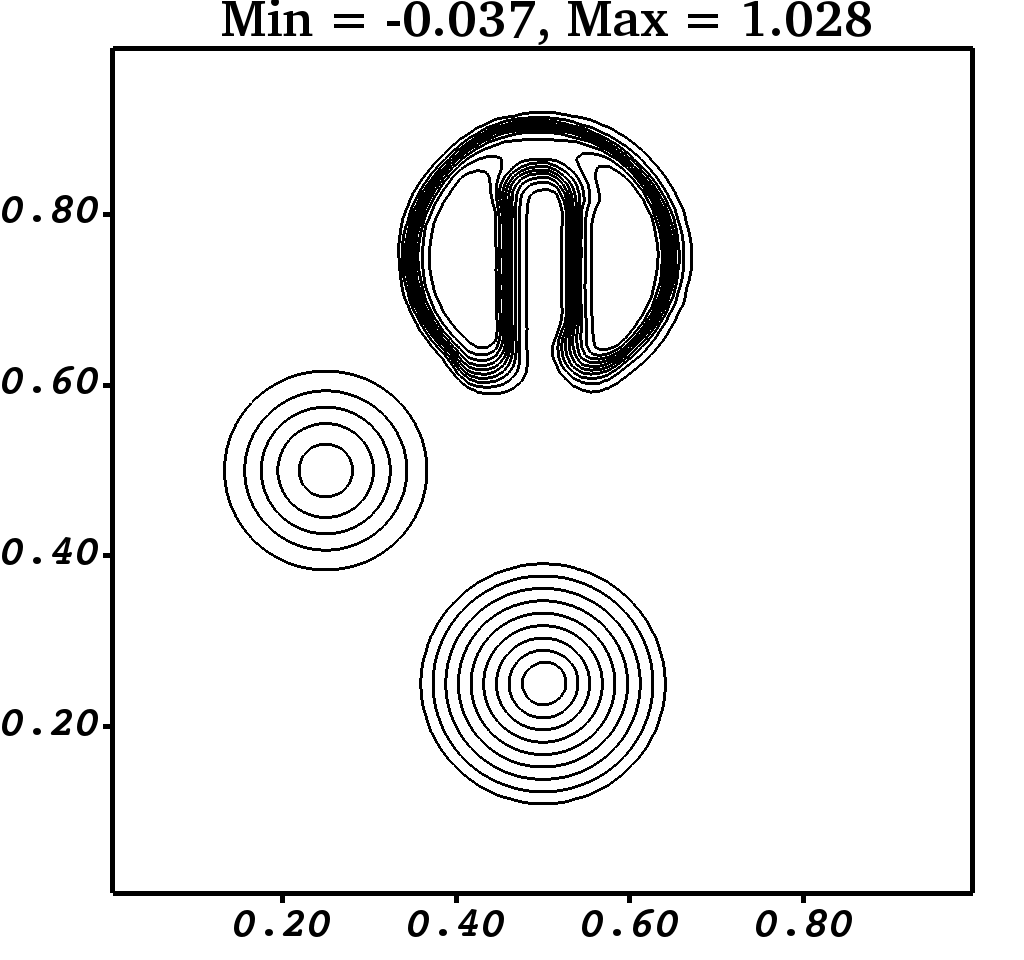}
\caption*{(c) LWFR}
\end{subfigure}
\caption*{}
\end{subfigure}
\begin{subfigure}{0.33\textwidth}
\begin{subfigure}{\textwidth}
\renewcommand\thesubfigure{\alph{subfigure}1}
\centering
\includegraphics[width=\textwidth]{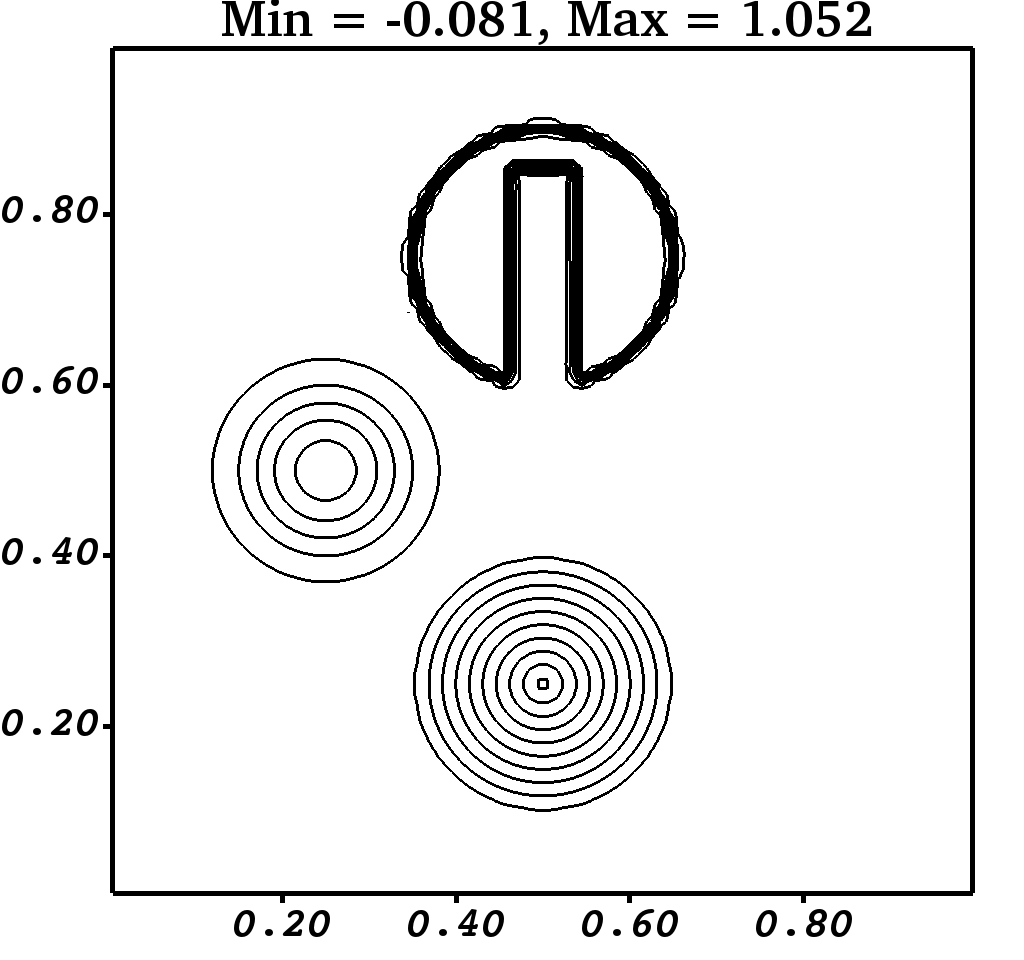}
\caption*{(d) RKFR}
\end{subfigure}
\begin{subfigure}{\textwidth}
\renewcommand\thesubfigure{\alph{subfigure}2}
\centering
\includegraphics[width=\textwidth]{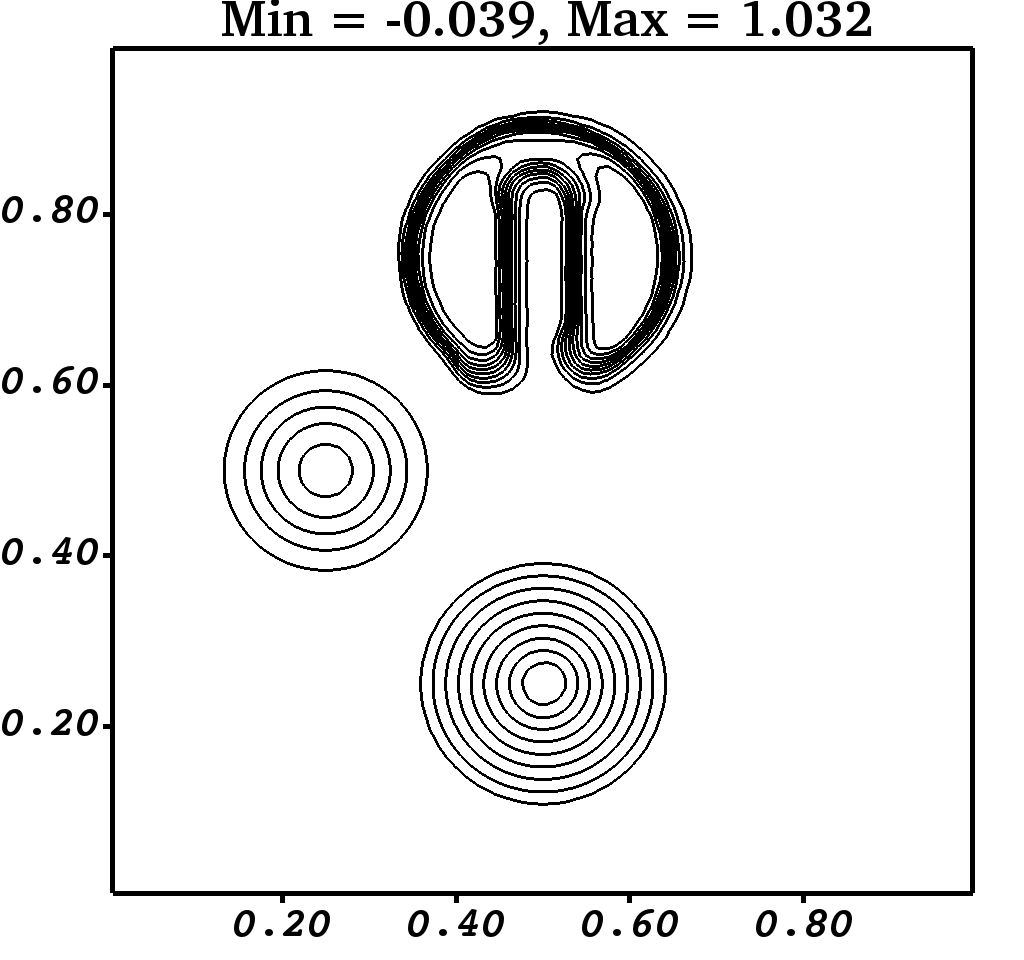}
\caption*{(e) RKFR}
\end{subfigure}
\caption*{}
\end{subfigure}
\caption{Numerical solutions of composite signal with velocity $\bm{a}=(\frac{1}{2}-y,x-\frac{1}{2})$ obtained for degree $N=3$ using Radau correction function and GL solution points. The solutions are plotted after 1 period of rotation on a mesh of $100 \times 100$ cells; No limiter is used in (b), (c) and TVB limiter $(M=100)$ is used in (d), (e).
}
\label{fig:interactions}
\end{figure}

In this example, we consider a classical test case~\cite{LeVeque1996} where the equation~\eqref{eq:2dvaradv} is solved with a divergence free velocity field $\bs a=(\frac{1}{2}-y,x-\frac{1}{2})$ and an initial condition which consists of a slotted disc, a cone and a smooth hump, given as follows
\begin{align*}
u(x,y,0) &= u_1(x,y)+u_2(x,y)+u_3( x,y),\quad  (x,y)\in [0,1]\times[0,1]\\
u_1(x,y) &= \frac{1}{4}(1+\cos (\pi q( x,y))),\quad q(x,y) = \min ( \sqrt{(x-\bar x)^2+(y-\bar y)^2},r_0 )/r_0,{(\bar x,\bar y)} = (0.25,0.5), r_0=0.15\\
u_2(x,y) &= \begin{cases}
1-\dfrac{1}{r_0} \sqrt{(x-\bar x)^2+(y-\bar y)^2} & \mbox{ if } (x-\bar x)^2+(y-\bar y)^2\le r_0^2\\
0 & \mbox{otherwise}
\end{cases}, \quad {(\bar x,\bar y)} =(0.5,0.25), r_0=0.15\\
u_3(x,y) &= \begin{cases}
1 & \mbox{ if } (x,y) \in \mathrm{C}\\
0 & \mbox{otherwise}
\end{cases}
\end{align*}
where $\mathrm{C}$ is a slotted disc with center at $(0.5,0.75)$ and radius of $0.15.$ The initial condition is shown in Figure~(\ref{fig:interactions}a). The numerical solutions of LWFR and RKFR after one rotation, without limiter, degree $N=3$ and $100 \times 100$ cells, are shown in Figures~(\ref{fig:interactions}b),~(\ref{fig:interactions}d), respectively. The same with a TVB limiter ($M=100$) are shown in Figures~(\ref{fig:interactions}c),~(\ref{fig:interactions}e). Without the limiter, the solution is captured well but there are some oscillations which take the solution outside the initial range of values. With the TVB limiter, the oscillations are reduced though it is not completely eliminated and  results  in increased numerical dissipation that smears the discontinuous profiles. However, in all cases, LWFR scheme performs comparably with RKFR scheme with the same limiter settings.

\subsection{Inviscid Burger's equation}
\begin{figure}
\centering
\begin{subfigure}[b]{0.45\textwidth}
\includegraphics[width=\textwidth]{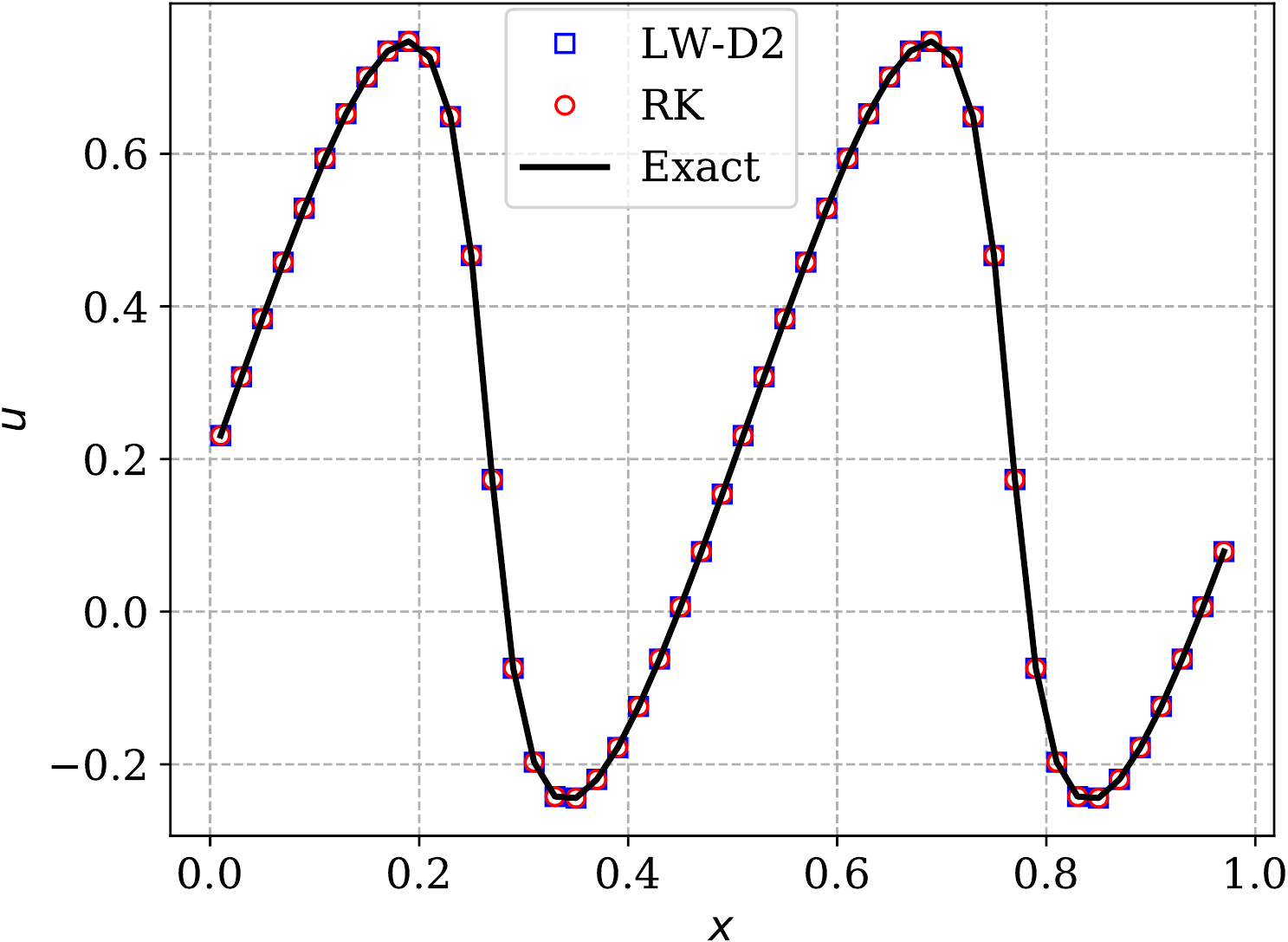}
\caption{$t=0.1$}
\end{subfigure}
\begin{subfigure}[b]{0.45\textwidth}
\includegraphics[width=\textwidth]{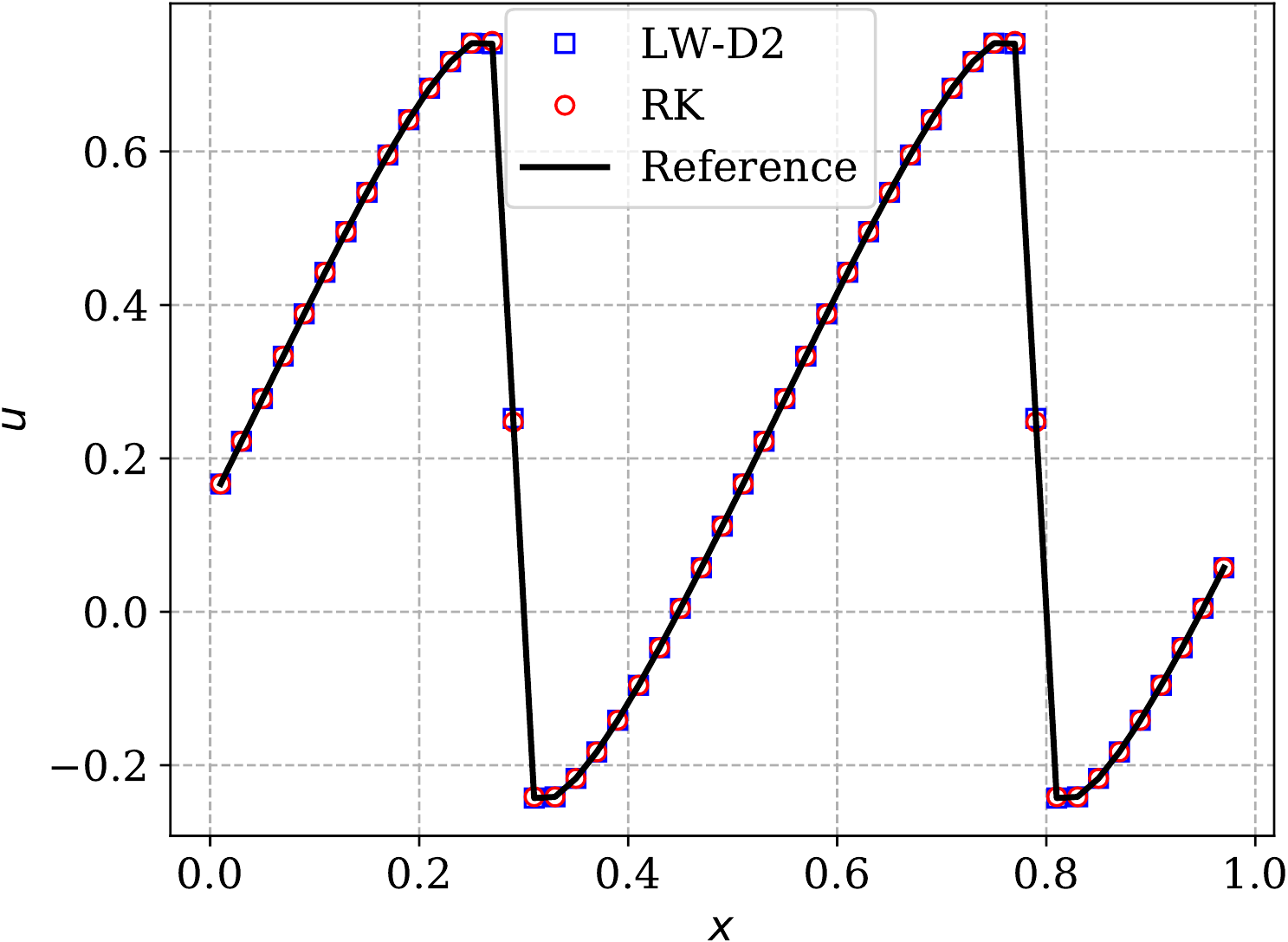}
\caption{$t=0.2$}
\end{subfigure}
\caption{Line plot across the diagonal of $[0,1]\times[0,1]$ of the solution of 2-D Burger's equation with $50 \times 50$ cells and degree $N=3$. The reference solution for $t=0.2$ is computed using RKFR scheme with degree $N=1$ on a mesh of $1000 \times 1000$ cells.}
\label{fig:lineplot_burg2d}
\end{figure}

\begin{figure}
\begin{center}
\begin{tabular}{cc}
\includegraphics[width=0.40\textwidth]{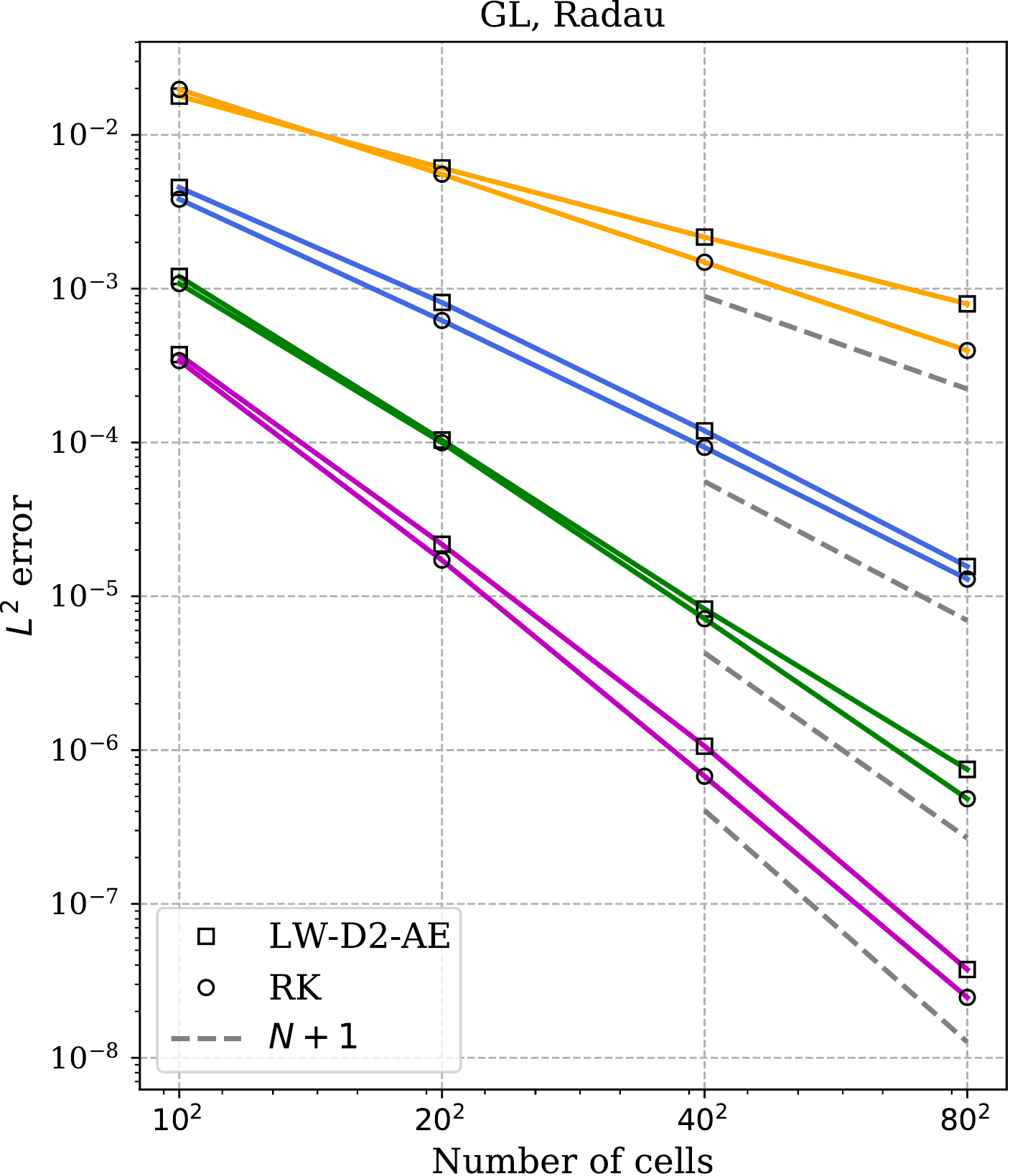} &
\includegraphics[width=0.40\textwidth]{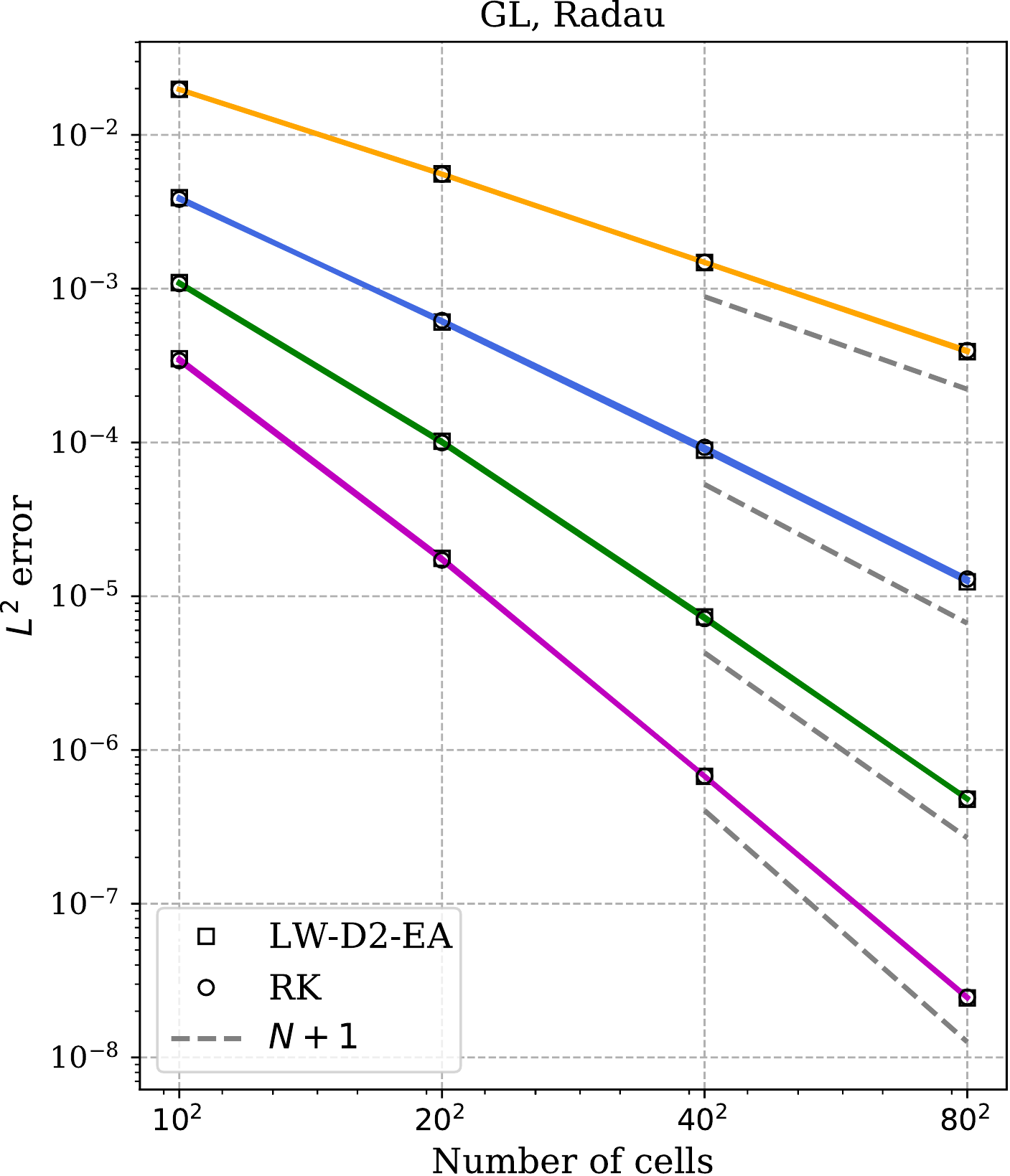} \\
(a) & (b)
\end{tabular}
\caption{Error convergence test for 2-D Burger's equation with initial condition $u(x,y,0)=\dfrac{1}{4}+\dfrac{1}{2}\sin(2\pi(x+y))$ in the domain $[0,1]\times[0,1]$ comparing the two boundary fluxes of LWFR (a) AE, (b) EA.  The errors are computed at $t=0.1$.}
\label{fig:conv_burger2d}
\end{center}
\end{figure}

We test the accuracy and robustness  of the LWFR scheme  for the two dimensional nonlinear scalar problem by considering a Burger-type equation~\cite{Qiu2005b}
\begin{equation}\label{eq:2dburger}
    u_t+\left(\dfrac{u^2}{2}\right)_{x}+\left(\dfrac{u^2}{2}\right)_{y}=0
\end{equation}
with an initial condition $u(x,y,0)=\dfrac{1}{4}+\dfrac{1}{2}\sin(2\pi(x+y))$ in the domain $\Omega = [0,1]\times[0,1]$. The  boundary conditions are set to be periodic in both directions. To test the  error convergence, the solutions are computed up to time $t=0.1$ as shown in Figure~(\ref{fig:lineplot_burg2d}a), when the solutions are still smooth and the exact solution is available. The error convergence results up to degree four are given in Figure~(\ref{fig:conv_burger2d}) using D2 dissipation model and Radau correction function.  Similar to that  in the 1-D case, the AE scheme shows optimal convergence rate for even polynomial degrees but sub-optimal convergence rates for odd degrees. The EA scheme on the other hand shows optimal convergence rates at all degrees and the error values are also comparable to those of RK scheme. In order to show the robustness of the LWFR scheme, we compute the numerical solution at time $t=0.2$ where the solution is discontinuous. The corresponding solution across the diagonal of the domain for mesh size $50\times 50$ with $N=3$ is shown in Figure~(\ref{fig:lineplot_burg2d}b) which shows that the shock is captured accurately and without spurious oscillations. In each case, when the interface fluxes are computed with EA scheme, the LWFR schemes perform at par with the RKFR schemes.

\section{Numerical results in 2-D: Euler equations}
We consider the two-dimensional Euler equations of gas dynamics given by
\begin{equation}\label{eq:2deuler}
\pd{}{t} \begin{pmatrix}
\rho \\
\rho u \\
\rho v \\
E
\end{pmatrix} +
\pd{}{x} \begin{pmatrix}
\rho u \\
p + \rho u^2 \\
\rho u v \\
(E+p)u
   \end{pmatrix} +
\pd{}{y} \begin{pmatrix}
\rho v \\
\rho u v \\
p + \rho v^2 \\
(E+p)v
\end{pmatrix}= 0
\end{equation}
where $\rho, p$ and $E$ denote the  density, pressure and total energy of the gas, respectively and $(u, v)$ are Cartesian components of the fluid velocity. For a polytropic gas, an equation of state $E=E(\rho, u, v, p)$ which leads to a closed system is given by
\begin{equation}\label{eq:2dstate}
E = E(\rho, u, v, p) = \frac{p}{\gamma -1}+\frac{1}{2} \rho (u^2 + v^2)
\end{equation}
where $\gamma > 1$ is the adiabatic constant, that will be taken as $1.4$ in the numerical tests, which is the typical value for air. 

We present results to test the accuracy and computational performance of the LW schemes for some 2-D problems and compare them to RK scheme. For each problem in this section, the corresponding CFL numbers of Lax-Wendroff schemes are chosen based on the Fourier stability analysis which are given in Table~\ref{tab:2Dcfl}. For a fair performance comparison, the Runge-Kutta schemes use their optimal CFL numbers~\cite{Gassner2011a}. The time averaged flux is always computed using the EA scheme. For RKFR, we use SSP Runge-Kutta time integration~\cite{Gottlieb2001} for degrees $N=1$ and 2,  the five stage SSP Runge-Kutta method of order four for  $N=3$~\cite{Spiteri2002}, and six-stage, fifth order Runge-Kutta (RK65) time integration for $N=4$~\cite{Tsitouras2011} implemented in {\tt DifferentialEquations.jl}~\cite{Rackauckas2017}. We make use of the HLLC flux with wave speeds from~\cite{Batten1997}. All the results in this section are produced using code written in Julia~\cite{Bezanson2017}.

\subsection{Isentropic vortex}

We next perform error convergence and Wall Clock Time (WCT) analysis using the isentropic vortex test case~\cite{Yee1999,Spiegel2016}. This problem consists of a vortex that advects at a constant velocity while the entropy is constant in both space and time. The initial condition is given by
\[
\rho = \left[ 1 - \frac{\beta^2 (\gamma - 1)}{8 \gamma \pi^2} \exp (1 -
r^2) \right]^{\frac{1}{\gamma - 1}}, \qquad u = M \cos \alpha -
\frac{\beta (y - y_c)}{2 \pi} \exp \left( \frac{1 - r^2}{2} \right)
\]
\[
v = M \sin \alpha + \frac{\beta (x - x_c)}{2 \pi} \exp \left( \frac{1 -
r^2}{2} \right), \qquad r^2 = (x - x_c)^2 + (y - y_c)^2
\]
and the pressure is given by $p = \rho^{\gamma}$. We choose the parameters $\beta = 5$, $M = 0.5$, $\alpha = 45^o$, $(x_c, y_c) = (0, 0)$ and the domain is taken to be $[- 10, 10] \times [- 10, 10]$ with periodic boundary conditions.  For this configuration, the vortex returns to its initial position after a time interval of $T=20\sqrt{2}/M$ units. We run the computations up to a  time  $t = T$ when the vortex has crossed the domain once in the diagonal direction.

\begin{figure}
\begin{center}
\begin{tabular}{cc}
\includegraphics[width=0.40\textwidth]{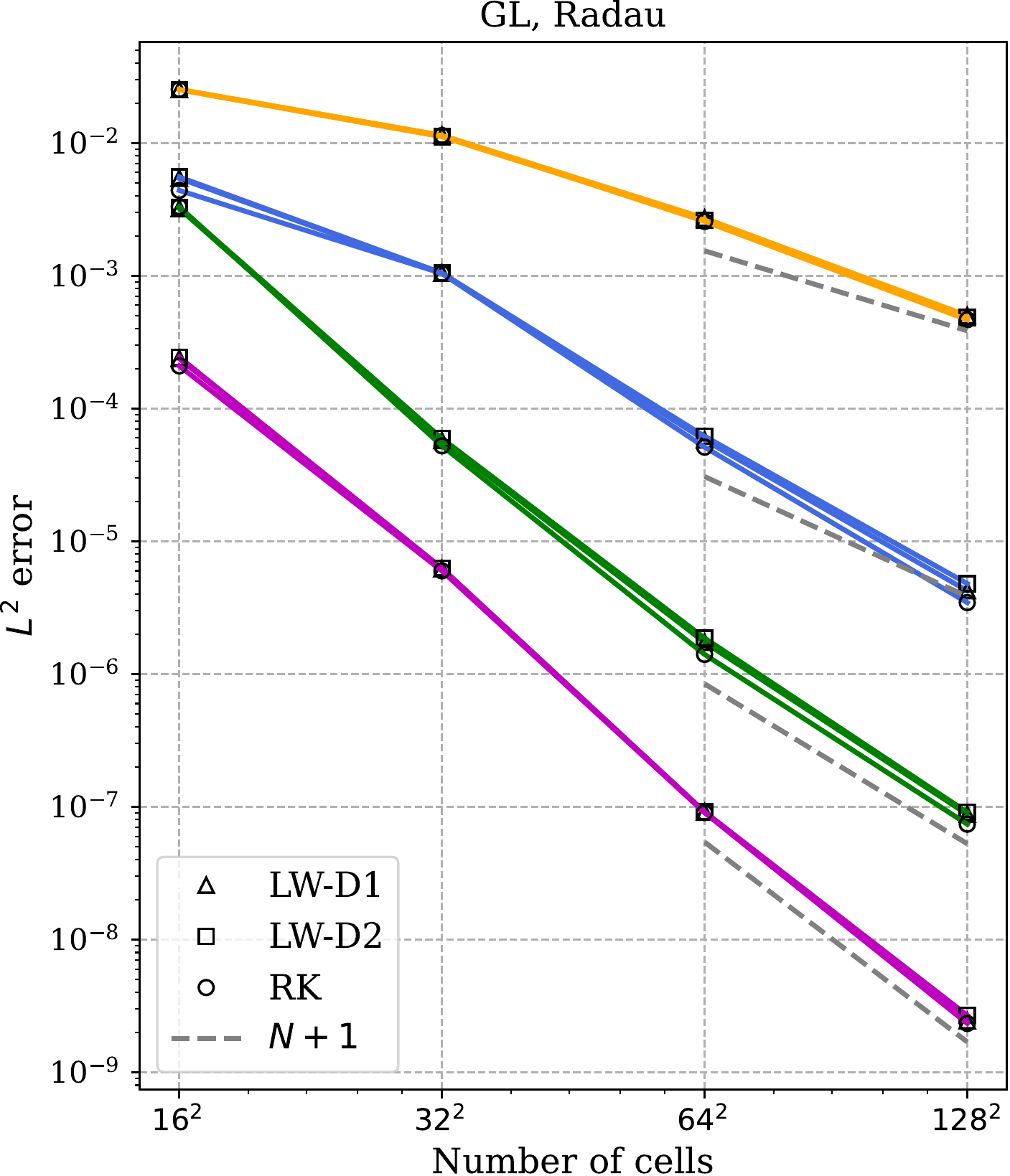} &
\includegraphics[width=0.40\textwidth]{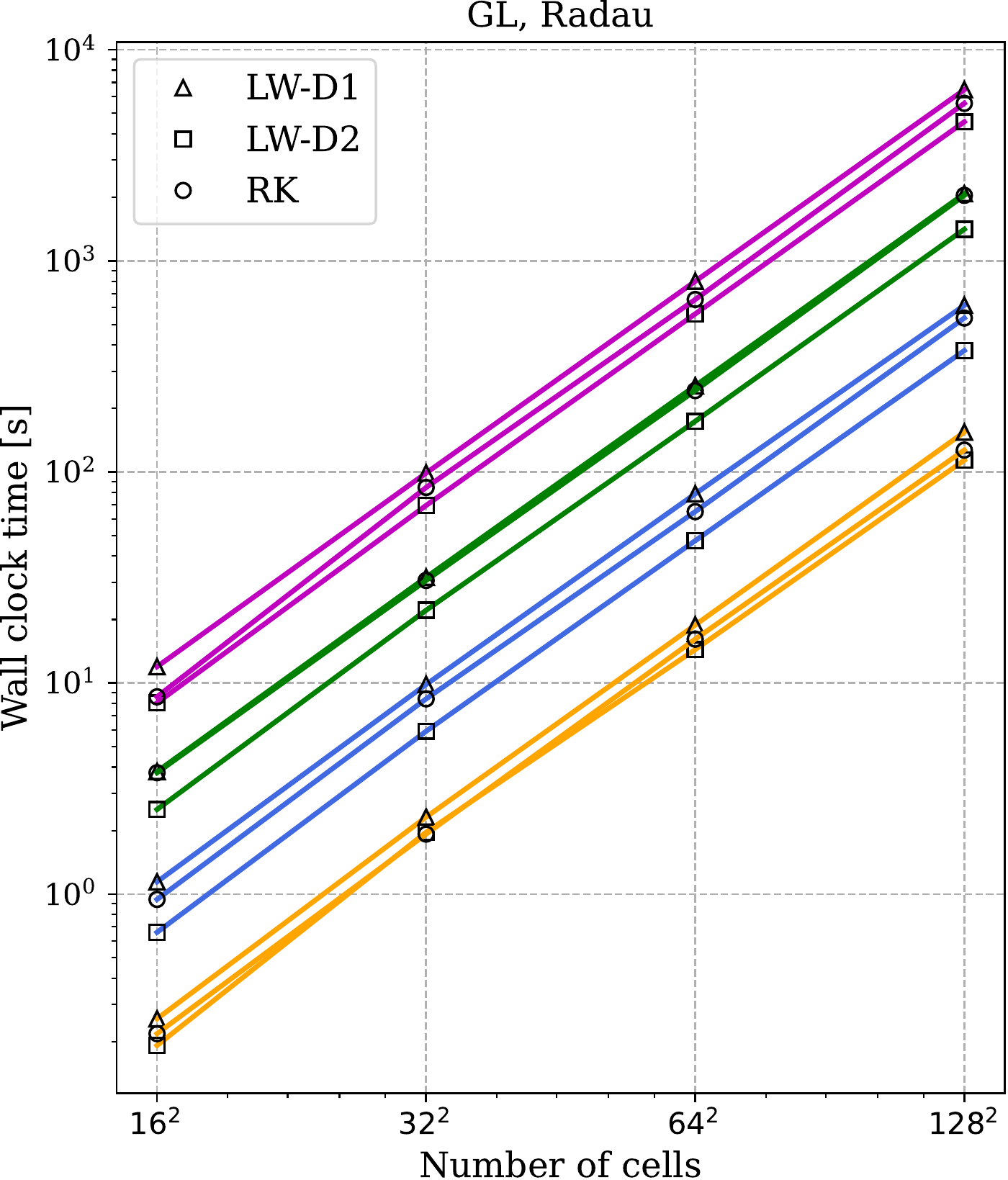}\\
(a) & (b)
\end{tabular}
\end{center}
\caption{$L^2$ error and Wall Clock Time (WCT) analysis of 2-D Euler equations (isentropic vortex) against grid resolution comparing LW-D1, LW-D2 and RK is shown in (a) and (b) respectively. The error is computed after one period. The time step size of each scheme is computed using its optimal CFL from Fourier stability analysis.}
\label{fig:isentropic.convergence}
\end{figure}

The $L^2$ error and Wall Clock Time (WCT) against grid resolution is shown in Figure~(\ref{fig:isentropic.convergence}). We observe optimal convergence rates for all  new LW schemes proposed in this paper. The WCT scales as the total number of cells to the power 1.5, which is the expected rate and the LW-D2 scheme shows smallest time as seen in Figure~(\ref{fig:isentropic.convergence}b). We denote by WCT(LW-D1) the Wall Clock Time corresponding to LW-D1 scheme, and similarly for other schemes.

\begin{figure}
\begin{center}
\begin{tabular}{c}
\includegraphics[width=0.65\textwidth]{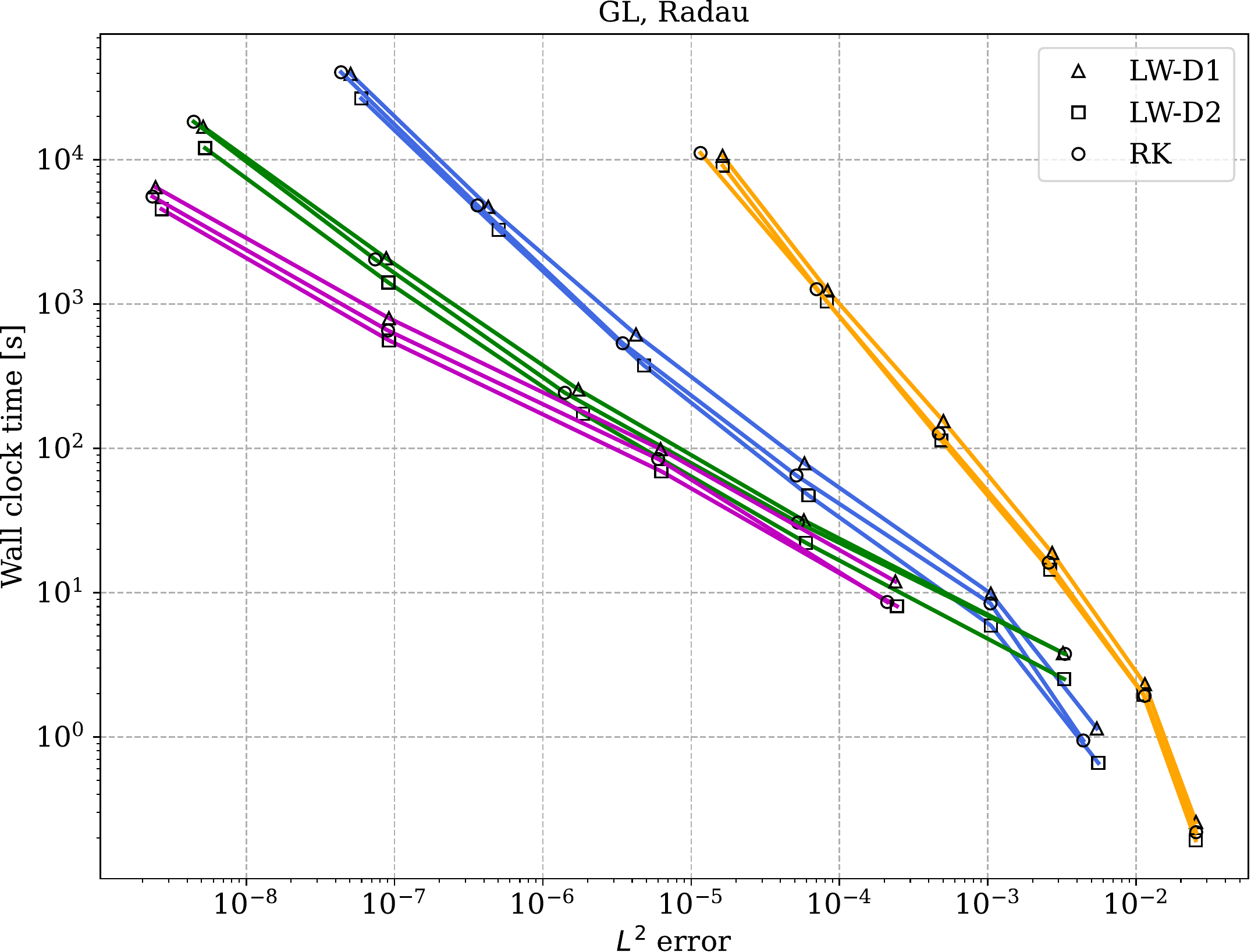}
\end{tabular}
\end{center}
\caption{Wall Clock Time (WCT) versus $L^2$ error for 2-D Euler equations (isentropic vortex) comparing LW-D1, LW-D2 and RK for degrees $N=1,2,3,4$. The different colors correspond to different degrees, with the degree increasing from right to left. The error is computed after one period. The time step size of each scheme is computed using its optimal CFL number from Fourier stability analysis.}
\label{fig:isentropic.time.vs.error}
\end{figure}

\begin{figure}
\begin{center}
\begin{tabular}{cc}
\includegraphics[width=0.40\textwidth]{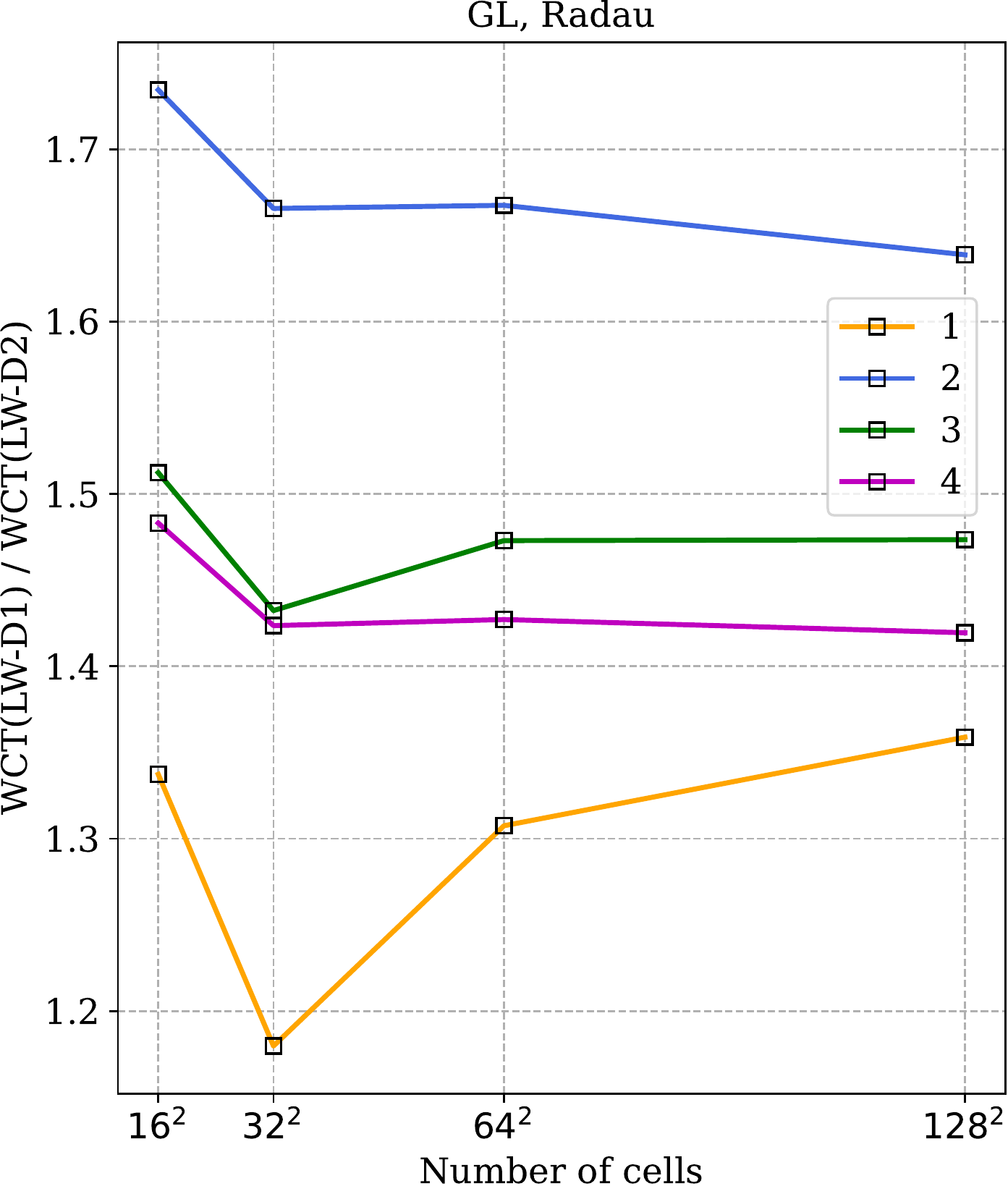} & 
\includegraphics[width=0.40\textwidth]{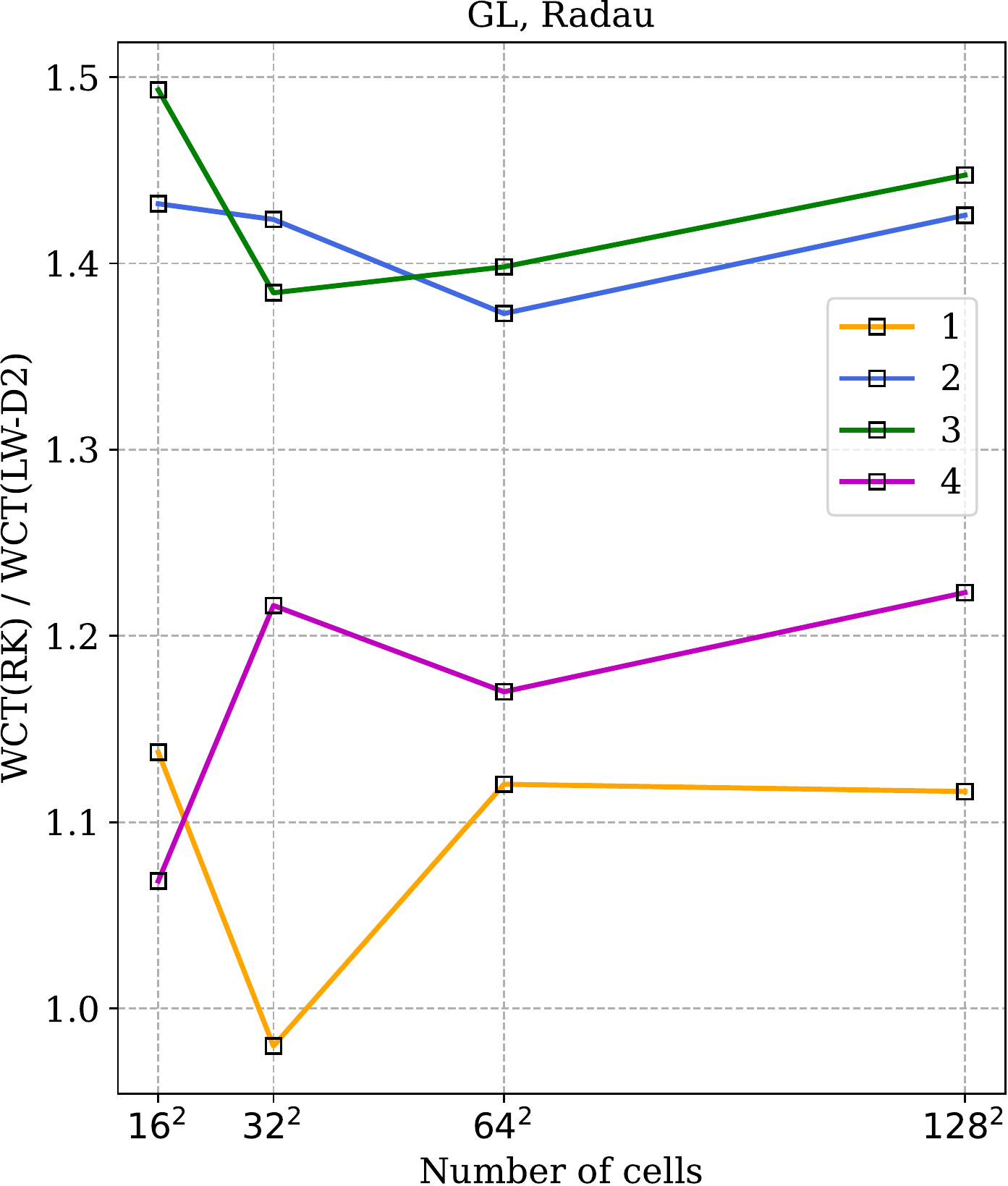} \\
(a) & (b)
\end{tabular}
\end{center}
\caption{Wall Clock Time (WCT) ratios versus grid resolution for 2-D Euler equations (isentropic vortex). (a) WCT ratio of LW-D1 and LW-D2, (b) WCT ratio of RK and LW-D2. The error is computed after one period. The time step size of each scheme is computed using its optimal CFL number from Fourier stability analysis.}
\label{fig:isentropic.ratios.vs.grid.1.2.3.4}
\end{figure}

The WCT versus $L^2$ error comparison has been made in  Figure~ (\ref{fig:isentropic.time.vs.error}) and the ratios WCT(LW-D1)/WCT(LW-D2) and WCT(RK)/WCT(LW-D2) are plotted against grid resolution in Figure~(\ref{fig:isentropic.ratios.vs.grid.1.2.3.4}a) and (\ref{fig:isentropic.ratios.vs.grid.1.2.3.4}b), respectively. We see that the newly proposed LW-D2 scheme is more efficient in comparison to the LW-D1 scheme since it can use a larger CFL number. As we expect from Table~\ref{tab:cfl} comparing the CFL ratios, the explicit time ratios of LW-D1 and LW-D2 are consistently in the range of 1.4 and 1.5 for $N>1$, as shown in Figure~(\ref{fig:isentropic.ratios.vs.grid.1.2.3.4}a).

Figure~(\ref{fig:isentropic.time.vs.error}) shows that LW-D2 has smaller Wall Clock Time than RK for all degrees and  Figure~(\ref{fig:isentropic.ratios.vs.grid.1.2.3.4}b) shows that the WCT ratio WCT(RK)/WCT(LW-D2) is close to $1.1, 1.4, 1.7$ for $N=1, 2, 3$ respectively. Thus, the ratio improves as we increase the degree up to 3. However, at $N=4$, the ratio deteriorates to approximately $1.2$.  The low CFL number of LW at $N=4$ relative to the RK scheme plays a role in this loss of performance. The amount of computations in assembling the right hand side per cell at $N=4$ is high in the LW scheme as compared to RK scheme, and it is possible that the performance of LW can be improved by a more careful implementation of the code to exploit modern CPU hardware. Figure~(\ref{fig:isentropic.time.vs.error}) shows that when the error levels are small, the higher order schemes are more efficient in terms of WCT than lower order methods.

\begin{figure}
\begin{center}
\begin{tabular}{cc}
\includegraphics[width=0.40\textwidth]{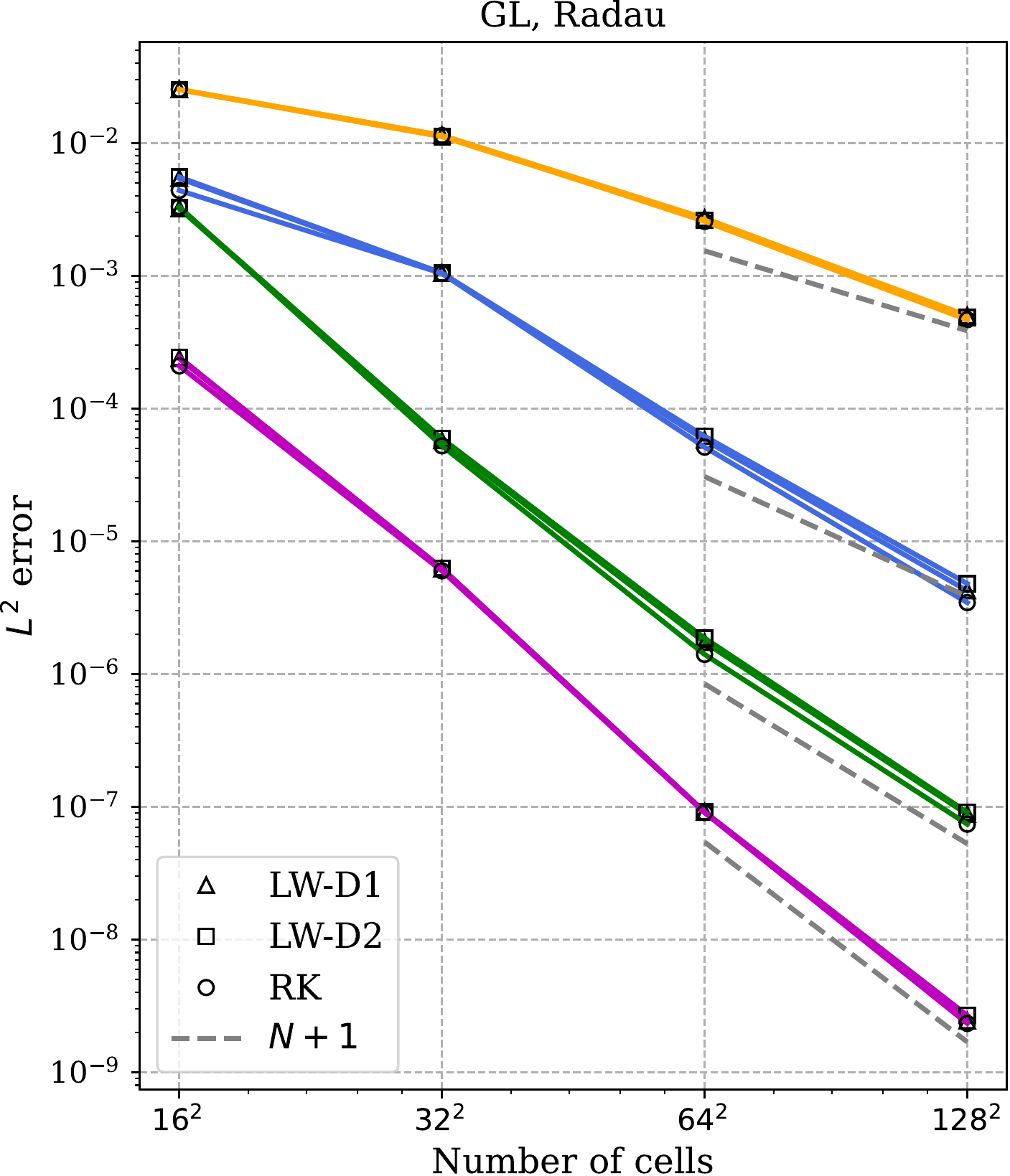} &
\includegraphics[width=0.40\textwidth]{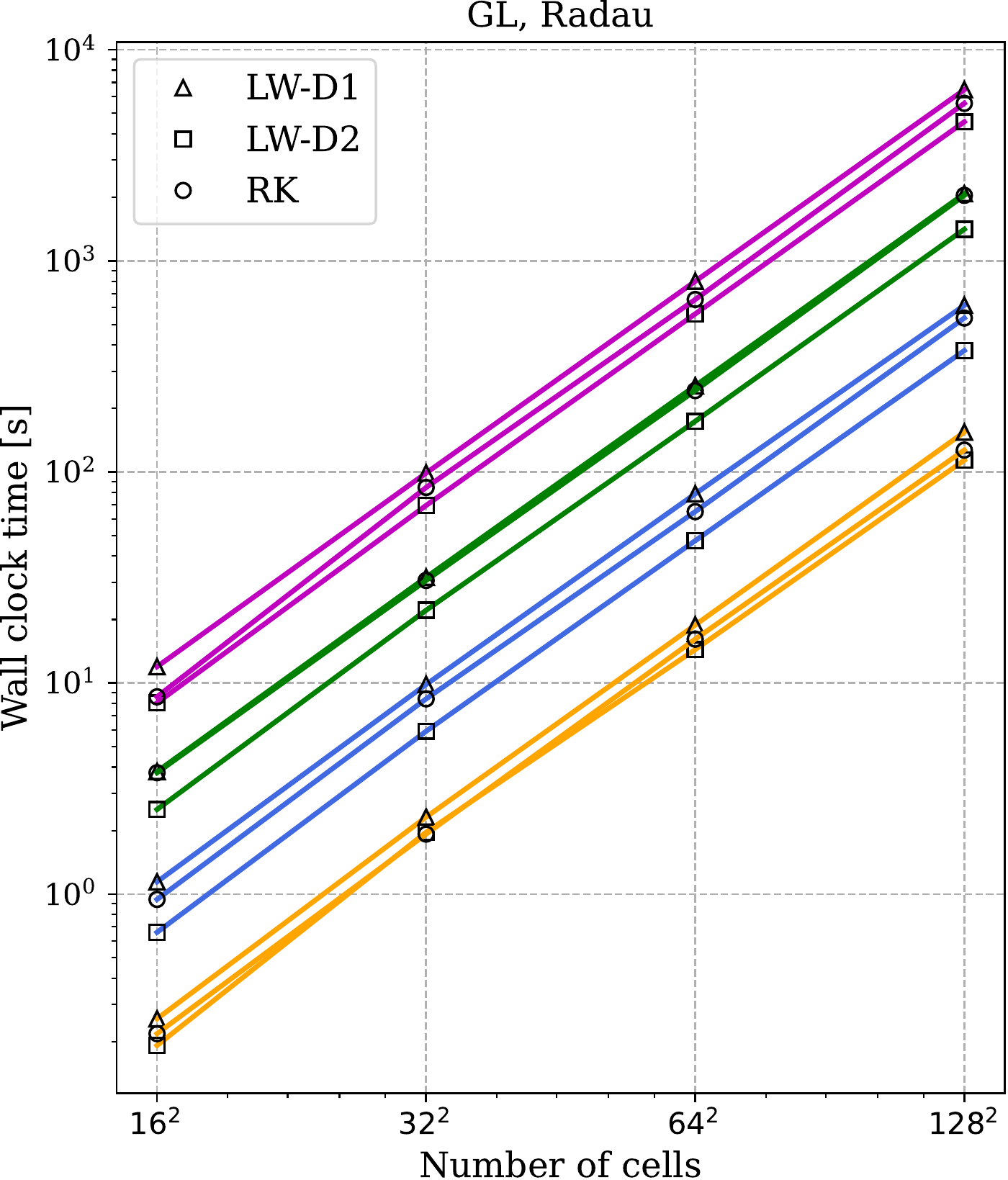}\\
(a) & (b)
\end{tabular}
\end{center}
\caption{$L^2$ error and Wall Clock Time (WCT) analysis of 2-D Euler equations (isentropic vortex) against grid resolution comparing LW-D1, LW-D2 and RK is shown in (a) and (b) respectively. The error is computed after one period. The time step size of each scheme is computed using its optimal CFL number from Fourier stability analysis.}
\label{fig:isentropic.convergence2}
\end{figure}

\subsection{Double Mach reflection}
\begin{figure}
\centering
\begin{tabular}{c}
\includegraphics[width=0.95\textwidth]{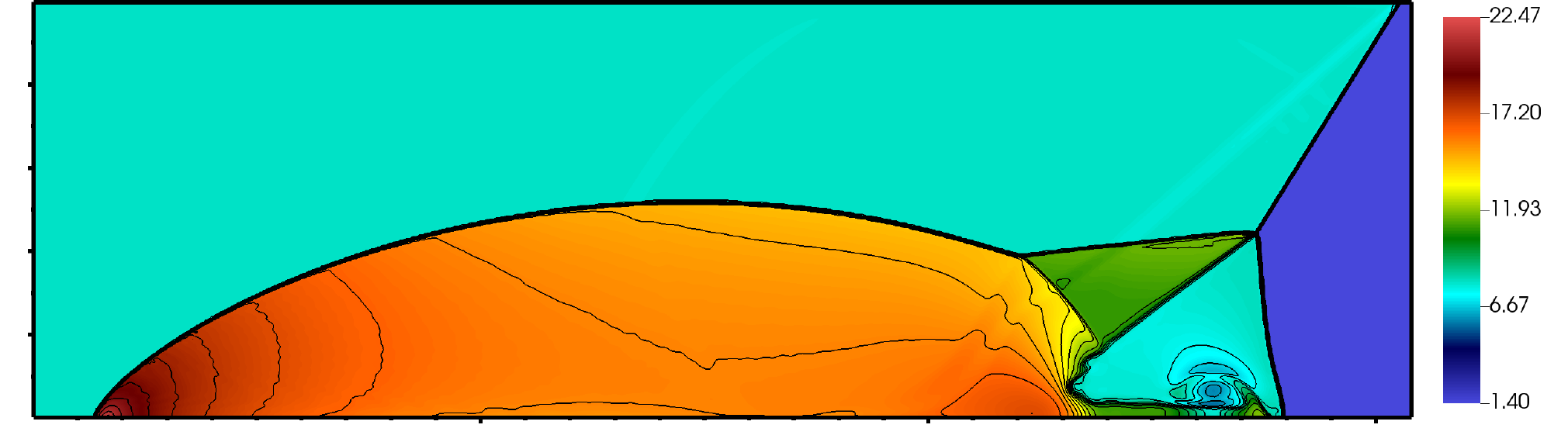} \\
(a) LWFR-D2 with EA scheme \\
\includegraphics[width=0.95\textwidth]{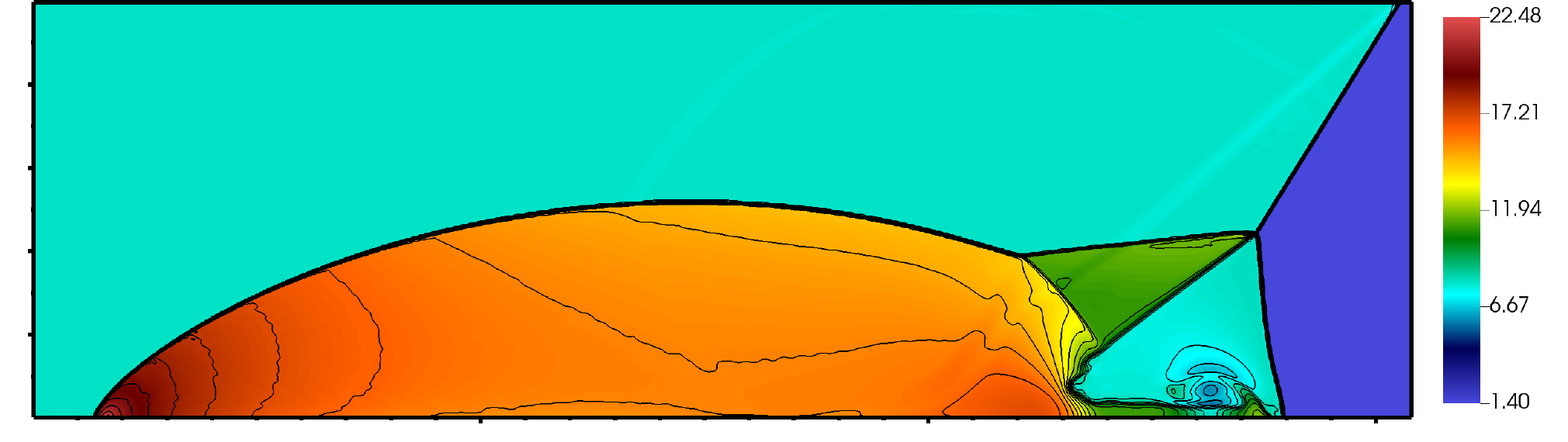} \\
(b) RKFR
\end{tabular}
\caption{Density profile  of numerical solutions of 2-D Euler equations (double Mach reflection problem) at $t=0.2$ for $N=2$, with $\Delta x = \Delta y = 1/240$. Contours of 30 steps from 1.4 to 22.5 are printed.}
\label{fig:dmr.plot}
\end{figure}

\begin{figure}
\centering
\begin{tabular}{cc}
\includegraphics[width=0.40\textwidth]{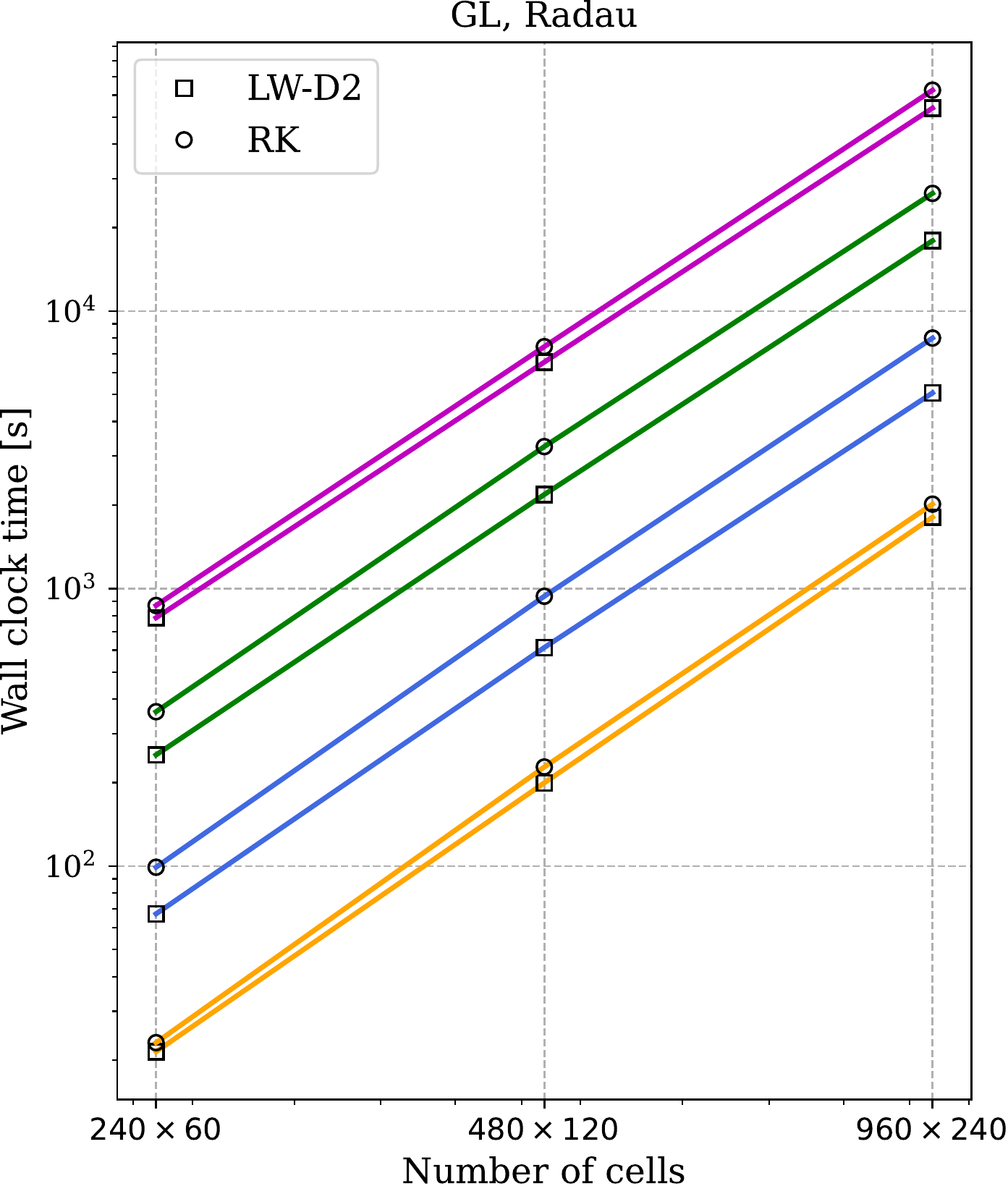} &
\includegraphics[width=0.40\textwidth]{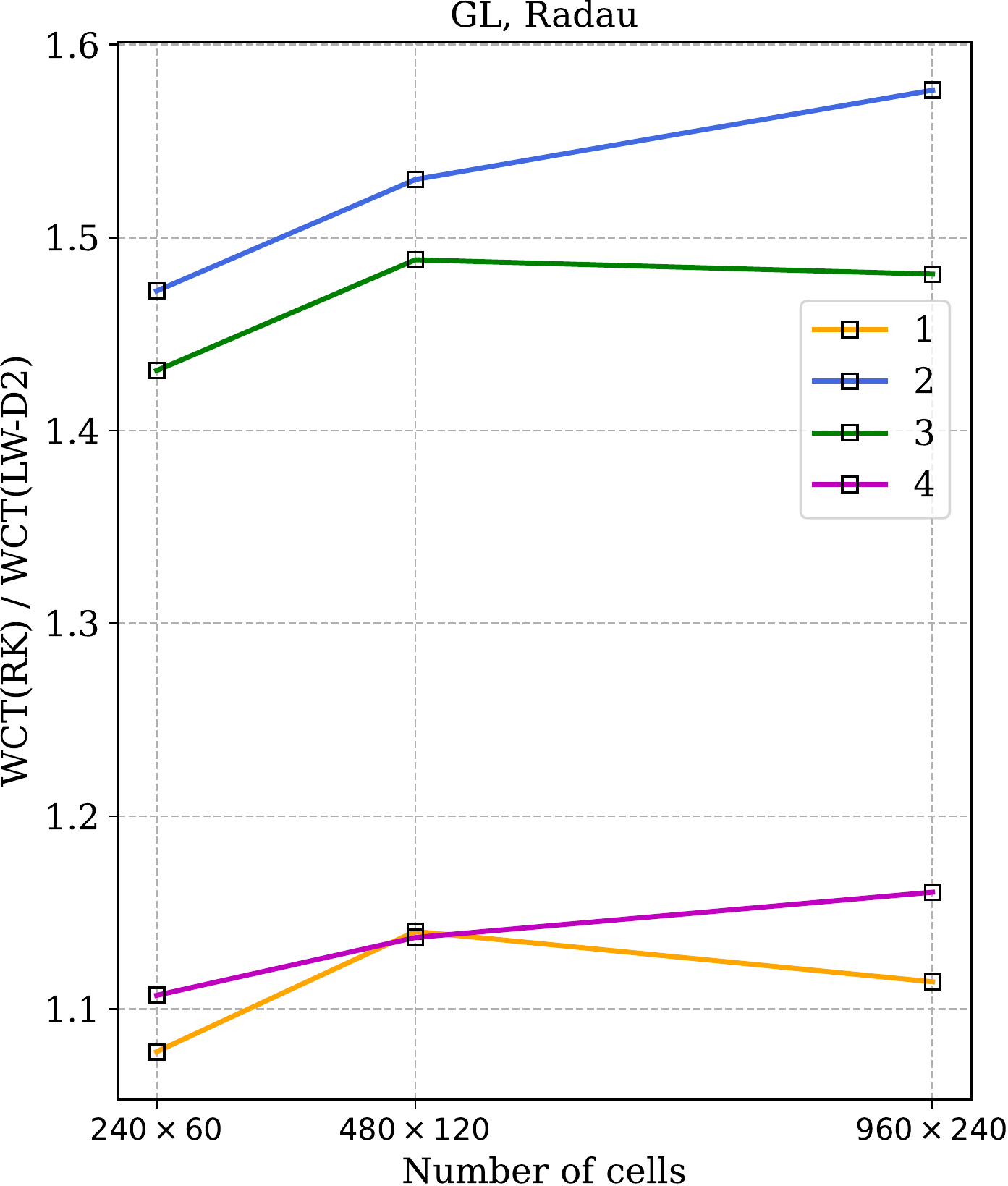} \\
(a) & (b)
\end{tabular}
\caption{Grid size versus WCT comparison of RK and LW schemes for 2-D Euler equations (double Mach reflection problem). Time step of each scheme has been chosen with its optimal CFL number from Fourier stability analysis.}
\label{fig:dmr.wct}
\end{figure}

We now test the double Mach reflection problem which was originally proposed by Woodward and Colella~\cite{Woodward1984}. The problem consists of a shock impinging on a wedge/ramp which is inclined by 30 degrees. The solution consists of a self similar shock structure with two triple points. The situation is simulated in the rectangular domain $\Omega = [0,4] \times [0,1],$ where the wedge/ramp is positioned at $x=1/6, y=0.$ Defining $\mathbf{u}_b = \mathbf{u}_b(x,y,t)$ with primitive variables given by

\begin{equation*}
(\rho,u,v,p)=\begin{cases}
(8, 8.25 \cos\left( \frac{\pi}{6} \right), -8.25 \sin\left( \frac{\pi}{6} \right), 116.5), & \mbox{ if } x < \frac{1}{6} + \frac{y + 20 t}{\sqrt{3}},\\
(1.4, 0, 0, 1), & \mbox{ if } x > \frac{1}{6} + \frac{y + 20 t}{\sqrt{3}},
\end{cases}
\end{equation*}
we define the initial condition to be $\mathbf{u}_0(x,y) = \mathbf{u}_b(x,y,0)$. With $\mathbf{u}_b$, we impose inflow boundary conditions at the left side $\{0\} \times [0,1]$, outflow boundary conditions both at $[0,1/6] \times \{0\}$ and $\{4\} \times [0,1]$, reflecting boundary conditions at $[1/6, 4] \times \{0\}$ and inflow boundary conditions at the upper side $[0,4] \times \{1\}$.
In Figure~(\ref{fig:dmr.plot}), we compare the density plots obtained using the LWFR and RKFR  schemes for $N=2$ at a resolution of $960 \times 240$ cells at $t=0.2.$ The non-linear TVB limiter is used with the parameter $M=100$~\cite{Qiu2005b}. The Lax-Wendroff solution is computed using D2 dissipation and EA scheme. We use GL points and Radau corrector in both LW and RK schemes. We observe similar resolution for both schemes; the similarity holds for other degrees also, which we have not shown to save space. In Figure~(\ref{fig:dmr.wct}a), we plot grid resolution against the Wall Clock Time for degrees $N=1,2,3,4$ which shows the expected dependence of time with grid size. Figure~(\ref{fig:dmr.wct}b) shows the ratio of WCT for RK and LW-D2 schemes, indicating better efficiency of LW scheme, and these results are similar to what is observed in~\cite{Qiu2005b}. 

\begin{figure}
\centering
\begin{tabular}{cc}
\includegraphics[width=0.45\textwidth]{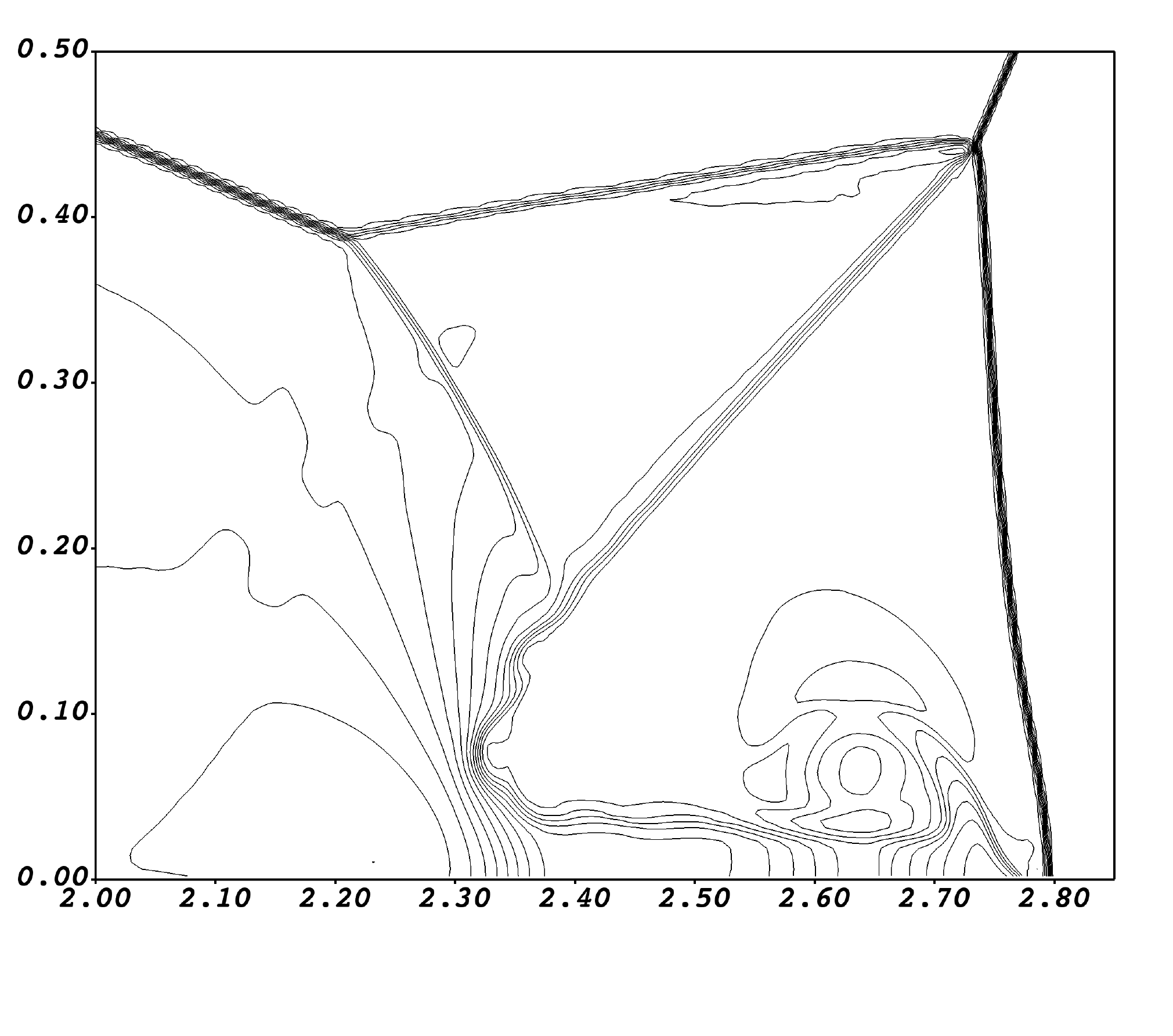} &
\includegraphics[width=0.45\textwidth]{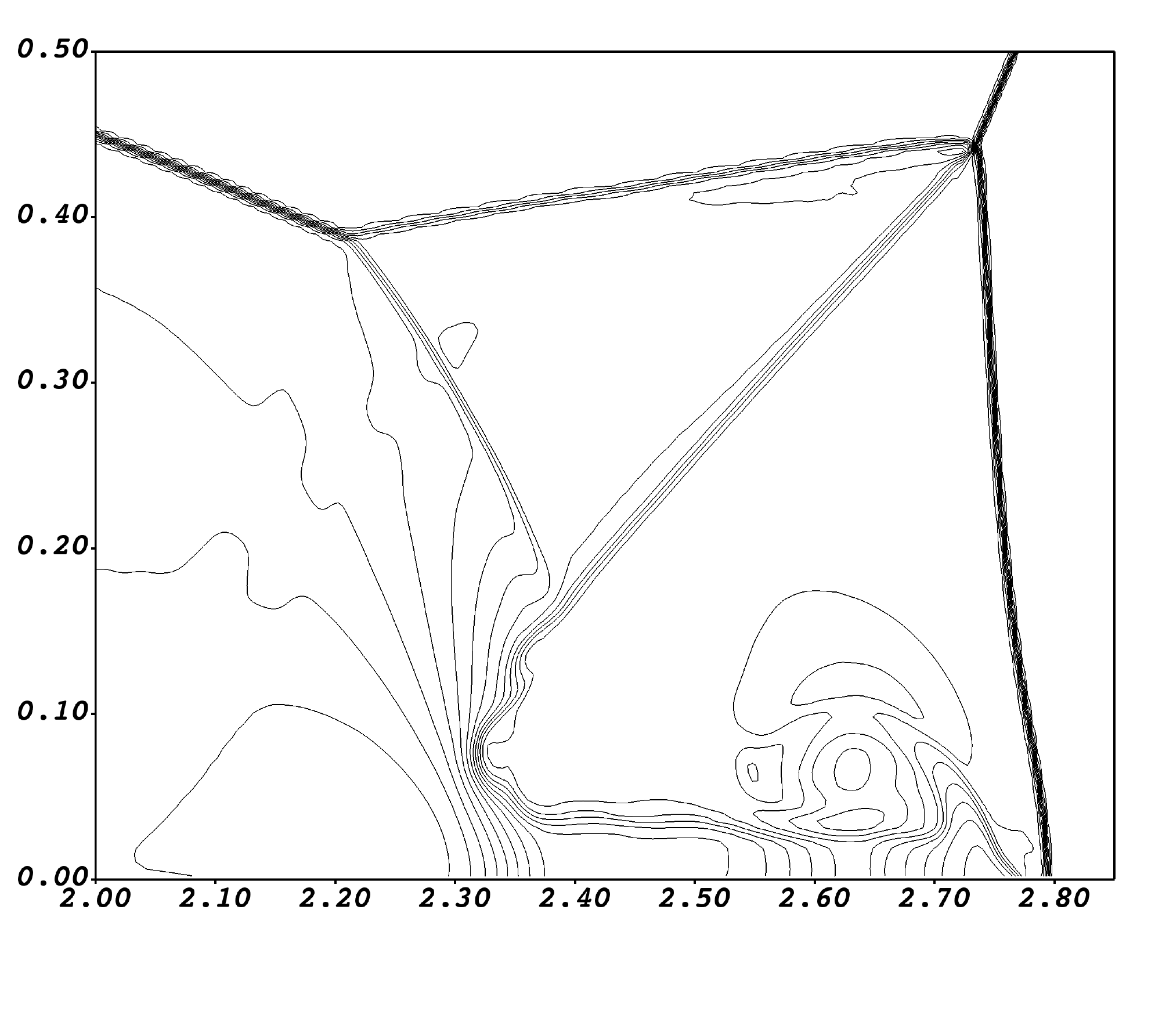} \\
(a) LWFR-D2 with EA scheme & (b) RKFR
\end{tabular}
\caption{Enlarged contours of density (2-D Euler equations, double Mach reflection problem)  at $t=0.2$ for $N=2$, with $\Delta x = \Delta y = 1/240$. Contours of 30 steps from 1.4 to 22.5 are printed.}
\label{fig:dmr.plot.zoom}
\end{figure}

\section{Summary and conclusions}\label{sec:sum}
A conservative, Jacobian-free and single step, explicit Lax-Wendroff method has been constructed in flux reconstruction context, and its implementation has been demonstrated for solving hyperbolic conservation laws in one and two dimensions. The Jacobian-free property is achieved by using a finite difference approach to compute time derivatives of the fluxes that are needed in the Taylor expansion. The method requires only the time average flux and its corresponding numerical fluxes. It is written in matrix-vector form that is useful for computer implementation. We have studied the effect of two commonly used correction functions and solution points. The stable CFL numbers are computed using Fourier stability analysis in one and two dimensions. The numerical fluxes are computed using both the time average flux and the time average solution which leads to improved CFL numbers compared to other existing methods which use the solution at previous time level to compute the dissipative part of the numerical flux. At fifth order ($N=4$), there is a mild linear instability for periodic problems, which seems to be present in other single step methods and also in RKDG schemes. For non-linear problems, we identify a loss of optimal convergence rate when a simple average-extrapolate (AE) approach is used to compute the central part of the numerical flux. We show that this can be improved to optimal rates by using an extrapolate-average procedure, and the resulting schemes perform comparably with RK schemes in terms of their error levels. The performance of the method is also demonstrated on 1-D and 2-D non-linear systems like Euler equations, where it is able to resolve all the waves at comparable accuracy to RK schemes. Many commonly used numerical fluxes based on approximate Riemann solvers and modeling even contact waves can be developed and used in these schemes. These studies show that the Radau correction function in combination with Gauss-Legendre solution points and the extrapolate-average (EA) technique leads to uniformly accurate LW scheme for non-linear problems.  The method has a simple structure which makes it easy to develop a general code that can be used to solve any conservation law; the user has to supply subroutines for the flux, numerical flux and maximum wave speed estimate used in the CFL condition.
\section*{Acknowledgments}
The work of Arpit Babbar and Praveen Chandrashekar is supported by the Department of Atomic Energy,  Government of India, under project no.~12-R\&D-TFR-5.01-0520. The work of Sudarshan Kumar Kenettinkara is supported by the  Science and Engineering Research Board, Government of India, under MATRICS project no.~MTR/2017/000649.
\appendix
\section{Numerical fluxes}\label{apendix:numfluxes}
We describe the procedure to compute the numerical flux for systems at one single face $\eph$. The numerical flux for LWFR is computed using the trace values of the solution $U_l = U_\eph^-$, $U_r = U_\eph^+$ and fluxes $F_l = F_\eph^-$, $F_r = F_\eph^+$. Here $U_l, U_r$ may be the solution values at time $t_n$ for dissipation model D1 or the time average value in case of dissipation model D2. Further, we use the cell average values at time $t=t_n$, $\bar U_l = \bar{u}_e^n$, $\bar U_r = \bar{u}_{e+1}^n$, to compute the dissipation coefficients. In the following sub-sections, we described different numerical fluxes which are functions of the quantities: $\bar U_l, \bar U_r, U_l, U_r, F_l, F_r$.
\subsection{Rusanov flux}
The Rusanov flux \cite{Rusanov1962} is a local version of the Lax-Friedrichs flux with the the wave speed being estimated locally. The flux approximation is given by
\[
F(U_l,U_r,F_l,F_r;\bar U_l, \bar U_r) = \half(F_l + F_r) - \half \lambda (U_r - U_l)
\]
where $\lambda$ is an estimate of the maximum wave speed in the two states
\[
\lambda = \max\{ \sigma(\bar U_l), \sigma(\bar U_r) \}
\]
and $\sigma$ denotes the spectral radius of the flux jacobian.
\subsection{Roe flux}
The Roe flux \cite{Roe1981} is built on a local linearization of the hyperbolic conservation law and solving the Riemann problem exactly. The Roe flux is given by
\[
F(U_l,U_r,F_l,F_r;\bar U_l, \bar U_r) = \half(F_l + F_r) - \half R |\Lambda| L (U_r - U_l)
\]
where $R, \Lambda, L$ are the right eigenvector matrix, diagonal matrix of eigenvalues and left eigenvector matrix corresponding to the flux Jacobian at the face, computed using the Roe average based on cell average values $\bar U_l$, $\bar U_r$.
\subsection{HLL flux}
The HLL Riemann solver~\cite{Harten1983a} models the solution of the Riemann problem using only the slowest and fastest waves with an intermediate state. Let the slowest and fastest speeds, denoted by $S_l < S_r$, be assumed to be known. We can determine the intermediate state and flux by writing the jump conditions across the two waves,
\[
F_* - F_l = S_l (U_* - U_l), \qquad F_r - F_* = S_r (U_r - U_*)
\]
whose solution is given by
\[
U_* = \frac{S_r U_r - S_l U_l - (F_r - F_l)}{S_r - S_l}, \qquad F_* = \frac{S_r F_l - S_l F_r + S_l S_r (U_r - U_l)}{S_r - S_l}
\]
The numerical flux is given by
\[
F(U_l,U_r,F_l,F_r;\bar U_l, \bar U_r)= \begin{cases}
F_l, & S_l > 0 \\
F_r, & S_r < 0 \\
F_*, & \textrm{otherwise}
\end{cases}
\]
The speeds $S_l, S_r$ are computed using the cell average values $\bar U_l$, $\bar U_r$ and there are various methods available~\cite{Einfeldt1988,Batten1997,Toro2009,Guermond2016,Toro2020}. In the numerical tests, we use the method from~\cite{Toro2009} to estimate the slowest and fastest speeds.

\subsection{HLLC flux}
The HLLC Riemann solver~\cite{Toro1994} uses a three wave model  with three wave speeds $S_l<S_*<S_r$ and two intermediate states $U_{* l}$ and $U_{* r}.$  The pressure and normal velocity are continuous across the contact wave, i.e.,
\[
p_{*l}=p_{*r}=p_*,\quad u_{*l}=u_{*r}=u_*
\]
and  the speed of the contact wave coincides with the intermediate velocity $S_*=u_*$. The jump condition across the $S_l$ and $S_r$ wave reads as
\[
F_{*\alpha}-F_\alpha=S_\alpha(U_{*\alpha}-U_\alpha), \quad \alpha=l,r.
\]
In the full form, the jump conditions are given by
\[
\begin{bmatrix}
\rho_{*\alpha} u_*\\
p_*+\rho_{*\alpha}u_*^2\\
(E_{*\alpha}+p_*)u_*
\end{bmatrix}
- S_\alpha\begin{bmatrix}
\rho_{*\alpha}\\
\rho_{*\alpha}u_*\\
E_{*\alpha}
\end{bmatrix}
=\begin{bmatrix}
F_\alpha^\rho \\
F_\alpha^m\\
F_\alpha^E
\end{bmatrix}
- S_\alpha\begin{bmatrix}
\rho_{\alpha}\\
m_\alpha\\
E_{\alpha}
\end{bmatrix}
\]
Using this expression we determine  the unknown variables $\rho_*,u_*,p_*$ and $E_*.$
From the first jump condition we get
\[
\rho_{*\alpha}=\frac{S_\alpha\rho_\alpha-F_{\alpha}^\rho}{S_\alpha-u_*}
\]
From the second equation we write the intermediate pressure
\begin{equation}\label{eq:pstar}
p_* = F_{\alpha}^m - S_\alpha m_\alpha + \rho_{*\alpha}u_*(S_\alpha-u_*) = F_{\alpha}^m - S_\alpha m_\alpha + u_* (S_\alpha\rho_\alpha-F_{\alpha}^\rho)
\end{equation}
We get two estimates of pressure $p_*$ from the $l,r$ states, and equating these two values
\begin{eqnarray*}
F_{l}^m - S_l m_l + u_* (S_l \rho_l - F_l^\rho) &=&
F_{r}^m - S_r m_r + u_* (S_r \rho_r - F_r^\rho)
\end{eqnarray*}
we obtain the intermediate velocity
\begin{equation*}
u_* = \frac{ (S_r m_r - F_{r}^m) - (S_l m_l  - F_{l}^m)}{
(S_r\rho_r-F_{r}^\rho) - (S_l\rho_l-F_{l}^\rho)}
\end{equation*}
The intermediate pressure can be computed from \eqref{eq:pstar} or from the following
\[
p_* = \frac{(S_r m_r - F_r^m)(S_l \rho_l - F_l^\rho) - (S_l m_l - F_l^m)(S_r \rho_r - F_r^\rho)}{(S_r\rho_r-F_{r}^\rho) - (S_l\rho_l-F_{l}^\rho)}
\]
From the last jump condition we obtain
\[
E_{*\alpha}=\frac{p_*u_* + S_\alpha E_\alpha - F_{\alpha}^E}{S_\alpha - u_*}
\]
The flux is now given by
\[
F(U_l,U_r,F_l,F_r;\bar U_l, \bar U_r)=\begin{cases}
F_l, & S_l > 0 \\
F_r, & S_r < 0 \\
F_{*l} = F_l+S_l(U_{*l}-U_l), & S_l<0<u_*\\
F_{*r} = F_r+S_r(U_{*r}-U_r), & u_*<0<S_r
\end{cases}
\]
where the wave speeds $S_l$ and $S_r$  are computed using the cell average values $\bar U_l,\bar U_r$.
\section{Fourier stability analysis in 2-D}\label{sec:fourier2d}
\begin{figure}
\begin{center}
\begin{tabular}{cccc}
\includegraphics[width=0.22\textwidth]{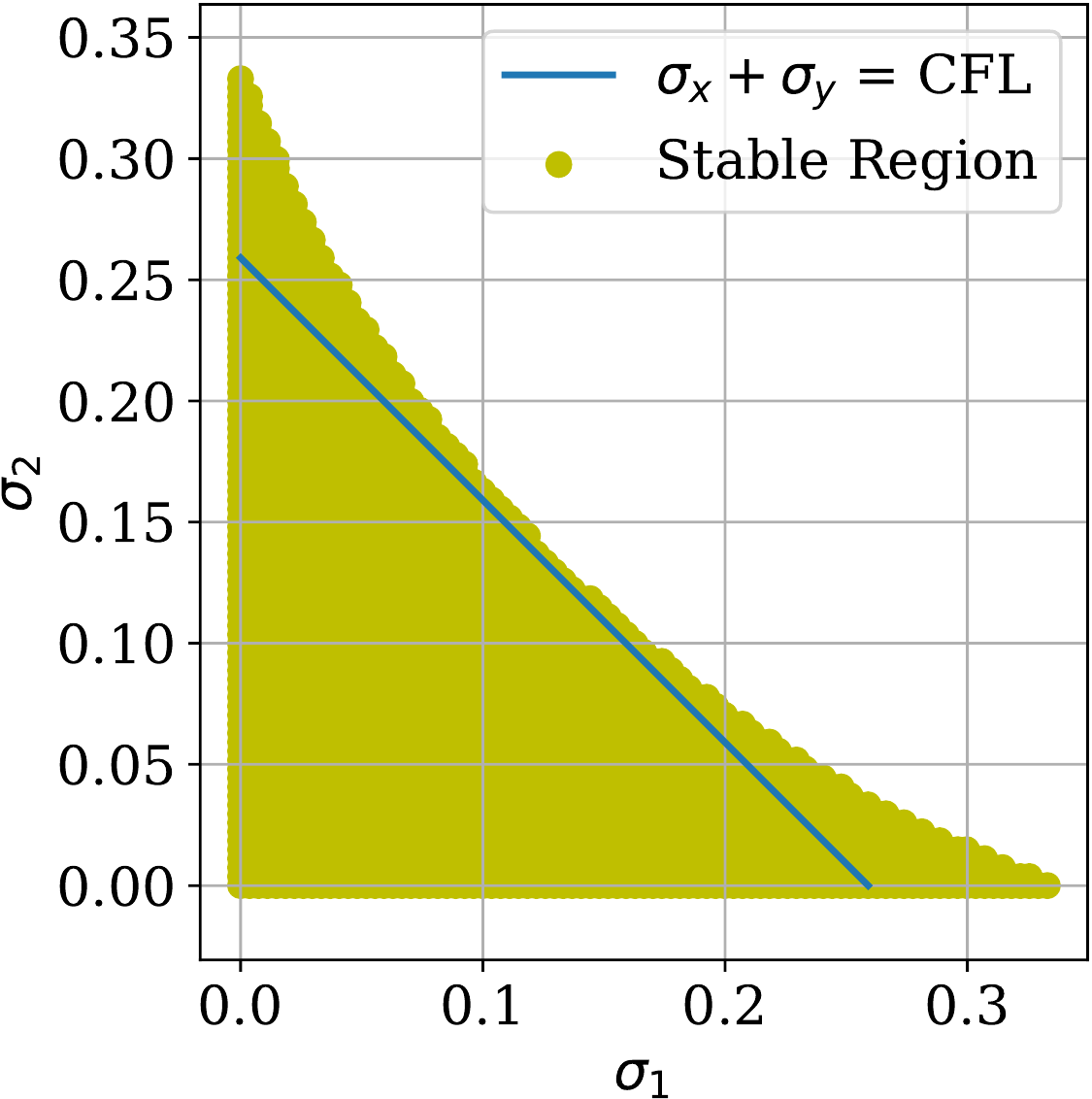} &
\includegraphics[width=0.22\textwidth]{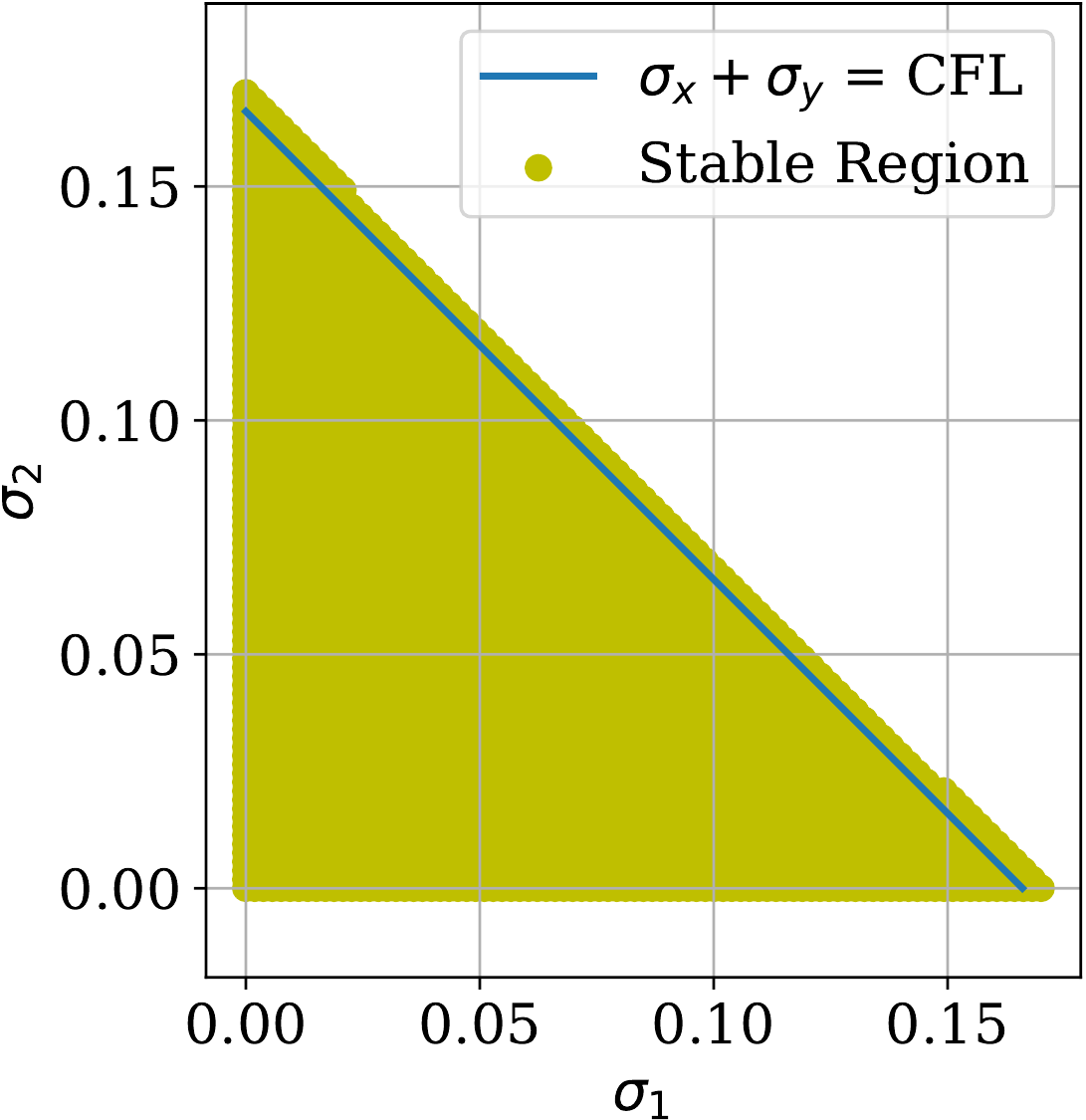} &
\includegraphics[width=0.22\textwidth]{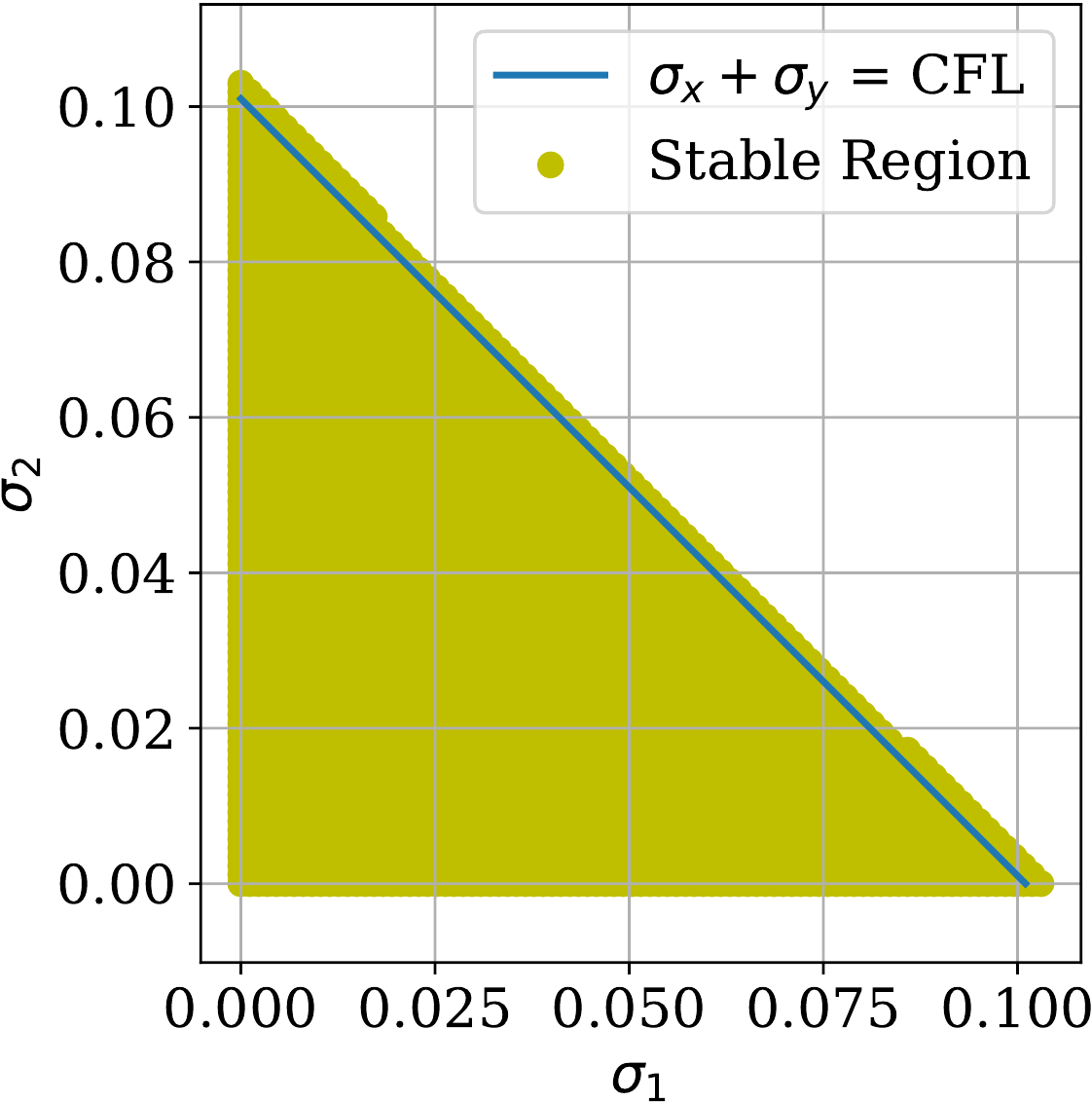} &
\includegraphics[width=0.22\textwidth]{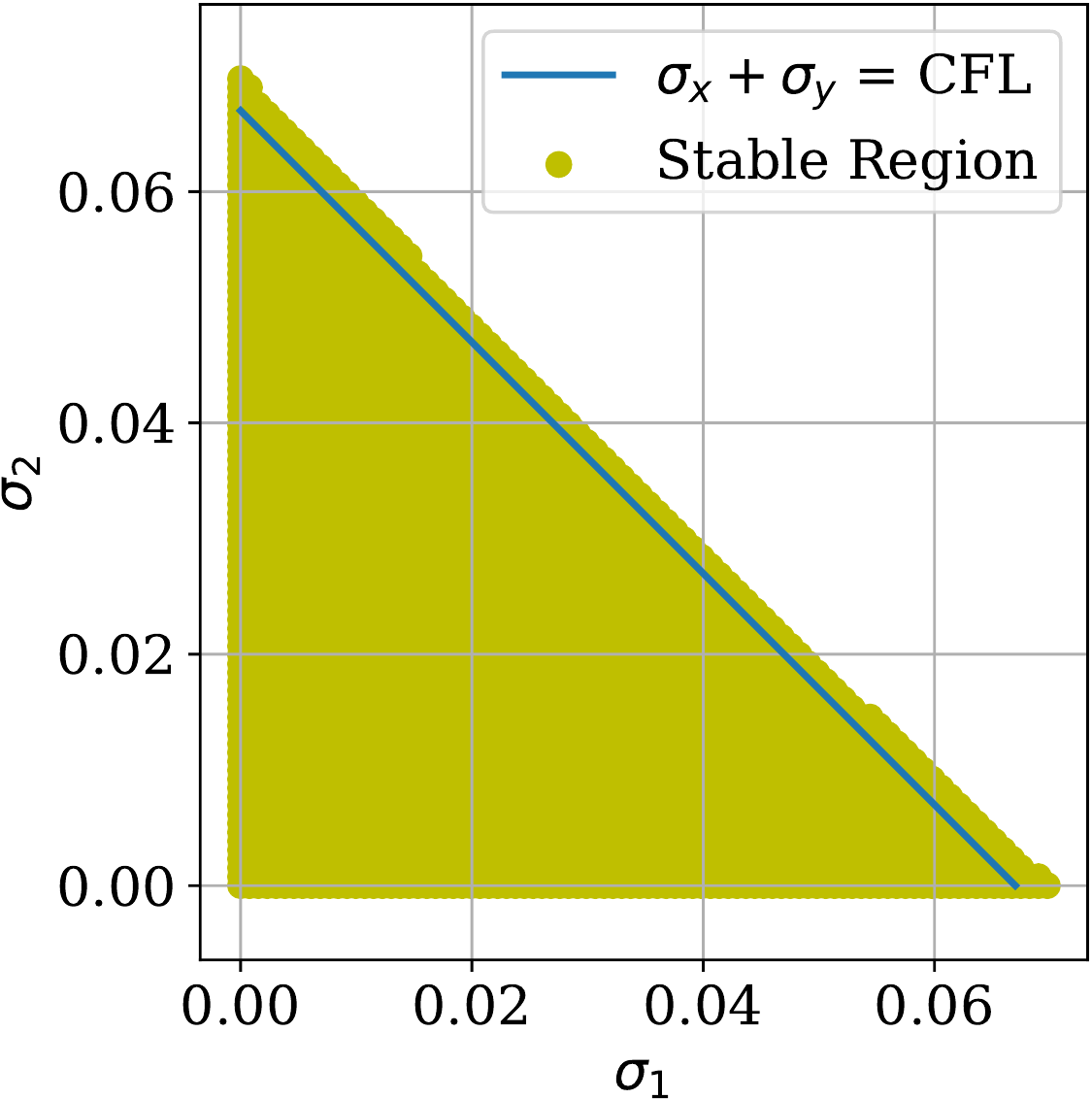} \\
(a) & (b) & (c) & (d)
\end{tabular}
\end{center}
\caption{Stability regions of LWFR scheme with  the  Radau correction function and D2 dissipation model in two dimensions. (a) $N=1$, (b) $N=2$, (c) $N=3$, (d) $N=4$.}
\label{fig:2dcfl_radau}
\end{figure}

\begin{figure}
\begin{center}
\begin{tabular}{cccc}
\includegraphics[width=0.22\textwidth]{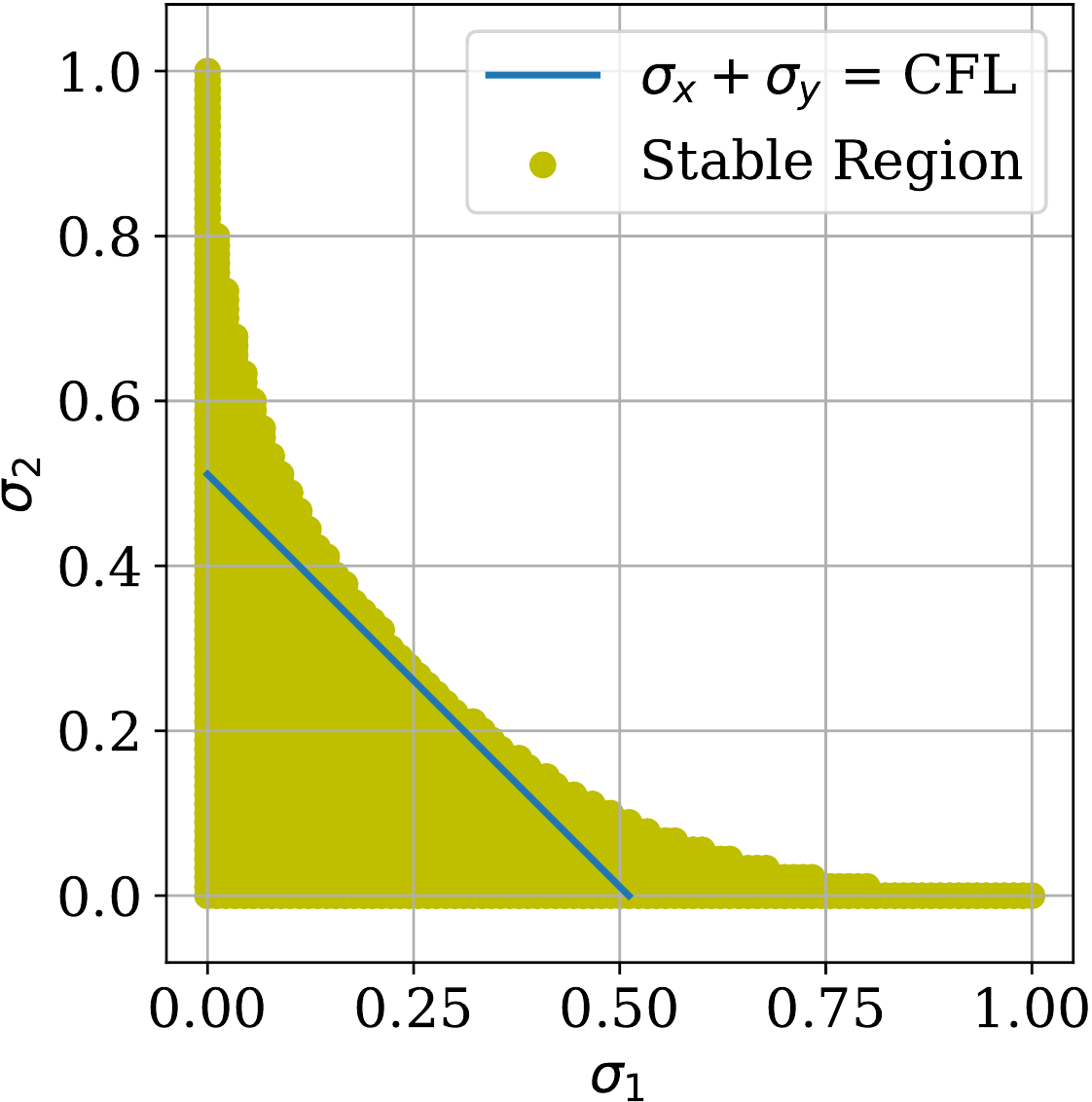} &
\includegraphics[width=0.22\textwidth]{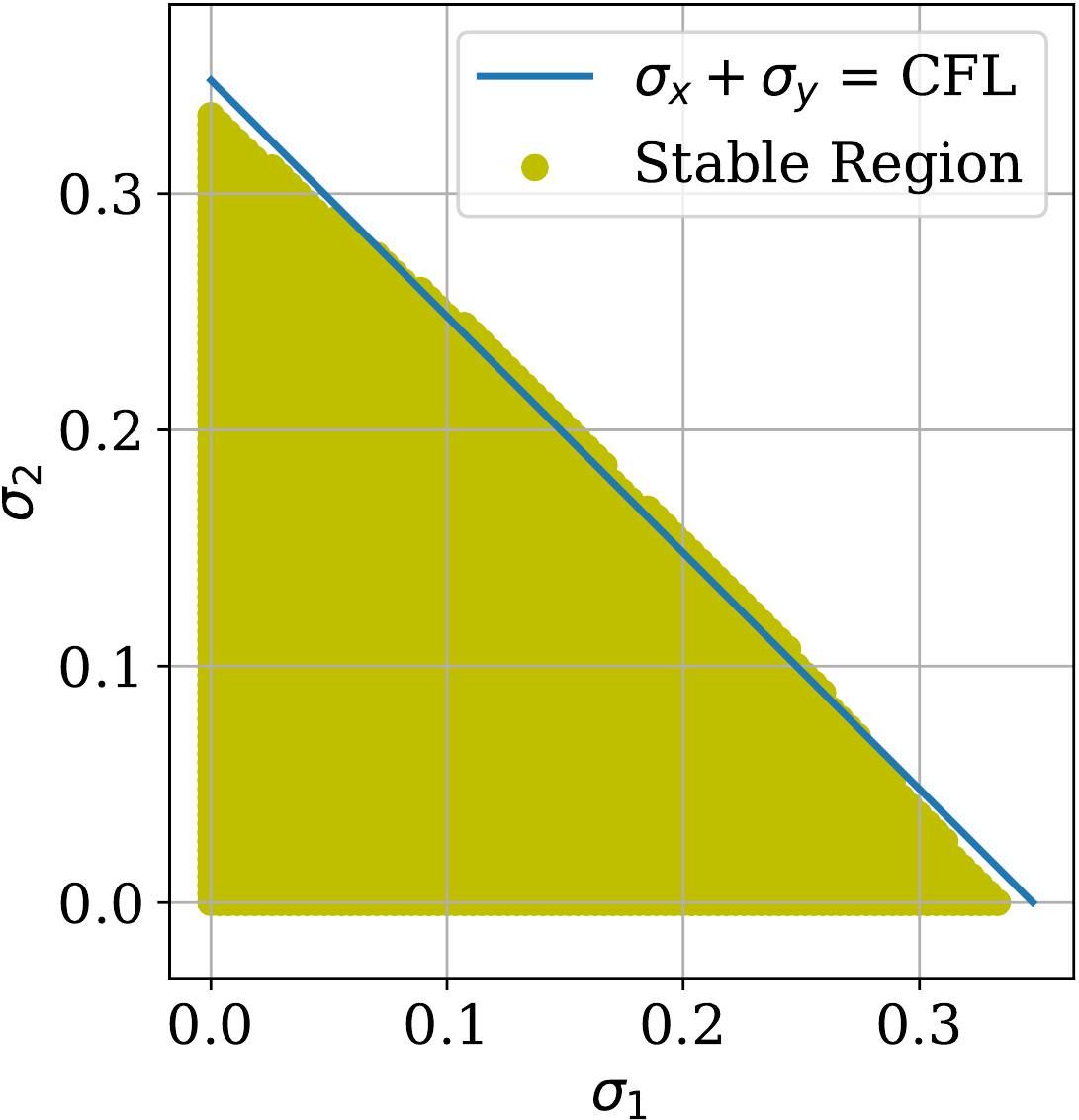} &
\includegraphics[width=0.22\textwidth]{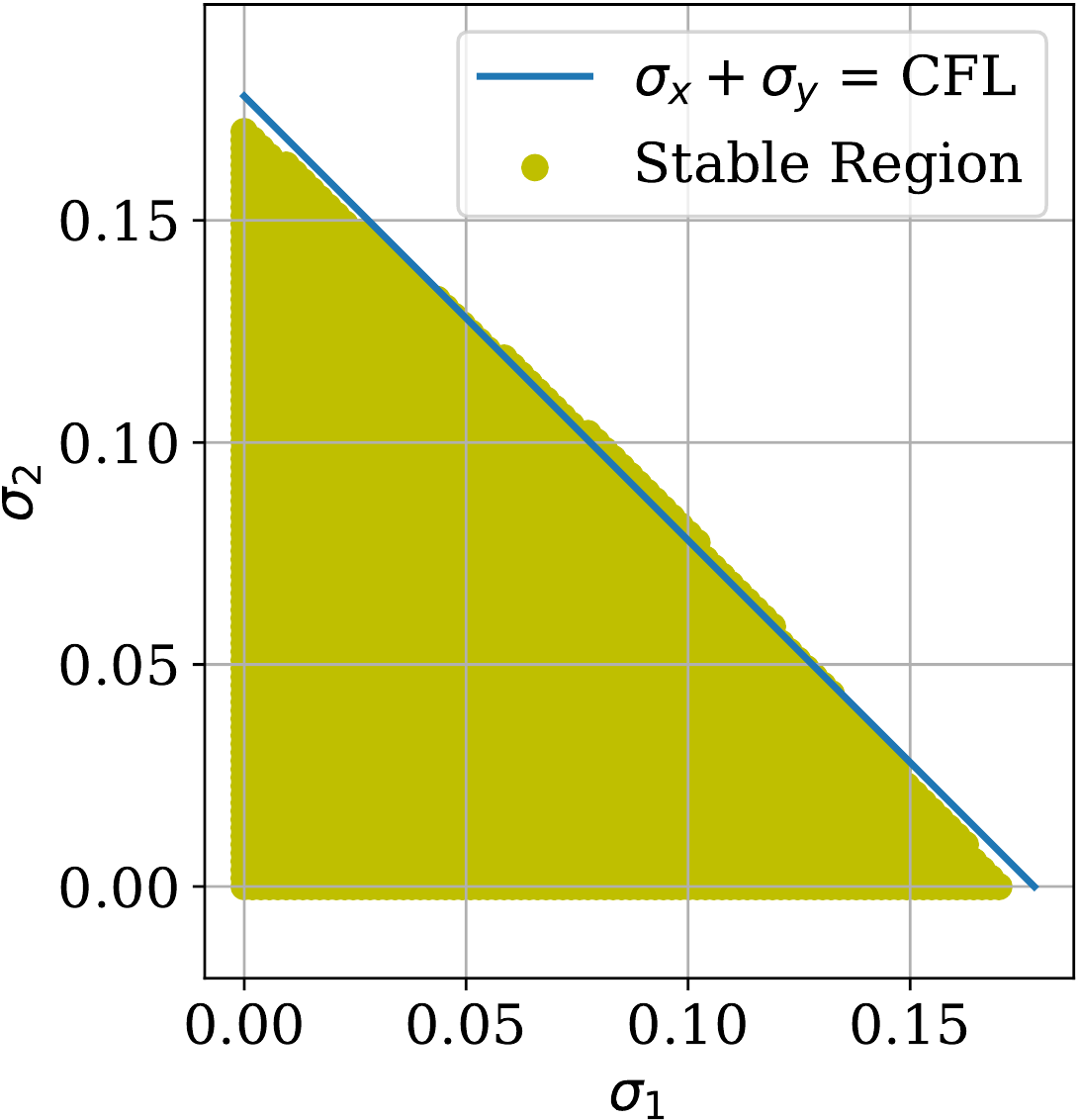} &
\includegraphics[width=0.22\textwidth]{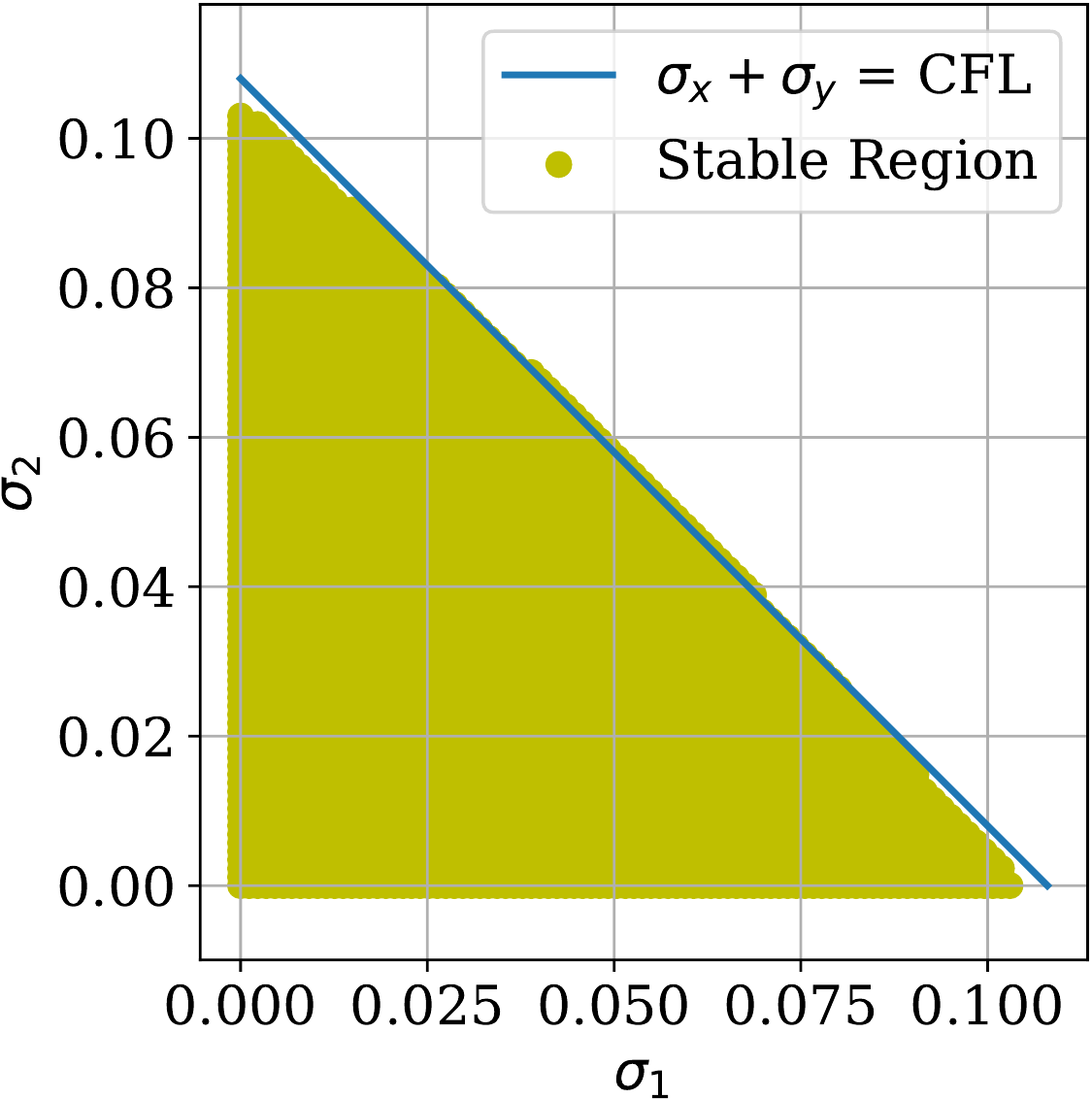} \\
(a) & (b) & (c) & (d)
\end{tabular}
\end{center}
\caption{Stability regions of LWFR scheme with  the  $g2$ correction function  and D2 dissipation model in two dimensions. (a) $N=1$, (b) $N=2$, (c) $N=3$, (d) $N=4$.}
\label{fig:2dcfl_g2}
\end{figure}

Consider the linear advection equation
\begin{equation}\label{eq:2dadv}
u_t + a_1 u_x + a_2 u_y = 0
\end{equation}
where $(a_1, a_2)$ is a constant velocity. We first write the LW scheme in matrix form which helps to derive the Fourier amplification term. Let us define the matrix of solution and flux values by
\[
\vu_e(i,j) = u_{ij}^e, \qquad \vf_e(i,j) = a_1 u_{ij}^e, \qquad \vg_e(i,j) = a_2 u_{ij}^e
\]
In the Lax-Wendroff procedure, the time derivative of the solution at all the solution points is given by
\[
\vu^{(m)}_e = - \sigma_1 \vD \vu^{(m-1)}_e - \sigma_2 \vu^{(m-1)}_e \vD^\top, \quad m=1,2,\ldots,N
\]
where $\sigma_1, \sigma_2$ are the CFL numbers along $x,y$ directions, respectively, which are given by
\[
\sigma_1 = \frac{a_1 \Delta t}{\Delta x_e}, \qquad
\sigma_2 = \frac{a_2 \Delta t}{\Delta y_e}
\]
Then the time average solution and fluxes are given by
\[
\vU_e = \sum_{m=0}^N \frac{\vu_e^{(m)}}{(m+1)!}, \qquad
\vF_e = a_1 \vU_e, \qquad
\vG_e = a_2 \vU_e
\]
To perform the Fourier analysis, we must write the scheme in matrix-vector form. To do this, let us renumber the two dimensional indices $(i,j)$ which denote solution points, into the one dimensional numbering by the following transformation
\[
k = i + (N+1)j, \qquad 0 \le i,j \le N
\]
Then $k$ takes the values $0$ to $M = (N+1)^2-1$. If $\phi_e \in \re^{(N+1) \times (N+1)}$ is some quantity defined at the solution points, we let $\rnum{\phi_e} \in \re^{M+1}$ denote the same renumbered as above. After renumbering, the matrix-matrix products become
\[
\rnum{\vA \phi_e} = R_1(\vA) \rnum{\phi_e}, \qquad \rnum{\phi_e \vA} = R_2(\vA) \rnum{\phi_e}
\]
where
\[
R_1(\vA) = \vI \otimes \vA, \qquad R_2(\vA) = \vA^\top \otimes \vI
\]
with $\otimes$ denoting the kronecker product. Then the renumbering of the solution time derivatives and time average solution and fluxes are given by
\[
\rnum{\vu^{(m)}_e} = \left( - \sigma_1 R_1(\vD) - \sigma_2 R_2(\vD^\top) \right) \rnum{\vu_e^{(m-1)}} =: \vH_1 \rnum{\vu_e^{(m-1)}}, \quad m=1,2,\ldots,N
\]
\[
\rnum{\vU_e} = \left( \sum_{m=0}^N \frac{\vH_1^{m}}{(m+1)!} \right) \rnum{\vu_e} =: \vT \rnum{\vu_e}, \qquad
\rnum{\vF_e} = a_1 \rnum{\vU_e}, \qquad
\rnum{\vG_e} = a_2 \rnum{\vU_e}
\]
Finally, the renumbered terms in the update equation~\eqref{eq:up2d} are given by
\[
\rnum{\vD_1 \vF_e} = R_1(\vD_1) \rnum{\vF_e} = a_1 R_1(\vD_1) \vT \rnum{\vu_e}, \qquad \rnum{\vG_e \vD_1^\top} = R_2(\vD_1^\top) \rnum{\vG_e} = a_2 R_2(\vD_1^\top) \vT \rnum{\vu_e}
\]
so that the cell terms can be written as
\[
\rnum{ \frac{\Delta t}{\Delta x_e} \vD_1 \vF_e + \frac{\Delta t}{\Delta y_e} \vG_e \vD_1^\top } = \left( \sigma_1 R_1(\vD_1) \vT + \sigma_2 R_2(\vD_1^\top) \vT \right) \rnum{\vu_e}
\]
For the terms involving the numerical flux, let us consider the case $a_1 \ge 0$, $a_2 \ge 0$. Let $\vu_{l}, \vu_{r}, \vu_{b}$ denote the solution in the elements to the left, right and bottom of the $e$'th element. Then, for the upwind flux which is obtained for dissipation model D2, we can renumber the terms involving the numerical flux as follows
\[
\frac{\Delta t}{\Delta x_e} \rnum{ \vb_L \vF_\emh^\top + \vb_R \vF_\eph^\top}=\frac{\Delta t}{\Delta x_e}\rnum{a_1 \vb_L \vV_R^\top \vU_l + a_1 \vb_R \vV_R^\top \vU_e}=\sigma_1 R_1(\vb_L \vV_R^\top) \vT \rnum{\vu_l} + \sigma_1 R_1(\vb_R \vV_R^\top) \vT\rnum{\vu_e}
\]
\[
\frac{\Delta t}{\Delta y_e}\rnum{\vG_\emh \vb_L^\top + \vG_\eph \vb_R^\top} =\frac{\Delta t}{\Delta y_e}\rnum{ a_2 \vU_b \vV_R \vb_L^\top + a_2 \vU_e \vV_R \vb_R^\top} = \sigma_2 R_2(\vV_R \vb_L^\top) \vT \rnum{\vu_b} + \sigma_2 R_2(\vV_R \vb_R^\top) \vT \rnum{\vu_e}
\]
The update equation can be written as
\[
\rnum{\vu_e^{n+1}} = \vA_l \rnum{\vu_{l}^n} + \vA_e \rnum{\vu_e^n} + \vA_b \rnum{\vu_{b}^n}
\]
where the coefficient matrices are given by
\[
\vA_l=-\sigma_1 R_1(\vb_L \vV_R^\top)\vT,\quad \quad \vA_b=-\sigma_2 R_2(\vV_R \vb_L^\top)\vT
\]
\[
\vA_e=I-\sigma_1 R_1(\vD_1) \vT - \sigma_2 R_2 (\vD_1^\top) \vT - \sigma_1 R_1 (\vb_R \vV_R^\top) \vT -\sigma_2 R_2 (\vV_R \vb_R^\top) \vT
\]
Assuming a solution of the form $\vu_e^n = \hat \vu_k^n \exp(\im(k_1 x_e + k_2 y_e))$, we get the amplification equation
\[
\rnum{\hat \vu_k^{n+1}} = (\vA_l \exp(-\im \kappa_1) + \vA_e + \vA_b \exp(-\im \kappa_2)) \rnum{\hat \vu_k^n} =: H(\sigma_1, \sigma_2; \kappa_1, \kappa_2) \rnum{\hat \vu_k^n}
\]
where $\kappa_1 = k_1 \Delta x$ and $\kappa_2 = k_2\Delta y$. For stability, it is required that the spectral radius of the matrix $H(\sigma_1, \sigma_2; \kappa_1, \kappa_2)$ is less than or equal to one for all wave numbers $\kappa_1, \kappa_2 \in[0,2\pi].$ Numerically, we compute the  region consisting of the pairs $(\sigma_1,\sigma_2)$ that ensures the stability. These regions for different degrees with dissipation model D2 are given in Figures~(\ref{fig:2dcfl_radau}) and (\ref{fig:2dcfl_g2}) for the Radau and g2 correction functions, respectively. We set $\mbox{CFL}=2c,$ where $c:=\max\{\sigma: (\sigma,\sigma) \mbox{ is a stable pair}\}$ which is the CFL limit when the advection velocity is in the direction $(1,1)$. These CFL numbers for different degrees are given in Table~\ref{tab:2Dcfl}. We see from the figures that the stable domain is bounded by a straight line except in case of degree $N=1$ so that this region is given by
\begin{equation}\label{eq:2dcfldom}
|\sigma_1| + |\sigma_2| \le \cfl
\end{equation}
If the advection velocity is along the $x$ or $y$ axis, the CFL corresponds to that of the 1-D case, but if the velocity is at an angle to the grid, then the allowed time step is reduced. This is because of the one dimensional numerical fluxes employed at the cell faces which couple each cell only to its left/right and bottom/top cells, without any coupling with the diagonal neighbours.

\begin{table}
    \centering
    \begin{tabular}{|c|c|c|c|c|}
    \hline
    $N$ & 1 & 2  & 3 & 4  \\
    \hline
    Radau & 0.259 & 0.166 & 0.101   & 0.067 \\
    \hline
    g2 & 0.511  & 0.348 & 0.178 & 0.108 \\
    \hline
    \end{tabular}
    \caption{Two dimensional CFL condition defined in equation~\eqref{eq:2dcfldom} for LWFR method with the dissipation model D2 and two correction functions}
    \label{tab:2Dcfl}
\end{table}

\section{Equivalence with direct flux reconstruction (DFR)}\label{sec:frdfr}
The direct flux reconstruction method does not require the choice of correction function. Following the ideas of~\cite{Romero2016}, we will prove that the LWFR scheme using Gauss-Legendre points and Radau correction function described in Section~\ref{sec:reconstruction} is equivalent to the LWDFR scheme described in Section~\ref{sec:DFR}, by showing that the $\vb_L, \vb_R, \vD_1$ are same for both. 
\paragraph{Equivalence of $\vb_L$.} We begin by proving the claim for $\vb_L$. For the FR scheme, we have
\[
\vb^\text{FR}_L = \begin{bmatrix}
g_L'(\xi_0) \\
\vdots \\
g_L'(\xi_N) \end{bmatrix}
\]
where $g_L$ is the Radau correction function and $\{\xi_i, 0 \le i \le N\}$ are Gauss-Legendre quadrature points on the interval $[0,1]$. For the DFR scheme, we have
\[
\vb^\text{DFR}_L = \begin{bmatrix}
\tilde\ell_{-1}'(\xi_0) \\
\vdots \\
\tilde\ell_{-1}'(\xi_N) \end{bmatrix}
\]
where $\tilde\ell_i$'s are Lagrange polynomials associated to the points $\{\xi_i, -1 \le i \le N+1\}$ where $\xi_{-1}=0$ and  $\xi_{N+1}=1.$  Since the $N+1$ zeros of $L_{N+1}$ are also zeros of $\tilde\ell_{-1}$ and $\tilde\ell_{-1}(0)=1$, $\tilde\ell_{-1}(1) = 0$, we must have
\[
\tilde\ell_{-1}(\xi) = (-1)^N (\xi-1) L_{N+1}(2\xi-1)
\]
To prove our claim, we need to prove
\[
\od{}{\xi}(g_L-\tilde\ell_{-1})(\xi_i)=0, \qquad i=0,1,\dots, N
\]
i.e.,
\[
L_N'(2\xi_i-1) - L_{N+1}(2\xi_i-1)-(2\xi_i-1)L_{N+1}'(2\xi_i-1)=0, \qquad i=0,1,\dots, N
\]
To work in $[-1,1]$ which is the natural domain of Legendre polynomials, we define the residual $R(\eta) = L_N'(\eta) - L_{N+1}(\eta)-\eta L_{N+1}'(\eta)$ so we have to show 
\[
R(\eta_i)=0,\qquad i=0,1,\dots,N
\]
where $\eta_i = 2\xi_i - 1$ are the Gauss-Legendre points in $[-1,+1]$. Using the recurrence relations
\begin{align*}
(1-\eta^2)L_{N+1}'(\eta) &= (N+1)[L_N(\eta)-\eta L_{N+1}(\eta)] \\
(1-\eta^2)L_N'(\eta) &= (N+1)[\eta L_N(\eta) - L_{N+1}(\eta)] 
\end{align*}
we get
\[
R(\eta) = -(N+2)L_{N+1}(\eta)
\]
proving that $R(\eta_i)=0$ for all $i=0,1,\dots,N$ since these $\eta_i$ are the zeros of $L_{N+1}$. Thus, $\vb^\text{FR}_L=\vb^\text{DFR}_L$. The claim for right correction follows analogously. 
\paragraph{Equivalence of $\vD_1$.} Writing $\vb_L = \vb^\text{FR}_L=\vb^\text{DFR}_L$ and $ \vb_R=\vb^\text{FR}_R=\vb^\text{DFR}_R$, proving that the $\vD_1$ matrices are same for both schemes is equivalent to showing that
\[
\vD=\vD_1^\text{DFR}+\vb_L\vV_L^\top+\vb_R\vV_R^\top,
\]
where $\vD$ is the differentiation matrix on Gauss-Legendre points. Further, to show that these two matrices are equal, it is enough to prove that their action on a set of $N+1$ linearly independent column vectors is the same. For this, we consider an arbitrary polynomial $p(\xi)$ of degree less than or equal to $N$, and let $\vp=[p(\xi_0),\cdots,p(\xi_N)]^\top$ and $\vp^\prime=[p^\prime(\xi_0),\cdots,p^\prime(\xi_N)]^\top = \vD \vp$. We have
\[
\vb_L\vV_L^\top \vp = \vb_L \sum_{i=0}^N p(\xi_i) \ell_i(0) = \vb_L p(0) = p(0)[\tilde \ell_{-1}'(\xi_0),\cdots,\tilde \ell_{-1}'(\xi_N)]^\top
\]
and
\[
\vb_R\vV_R^\top \vp = \vb_R \sum_{i=0}^N p(\xi_i) \ell_i(1) = \vb_R p(1) = p(1)[\tilde \ell_{N+1}'(\xi_0),\cdots,\tilde \ell_{N+1}'(\xi_N)]^\top
\]
As $p$ is a polynomial of degree less than or equal to $N$, we can write 
$$
p(\xi)=\sum_{i=-1}^{N+1}p(\xi_i)\tilde\ell_{i}(\xi), \qquad p^\prime(\xi)=\sum_{i=-1}^{N+1}p(\xi_i)\tilde\ell^\prime_{i}(\xi)
$$
We get 
\[
(\vD_1^\text{DFR}+\vb_L\vV_L^\top+\vb_R\vV_R^\top)\vp=\begin{bmatrix}
\sum_{j=0}^{N}p(\xi_0)\tilde \ell^\prime_j(\xi_0)\\
\vdots\\
\sum_{j=0}^{N}p(\xi_N)\tilde \ell^\prime_j(\xi_N)
\end{bmatrix}+\begin{bmatrix}
p(0)\tilde \ell_{-1}'(\xi_0)\\
\vdots\\
p(0)\tilde \ell_{-1}'(\xi_N)
\end{bmatrix}+
\begin{bmatrix}
p(1)\tilde \ell_{N+1}'(\xi_0)\\
\vdots\\
p(1)\tilde \ell_{N+1}'(\xi_N)
\end{bmatrix}=\vp^\prime = \vD \vp
\]
for all $\vp \in \re^{N+1}$,  which proves the claim.

\bibliographystyle{siam}
\bibliography{references}
\end{document}